\documentclass[11pt]{article}
 \usepackage{amsfonts,amssymb,enumerate,epsf,amsmath,bbm}
 \usepackage[mathscr]{eucal}
 \usepackage{graphicx}
 \usepackage{xcolor}
 \usepackage{standalone}
\date{}
\usepackage{hyperref}
\usepackage{verbatim}
\usepackage{tikz}
\usepackage{pgfplots}
\usepackage{pgfplotstable}
\usepgfplotslibrary{colorbrewer}
\usepgfplotslibrary{groupplots}
\usetikzlibrary{patterns}
\usetikzlibrary{arrows}
\usepackage{authblk}
\usepackage{mathrsfs}

\definecolor{grigio}{RGB}{232,232,232}

 \usepackage[english]{babel}
 \def\botcaption#1#2{\medskip\centerline{{\scshape #1.}\kern8pt
 {\rm #2}}\bigskip}

  \newtheorem{ittheorem}{Theorem}[section]
 \newtheorem{itlemma}[ittheorem]{Lemma}
 \newtheorem{itproposition}[ittheorem]{Proposition}
 \newtheorem{itdefinition}[ittheorem]{Definition}
 \newtheorem{itremark}[ittheorem]{Remark}
 \newtheorem{itclaim}[ittheorem]{Claim}
 \newtheorem{itcorollary}[ittheorem]{Corollary}
 \numberwithin{equation}{section}

 \newenvironment{theorem}{
 \begin{ittheorem}}{\end{ittheorem}}

 \newenvironment{lemma}{
 \begin{itlemma}}{\end{itlemma}}

 \newenvironment{proposition}{
 \begin{itproposition}}{\end{itproposition}}

 \newenvironment{definition}{
 \begin{itdefinition}}{\end{itdefinition}}

 \newenvironment{remark}{
 \begin{itremark}}{\end{itremark}}

 \newenvironment{claim}{
 \begin{itclaim}}{\end{itclaim}}

\newenvironment{corollary}{
	\begin{itcorollary}}{\end{itcorollary}}

 \newenvironment{proof}{\noindent {\bf Proof.\,}
 }{\hspace*{\fill}$\qed$\medskip}

 \newenvironment{proof*}{\noindent {\bf Proof\,}
}{\hspace*{\fill}$\qed$\medskip}

 \newcommand{\be}[1]{\begin{equation}\label{#1}}
 \newcommand{\ee}{\end{equation}}

 \newcommand{\bl}[1]{\begin{lemma}\label{#1}}
 \newcommand{\el}{\end{lemma}}

 \newcommand{\br}[1]{\begin{remark}\label{#1}}
 \newcommand{\er}{\end{remark}}

 \newcommand{\bt}[1]{\begin{theorem}\label{#1}}
 \newcommand{\et}{\end{theorem}}

 \newcommand{\bd}[1]{\begin{definition}\label{#1}}
 \newcommand{\ed}{\end{definition}}

 \newcommand{\bcl}[1]{\begin{claim}\label{#1}}
 \newcommand{\ecl}{\end{claim}}

 \newcommand{\bp}[1]{\begin{proposition}\label{#1}}
 \newcommand{\ep}{\end{proposition}}

 \newcommand{\bc}[1]{\begin{corollary}\label{#1}}
 \newcommand{\ec}{\end{corollary}}

 \newcommand{\bpr}{\begin{proof}}
 \newcommand{\epr}{\end{proof}}

 \newcommand{\bi}{\begin{itemize}}
 \newcommand{\ei}{\end{itemize}}

 \newcommand{\ben}{\begin{enumerate}}
 \newcommand{\een}{\end{enumerate}}

\def\ss{{\cX^s}}
\def\sm{{\cX^m}}

 \def \ba {\begin{array}}
 \def \ea {\end{array}}

 \def \qed {{\square\hfill}}
 \def \Z {{\mathbb Z}}
 \def \R {{\mathbb R}}
 
 \def \N {{\mathbb N}}
 \def \P {{\mathbb P}}
 \def \E {{\mathbb E}}

 \def \ra {\rightarrow}

 \def \cN {{\cal N}}
 \def \cS {{\cal S}}
 \def \cE {{\cal E}}
 \def \cF {{\cal F}}
 \def \cA {{\cal A}}
 \def \cR {{\cal R}}
  \def \cG {{\cal G}}
\def \cI {{\cal I}}
\def \cL{{\cal L}}
 \def \cQ {{\cal Q}}
 \def \cC {{\cal C}}
 \def \cD {{\cal D}}
 \def \cX {{\cal X}}
 \def \cB {{\cal B}}
 \def \cM {{\cal M}}
 \def \cW {{\cal W}}
 
 \def \cP {{\cal P}}
 \def \cL {{\cal L}}
 \def \cV {{\cal V}}
 \def \cZ {{\cal Z}}
 \def \cJ {{\cal J}}

 \def \G {{\Gamma}}
 \def \L {{\Lambda}}
 
 \def \a {{\alpha}}
 \def \b {{\beta}}
 \def \e {{\varepsilon}}
 \def \D {{\Delta}}
 \def \r {{\rho}}
\def \m {{\mu}}

 \def \h {{\eta}}
 \def \s {{\sigma}}
 \def \z {{\zeta}}
 \def \g {{\gamma}}
 \def \t {{\tau}}
 \def \o {{\omega}}
 \def \d {{\delta}}
\def \p {{\pi}}
\def \x {{\xi}}

 \def \CAPA {{\hbox{\footnotesize\rm CAP}}}

\def\pieno{{\blacksquare}}
\def\vuoto{{\square}}

\def\cigeo{{\cal C}_{is}^*}
\def\cwgeo{{\cal C}_{wa}^*}

\def\gi{\G_{is}^*}

\def\gw{\G_{wa}^*}

\begin{document}

\title{Critical Droplets and sharp asymptotics for Kawasaki dynamics \\ with weakly anisotropic interactions. Extended version}
	
	\author[
	{}\hspace{0.5pt}\protect\hyperlink{hyp:email1}{1},\protect\hyperlink{hyp:affil1}{a} 
	]
	{\protect\hypertarget{hyp:author1}{Simone Baldassarri}}
	
	\author[
	{}\hspace{0.5pt}\protect\hyperlink{hyp:email2}{2},\protect\hyperlink{hyp:affil1}{a,b}
	]
	{\protect\hypertarget{hyp:author2}{Francesca R.\ Nardi}}

	\affil[ ]{\centering
		\small\parbox{365pt}{\centering
			\parbox{5pt}{\textsuperscript{\protect\hypertarget{hyp:affil2}{a}}}Dipartimento di Matematica e Informatica ``Ulisse Dini", Universit\`{a} degli Studi di Firenze, Firenze, Italy.
		}
	}
	
	\affil[ ]{\centering
		\small\parbox{365pt}{\centering
			\parbox{5pt}{\textsuperscript{\protect\hypertarget{hyp:affil2}{b}}}Department of Mathematics and Computer Science, Eindhoven University of Technology, Eindhoven, the Netherlands.
		}
	}
	
	\affil[ ]{\centering
		\small\parbox{365pt}{\centering
			\parbox{5pt}{\textsuperscript{\protect\hypertarget{hyp:email1}{1}}}\texttt{\footnotesize\href{mailto:simone.baldassarri@unifi.it}{simone.baldassarri@unifi.it}},
			\parbox{5pt}{\textsuperscript{\protect\hypertarget{hyp:email2}{2}}}\texttt{\footnotesize\href{mailto:francescaromana.nardi@unifi.it}{francescaromana.nardi@unifi.it}},
		}
	}

	\maketitle
	
		\vspace{-1.2cm}

		\begin{center}
		{\it Unfortunately my coauthor Francesca Nardi passed away on 21 October 2021 \\ before the review process of the paper. I wish to thank her for the bright \\ person and talented mathematician she was.}
		
	\end{center}
	
	\begin{abstract}
	In this paper we analyze metastability and nucleation in the context of the Kawasaki dynamics for the two-dimensional Ising lattice gas at very low temperature with periodic boundary conditions. Let $\b>0$ be the inverse temperature and let $\L\subset\L^\b\subset\Z^2$ be two boxes. We consider the asymptotic regime corresponding to the limit as $\beta\ra\infty$ for finite $\L$ and $\lim_{\b\ra\infty}\frac{1}{\b}\log|\L^\b|=\infty$. We study the simplified model, in which particles perform independent random walks on $\L^\b\setminus\L$, while inside $\L$ particles perform simple exclusion, but when they occupy neighboring sites they feel a binding energy $-U_1<0$ in the horizontal direction and $-U_2<0$ in the vertical one. Thus the Kawasaki dynamics is conservative inside $\Lambda^\b$. The initial configuration is chosen such that $\L$ is empty and $\rho|\L^\b|$ particles are distributed randomly over $\L^\b\setminus\L$. Our results will use a deep analysis of a local model, i.e., particles perform Kawasaki dynamics inside $\L$ and along each bond touching the boundary of $\Lambda$ from the outside to the inside, particles are created with rate $\rho=e^{-\Delta\beta}$, while along each bond from the inside to the outside, particles are annihilated with rate $1$, where $\Delta>0$ is an activity parameter. Thus, in the local model the boundary of $\Lambda$ plays the role of an infinite gas reservoir with density $\rho$. We take $\Delta\in{(U_1,U_1+U_2)}$, so that the empty (resp.\ full) configuration is a metastable (resp.\ stable)  pair of configurations. We investigate how the transition from empty to full takes place in the local model with particular attention to the critical configurations that asymptotically have to be crossed with probability 1. To this end, we provide a \emph{model-independent strategy} to identify unessential saddles (that are not in the union of minimal gates) for the transition from the metastable (or stable) to the stable states and we apply this method to the local model. The derivation of further geometrical properties of the saddles allows us to use this strategy and to identify the union of all the minimal gates for the nucleation in the isotropic ($U_1=U_2$) and weakly anisotropic ($U_1<2U_2$) cases. More precisely, in the weakly anisotropic case we are able to identify the full geometry of the minimal gates and their boundaries. We observe very different behaviors compared to the strongly anisotropic case ($U_1>2U_2$). Moreover, we investigate the asymptotics, mixing time and spectral gap for isotropic and weakly anisotropic cases.

		\medskip
		{\it AMS} 2020 {\it subject classifications.} 60J10; 60K35; 82C20; 82C22; 82C26
		
		
		\medskip
		{\it Acknowledgment.} F.R.N. was partially supported by the Netherlands Organisation for Scientific Research (NWO) [Gravitation Grant number 024.002.003--NETWORKS]. This work was supported by the ``Gruppo Nazionale per l'Analisi Matematica e le loro Applicazioni" (GNAMPA-INdAM). The authors would like to thank Anna Gallo, Vanessa Jacquier and Cristian Spitoni for useful and fruitful discussions.

	\end{abstract}
	
	\tableofcontents

	\section{Introduction}
	{\bf Background.} Metastability is a dynamical phenomenon that occurs when a thermodynamic system is close to a first order phase transition, that takes place when some physical parameter such as the temperature, pressure or magnetic field abruptly changes. The phenomenon of metastability is characterized by the tendency of the system to remain for a long time in a state (the metastable state $m$) different from the stable states denoted by $\cX^s$. Moreover, the system leaves this apparent equilibrium at some random time performing a sudden transition to the stable state. This transition is called {\it metastability} or {\it metastable behavior}. In the study of metastablity there are three main issues that are tipically investigated. The first one is the study of the {\it typical transition time} from the metastable to the stable state. The second and third issues, that are physically more interesting, concern the geometrical description of the {\it gate configurations} (also called {\it critical configurations}) and the {\it tube of typical trajectories}, that we will discuss in the sequel. A central role in these descriptions is played by the {\it optimal paths}, i.e., the set of paths realizing the minimal value among all the paths going from $m$ to $\cX^s$ of the maximal energy reached in a single path. A basic notion for the second issue is the set of {\it saddles} $\cS(m,\cX^s)$ defined as the set of all maxima in the optimal paths between $m$ and $\cX^s$. Since we want to focus on the subsets of saddles that are typically visited during the last excursion from $m$ to $\cX^s$, we introduce the {\it gates} $\cW(m,\cX^s)$ from $m$ to $\cX^s$, defined as the subsets of $\cS(m,\cX^s)$ that are visited by all the optimal paths. The process, while performing the transition, is likely (with probability that tends to 1 in the chosen asymptotic regime) to cross a subset of the optimal paths. Therefore, gates are subsets of $\cS(m,\cX^s)$ that are typically visited during the transition. A {\it minimal gate} has the physical meaning of ``minimal set of critical configurations" and it is defined as a gate such that, after removing any configuration from it, the new set has not the gate property. For this purpose the characterization of the set $\cG(m,\cX^s)\subset \cS(m,\cX^s)$, defined as the union of minimal gates (see (\ref{defg})), is important. The third issue, the so-called tube of typical trajectories, is the set of typical paths followed by the system during the transition from the metastable to the stable state. We note that the hypotheses needed to discuss the gates are weaker than the ones necessary to completely characterize the tube of typical paths. The geometrical characterization of the union of minimal gates $\cG(m,\cX^s)$ is a central issue both from a probabilistic and from a physical point of view and it is a crucial step in the description of the typical trajectories. Metastability is an ubiquitous phenomenon with many examples from physical systems such as supersatured vapour, superheated and supercooled water, magnetic hysteresis loop, and from wireless networks. We remark that in several models proposed to describe ferromagnetic systems and analyzed in the literature in the context of F-W Markov chains evolving under Glauber dynamics, the minimal gate is unique but, in general, there may exist many minimal sets with the gate property, either distinct or overlapping. To model mathematically phenomena such as superheated or supercooled water, it is often proposed to use lattice gas models evolving according to Kawasaki dynamics since the dynamics conserves the number of particles.

	{\bf Anisotropic lattice gas.} In this paper we consider the metastable behavior of the two-dimensional Ising lattice gas at very low temperature and low density that evolves under Kawasaki dynamics, i.e., a discrete time Markov chain defined by the Metropolis algorithm with transition probabilities given in (\ref{defkaw}). Let $\b>0$ be the inverse temperature and let $\L\subset\L^\b\subset\Z^2$ be two boxes.  We consider periodic boundary conditions, fix the density $\rho$ of particles in $\L^\b$ and assume that $\lim_{\b\ra\infty}\frac{1}{\b}\log|\L^\b|=\infty$. Particles live and evolve in a conservative way in $\L^\b$, but when they occupy neighboring sites they feel a binding energy $-U_1<0$ in the horizontal direction and $-U_2<0$ in the vertical direction. Without loss of generality we may assume $U_1\geq U_2$. Unfortunately, we are unable to handle this model. Thus we consider a simplified model (see Section \ref{simplemodel} for more details) obtained by removing all the interactions outside $\L\setminus\partial^-\L$ and the exclusion outside $\L$, i.e., the dynamics of the gas in $\L^\b\setminus\L$ is that of independent random walks.

	For this simplified model, our goal is to identify the full geometry of the union of the minimal gates for the isotropic case (see Theorem \ref{simpleiso}) and we do this also for the weakly anisotropic case together with the characterization of a gate, the estimate in probability of the transition time from the metastable to the stable state and the description of the subcritical and supercritical rectangles (see Theorem \ref{simpleweak}). In order to obtain these results, we rely on the study of the local version of the model inside the box $\L$ with open boundary conditions (see Section \ref{S1.1} for more details). In particular, along each bond touching the boundary of $\Lambda$ from the outside to the inside, particles are created with rate $\rho=e^{-\Delta\beta}$, while along each bond from the inside to the outside, particles are annihilated with rate $1$, where $\Delta>0$ is an activity parameter. Thus, in the local model the boundary of $\Lambda$ plays the role of an infinite gas reservoir with density $\rho$. We take $\Delta\in{(U_1,U_1+U_2)}$, so that the empty (respectively full) configuration can be naturally related to metastability (respectively stability). We consider the asymptotic regime corresponding to finite volume $\L$ in the limit of large inverse temperature $\beta$. We investigate how the system {\it nucleates}, i.e., how it reaches $\pieno$ (box full of particles) starting from $\vuoto$ (empty box). We fix the parameters $U_1$, $U_2$ and $\D$ such that $\e:=U_1+U_2-\D>0$ is sufficiently small (we will consider in the sequel $0<\e\ll U_2$) and $U_1=U_2=U$ in what we call the isotropic case and $U_1<2U_2-2\e$ in the weakly anisotropic case.

	{\bf Goals.} For this local model, one of our goals of the paper is to investigate the isotropic and weakly anisotropic cases, giving the geometrical description of the set $\cG(\vuoto,\pieno)$ in Theorems \ref{giso} and \ref{gweak}. We will prove that in these different cases there are many distinct minimal gates, which we will geometrically characterize together with their union. Let us explain the strategy we adopt in our paper. In \cite[Theorem 5.1]{MNOS} the authors characterized in general the set $\cG(m,\cX^s)$ as the union of all the {\it essential saddles} for the transition from $m$ to $\cX^s$.  These are defined as the configurations in $\cS(m,\cX^s)$ that cannot be avoided with a short-cut of an optimal path (where a short-cut of a path $\o$ is a path $\o'$ whose set of maxima is a subset of the set of maxima in $\o$). Thanks to this equivalence (see Section \ref{modinddef} point 4 for the definition of essential saddle), we reduce our study to the identification of the set of all the essential saddles that has to be crossed during the transition between the metastable state $\vuoto$ and the stable state $\pieno$ in these two cases. First, we give a {\it model-independent strategy} that is useful to eliminate some unessential saddles, i.e., those that are not essential. More precisely, we require some model-dependent inputs (see Subsection \ref{modindep} for more details) in order to state Propositions \ref{selle1} and \ref{selle2}, in which we prove that two types of saddles are unessential and therefore are not in $\cG(m,\cX^s)$. In the sequel we apply this strategy to the isotropic and weakly anisotropic cases for two-dimensional models that evolve under Kawasaki dynamics where $m=\vuoto$ and $\cX^s=\{\pieno\}$. In order to do this, we need to verify that the required model-dependent inputs are valid for our model in the two cases. This study, together with the characterization of the essential saddles, relies on a detailed analysis of the motion of particles along the border of the droplet, which is a typical feature of the Kawasaki dynamics. Indeed in the metastable regime, particles move along the border of a droplet more rapidly than they arrive from the boundary of the box. More precisely, before the arrival of the next particle, we have that single particles attached to one side of a droplet tipycally detach (because $e^{U_1\b}\ll e^{\D\b}$ and $e^{U_2\b}\ll e^{\D\b}$), while bars of two or more particles tipycally do not detach (because $e^{\D\b}\ll e^{(U_1+U_2)\b}$). Roughly speaking, we will investigate the saddles that are crossed ``just before visiting" and ``just after visiting" the gate of the transition $\cG(\vuoto,\pieno)$. For the isotropic case, in \cite[Theorem 1.4.3(ii)]{BHN} the authors give a description of a gate, but they state that the geometrical characterization of the set $\cG(\vuoto,\pieno)$ is an open problem. We fill this gap in Theorem \ref{giso}. This suggests that an analogue detailed description for the weakly anistropic case is needed and we give it in Theorems \ref{thgateweak} and \ref{gweak}. Additionally, we prove sharp asymptotics for the transition time in Theorem \ref{sharptimeweak} and for the uniform entrance distribution in Theorem \ref{entrunif} for the weakly anisotropic case, while we refer to \cite[Theorems 1.4.4 and 1.4.3(ii)]{BHN} for the isotropic case. Moreover, in Theorem \ref{gapspettrale} we investigate the spectral gap and mixing time in both cases.

The corresponding analysis for the strongly anisotropic case, i.e., $U_1>2U_2$, is given in \cite{BN3} and we discuss the differences and similarities below.
Despite the fact that the structure of the gate is similar for the three cases, we emphasize that the entrance in them is very different. In particular, for the strongly anisotropic case there are two different mechanisms to enter the gate (see \cite[Lemma 6.17]{BN3}), while for the other two cases there is a unique one (see Lemma \ref{entrataweak}). On the other hand, it is clear that the properties that are related to the horizontal and vertical interactions are the same for both weakly and strongly anisotropic cases, while some properties that involve the motion of particles along the border of the droplet are very different. Intuitively, we can think of the weakly anisotropic case as a ``sort of interpolation" between the isotropic and strongly anisotropic ones, indeed it has some properties similar to the first, others to the latter. This specific difference between these cases,together with applications, motivates the rigorous investigation of the anisotropic cases. Moreover, we highlight this difference in the description of the set $\cG(\vuoto,\pieno)$. Indeed, for the isotropic case more motions along the border are allowed and thus a totally explicit geometric description of the set is more difficult (see Theorem \ref{giso}), but for the anisotropic cases we fully obtain it, since the condition $U_1\neq U_2$ makes more difficult the sliding of particles along the border of the droplet. Among the anisotropic cases, by Theorem \ref{gweak} and \cite[Theorem 4.10]{BN3} it is clear that the structure of the set $\cG(\vuoto,\pieno)$ ``strongly" depends on ``how large" is $U_1$ with respect to $U_2$, indeed in the case $U_1>2U_2$ less slidings along the border are allowed and thus the structure of the union of minimal gates is less rich than the weakly anisotropic case.

{\bf Literature.} Some properties of the metastable behavior for these three cases have already been derived in the literature. In particular, for isotropic interactions in \cite{HOS} the authors investigated the simplified and local models. In particular, they studied the asymptotic properties of the transition time together with an intrinsic description of a gate (see Section \ref{modinddef} point 4 for the precise definition). This paper initiated the study of the local model that we describe in the following discussion. In \cite{BHN} a geometric characterization of a subset of $\cG(\vuoto,\pieno)$ is given. We further improve this result in Theorem \ref{giso} by determining the minimal gates and their union $\cG(\vuoto,\pieno)$. For the three-dimensional lattice gas we refer to \cite{HNOS}, where the authors investigated the asymptotic properties of the transition time and an intrinsic description of a gate. Moreover, for both two and three-dimensional isotropic case, using the {\it potential theoretic approach} the authors investigated in \cite{BHN} the sharp asymptotics of the mean transition time, which includes the so-called {\it pre-factor}. They proved that it is a constant that asymptotically depends only on the size of the box and the cardinality of the gate that they identified, but not on the parameter $\b$. In the framework of the {\it pathwise approach} it is natural to study the third issue of metastability, namely the tube of typical trajectories realizing the transition between $\vuoto$ and $\pieno$.
This has been analyzed only in \cite{GOS} for two dimensions. Indeed, actually there are no known results about the tube for the three-dimensional isotropic case and for the anisotropic one. Concerning the anisotropic case, the asymptotic behavior of the transition time in probability, law and expectation has been derived for the weakly anisotropic interactions (see \cite{NOS}) and for the strongly anisotropic interactions (see \cite{BN}). In those papers there is also the geometric description of a gate, but with less control over the geometry of the minimal gates, their union and the entrance in it, with respect to what we give in Theorems \ref{thgateweak} and \ref{gweak} for the weakly anisotropic case.

	{\bf State of the art.} The first dynamical approach, known as \emph{pathwise approach}, was initiated in 1984 in \cite{CGOV}, developed in \cite{OS,OS2} and summerized in the monograph \cite{OV}. For Metropolis chains associated with statistical mechanics systems, metastability has been described by this approach in an elegant way in terms of the energy landscape associated to the Hamiltonian of the system. This approach focuses on the {\it dynamics} of the transition from metastable to stable states and it is so flexible that has been later developed to treat the tunnelling, namely the transition from a stable state to another stable state or stable states. Independently, a graphical approach was introduced in \cite{CC} and later used for Ising-like models \cite{CaTr}. Using the pathwise approach it is possible to obtain a detailed description of metastable behavior of the system and it made possible to answer all the three questions of metastability. A modern version of the pathwise approach can be found in \cite{{MNOS},{CNbc},{CNS2},{BNZ}}. In particular, in \cite{MNOS}, for the Metropolis markov chains, there are model-independent results concerning the transition time in probability, expectation and distribution, and concerning minimal gates and their union disentangled with respect to the tube of typical trajectories. In \cite{MNOS} the results on hitting times are obtained with minimal model-dependent knowledge,
	i.e., find all the metastable states and the minimal energy barrier which separates them from the stable states. In \cite[Sections 2-3]{CNbc} the authors prove model-independent results to treat systems with multiple metastable states and give a sufficient condition to identify them. In \cite{CNS2} the authors extend the results of \cite{MNOS} to general Markov chains (reversible and non reversible) with rare transitions setup, also called Freidlin-Wentzel Markov chains. These results are a useful tool to approach metastability for non-Metropolis systems such as Probabilistic Cellular Automata. In \cite[Section 3]{BNZ} the authors extended the model-independent framework of \cite{MNOS} to study the first hitting times from any starting configuration (not necessarily metastable) to any target subset of configurations (not necessarily the set of stable configurations). This approach developed over the years has been extensively applied to study metastability in Statistical Mechanics lattice models. In this context, this approach and the one that follows (\cite{BEGK,MNOS,OV}) have been developed with the aim of finding answers valid with maximal generality and to reduce as much as possible the number of model dependent inputs necessary to describe the metastable behavior of any given system.
	
	Another approach is the \emph{potential-theoretic approach}, initiated in \cite{BEGK}. We refer to \cite{BH} for an extensive discussion and applications to different models. In this approach, the metastability phenomenon is interpreted as a sequence of visits of the path to different metastable sets. This method focuses on a precise analysis of hitting times of these sets with the help of \emph{potential theory}. In the potential-theoretic approach the mean transition time is given in terms of the so-called \emph{capacities} between two sets. Crucially capacities can be estimated by exploiting powerful variational principles. This means that the estimates of the average crossover time that can be derived are much sharper than those obtained via the pathwise approach. The quantitative success of the potential-theoretic approach is however limited to the case of reversible Markov processes. 
	
	These mathematical approaches, however, are not equivalent as they rely on different definitions of metastable states (see \cite[Section 3]{CNbc} for a comparison) and thus involve different properties of hitting and transition times. The situation is particularly delicate for evolutions of infinite-volume systems, for irreversible systems, and degenerate systems, i.e., systems where the energy landscape has configurations with the same energy (as discussed in \cite{CNbc,CNS2,CNS2017}). More recent approaches are developed in \cite{BL1,BL2,BiGa}.
	
	Statistical mechanical models for magnets deal with dynamics that do not conserve the total number of particles or the total magnetization. They include single spin-flip Glauber dynamics and probabilistic cellular automata (PCA), that is, a parallel dynamics. The pathwise approach was applied in finite volume at low temperature in \cite{CGOV,NevSchbehavdrop,CaTr,KOd,KOs,CO,NO,CL,AJNT,BGNneg2021,BGNpos2021} for single-spin-flip Glauber dynamics and in \cite{CN,CNSp,CNS22,CNS2016} for parallel dynamics. The potential theoretic approach was applied to models at finite volume and at low temperature in \cite{BM,BHN,HNT1,HNT2,NS1,HNTA2018,BJN}. The more involved infinite volume limit at low temperature or vanishing magnetic field was studied in \cite{DS1,DS2,S2,SS,MO1,MO2,HOS,GHNOS,GN,BHS,CeMa,GMV} for Ising-like models under single-spin-flip Glauber and Kawasaki dynamics.

	{\bf Outline.} The outline of the paper is as follows. In Section \ref{S2} we define the simplified model with periodic boundary conditions and the local model with open boundary conditions, and the Kawasaki dynamics. In Section \ref{S3} we give some model-independent definitions in order to state our main model-independent results in Propositions \ref{selle1} and \ref{selle2}. In Section \ref{S4} we give some geometric definitions valid for the two models (see Section \ref{moddepdef}). We state our main results concerning the gates for the isotropic case in Section \ref{iso} and for the weakly anisotropic in Section \ref{wani}. The main results about the sharp asymptotics are given in Section \ref{sharpestimates}, while the results for the simplified model are given in Section \ref{simpletheorem}. In Section \ref{modproof} we prove the model-independent results that we apply to our model in Section \ref{dependentdef}. In particular, in Section \ref{sitigood} we give some model-dependent definitions, in Section \ref{lemmi3modelli} some tools that are useful in Section \ref{S6.4} for our model-dependent strategy. In Section \ref{proofiso} we give the proof of the main result for the isotropic interactions regarding the identification of the union of all the minimal gates (see Theorem \ref{giso}). In Section \ref{proofweak} we give the proof of the main results for the weakly anisotropic interactions regarding the description of the gate (see Theorem \ref{thgateweak}) and the geometric characterization of the union of all the minimal gates (see Theorem \ref{gweak}). In Section \ref{sharpasymptotics} we give the proof of the main theorems about the sharp asymptotics (see Theorems \ref{sharptimeweak}, \ref{entrunif} and \ref{gapspettrale}) and in Section \ref{simpleproof} we give the proof of the results for the simplified model (see Theorems \ref{simpleiso} and \ref{simpleweak}). In the Appendix \ref{appendice} we give additional explicit proofs and computations.

	\section{Definition of the model}
	\label{S2}
	
	\subsection{The model with open boundary conditions}
	\label{S1.1}
	
	Let $\Lambda=\{0,..,L\}^2\subset \Z^2$ be a finite box centered at the
	origin. The side length $L$ is fixed, but arbitrary, and later we will require $L$ to be sufficiently large. Let
	\be{inbd}
	\partial^- \L:= \{x\in\L: \exists\; y \notin\L\: |y-x|=1\},
	\ee
	
	\noindent be   the interior  boundary of $ \Lambda$ and let
	$ \Lambda_0:= \Lambda\setminus\partial^- \Lambda$ be the interior of $\L$.
	With each $x\in \Lambda$ we associate an occupation variable
	$\eta(x)$, assuming the values 0 or 1. A lattice configuration is
	denoted by $\eta\in {\cal X} =\{ 0,1\} ^{ \Lambda }$. Each configuration $\h\in \cX$ has an energy given by the following Hamiltonian:
	
	\be{hamilt} H(\eta):= -U_1 \sum_{(x,y)\in   \Lambda_{0,h}^{*}}
	\eta(x)\eta(y) -U_2\sum_{(x,y)\in   \Lambda_{0,v}^{*}} \eta(x)
	\eta(y)+ \D \sum _{x\in \L} \eta (x), \ee
	
	\noindent
	where $ \Lambda_{0,h}^{*}$ (resp.\ $ \Lambda_{0,v}^{*}$) is the set of
	the horizontal (resp.\ vertical) unoriented  bonds joining nearest-neighbors points in
	$ \Lambda_0$. Thus the interaction is acting only inside $
	\Lambda_0$; the binding energy associated to a horizontal
	(resp.\ vertical) bond is $-U_1<0$ (resp.\ $-U_2<0$). We may assume without
	loss of generality that $U_1\ge U_2$. 
	
	The grand-canonical Gibbs measure associated with $H$ is
	\be{misura} \m(\eta):= {  e^{- \b H(\eta) }\over Z} \qquad \h\in
	\cX, \ee
	
	\noindent
	where
	\be{partfunc} Z:=\sum_{\eta\in {\cal X}}e^{-\b H(\eta)}
	\ee
	
	\noindent
	is the so-called {\it partition function}.
	
	\subsection{Local Kawasaki dynamics}
	\label{S1.2}
	
	Next we define Kawasaki dynamics on $\L$ with boundary
	conditions that mimic the effect of an infinite gas reservoir
	outside $\L$ with density $ \r = e^{-\D\b}.$ Let $b=(x \to y)$ be
	an oriented bond, i.e., an {\it ordered} pair of nearest neighbour
	sites, and define
	
	\be{Loutindef}
	\ba{lll}
	\partial^* \L^{out} &:=& \{b=(x \to y): x\in\partial^- \L,
	y\not\in\L\},\\
	\partial^* \L^{in}  &:=& \{b=(x \to y): x\not\in
	\L, y\in\partial^-\L\},\\
	\L^{*, orie} &:=& \{b=(x \to y): x,y\in\L\},
	\ea
	\ee
	
	\noindent and put $ \bar\L^{*, orie}:=\partial^* \L ^{out}\cup
	\partial^* \L ^{in}\cup\L^{*,\;orie}$.
	Two configurations $  \eta,
	\eta'\in {\cal X}$ with $ \eta\ne \eta'$ are said to be {\it
		communicating states} if there exists a bond
	$b\in  \bar\L^{*,orie}$ such that $ \eta' = T_b \eta$, where $T_b   \eta$ is the
	configuration obtained from $ \eta$ in any of these ways:
	
	\begin{itemize}
		
		\item
		For $b=(x \to y)\in\L^{*,\;orie}$, $T_b \eta$ denotes the
		configuration obtained from $ \eta$ by interchanging particles
		along $b$:
		\be{Tint}
		T_b \h(z) =
		\left\{\ba{ll}
		\h(z) &\mbox{if } z \ne x,y,\\
		\h(x) &\mbox{if } z = y,\\
		\h(y) &\mbox{if } z = x.
		\ea
		\right.
		\ee
		
		\item
		For  $b=(x \to y)\in\partial^*\L^{out}$ we set:
		\be{Texit}
		T_b \h(z) =
		\left\{\ba{ll}
		\h(z) &\mbox{if } z \ne x,\\
		0     &\mbox{if } z = x.
		\ea
		\right.
		\ee
		
		\noindent
		This describes the annihilation of a particle along the border;
		
		\item
		For  $b=(x  \to y)\in\partial^*\L^{in}$ we set:
		\be{Tenter}
		T_b \h(z) =
		\left\{\ba{ll}
		\h(z) &\mbox{if } z \ne y,\\
		1     &\mbox{if } z=y.
		\ea
		\right.
		\ee
		
		\noindent
		This describes the creation of a particle along the border.
		
	\end{itemize}
	
	\noindent
	The Kawasaki dynamics is  the discrete time Markov chain
	$(\eta_t)_{t\in \mathbb{N}}$ on state space $ {\cal X} $ given by
	the following transition  probabilities: for  $  \eta\not= \eta'$:
	\be{defkaw}
	P( \eta,  \eta'):=\left\{\ba{ll}
	{ |\bar\L^{*,\;orie}|}^{-1} e^{-\b[H( \eta') - H( \eta)]_+}
	&\mbox{if }  \exists b\in \bar\L ^{*, orie}: \eta' =T_b \eta,  \\
	0   &\mbox{ otherwise, }  \ea \right.
	\ee
	
	\noindent
	where $[a]_+ =\max\{a,0\}$ and $P(\h,\h):=1-\sum_{\h'\neq\h}P(\h,\h')$.
	This describes a standard Metropolis dynamics with open boundary conditions: along each
	bond touching $\partial^-\L$ from the outside, particles are created with
	rate $\rho=e^{-\D\b}$ and are annihilated with rate 1, while inside $\L_0$ particles
	are conserved.
	Note that an exchange of occupation
	numbers $\h(x)$ for any $x$ inside the ring $ \L\setminus  \L_0$
	does not involve any change in energy.\par
	
	\br{p1}
	The stochastic dynamics defined by
	(\ref{defkaw}) is reversible w.r.t. Gibbs measure (\ref{misura}) corresponding
	to $ H$.
	\er
	
	\subsection{The model with periodic boundary conditions and its simplified version}
	\label{simplemodel}
	In this Section we consider a lattice gas model to study the metastable behavior of conservative systems. Let $\L^{\b}\subset\Z^2$ be a large finite box centered at the origin, with periodic boundary conditions. With each $x\in\L^\b$ we associate an occupation variable $\h^\b(x)$, assuming the values 0 and 1. A lattice gas configuration is denoted by $\h^\b\in\cX^\b=\{0,1\}^{\L^\b}$. We consider the interaction defined by the following Hamiltonian:
	\be{Hbeta}
	H(\h^\b):=-U_1 \sum_{(x,y)\in   \Lambda_{\b,h}^{*}}
	\eta^\b(x)\eta^\b(y) -U_2\sum_{(x,y)\in   \Lambda_{\b,v}^{*}} \eta^\b(x)
	\eta^\b(y)+ \D \sum _{x\in \L^\b} \eta^\b(x),
	\ee
	
	\noindent
	where $ \Lambda_{\b,h}^{*}$ (resp.\ $ \Lambda_{\b,v}^{*}$) is the set of the horizontal (resp.\ vertical) unoriented  bonds joining nearest-neighbors points in $\Lambda^\b$. Thus the binding energy associated to a horizontal (resp.\ vertical) bond is $-U_1<0$ (resp.\ $-U_2<0$). We may assume without loss of generality that $U_1\ge U_2$. We fix the particle density in $\L^\b$ at
	\be{density}
	\frac{1}{|\L^\b|}\sum_{x\in\L^\b}\h^\b(x)=\r=e^{-\D\b},
	\ee
	
	\noindent
	where $\D>0$ is an activity parameter. This corresponds to a total number of particles in $\L^\b$ equal to $N=\r|\L^\b|$. From (\ref{density}) we see that in order to have particles at all we must pick $|\L^\b|$ at least exponentially large in $\b$. This means that the regime where $\L^\b$ is fixed, as in the non conservative case considered in Section \ref{S1.2}, has no relevance here. We will in fact be interested in the regime
	\be{Hbeta}
	\D\in(U_1,U_1+U_2), \quad \b\ra\infty, \quad \displaystyle\lim_{\b\ra\infty}\frac{1}{\b}\log|\L^\b|=\infty.
	\ee
	
	\noindent
	We can not define the dynamics in $\L^\b$ as in Section \ref{S1.2} due to the presence of periodic boundary conditions, but if we change the open boundary conditions, then the model is too difficult to treat due to the behavior of the gas in $\L^\b\setminus\L$. Indeed in the conservative case the dynamics is not really local: particles must arrive from or return to the gas, which acts as a reservoir. It is therefore not possible to decouple the dynamics of the particles inside $\L$ from the dynamics of the gas in $\L^\b\setminus\L$. For this reason we consider a {\it simplified model}, in which we remove the interactions outside the box $\L_0=\L\setminus\partial^-\L$, where $\partial^-\L$ is defined in (\ref{inbd}). Moreover, we also remove the exclusion outside $\L$. Thus the dynamics of the gas in $\L^\b\setminus\L$ is that of independent random walks, and we replace the Hamiltonian defined in (\ref{Hbeta}) with 
	\be{Hbeta'}
	H'(\h^\b):=-U_1 \sum_{(x,y)\in   \Lambda_{0,h}^{*}}
	\eta^\b(x)\eta^\b(y) -U_2\sum_{(x,y)\in   \Lambda_{0,v}^{*}} \eta^\b(x)
	\eta^\b(y)+ \D \sum _{x\in \L^\b} \eta^\b(x).
	\ee
	
	\noindent
	We will use the results for the local model that concern the Hamiltonian defined in (\ref{hamilt}) to derive results for the simplified model, as will be come clear in Section \ref{simpleproof}. In this setting, we can view the model described in Section \ref{S1.1} as the local version of the one described in this Section, where the effect on $\L$ of the gas in $\L^\b\setminus\L$ may be described in terms of the creation of new particles with rate $\r=e^{-\D\b}$ at sites on the interior boundary of $\L$ and the annihilation of particles with rate 1 at sites on the exterior boundary of $\L$.
 	
	Let $\cV_N=\{\h^\b\in\cX^\b: N_{\L^\b}=N\}$ denote the set of configurations with $N$ particles. For $\h\in\cX=\{0,1\}^\L$, let $\nu_\h$ denote the canonical Gibbs measure conditioned on the configuration inside $\L$ being $\h$, i.e.,
	\be{Gibbscond}
	\nu_\h(\h^\b)=\frac{\nu_N(\h^\b)\mathbbm{1}_{\cJ_\b(\h)}(\h^\b)}{\nu_N(\cJ_\b(\h))}, \quad \h^\b\in\cX^\b,
	\ee
	
	\noindent
	where $\cJ^\b(\h)=\{\h^\b\in\cX^\b: {\h^\b}_{|\L}=\h\}$, with ${\h^\b}_{|\L}$ the restriction of $\h^\b$ to $\L$, and $\nu_N$ the canonical Gibbs measure defined as
	\be{canGibbs}
	\nu_N(\h^\b)=\frac{e^{-\b H(\h^\b)}\mathbbm{1}_{\cV_N(\h^\b)}}{Z_N}, \quad \h^\b\in\cX^\b,
	\ee 
	
	\noindent
	where
	\be{}
	Z_N=\displaystyle\sum_{\h^\b\in\cV_N}e^{-\b H(\h^\b)}
	\ee
	
	\noindent
	For $\h\in\cX$, write $\P_{\nu_\h}$ and $\E_{\nu_\h}$ to denote respectively the probability law and expectation for the Markov process $(\h_t)_{t\geq0}$ on $\cX^\b$ following the Kawasaki dynamics with Hamiltonian 
	\be{}
	H''(\h^\b)=-U_1 \sum_{(x,y)\in   \Lambda_{0,h}^{*}}
	\eta^\b(x)\eta^\b(y) -U_2\sum_{(x,y)\in   \Lambda_{0,v}^{*}} \eta^\b(x)
	\eta^\b(y),
	\ee
	
	\noindent
	given that $\h_0$ is chosen according to $\nu_\h$. Write $\vuoto$ to denote the empty configuration in $\Lambda$, i.e., $\vuoto^\b=\cJ^\b(\vuoto)$.

\section{Model-independent definitions and results}
\label{S3}
We will use italic capital letters for subsets of $\L$, script
capital letters for subsets of $\cX$, and boldface capital letters for events
under the Kawasaki dynamics. We use this convention in order to keep the
various notations apart. We will denote by $\P_{\h_0}$ the probability law of the Markov
process  $(\eta_t)_{t\geq 0}$ starting at $\h_0$ and by
$\E_{\h_0}$ the corresponding expectation.

\subsection{Model-independent definitions}
\label{modinddef}
\noindent
{\bf 1. Paths and hitting times.}
\bi

\item
A {\it path\/} $\o$ is a sequence $\o=(\o_1,\dots,\o_k)$, with
$k\in\N$, $\o_i\in\cX$ and $P(\o_i,\o_{i+1})>0$ for $i=1,\dots,k-1$.
We write $\o\colon\;\h\to\h'$ to denote a path from $\h$ to $\h'$,
namely with $\o_1=\h,$ $\o_k=\h'$. A set
$\cA\subset\cX$ with $|\cA|>1$ is {\it connected\/} if and only if for all
$\h,\h'\in\cA$ there exists a path $\o:\h\to\h'$ such that $\o_i\in\cA$
for all $i$. We indicate with $\o_1\circ\o_2$ the composition of two paths $\o_1$ and $\o_2$, namely if $\o_1=(\o_1^1,...,\o_k^1)$ and $\o_2=(\o_1^2,...,\o_m^2)$ then $\o_1\circ\o_2=(\o_1^1,...,\o_k^1,\o_1^2,...,\o_m^2)$.

\item[$\bullet$]
Given a non-empty
set $\cA\subset\cX$, define the {\it first-hitting time of} $\cA$
as
\be{tempo}
\t_{\cA}:=\min \{t\geq 0:  \eta_t \in \cA \}.
\ee
\ei
\noindent

\medskip
\noindent
{\bf 2. Min-max and communication height}
\bi

\item Given a function $f:\cX\ra\R$ and a subset $\cA\subseteq\cX$, we denote by 
\be{defargmax}
\arg \hbox{max}_{\cA}f:=\{\h\in\cA: f(\h)=\max_{\z\in\cA}f(\z)\}
\ee

\noindent
the set of points where the maximum of $f$ in $\cA$ is reached. If $\o=(\o_1,...,\o_k)$ is a path, in the sequel we will write $\arg \max_{\o}H$ to indicate $\arg \max_{\cA}H$, with $\cA=\{\o_1,...,\o_k\}$ and $H$ the Hamiltonian.

\item 
The {\it bottom} $\cF(\cA)$ of a  non-empty
set $\cA\subset\cX$ is the
set of {\it global minima} of the Hamiltonian $H$ in  $\cA$:
\be{Fdef}
\cF(\cA):=\arg \hbox{min}_{\cA}H=\{\h\in\cA: H(\h)=\min_{\z\in\cA} H(\z)\}.
\ee
For a set $\cA\subset\cX$ such that all the configurations have the same energy, with an abuse of notation we denote this energy by $H(\cA)$.

\item 
The {\it communication height} between a pair $\h$, $\h'\in\cX$ is
\be{}
\Phi(\h,\h'):= \min_{\o:\h\ra\h'}\max_{\z\in\o} H(\z).
\ee
\noindent
Given $\cA\subset\cX$, we define the {\it restricted communication height} between $\h,\h'\in\cA$ as
\be{}
\Phi_{|\cA}(\h,\h'):= \min_{\o:\h\ra\h'\atop\o\subseteq\cA}
\max_{\z\in\o} H(\z),
\ee
\noindent
where $(\o_1,...,\o_k)=\o\subseteq\cA$ means $\o_i\in\cA$ for every $i$.
\ei

\medskip
\noindent
{\bf 3. Stability level, stable and metastable states}
\bi

\item
We call
{\it stability level} of
a state $\z \in \cX$
the energy barrier
\be{stab}
V_{\z} :=
\Phi(\z,\cI_{\z}) - H(\z),
\ee

\noindent where $\cI_{\z}$ is the set of states with
energy below $H(\z)$:
\be{iz} \cI_{\z}:=\{\eta \in \cX: H(\eta)<H(\z)\}. \ee

\noindent 
We set $V_\z:=\infty$ if $\cI_\z$ is empty.

\item
We call {\it $V$-irreducible states}
the set of all states with stability level larger than  $V$:
\be{xv} 
\cX_V:=\{\h\in\cX: V_{\h}>V\}. 
\ee

\item
The set of {\it stable states} is the set of the global minima of
the Hamiltonian:
\be{st.st.}
\cX^s:=\cF(\cX). 
\ee

\item
The set of {\it metastable states} is given by
\be{st.metast.} \sm:=\{\h\in\cX:
V_{\h}=\max_{\z\in\cX\backslash \ss}V_{\z}\}. \ee
\noindent
We denote by $\G_m$ the stability level of the states in $\cX^m$.

\ei

\medskip
\noindent
{\bf 4. Optimal paths, saddles and gates}
\bi

\item 
We denote by $(\h\to\h')_{opt} $ the {\it set of optimal paths\/} as the set of all
paths from $\h$ to $\h'$ realizing the min-max in $\cX$, i.e.,
\be{optpath}
(\h\to\h')_{opt}:=\{\o:\h\to\h'\; \hbox{such that} \; \max_{\xi\in\o} H(\xi)=  \Phi(\h,\h') \}.
\ee

\item
The set of {\it minimal saddles\/} between
$\h,\h'\in\cX$
is defined as
\be{minsad}
\cS(\h,\h'):= \{\z\in\cX\colon\;\; \exists\o\in (\h\to\h')_{opt},
\ \o\ni\z \hbox{ such that } \max_{\xi\in\o} H(\xi)= H(\z)\}.
\ee

\item A saddle $\x\in\cS(\h,\h')$ is called {\it unessential} if for any $\o\in(\h\ra\h')_{opt}$ such that $\o\cap\x\neq\emptyset$ we have $\{\arg\max_{\o}H\}\setminus\{\x\}\neq\emptyset$ and there exists $\o'\in(\h\ra\h')_{opt}$ such that $\{\arg\max_{\o'}H\}\subseteq\{\arg\max_{\o}H\}\setminus\{\x\}$.

\item A saddle $\x\in\cS(\h,\h')$ is called {\it essential} if it is not unessential, i.e., if either

\bi
\item[(i)] there exists $\o\in(\h\ra\h')_{opt}$ such that $\{\hbox{arg max}_{\o}H\}=\{\x\}$ or
\item[(ii)] there exists $\o\in(\h\ra\h')_{opt}$ such that $\{\hbox{arg max}_{\o}H\}\supset\{\x\}$ and $\{\hbox{arg max}_{\o'}H\}\nsubseteq\{\hbox{arg max}_{\o}H\}\setminus\{\x\}$ for all $\o'\in(\h\ra\h')_{opt}$.
\ei

\item
Given a pair $\h,\h'\in\cX$,
we say that $\cW\equiv\cW(\h,\h')$ is a {\it gate\/}
for the transition $\h\to\h'$ if $\cW(\h,\h')\subseteq\cS(\h,\h')$
and $\o\cap\cW\neq\emptyset$ for all $\o\in (\h\to\h')_{opt}$. In words, a gate is a subset of $\cS(\h,\h')$ that is visited by all optimal paths.

\item
We say that $\cW(\h,\h')$ is a {\it minimal gate\/} for the transition
$\h\to\h'$ if it is a gate and for any $\cW'\subsetneq \cW(\h,\h')$ there
exists $\o'\in (\h\to \h')_{opt}$ such that $\o'\cap\cW'
=\emptyset$. In words, a minimal gate is a minimal subset of $\cS(\h,\h')$ by inclusion that is visited by all optimal paths.

\item
For a given pair of configurations $\h,\h'$, we denote by $\cG(\h,\h')$ the union of all minimal gates:
\be{defg}
\cG(\h,\h'):=\displaystyle\bigcup_{\cW(\h,\h') \hbox{ minimal gate}} \cW(\h,\h')
\ee

\ei

\subsection{Model-independent strategy in finite volume}
\label{genstrategy}
In this Section we give a general strategy to analyze the geometry of the set $\cG(m,\cX^s)$, where either $m\in\cX^m$ is a metastable state if we analyze metastability or $m\in\cX^s$ is a stable state if we analyze tunnelling between two stable states. Assume that we are in finite volume and $\cW(m,\cX^s)$ is a set of configurations that has been proven to be a gate. The following strategy is useful to eliminate some unessential saddles from the set $\cS(m,\cX^s)\setminus\cW(m,\cX^s)$ in order to determine the set $\cG(m,\cX^s)$. This strategy is more efficient if the gate proposed is minimal or union of minimal gates.

\subsubsection{Model-independent results for unessential saddles}
\label{modindep}

In order to state our results concerning the unessential saddles we need the following definitions.

\bi

\item A nonempty set $\cA\subset\cX$ is a cycle if it is either a singleton or it verifies the relation
\be{}
\max_{x,y\in\cA}\Phi(x,y)<\Phi(\cA,\cX\setminus\cA).
\ee

\noindent
See \cite[equation (3.40)]{CNS2}. In the case of Metropolis dynamics, this definition coincides with \cite[equation (2.7)]{MNOS}. 

\item Given $\s\in\cX$, $\G>0$ and $\cA$ a set of target configurations, we say that the {\it initial cycle for the transition from $\s$ to $\cA$ with depth $\G$} is
\be{ciclo}
\cC_{\cA}^{\s}(\G):=\s\cup\{\h\in\cX: \ \Phi(\s,\h)-H(\s)<\G=\Phi(\s,\cA)-H(\s)\}.
\ee

\noindent
Note that in definition (\ref{ciclo}) we emphasize the dependence on $\s$ and $\cA$ and that $\G$ is identified by them. Note that this definition of $\cC_{\cA}^{\s}(\G)$ concides with $\cC_{\cA}(\s)$ defined in \cite[equation (2.25)]{MNOS}.

\ei

\begin{figure}[h!]
\begin{center}
		\adjustbox{max width=\textwidth}{\includestandalone[mode=image|tex]{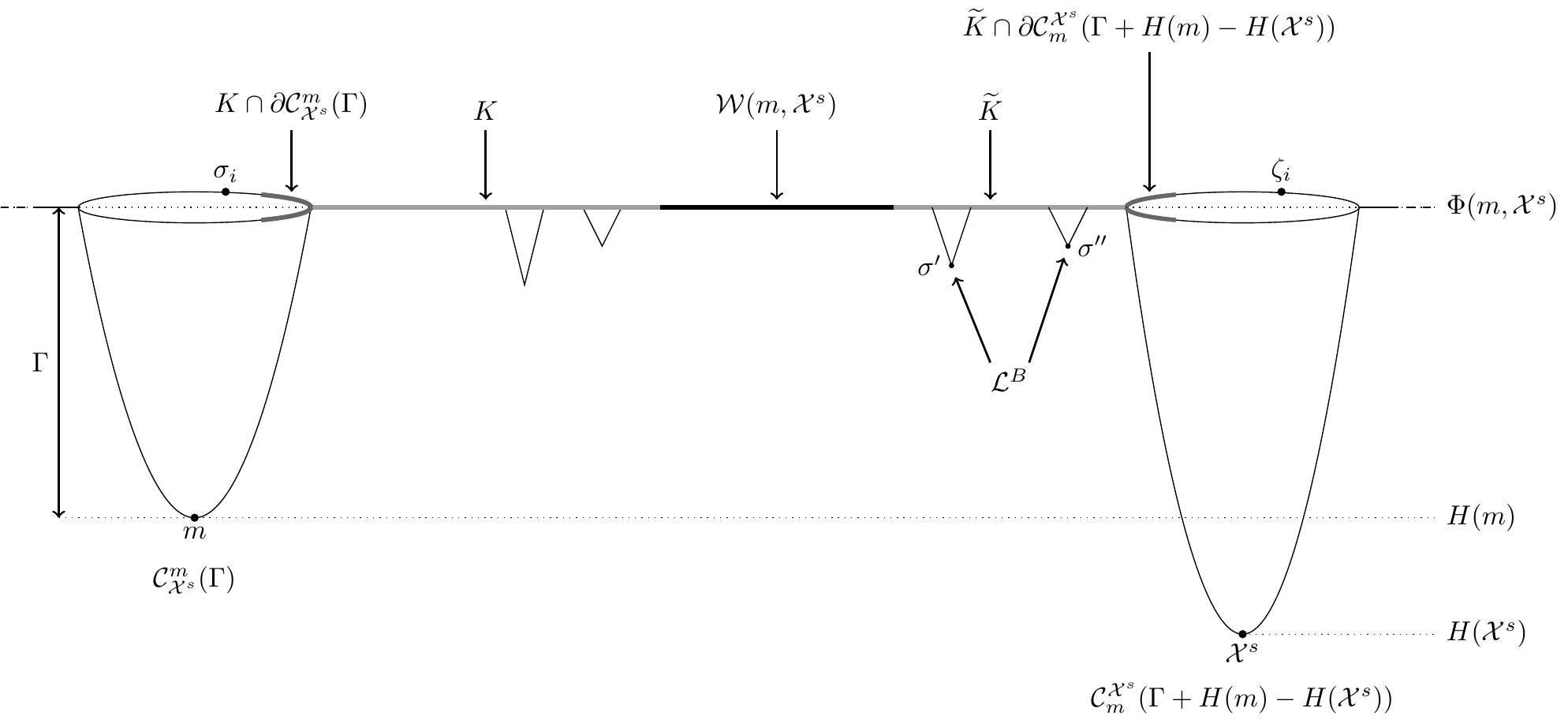}}
	\end{center}
	\vskip -0.75 cm
	\caption{We depict an example of the energy landscape for the transition between the metastable state $m$ and the stable states $\cX^s$. We depict on the left the cycle of the metastable state $\cC^m_{\cX^s}(\G)$ and on the right the cycle of the stable states $\cC^{\cX^s}_m(\G+H(m)-H(\cX^s))$. We indicate in black $W(m,\cX^s)$, in light grey $K$ and $\widetilde{K}$, emphasizing with dark grey the part of $K$ and $\widetilde{K}$ that intersect the boundaries of the two previous cycles. We give an example of two configurations $\s'$ and $\s''$ that are in $\cL^B$.}
	\label{fig:land}
\end{figure}

In order to apply this strategy to a concrete model, we require the following model-dependent inputs (we encourage the reader to inspect Figure \ref{fig:land}):

\bi
\item[(i)] Identify the sets $\cX^m$ and $\cX^s=\{\h_1^s,...,\h_k^s\}$, where $\h_1^s,...,\h_k^s$ have to be in $\cC^{\cX^s}_{m}(\G+H(m)-H(\cX^s))$. For a given $m\in\cX^m$, compute $\Phi(m,\cX^s)$ and set $\G:=\Phi(m,\cX^s)-H(m)$, the energy barrier between $m$ and $\cX^s$.
\item[(ii)] Find a set $\cW(m,\cX^s)$ and prove that it is a gate for the transition $m\ra\cX^s$.
\item[(iii)] Find two sets of configurations $\cL^G$ and $\cL^B$ and prove the following conditions for any $\h\in\cW(m,\cX^s)$:
\bi
\item[(a)] there exist a path $\o_1^G:\h\ra\cL^G$ such that $\max_{\s\in\o_1^G}H(\s)\leq\G+H(m)$ and a path $\o_2^G:\cL^G\ra\cX^s$ such that $\max_{\s\in\o_2^G}H(\s)<\G+H(m)$;
\item[(b)] there exists a path $\o_1^B:\h\ra\cL^B$ such that $\max_{\s\in\o_1^B}H(\s)\leq\G+H(m)$ and $\nexists \ \o_2^B:\cL^B\ra\cX^s$ and $\nexists \ \o_2^B:\cL^B\ra m$ such that $\max_{\s\in\o_2^B}H(\s)<\G+H(m)$;
\ei
\item[(iv)] Identify the subset $K$ (resp.\ $\widetilde{K}$) of the saddles that are visited by the optimal paths ``just before entering" (resp.\ ``just after visiting") $\cW(m,\cX^s)$. More precisely,
\be{defK}
\ba{ll}
K:=\{\bar\h\in\cS(m,\cX^s)\setminus\cW(m,\cX^s): \exists \ \h\in\cW(m,\cX^s) \text{ and } \o=\o_1 \circ \o_2,\\
\quad \qquad \text{ with } \o_1:\h\ra\bar\h  \text{ s.t. } \o_1\cap\cW(m,\cX^s)=\{\h\}, \ \o_1\cap\cC_{\cX^s}^{m}(\G)=\emptyset\\
\quad \qquad \text{ and } \o_2:\bar\h\ra m \text{ s.t. } \o_2\cap\cW(m,\cX^s)=\emptyset, \ \max_{\s\in\o}H(\s)\leq\G+H(m)\} 
\ea
\ee

\noindent
and
\be{defKtilde}
\ba{lll}
\widetilde{K}:=\{\bar\h\in\cS(m,\cX^s)\setminus\cW(m,\cX^s): \exists \ \h\in\cW(m,\cX^s) \text{ and } \o=\o_1 \circ \o_2, \text{ with}\\
\quad \qquad  \o_1:\h\ra\bar\h \text{ s.t. } \o_1\cap\cW(m,\cX^s)=\{\h\}, \o_1\cap\cC_{m}^{\cX^s}(\G+H(m)-H(\cX^s))=\emptyset\\
\qquad \quad  \text{and } \o_2:\bar\h\ra\cX^s \text{ s.t. } \o_2\cap\cW(m,\cX^s)=\emptyset, \ \max_{\s\in\o}H(\s)\leq\G+H(m)\}.
\ea
\ee
\ei

\noindent
If $\cX^s$ is a singleton, then it belongs to $\cC^{\cX^s}_{m}(\G+H(m)-H(\cX^s))$. Conditions (iii)-(a) and (iii)-(b) guarantee that when the dynamics reaches $\cL^G$ it has gone ``over the hill", while when it reaches $\cL^B$ the energy has to increase again to the level $\G+H(m)$ to visit $m$ or $\cX^s$. In particular, this implies that $\cL^G\subset\cC_{m}^{\cX^s}(\G+H(m)-H(\cX^s))$ and $\cL^B\nsubseteq\cC_{m}^{\cX^s}(\G+H(m)-H(\cX^s))$. We will show in Section \ref{S6.4} how the model-dependent inputs (iii)-(a) and (iii)-(b) apply to isotropic and weakly anisotropic models evolving under Kawasaki dynamics. For the strongly anisotropic case, see \cite[Section 5.4]{BN3}. In Section \ref{sharpasymptotics} we will refer to $\cC_{\cX^s}^{m}(\G)$ as $\cX^{meta}$ and to $\partial\cC_{m}^{\cX^s}(\G+H(m)-H(\cX^s))$ as $\cX^{stab}$.

\br{}
Note that it is possible that $K=\emptyset$ and/or $\widetilde{K}=\emptyset$, since the gate $\cW(m,\cX^s)$ could contain all the configurations with such properties. Indeed in \cite[Lemma 7.4(a)]{BGNneg2021} it is proved that this is the case for the $q$-state Potts model with negative external magnetic field evolving under Glauber dynamics. In particular, assuming $q=2$, the same result holds for the Ising model.
\er

\noindent
See Figure \ref{fig:land} for Propositions \ref{selle1} and \ref{selle2}.

\bd{sellesigma}
A saddle $\s$ is of the first type if it is not in $\cW(m,\cX^s)\cup K$ and belongs to the boundary of the cycle $\cC_{\cX^s}^{m}(\G)$, i.e., $\s\in\partial\cC_{\cX^s}^{m}(\G)\cap(\cS(m,\cX^s)\setminus(\cW(m,\cX^s)\cup K))$, where $\cC_{\cX^s}^{m}(\G)$ is defined in (\ref{ciclo}).
\ed

\bp{selle1}
Any saddle $\s$ of the first type is unessential and therefore it is not in $\cG(m,\cX^s)$.
\ep

\noindent
We refer to Section \ref{proofdep1} for the proof of Proposition \ref{selle1}. As we can see in the proof, it will be clear that this result is guaranteed only by the model-dependent inputs (i), (ii) and (iv).

\bd{sellezeta}
A saddle $\z$ is of the second type if it is not in $\cW(m,\cX^s)\cup\widetilde K$ and belongs to the boundary of the cycle $\cC_{m}^{\cX^s}(\G+H(m)-H(\cX^s))$, i.e., $\z\in\partial\cC_{m}^{\cX^s}(\G+H(m)-H(\cX^s))\cap(\cS(m,\cX^s)\setminus(\cW(m,\cX^s)\cup\widetilde{K}))$, where $\cC_{m}^{\cX^s}(\G+H(m)-H(\cX^s))$ is defined in (\ref{ciclo}).
\ed

\bp{selle2}
Any saddle $\z$ of the second type is unessential and therefore it is not in $\cG(m,\cX^s)$.
\ep

\noindent
We refer to Section \ref{proofdep2} for the proof of Proposition \ref{selle2}. For this result all of the four model-dependent inputs are necessary.

\br{}
This strategy can be applied also in the tunnelling scenario, i.e., the transition between two stable states, which corresponds to selecting the starting state $m\in\cX^s$. The model-dependent input (i) has to be modified by requiring that the configurations in $\cX^s\setminus\{m\}$ are in the same cycle that does not contain $m$, while the inputs (ii)-(iv) remain the same. Thus Propositions \ref{selle1} and \ref{selle2} still hold after replacing the set $\cX^s$ by $\cX^s\setminus\{m\}$. In this case, since $H(m)=H(\cX^s)$, note that the cycles $\cC_{\cX^s\setminus\{m\}}^m(\G)$ and $\cC^{\cX^s\setminus\{m\}}_m(\G+H(m)-H(\cX^s))$ have the same depth. The idea of this strategy can be applied also in the tunnelling scenario in which the configurations in $\cX^s\setminus\{m\}$ are not in the same cycle, but this requires an extension of this strategy. This occurs in the $q$-state Potts model with $q$ possible spins and zero external magnetic field \cite{BGN}, where the stable states are the configurations with all spins of the same type. In \cite[Theorem 3.4]{BGN} the authors study the gates relevant for the tunnelling between one stable state $m$ to the set of the other stable states $\cX^s\setminus\{m\}$. For the proof of this theorem they use \cite[Theorem 3.2]{BGN}, in which they identify all the unessential saddles for the transition between the selected stable state $m$ to one of the other stable states $s\in\cX^s\setminus\{m\}$ when the dynamics is restricted only to the optimal paths that do not visit the rest of the stable states $\cX^s\setminus\{m,s\}$. The proof of \cite[Theorem 3.2]{BGN} uses, in the specific model, the ideas presented in this Section and the symmetry of the energy landscape for $q$-state Potts model with zero external magnetic field.
\er

\noindent
These model-independent propositions will be applied to the isotropic and weakly anisotropic models evolving under Kawasaki dynamics (in Section \ref{S6.4}) to identify the set $\cG(m,\cX^s)$. For the application to the strongly anisotropic model, see \cite[Section 5.4]{BN3}.

\section{Main results: the gates for the local model}
\label{S4}
In this Section we state our main results. In Section \ref{iso} (resp.\ Section \ref{wani}) we obtain the geometrical characterization of the union of all minimal gates for the isotropic (resp.\ weakly anisotropic) case. In order to do this we need some model-dependent definitions for the Kawasaki dynamics (see Section \ref{moddepdef}) and some specific definitions for the weakly anisotropic case (see Section \ref{wani}). In Section \ref{sharpestimates} we derive sharp estimates for the asymptotic transition time in the weakly anisotropic case. Moreover, we derive the mixing time and spectral gap in the isotropic and weakly anisotropic cases. For the corresponding results obtained in the strongly anisotropic case, i.e., in the parameter regime $U_1>2U_2$, we refer to \cite[Section 4.2]{BN3} for results concerning the gates and union of minimal gates and to \cite[Section 4.3]{BN3} for results concerning the asymptotic transition time, mixing time and spectral gap.

\subsection{Geometric definitions for Kawasaki dynamics}
\label{moddepdef}
We give some model-dependent definitions and notations in order to state our main theorems.

\medskip
\noindent
{\bf 1. Free particles and clusters}

\bi

\item[$\bullet$]
For $x\in\L_0$, let $ \hbox{nn}(x):=\{ y\in \L_0\colon\;d(y,x)=1\}$ be the set of
nearest-neighbor sites of $x$ in $\L_0$, where $d$ in the entire paper denotes the lattice distance.
\item[$\bullet$]
A {\it free particle\/} in $\h\in\cX$ is a site $x$, with $\h(x)=1$, such that either $x\in\partial^-\L$, or $x\in\L_0$ and $\sum_{y\in nn(x)\cap\L_0}\eta(y)$ $=0$. We denote by $\h_{fp}$ the union of free particles in $\partial^-
\L$ and  free particles in $\L_0$. We denote by $n(\h)$ the number of free particles in $\h$.

We denote by $\h_{cl}$ the clusterized
part of the occupied sites of $\h$:
\be{hcl}\h_{cl} :=\{x\in\L_0: \ \h(x)=1\}\setminus\h_{fp}.
\ee

\item[$\bullet$]
We denote by $\h^{fp}$ the addition of a free particle anywhere in $\L$ to the configuration $\h$.

\item[$\bullet$]
Given a configuration $\h \in\cX$, consider the subset $C(\h_{cl})$ of $\R^2$ defined as the union of the $1\times 1$ closed squares centered at the occupied sites of $ \h_{cl}$ in $ \L_0$ and call the maximal connected components of this set the clusters of $\h_{cl}$. 

\item[$\bullet$]
Given a set $A\subset\R^2$, we define the number of $1\times1$ closed occupied squares in $A$ as 
\be{defsbarra}
|A|:=|A\cap C(\h_{cl})|
\ee

\noindent
and as $||A||$ the numbers of $1\times1$ closed squares in $A$. Note that $||\cdot||$ takes into account the possibility that the squares are occupied or not.

\ei

\medskip
\noindent
{\bf 2. Projections, semi-perimeter and vacancies}

\bi

\item[$\bullet$]
For $\h\in\cX$, we denote by $g_1(\eta)$ (resp.\ $g_2(\eta)$) one half of the horizontal (resp.\ vertical) length of the Euclidean boundary of $C(\h_{cl})$. Then the energy associated with $\h$ is given by
\be{Hcont}
H(\h) = - (U_1+U_2-\D)|C(\eta_{cl})| + {U_1}
g_2(\eta)+
{U_2} g_1(\eta)
+\D n(\h).
\ee

\item[$\bullet$]
Let $p_1(\h)$
and $p_2(\h)$ be the total lengths of  horizontal and vertical
projections of $ C(\h_{cl})$ respectively. More precisely, let
$r_{j,1}=\{x \in \Z^2:(x)_1=j\}$ be the $j$-th column and
$r_{j,2}=\{x \in \Z^2:(x)_2=j\}$ be the $j$-th row, where
$(x)_1$ or $(x)_2$ denote the first or second component of $x$. Let

\be{proie1} \p_1(\h):=\{j \in \Z:\, r_{j,1}\cap
C(\h_{cl})\not=\emptyset\} \ee

\noindent and $p_1(\h):=|\p_1(\h)|$. In a similar way we define
the vertical projection $\p_2(\h)$ and $p_2(\h)$.

\item[$\bullet$]
We define $g'_i(\h):= g_i(\h)- p_i(\h)\ge 0$; we call {\it monotone} a configuration such that $g_i(\h)= p_i(\h)$ for $i=1,2$.

\item[$\bullet$]
We define the {\it semi-perimeter} $s(\h)$ and the {\it vacancies} $v(\h)$ as
\be{defdep1}
\ba{lll}
s(\h)&:=& p_1(\h)+p_2(\h),\\
v(\h)&:=& p_1(\h)p_2(\h)- |C(\h_{cl})|.
\ea
\ee

\ei

\medskip
\noindent
{\bf 3. $\boldsymbol{n}$-manifold, rectangles and corners}

\bi 
\item[$\bullet$]
The configuration space $\cX$ can be partitioned as
\be{}
\cX=\displaystyle\bigcup_{n}\cV_n,
\ee
\noindent
where $\cV_n:=\{\h\in\cX: |C(\h_{cl})|+n(\h)=n\}$ is the set of configurations with $n$ particles, called the {\it $n$-manifold}.

\item[$\bullet$]
We denote by $\cR(l_1,l_2)$ the set of configurations that have no free particle and a single cluster such that $C(\h_{cl})$ is a rectangle $R(l_1,l_2)$, with $l_1,l_2\in\N$. For any $\h,\h'\in\cR (l_1,l_2)$ we have immediately:
\be{enerettan} H(\h)=H(\h')=H(\cR(l_1,l_2))={U_1}l_2 +{U_2} l_1
-\varepsilon l_1 l_2, 
\ee
\noindent
where
\be{defepsilon}
\e:= U_1+U_2 -\D.
\ee

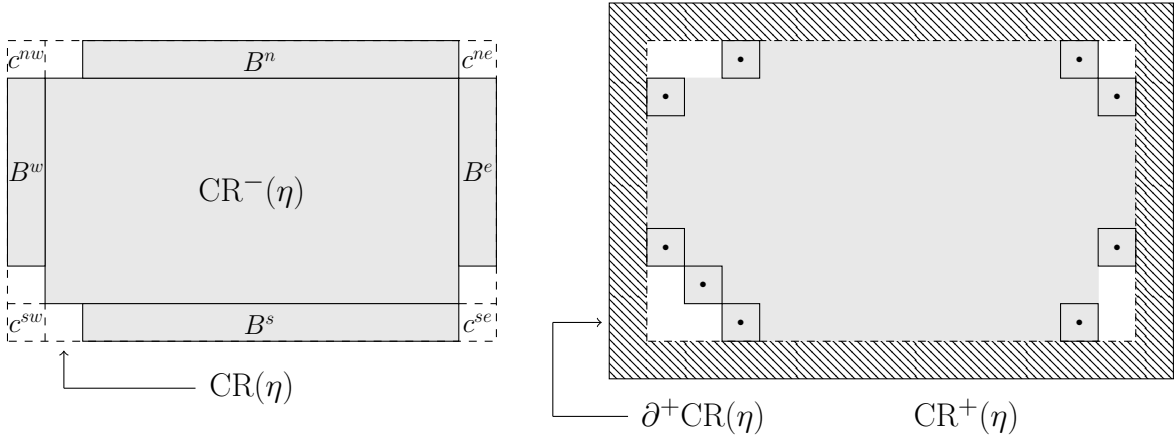
\begin{figure}
	\centering
	\begin{tikzpicture}[scale=0.5,transform shape]

		\draw [fill=grigio] (6,3) rectangle (17,9);
		\draw [fill=grigio] (17,4) rectangle (18,9);
		\node at (17.5,6.5){\huge{$B^e$}};
		\node at (11.5,6){\Huge{$\hbox{CR}^-(\h)$}};
		\draw [fill=grigio] (7,9) rectangle (17,10);
		\node at (11.7,9.45){\huge{$B^n$}};
		\draw [fill=grigio] (5,4) rectangle (6,9);
		\node at (5.5,6.5){\huge{$B^w$}};
		\draw [fill=grigio] (7,3) rectangle (17,2);
		\node at (11.7,2.45){\huge{$B^s$}};
		\draw [dashed] (5,2) rectangle (18,10);
		\draw [dashed] (5,3)--(6,3);
		\draw [dashed] (6,2)--(6,3);
		\node at (5.5,2.5){\huge{$c^{sw}$}};
		\draw [dashed] (17,3)--(18,3);
		\node at (17.5,2.5){\huge{$c^{se}$}};
		\node at (17.5,9.5){\huge{$c^{ne}$}};
		\draw [dashed] (6,9)--(6,10);
		\node at (5.5,9.5){\huge{$c^{nw}$}};
		\node at (11.5,0.75){\Huge{$\hbox{CR}(\h)$}};
		\draw (10,0.75) -- (6.5,0.75);
		\draw[->](6.5,0.75)--(6.5,1.75);

		\draw [color=grigio,fill=grigio] (23,3) rectangle (34,9);
		\draw [color=grigio,fill=grigio] (34,4) rectangle (35,9);
		\draw [color=grigio,fill=grigio] (24,9) rectangle (34,10);
		\draw [color=grigio,fill=grigio] (22,4) rectangle (23,9);
		\draw [color=grigio,fill=grigio] (24,3) rectangle (34,2);
		\draw [dashed] (22,2) rectangle (35,10);
		\draw (21,1) rectangle (36,11);
		\draw (22,4) rectangle (23,5);
		\node at (22.5,4.5) {$\bullet$};
		\draw (23,3) rectangle (24,4);
		\node at (23.5,3.5) {$\bullet$};
		\draw (24,2) rectangle (25,3);
		\node at (24.5,2.5) {$\bullet$};
		\draw (33,2) rectangle (34,3);
		\node at (33.5,2.5) {$\bullet$};
		\draw (34,4) rectangle (35,5);
		\node at (34.5,4.5) {$\bullet$};
		\draw (34,8) rectangle (35,9);
		\node at (34.5,8.5) {$\bullet$};
		\draw (33,9) rectangle (34,10);
		\node at (33.5,9.5) {$\bullet$};
		\draw (24,9) rectangle (25,10);
		\node at (24.5,9.5) {$\bullet$};
		\draw (22,8) rectangle (23,9);
		\node at (22.5,8.5) {$\bullet$};
		\fill[pattern=north west lines] (21,1) rectangle (22,11);
		\fill[pattern=north west lines] (22,1) rectangle (36,2);
		\fill[pattern=north west lines] (35,2) rectangle (36,11);
		\fill[pattern= north west lines] (22,10) rectangle (35,11);
		\node at (30.5,0){\Huge{$\hbox{CR}^+(\h)$}};
		\node at (23.5,0){\Huge{$\partial^+\hbox{CR}(\h)$}};
		\draw (21.5,0) -- (19.5,0);
		\draw (19.5,0)--(19.5,2.5);
		\draw[->] (19.5,2.5)--(20.8,2.5);

	\end{tikzpicture}
	
	\vskip 0 cm
	\caption{Here we depict the same configuration $\h$ on the left and on the right to emphasize different geometrical definitions. The grey area in both pictures represents $C(\h_{cl})$. In particular, on the left-hand side we stress the frame-angles $c^{\a\a'}(\h )$, the bars $B^\a(\h)$, $\hbox{CR}^-(\h)$ and the circumscribing rectangle $\hbox{CR}(\h)$ (respresented with a dashed line). While on the right-hand side we stress the sites that are in a corner (represented with a dot), $\hbox{CR}^+(\h)$ and the external frame $\partial^+\hbox{CR}(\h)$ (the dashed area).}
	\label{fig:figesempio}
\end{figure}

\item[$\bullet$]
A {\it corner} in $\h\in\cX$ is a closed $1\times1$ square centered in an occupied site $x\in\L_0$ such that, if we order clockwise its four nearest neighbors $x_1,x_2,x_3,x_4$, then $\sum_{y\in \hbox{nn}(x)}\h(y)=2$, with $\h(x_i)=\h(x_{i+1})=1$, with $i=1,...,4$ and the convention that $x_5=x_1$ (see Figure \ref{fig:figesempio} on the right-hand side).
\ei

\medskip
\noindent
{\bf 4. Circumscribed rectangle, frames and bars}
\bi

\item[$\bullet$]
If $\h$ is a configuration with a single cluster then we denote by CR$(\h)$ the rectangle {\it circumscribing} $C(\h_{cl})$, i.e., the smallest rectangle containing $\h$. 

We denote $\partial^+\hbox{CR}(\h)$ the {\it external frame of} $\hbox{CR}(\h)$ as the union of squares $1\times1$ centered at sites that are not contained in $\hbox{CR}(\h)$ such that those sites have Euclidean distance with sites in $\hbox{CR}(\h)$ less or equal than $\sqrt{2}$ (see Figure \ref{fig:figesempio} on the right-hand side). Note that the external frame of $\hbox{CR}(\h)$ contains only non occupied sites.

We denote $\partial^-\hbox{CR}(\h)$ the {\it internal frame of} $\hbox{CR}(\h)$ as the union of squares $1\times1$ centered at sites that are contained in $\hbox{CR}(\h)$ such that those sites have Euclidean distance with sites not in $\hbox{CR}(\h)$ less or equal than $\sqrt{2}$. If this distance is equal to $\sqrt{2}$, we say that the unit square is a {\it frame-angle} $c^{\alpha\alpha'}(\h)$ in $\partial^-\hbox{CR}(\h)$, where $\alpha\alpha'\in\{ne,nw,se,sw\}$, with $n=\hbox{north}$, $s=\hbox{south}$, etc. Note that the internal frame of $\hbox{CR}(\h)$ is a geometrical object contained in $\R^2$ that can contain both occupied and non occupied sites (see Figure \ref{fig:figesempio} on the left-hand side). We partition the set $\partial^-\hbox{CR}(\h)$ without frame-angles in {\it two horizontal and two vertical rows} $r^{\alpha}(\h)$, with $\alpha\in\{n,w,e,s\}$.

Moreover, we set
\be{cr}
\ba{ll}
\hbox{CR}^-(\h)=\hbox{CR}(\h)\setminus\partial^-\hbox{CR}(\h),\\
\hbox{CR}^+(\h)=\hbox{CR}(\h)\cup\partial^+\hbox{CR}(\h).
\ea
\ee

See Figure \ref{fig:figesempio} for an example.

\br{}
Note that, for example, the frame-angles $c^{ne}(\h)$ and $c^{en}(\h)$ are the same, but this distinction will be useful in Definitions \ref{movepart} and \ref{trenino}.
\er

\item[$\bullet$]
A {\it vertical} (respectively {\it horizontal}) {\it bar} $B^{\a}(\h)$ of a single cluster of $\h$ with length $k$ is a $1\times k$ (respectively $k\times1$) rectangle contained in $C(\h_{cl})$, with $\alpha\in\{n,w,e,s\}, \ k\geq1$, such that each square $1\times1$ of the bar is attached only to one square of $C(\h_{cl})\setminus B^\a(\h)$ (see Figure \ref{fig:figesempio} on the left-hand side). In the cases in which it is not specified if the bar is vertical or horizontal we call it simply {\it bar}. If $k=1$ we say that the bar is a {\it protuberance}.

\br{sommabarre}
Note that two bars $B^{\a}(\h)$ and $B^{\a'}(\h)$, with $\a,\a'\in\{n,s,w,e\}$, can possibly intersect in the frame-angle $c^{\a\a'}(\h)$. If this is the case, we get $|B^{\a}(\h)\cup B^{\a'}(\h)|=|B^{\a}(\h)|+|B^{\a'}(\h)|-1$.
\er

\ei

\medskip
\noindent
{\bf 5. Motions along the border}

\noindent
Recall definitions of $|\cdot|$ and $||\cdot||$ in (\ref{defsbarra}) and below. In the following, we give the precise notion of translation by 1 of a bar, for example to the left or to the right, while keeping all the squares of the bar attached to the cluster below.

\bd{translation}
Given $\h$ and a bar $B^\a(\h)$ of length $k$, with $\a\in\{n,s,e,w\}$, we say that it is possible to \emph{translate the bar $B^\a(\h)$} if 
\be{condtrasl}
k=|B^\a(\h)|<|\partial^+ B^\a(\h)|.
\ee

\noindent
We define the \emph{$1$-translation of a bar $B^\a(\h)$} of length $k$ as a sequence of configurations $(\h_1,...,\h_k)$ such that $\h_1=\h$ and $\h_i$ is obtained from $\h_{i-1}$  translating by $1$ a unit square along the rectangle $\partial^+ B^\a(\h)\cap C(\h_{cl})$ for any $2\leq i\leq k$.
\ed

\setlength{\unitlength}{1.1pt}
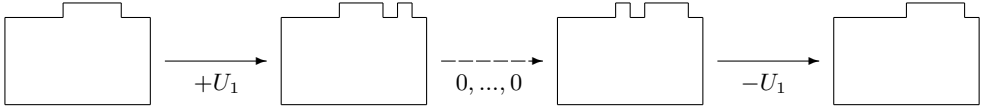
\begin{figure}
	\centering
	\begin{picture}(400,40)(0,40)
		\thinlines
		\put(20,50){\line(1,0){50}}
		\put(20,80){\line(1,0){20}}
		\put(40,80){\line(0,1){5}}
		\put(40,85){\line(1,0){20}}
		\put(60,85){\line(0,-1){5}}
		\put(60,80){\line(1,0){10}}
		\put(20,50){\line(0,1){30}}
		\put(70,50){\line(0,1){30}}
		\put(75,65){\vector(1,0){35}}
		\begin{footnotesize} \put(85,55){$+U_1$} \end{footnotesize}
		\thinlines
		\put(115,50){\line(1,0){50}}
		\put(115,80){\line(1,0){20}}
		\put(135,80){\line(0,1){5}}
		\put(135,85){\line(1,0){15}}
		\put(150,85){\line(0,-1){5}}
		\put(150,80){\line(1,0){5}}
		\put(155,80){\line(0,1){5}}
		\put(155,85){\line(1,0){5}}
		\put(160,85){\line(0,-1){5}}
		\put(160,80){\line(1,0){5}}
		\put(115,50){\line(0,1){30}}
		\put(165,50){\line(0,1){30}}
		\put(170,65){\line(1,0){5}}
		\put(176,65){\line(1,0){5}}
		\put(182,65){\line(1,0){5}}
		\put(188,65){\line(1,0){5}}
		\put(194,65){\line(1,0){5}}
		\put(200,65){\vector(1,0){5}}
		
		\begin{footnotesize}
			\put(175,55){$0,...,0$}
		\end{footnotesize}
		
		\thinlines
		\put(210,50){\line(1,0){50}}
		\put(210,80){\line(1,0){20}}
		\put(230,80){\line(0,1){5}}
		\put(230,85){\line(1,0){5}}
		\put(235,85){\line(0,-1){5}}
		\put(235,80){\line(1,0){5}}
		\put(240,80){\line(0,1){5}}
		\put(240,85){\line(1,0){15}}
		\put(255,85){\line(0,-1){5}}
		\put(255,80){\line(1,0){5}}
		\put(210,50){\line(0,1){30}}
		\put(260,50){\line(0,1){30}}
		\put(265,65){\vector(1,0){35}}
		\begin{footnotesize}
			\put(273,55){$-U_1$}
		\end{footnotesize}
		
		\thinlines
		\put(305,50){\line(1,0){50}}
		\put(305,80){\line(1,0){25}}
		\put(330,80){\line(0,1){5}}
		\put(330,85){\line(1,0){20}}
		\put(350,85){\line(0,-1){5}}
		\put(350,80){\line(1,0){5}}
		\put(305,50){\line(0,1){30}}
		\put(355,50){\line(0,1){30}}
	\end{picture}
	\vskip -0.5 cm
	\caption{$1$-translation of the horizontal bar $B^n(\h)$ at cost $U_1$.}
	\label{fig:traslazione1}
\end{figure}

\setlength{\unitlength}{1.1pt}
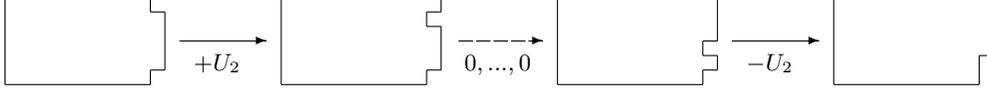
\begin{figure}
	\centering
	\begin{picture}(400,40)(0,40)
		\thinlines
		\put(20,50){\line(1,0){50}}
		\put(20,80){\line(1,0){50}}
		\put(70,75){\line(1,0){5}}
		\put(75,75){\line(0,-1){20}}
		\put(75,55){\line(-1,0){5}}
		\put(60,80){\line(1,0){10}}
		\put(20,50){\line(0,1){30}}
		\put(70,50){\line(0,1){5}}
		\put(70,75){\line(0,1){5}}
		\put(80,65){\vector(1,0){30}}
		\begin{footnotesize} \put(85,55){$+U_2$} \end{footnotesize}
		\thinlines
		\put(115,50){\line(1,0){50}}
		\put(115,80){\line(1,0){55}}
		\put(170,80){\line(0,-1){5}}
		\put(170,75){\line(-1,0){5}}
		\put(165,75){\line(0,-1){5}}
		\put(165,70){\line(1,0){5}}
		\put(170,70){\line(0,-1){15}}
		\put(170,55){\line(-1,0){5}}
		\put(160,80){\line(1,0){5}}
		\put(115,50){\line(0,1){30}}
		\put(165,50){\line(0,1){5}}
		\put(176,65){\line(1,0){5}}
		\put(182,65){\line(1,0){5}}
		\put(188,65){\line(1,0){5}}
		\put(194,65){\line(1,0){5}}
		\put(200,65){\vector(1,0){5}}
		
		\begin{footnotesize}
			\put(178,55){$0,...,0$}
		\end{footnotesize}
		
		\thinlines
		\put(210,50){\line(1,0){50}}
		\put(210,80){\line(1,0){55}}
		\put(265,80){\line(0,-1){15}}
		\put(265,65){\line(-1,0){5}}
		\put(260,65){\line(0,-1){5}}
		\put(260,60){\line(1,0){5}}
		\put(265,60){\line(0,-1){5}}
		\put(265,55){\line(-1,0){5}}
		\put(255,80){\line(1,0){5}}
		\put(210,50){\line(0,1){30}}
		\put(260,50){\line(0,1){5}}
		\put(270,65){\vector(1,0){30}}
		\begin{footnotesize}
			\put(275,55){$-U_2$}
		\end{footnotesize}
		
		\thinlines
		\put(305,50){\line(1,0){50}}
		\put(305,80){\line(1,0){55}}
		\put(360,80){\line(0,-1){20}}
		\put(355,60){\line(1,0){5}}
		\put(305,50){\line(0,1){30}}
		\put(355,50){\line(0,1){10}}
	\end{picture}
	\vskip -0.5 cm
	\caption{$1$-translation of the vertical bar $B^e(\h)$ at cost $U_2$.}
	\label{fig:traslazione2}
\end{figure}

\noindent
In Figure \ref{fig:traslazione1} (resp.\ Figure \ref{fig:traslazione2})we depict a $1$-translation of a horizontal (resp.\ vertical) bar at cost $U_1$ (resp.\ $U_2$).

In the following, we give the precise notion of sliding a unit square from row $r^{\a}(\h)$ to $r^{\a'}(\h)$ passing through the frame angle $c^{\a\a'}(\h)$.

\setlength{\unitlength}{0.92pt}
\begin{figure}
	\begin{picture}(400,40)(0,40)
		\thinlines
		\put(20,50){\line(1,0){40}}
		\put(20,80){\line(1,0){25}}
		\put(45,80){\line(0,1){5}}
		\put(45,85){\line(1,0){15}}
		\put(60,85){\line(0,-1){5}}
		\put(60,80){\line(1,0){5}}
		\put(65,80){\line(0,-1){20}}
		\put(20,50){\line(0,1){30}}
		\put(60,50){\line(0,1){10}}
		\put(60,60){\line(1,0){5}}
		\put(70,65){\vector(1,0){15}}
		\begin{scriptsize} \put(68,55){$+U_1$} \end{scriptsize}
		\thinlines
		\put(90,50){\line(1,0){40}}
		\put(90,80){\line(1,0){25}}
		\put(115,80){\line(0,1){5}}
		\put(115,85){\line(1,0){10}}
		\put(125,85){\line(0,-1){5}}
		\put(125,80){\line(1,0){5}}
		\put(130,80){\line(0,1){5}}
		\put(130,85){\line(1,0){5}}
		\put(135,85){\line(0,-1){25}}
		\put(135,60){\line(-1,0){5}}
		\put(90,50){\line(0,1){30}}
		\put(130,50){\line(0,1){10}}
		\put(140,65){\vector(1,0){15}}
		
		\begin{scriptsize}
			\put(145,55){$0$}
		\end{scriptsize}
		
		\thinlines
		\put(160,50){\line(1,0){40}}
		\put(160,80){\line(1,0){25}}
		\put(185,80){\line(0,1){5}}
		\put(185,85){\line(1,0){5}}
		\put(190,85){\line(0,-1){5}}
		\put(190,80){\line(1,0){5}}
		\put(195,80){\line(0,1){5}}
		\put(195,85){\line(1,0){10}}
		\put(205,85){\line(0,-1){25}}
		\put(205,60){\line(-1,0){5}}
		\put(160,50){\line(0,1){30}}
		\put(200,50){\line(0,1){10}}
		\put(210,65){\vector(1,0){15}}
		\begin{scriptsize}
			\put(208,55){$-U_1$}
		\end{scriptsize}
		
		\thinlines
		\put(230,50){\line(1,0){40}}
		\put(230,80){\line(1,0){30}}
		\put(260,80){\line(0,1){5}}
		\put(260,85){\line(1,0){5}}
		\put(265,85){\line(1,0){10}}
		\put(275,85){\line(0,-1){25}}
		\put(275,60){\line(-1,0){5}}
		\put(230,50){\line(0,1){30}}
		\put(270,50){\line(0,1){10}}
		\put(280,65){\vector(1,0){15}}
		\begin{scriptsize}
			\put(278,55){$+U_2$}
		\end{scriptsize}

		\thinlines
		\put(300,50){\line(1,0){40}}
		\put(300,80){\line(1,0){30}}
		\put(330,80){\line(0,1){5}}
		\put(330,85){\line(1,0){5}}
		\put(335,85){\line(1,0){10}}
		\put(345,85){\line(0,-1){20}}
		\put(345,65){\line(-1,0){5}}
		\put(340,65){\line(0,-1){5}}
		\put(340,60){\line(1,0){5}}
		\put(345,60){\line(0,-1){5}}
		\put(345,55){\line(-1,0){5}}
		\put(300,50){\line(0,1){30}}
		\put(340,50){\line(0,1){5}}
		\put(350,65){\line(1,0){2}}
		\put(353,65){\line(1,0){2}}
		\put(356,65){\line(1,0){2}}
		\put(359,65){\line(1,0){2}}
		\put(362,65){\vector(1,0){3}}
		\begin{scriptsize}
			\put(355,55){$0$}
		\end{scriptsize}
		
		\thinlines
		\put(370,50){\line(1,0){40}}
		\put(370,80){\line(1,0){30}}
		\put(400,80){\line(0,1){5}}
		\put(400,85){\line(1,0){5}}
		\put(405,85){\line(1,0){10}}
		\put(415,85){\line(0,-1){5}}
		\put(415,80){\line(-1,0){5}}
		\put(410,80){\line(0,-1){5}}
		\put(410,75){\line(1,0){5}}
		\put(415,75){\line(0,-1){20}}
		\put(415,55){\line(-1,0){5}}
		\put(370,50){\line(0,1){30}}
		\put(410,50){\line(0,1){5}}
		\put(420,65){\vector(1,0){15}}
		\begin{scriptsize}
			\put(418,55){$-U_2$}
		\end{scriptsize}
		
		\put(440,50){\line(1,0){40}}
		\put(440,80){\line(1,0){30}}
		\put(470,80){\line(0,1){5}}
		\put(470,85){\line(1,0){5}}
		\put(475,85){\line(1,0){5}}
		\put(480,85){\line(0,-1){5}}
		\put(480,80){\line(1,0){5}}
		\put(485,80){\line(0,-1){5}}
		\put(485,75){\line(0,-1){20}}
		\put(485,55){\line(-1,0){5}}
		\put(440,50){\line(0,1){30}}
		\put(480,50){\line(0,1){5}}

	\end{picture}
	\vskip -0.5 cm
	\caption{Sliding of a unit square around the frame angle $c^{ne}(\h)$ at cost $U_1$. In this case $\a=n$, $\a'=e$, $\a''=w$ and $\a'''=s$.}
	\label{fig:trenino}
\end{figure}
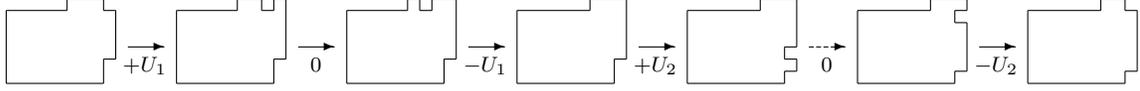

\bd{movepart}
Given $\h$, let $\alpha\alpha'$ such that $c^{\a\a'}(\h)$ is a frame-angle. We say that it is possible to {\rm slide a unit square around a frame-angle $c^{\alpha\alpha'}(\h)\subseteq\partial^-\hbox{CR}(\h)$} from a row $r^{\alpha}(\h)\subseteq\partial^-\hbox{CR}(\h)$ to a row $r^{\alpha'}(\h)\subseteq\partial^-\hbox{CR}(\h)$ via a frame-angle $c^{\alpha\alpha'}(\h)$ if
\be{condcorner}
|c^{\alpha\alpha'}(\h)|=0, \quad |r^{\alpha}(\h)|\geq1, \quad 1\leq|r^{\alpha'}(\h)|<||r^{\alpha'}(\h)||+1.
\ee

\noindent
Let $\a''\neq\a'$ such that $c^{\a\a''}(\h)$ is a frame-angle. See Figure \ref{fig:trenino} for an example. We define a {\rm sliding of a unit square around a frame-angle $c^{\alpha\alpha'}(\h)\subseteq\partial^-\hbox{CR}(\h)$} as the composition of a sequence of $1$-translations of the bar $B^\a(\h)$ from $r^{\alpha}(\h)\cup c^{\a\a''}(\h)$ to $r^{\alpha}(\h)\cup c^{\alpha\alpha'}(\h)$, namely $(\h^1,...,\h^k)$, and the $1$-translation of a bar $B^{\a'}(\h)=C(\h_{cl}^k)\cap (r^{\alpha'}(\h)\cup c^{\alpha\alpha'}(\h))$ from $r^{\alpha'}(\h)\cup c^{\alpha\alpha'}(\h)$ to $r^{\alpha'}(\h)\cup c^{\a'\a'''}(\h)$, where $\a'''\neq\a$ is such that $c^{\a'\a'''}(\h)$ is a frame-angle.
\ed

The definition above is used only to define the following sliding of a bar from row $r^{\a}(\h)$ to $r^{\a'}(\h)$ passing through the frame angle $c^{\a\a'}(\h)$, that corresponds to iteratively apply the sliding of a unit square around a frame-angle.

\bd{trenino}
Given $\h$, let $\alpha\alpha'$ such that $c^{\a\a'}(\h)$ is a frame-angle. Before sliding a bar around a frame-angle, we translate the bars $B^{\a}(\h)$ and $B^{\a'}(\h)$ at distance 1 to the frame-angle $c^{\a\a'}(\h)$ obtaining a configuration $\h'$. We say that it is possible to \emph{slide a bar $B^\a(\h')$ around a frame-angle $c^{\alpha\alpha'}(\h')\subseteq\partial^-\hbox{CR}(\h')$} if it is possible to move all the unit squares in $B^\a(\h')$ around a frame-angle $c^{\alpha\alpha'}(\h')$ from a row $r^{\alpha}(\h')\cup c^{\alpha\alpha''}(\h')$ to a row $r^{\alpha'}(\h')\cup c^{\a'\a'''}(\h')$, where $\a''\neq\a'$ and $\a'''\neq\a$ are such that $c^{\a\a''}(\h')$ and $c^{\a'\a'''}(\h')$ are frame-angles. Namely,
\be{condtrenino}
|B^\a(\h')|+|r^{\alpha'}(\h')|\leq||r^{\alpha'}(\h')||+1.
\ee

\noindent
Moreover, we define a \emph{sliding of a bar $B^\a(\h')$ around a frame-angle $c^{\alpha\alpha'}(\h')$} as the sequence of $|B^\a(\h')|$ slidings of unit squares around a frame-angle $c^{\alpha\alpha'}(\h')$.
\ed

\noindent
See the path described in Figure \ref{fig:columntorow}, that connects the configuration $\h$ to the configuration $(12)$ for an example of sliding of the bar $B^e(\h)$ around the frame-angle $c^{en}(\h)$, with $\h$ as the configuration $(3)$.

\setlength{\unitlength}{1.1pt}
\begin{figure}
	\begin{picture}(400,50)(0,30)
		\thinlines
		\put(20,50){\line(1,0){50}}
		\put(20,80){\line(1,0){50}}
		\put(20,50){\line(0,1){30}}
		\put(70,50){\line(0,1){30}}
		\thinlines  
		\put(40,60){$(1)$}
		\begin{Huge}
			\put(75,60){$\longrightarrow$} \end{Huge}
		\begin{footnotesize} \put(85,55){$\D$} \end{footnotesize}
		\thinlines
		\put(115,50){\line(1,0){50}}
		\put(115,80){\line(1,0){50}}
		\put(115,50){\line(0,1){30}}
		\put(165,50){\line(0,1){30}}
		\put(160,85){\line(0,1){5}}
		\put(155,85){\line(1,0){5}}
		\put(155,85){\line(0,1){5}}
		\put(155,90){\line(1,0){5}}
		\thinlines 
		\put(135,60){$(2)$}
		\begin{Huge}
			\put(170,60){$\longrightarrow$}
		\end{Huge}
		\begin{footnotesize}
			\put(180,55){$-U_2$}
		\end{footnotesize}
		
		\thinlines
		\put(210,50){\line(1,0){50}}
		\put(260,50){\line(0,1){30}}
		\put(210,50){\line(0,1){30}}
		\put(210,80){\line(1,0){40}}
		\put(255,80){\line(0,1){5}}
		\put(255,80){\line(1,0){5}}
		\put(250,80){\line(0,1){5}}
		\put(250,85){\line(1,0){5}}
		\thinlines 
		\put(230,60){$(3)$}
		\begin{Huge}
			\put(265,60){$\longrightarrow$}
		\end{Huge}
		\begin{footnotesize}
			\put(275,55){$+U_2$}
		\end{footnotesize}
		
		\thinlines
		\put(305,50){\line(1,0){50}}
		\put(305,50){\line(0,1){30}}
		\put(305,80){\line(1,0){40}}
		\put(345,80){\line(0,1){5}}
		\put(345,85){\line(1,0){10}}
		\put(355,50){\line(0,1){25}}
		\put(350,75){\line(1,0){5}}
		\put(350,75){\line(0,1){5}}
		\put(350,80){\line(1,0){5}}
		\put(355,80){\line(0,1){5}}
		\thinlines 
		\put(325,60){$(4)$}
		\begin{Huge}
			\put(360,60){$\dashrightarrow$}
		\end{Huge}
		\begin{footnotesize}
			\put(365,55){$0,..,0$}
		\end{footnotesize}

		\thinlines
		\put(20,0){\line(1,0){50}}
		\put(20,0){\line(0,1){30}}
		\put(70,0){\line(0,1){5}}
		\put(20,30){\line(1,0){40}}
		\put(60,30){\line(0,1){5}}
		\put(60,35){\line(1,0){10}}
		\put(65,5){\line(1,0){5}}
		\put(65,5){\line(0,1){5}}
		\put(65,10){\line(1,0){5}}
		\put(70,10){\line(0,1){25}}
		\thinlines 
		\put(37,13){$(5)$}
		\begin{Huge}
			\put(75,10){$\longrightarrow$}
		\end{Huge}
		\begin{footnotesize}
			\put(83,5){$-U_2$}
		\end{footnotesize}

		\thinlines
		\put(115,0){\line(1,0){45}}
		\put(115,0){\line(0,1){30}}
		\put(160,0){\line(0,1){5}}
		\put(115,30){\line(1,0){40}}
		\put(160,5){\line(1,0){5}}
		\put(165,5){\line(0,1){30}}
		\put(155,35){\line(1,0){10}}
		\put(155,30){\line(0,1){5}}
		\thinlines 
		\put(135,13){$(6)$}
		\begin{Huge}
			\put(170,10){$\longrightarrow$}
		\end{Huge}
		\begin{footnotesize}
			\put(180,5){$+U_1$}
		\end{footnotesize}
		
		\thinlines
		\put(210,0){\line(1,0){45}}
		\put(210,0){\line(0,1){30}}
		\put(210,30){\line(1,0){35}}
		\put(245,30){\line(0,1){5}}
		\put(245,35){\line(1,0){5}}
		\put(250,30){\line(0,1){5}}
		\put(250,30){\line(1,0){5}}
		\put(255,30){\line(0,1){5}}
		\put(255,35){\line(1,0){5}}
		\put(255,0){\line(0,1){5}}
		\put(255,5){\line(1,0){5}}
		\put(260,5){\line(0,1){30}}
		\thinlines 
		\put(230,13){$(7)$}
		\begin{Huge}
			\put(265,10){$\longrightarrow$}
		\end{Huge}
		\begin{footnotesize}
			\put(275,5){$-U_1$}
		\end{footnotesize}

		\thinlines
		\put(305,0){\line(1,0){45}}
		\put(305,0){\line(0,1){30}}
		\put(305,30){\line(1,0){35}}
		\put(340,30){\line(0,1){5}}
		\put(340,35){\line(1,0){10}}
		\put(350,30){\line(0,1){5}}
		\put(350,30){\line(1,0){5}}
		\put(350,0){\line(0,1){5}}
		\put(350,5){\line(1,0){5}}
		\put(355,5){\line(0,1){25}}
		\thinlines 
		\put(325,13){$(8)$}
		\begin{Huge}
			\put(360,10){$\longrightarrow$}
		\end{Huge}
		\begin{footnotesize}
			\put(370,5){$+U_2$}
		\end{footnotesize}

		\thinlines
		\put(20,-50){\line(1,0){45}}
		\put(20,-50){\line(0,1){30}}
		\put(20,-20){\line(1,0){35}}
		\put(55,-20){\line(0,1){5}}
		\put(55,-15){\line(1,0){15}}
		\put(65,-50){\line(0,1){5}}
		\put(65,-45){\line(1,0){5}}
		\put(70,-45){\line(0,1){20}}
		\put(65,-25){\line(1,0){5}}
		\put(65,-25){\line(0,1){5}}
		\put(65,-20){\line(1,0){5}}
		\put(70,-20){\line(0,1){5}}
		\thinlines 
		\put(37,-37){$(9)$}
		\begin{Huge}
			\put(75,-40){$\dashrightarrow$}
		\end{Huge}
		\begin{tiny}
			\put(75,-45){$0,..,0,-U_2$}
		\end{tiny}

		\thinlines
		\put(115,-50){\line(1,0){45}}
		\put(115,-50){\line(0,1){30}}
		\put(115,-20){\line(1,0){35}}
		\put(150,-20){\line(0,1){5}}
		\put(150,-15){\line(1,0){15}}
		\put(160,-50){\line(0,1){10}}
		\put(160,-40){\line(1,0){5}}
		\put(165,-40){\line(0,1){25}}
		\thinlines 
		\put(130,-37){$(10)$}
		\begin{Huge}
			\put(170,-40){$\dashrightarrow$}
		\end{Huge}
		\begin{footnotesize}
			\put(177,-45){$+U_1$}
		\end{footnotesize}

		\thinlines
		\put(210,-50){\line(1,0){45}}
		\put(210,-50){\line(0,1){30}}
		\put(210,-20){\line(1,0){15}}
		\put(225,-20){\line(0,1){5}}
		\put(225,-15){\line(1,0){5}}
		\put(230,-15){\line(0,-1){5}}
		\put(230,-20){\line(1,0){5}}
		\put(235,-20){\line(0,1){5}}
		\put(235,-15){\line(1,0){25}}
		\put(255,-50){\line(0,1){30}}
		\put(255,-20){\line(1,0){5}}
		\put(260,-20){\line(0,1){5}}
		
		\thinlines 
		\put(225,-35){$(11)$}
		\begin{Huge}
			\put(265,-40){$\longrightarrow$}
		\end{Huge}
		\begin{footnotesize}
			\put(272,-45){$-U_1$}
		\end{footnotesize}

		\thinlines
		\put(305,-50){\line(1,0){45}}
		\put(305,-50){\line(0,1){30}}
		\put(305,-20){\line(1,0){25}}
		\put(330,-20){\line(0,1){5}}
		\put(330,-15){\line(1,0){20}}
		\put(350,-50){\line(0,1){30}}
		\put(350,-20){\line(1,0){5}}
		\put(350,-15){\line(1,0){5}}
		\put(355,-20){\line(0,1){5}}
		\put(357,-18){\line(1,0){2}}
		\put(360,-18){\line(1,0){2}}
		\put(363,-18){\line(1,0){2}}
		\put(366,-18){\line(1,0){2}}
		\put(370,-18){\line(1,0){2}}
		\put(372,-17){\line(0,1){2}}
		\put(372,-13){\line(0,1){2}}
		\put(372,-10){\line(0,1){2}}
		\put(370,-8){\line(-1,0){2}}
		\put(367,-8){\line(-1,0){2}}
		\put(364,-8){\line(-1,0){2}}
		\put(361,-8){\line(-1,0){2}}
		\put(358,-8){\line(-1,0){2}}
		\put(355,-8){\line(-1,0){2}}
		\put(352,-8){\line(-1,0){2}}
		\put(349,-8){\line(-1,0){2}}
		\put(346,-8){\line(-1,0){2}}
		\put(343,-8){\line(-1,0){2}}
		\put(340,-8){\line(-1,0){2}}
		\put(337,-8){\line(-1,0){2}}
		\put(334,-8){\line(-1,0){2}}
		\put(331,-8){\line(-1,0){2}}
		\put(328,-8){\line(-1,0){2}}
		\put(326,-9){\line(0,-1){2}}
		\put(326,-11){\vector(0,-1){5}}
		\thinlines 
		\put(320,-35){$(12)$}
		
	\end{picture}
	\vskip 3.5 cm
	\caption{Sliding of the bar $B^e(\h)$ around the frame-angle $c^{en}(\h)$ that connects the configuration $\h$ to the configuration $(12)$, with $\h$ as the configuration $(3)$.}
	\label{fig:columntorow}
\end{figure}
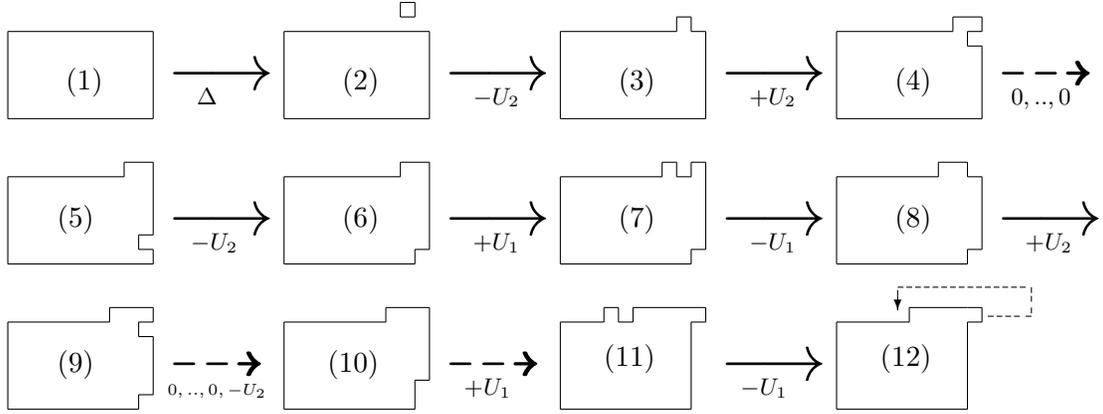

\subsection{Gate for isotropic interactions}
\label{iso}
In this Section we impose $U_1=U_2=U$ in (\ref{hamilt}), i.e., we consider isotropic interactions between nearest-neighbors sites. Recall (\ref{defepsilon}) for the definition of $\e$. We will consider $0<\e\ll U$, where $\ll$ means sufficiently smaller; for instance $\e\leq\frac{U}{100}$ is enough. Many interesting quantities that follow have lower index $is$ to remind that they refer to isotropic interactions.

In order to state our main result we recall some important definitions that are given in \cite{BHN}.
\bd{critsets}
(a) Let 
\be{deflc}
l_c:=\left\lceil {U\over 2U-\D}\right\rceil
\ee

\noindent
the critical length, where $\lceil \;\rceil$ denotes the integer part plus 1. \\
(b) Let $\cQ_{is}$ denote the set of configurations having one cluster consisting 
of an $(l_c-1)\times l_c$ quasi-square anywhere in $\Lambda_0$ with a 
single particle attached anywhere to one of its sides.\\
(c) Let $\cD_{is}$ denote the set of configurations that can be reached from some 
configuration in $\cQ_{is}$ via a $U$-path, i.e.,
\be{cRdef}
\cD_{is}:= \Big\{\h'\in\cV_{n_{is}^c}\colon\, \exists\,\h\in\cQ_{is}\colon\,
H(\h) = H(\h'),\,\Phi_{|\cV_{n_{is}^c}}(\h,\h') \leq H(\h) + U \Big\},
\ee
where $n_{is}^c:=l_c(l_c-1)+1$ is the volume of the clusters in $\cQ_{is}$.\\
(d) Let $\cigeo:=\cD_{is}^{fp}$. \\
(e) Let
\be{Gdef}
\ba{lll}
\gi &:=& H(\cigeo) = H(\cD_{is}^{fp}) = H(\cD_{is})+\D= H(\cQ_{is})+\D\\ 
&=&-U[(l_c-1)^2+l_c(l_c-2)+1]+\D[l_c(l_c-1)+2]\\
&=&2U[l_c+1] - (2U-\D)[l_c(l_c-1)+2]
\ea
\ee

\noindent
denote the energy of the configurations in $\cC_{is}^*$. 
\ed

By \cite[Theorem 1.4.1]{BHN} we obtain the geometric description of the set $\cD_{is}$ as $\cD_{is}=\bar\cD_{is}\cup\widetilde\cD_{is}$ that will be useful later on. Roughly speaking, one can think of $\cD_{is}$ as the set of configurations consisting of a rectangular cluster with four bars attached to its four sides, whose lengths satisfy precise conditions. See Figure \ref{fig:Upath} for an example of a configuration in $\bar\cD_{is}$ that is obtained via a $U$-path from a configuration in $\bar\cQ_{is}$.

\setlength{\unitlength}{0.25cm}

\begin{figure}
	\begin{picture}(42,15)(-13,0)

		
		\thinlines 
		\qbezier[51](0,0)(7,0)(14,0)
		\thinlines
		\qbezier[51](0,0)(0,7)(0,14)
		\qbezier[51](14,0)(14,7)(14,14)
		
		\qbezier[51](1,13)(6,13)(13,13) 
		\qbezier[51](1,1)(1,6)(1,13)
		\qbezier[51](13,1)(13,6)(13,13)
		
		\put(0,1){\line(1,0){4}}
		\put(4,0){\line(1,0){1}}
		\put(5,1){\line(1,0){9}}
		\put(0,14){\line(1,0){14}}
		
		\put(4,0){\line(0,1){1}}
		\put(5,0){\line(0,1){1}}
		\put(14,1){\line(0,1){13}}
		\put(0,1){\line(0,1){13}}
		\put(4.5,6.5){$12 \times 12$} 
		\put(6,-4){$\bar {\cal Q}_{is}$}
		
		\put(17,7){\vector(1,0){5}} 
		\put(17,5){$U$-path}
		
		
		\thinlines
		\qbezier[51](26,0)(37,0)(40,0)
		\qbezier[51](40,0)(40,7)(40,14) 
		\qbezier[51](26,14)(37,14)(40,14)
		\qbezier[51](26,0)(26,6)(26,14)
		
		\thinlines
		\qbezier[51](27,1)(31,1)(39,1)
		\qbezier[51](27,13)(31,13)(39,13)
		\qbezier[51](27,1)(27,7)(27,13)
		\qbezier[51](39,1)(39,7)(39,13)
		
		\put(26,1){\line(1,0){4}}
		\put(30,0){\line(1,0){8}}
		\put(38,1){\line(1,0){1}}
		\put(39,2){\line(1,0){1}}
		\put(38,13){\line(1,0){2}}
		\put(28,14){\line(1,0){10}}
		\put(27,13){\line(1,0){1}}
		\put(26,11){\line(1,0){1}}
		
		\put(30,0){\line(0,1){1}}
		\put(38,0){\line(0,1){1}}
		\put(39,1){\line(0,1){1}}
		\put(40,2){\line(0,1){11}}
		\put(38,13){\line(0,1){1}}
		\put(28,13){\line(0,1){1}}
		\put(38,13){\line(0,1){1}}
		\put(27,11){\line(0,1){2}}
		\put(26,1){\line(0,1){10}} 
		\put(30.5,6.5){$12 \times 12$}
		\put(32,-4){$\bar {\cal D}_{is}$}

	\end{picture}
	\vskip 1 cm
	\caption{Configurations in $\bar\cQ_{is}$ and $\bar\cD_{is}$ for $l_c=14$. A similar picture applies for $\widetilde\cQ_{is}$ and $\widetilde\cD_{is}$ with a $11 \times 13$ rectangle in the center. Note that in the anisotropic regime the cost is $\max\{U_1,U_2\}=U_1$.}
	\label{fig:Upath}
\end{figure}
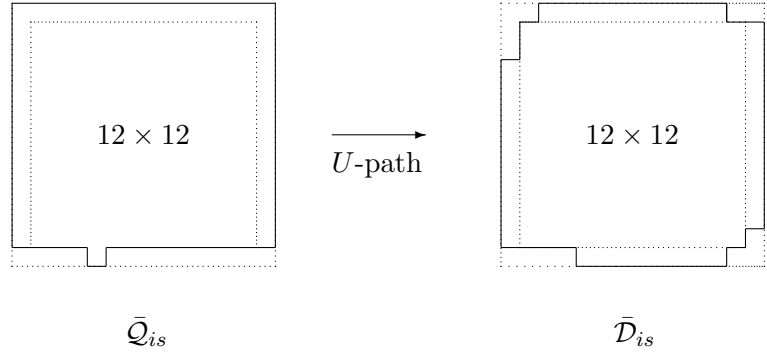

Next we define two types of sets that will be useful in order to characterize the set $\cG_{is}(\vuoto,\pieno)$. For any $i=0,1,2,3$ we define $\cP_{is,i}\subseteq\cS_{is}(\vuoto,\pieno)$ to consist of configurations with a single cluster and no free particle, and a fixed number of vacancies that is not monotone with circumscribed rectangle obtained from the one of the configurations in $\cD_{is}$ via increasing and/or decreasing the horizontal or vertical length. More precisely,
\be{defpi}
\ba{ll}
\cP_{is,i}:= \{\h&:\, n(\h)= 0,\, v(\h)=2l_c-i(i+1)-2,\,
\h_{cl} \hbox{ is connected},\, g_1'(\h)+g_2'(\h)=1,\\
&\hbox{ with a } (l_c+i+1)\times(l_c-i) \hbox{ circumscribed rectangle}\}, \ i=0,1,2,3.
\ea
\ee

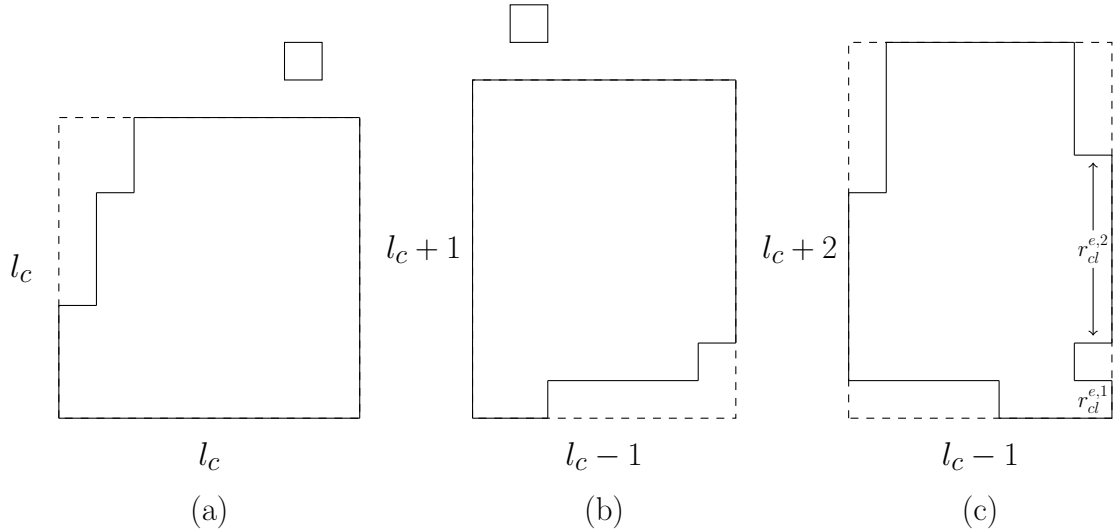
\begin{figure}[htp]
	\centering
	\begin{tikzpicture}[scale=0.5,transform shape]
		
		\draw[dashed] (-1,0) rectangle (7,8);
		\draw (-1,0)--(7,0);
		\draw (7,0)--(7,8);
		\draw (1,8)--(7,8);
		\draw (1,8)--(1,6);
		\draw (0,6)--(1,6);
		\draw(-1,0)--(-1,3);
		\draw (-1,3)--(0,3);
		\draw(0,3)--(0,6);
		\node at (3,-1){\Huge{$l_c$}};
		\node at (-2,4){\Huge{$l_c$}};
		\node at (3,-2.5){\Huge{(a)}};
		\draw (5,9) rectangle (6,10);
		
		\draw[dashed] (10,0) rectangle (17,9);
		\node at (13.5,-1){\Huge{$l_c-1$}};
		\node at (13.5,-2.5){\Huge{(b)}};
		\node at (8.7,4.5){\Huge{$l_c+1$}};
		\draw (10,0)--(10,9);
		\draw (10,9)--(17,9);
		\draw(17,9)--(17,2);
		\draw(17,2)--(16,2);
		\draw(16,2)--(16,1);
		\draw (16,1)--(12,1);
		\draw (12,1)--(12,0);
		\draw(12,0)--(10,0);
		\draw (11,10) rectangle (12,11);

		\draw[dashed] (20,0) rectangle (27,10);
		\node at (23.5,-1){\Huge{$l_c-1$}};
		\node at (23.5,-2.5){\Huge{(c)}};
		\node at (18.7,4.5){\Huge{$l_c+2$}};
		\draw (24,0)--(27,0);
		\draw (27,0)--(27,1);
		\draw(26,1)--(27,1);
		\draw(26,1)--(26,2);
		\draw(27,2)--(26,2);
		\draw(27,2)--(27,3);
		\draw(27,3)--(27,7);
		\draw(27,7)--(26,7);
		\draw(26,7)--(26,10);
		\draw(26,10)--(24,10);
		\draw(24,10)--(21,10);
		\draw(21,10)--(21,6);
		\draw(21,6)--(20,6);
		\draw (20,6)--(20,1);
		\draw(20,1)--(24,1);
		\draw(24,1)--(24,0);
		\node at (26.5,0.5){\LARGE{$r_{cl}^{e,1}$}};
		\draw[->](26.5,4)--(26.5,2.2);
		\draw[->](26.5,5)--(26.5,6.8);
		\node at (26.5,4.5){\LARGE{$r_{cl}^{e,2}$}};

	\end{tikzpicture}
	
	\vskip 0 cm
	\caption{Critical configurations in the isotropic case: in (a) is represented a configuration in $\cP_{is,-1}^{fp}\setminus\cigeo$, in (b) a configuration in $\cP_{is,0}^{fp}\cap\cigeo$ and in (c) a configuration in $\cP_{is,1}$ and not in $\cigeo$.}
	\label{fig:figiso}
\end{figure}

\noindent
See Figure \ref{fig:figiso}(c) for an example of configurations in $\cP_{is,1}$ and not in $\cigeo$.

For any $ i=-1,0,1,2$ we define $\cP_{is,i}^{fp}\subseteq\cS_{is}(\vuoto,\pieno)$ to consist of configurations with a single cluster and one free particle, and a fixed number of vacancies that is monotone with circumscribed rectangle obtained from the one of the configurations in $\cP_{is,i}$ via decreasing by one the shortest length. More precisely,
\be{defpifp}
\ba{ll}
\cP_{is,i}^{fp}:=& \{\h:\, n(\h)= 1,\, v(\h)=l_c-i(i+2)-2,\,
\h_{cl} \hbox{ is connected},\, g_1'(\h)+g_2'(\h)=0,\\
&\hbox{with a }(l_c+i+1)\times(l_c-i-1) \hbox{ circumscribed rectangle}\}, \ i=-1,0,1,2.
\ea
\ee

\noindent
See Figure \ref{fig:figiso} for an example of configuration in $\cP_{is,-1}^{fp}\setminus\cigeo$ (in (a)) and in $\cP_{is,0}^{fp}$ and in $\cigeo$ (in (b)). Note that other examples of configurations in $\cigeo$ can be obtained by those in Figure \ref{fig:Upath} by adding a free particle.

\br{}
Note that $\cigeo\subsetneq\cP_{is,-1}^{fp}\cup\cP_{is,0}^{fp}$. Indeed, there exist configurations that are in $\cP_{is,-1}^{fp}\cup\cP_{is,0}^{fp}$ and not in $\cigeo$: for an example see Figure \ref{fig:figiso}(a).
\er

The set $\cG_{is}(\vuoto,\pieno)$ contains all the configurations that are in the sets defined in (\ref{defpi}) and (\ref{defpifp}) with the following further conditions. First, for any $i=0,1,2,3$ we define the subset $\mathscr{I}_i^{\a}$ of the saddles in $\cP_{is,i}$ with the condition that the configurations have only one occupied unit square in either a row or in one of its two adjacent frame-angles. More precisely,
\be{I1}
\mathscr{I}_i^{\a}:=\{\h\in\cP_{is,i}: |r^{\a}(\h)\cup c^{\a\a'}(\h)\cup c^{\a\a''}(\h)|=1\}, \quad i=0,1,2,3,
\ee

\noindent
with $\a,\a',\a''\in\{n,s,w,e\}$ such that $c^{\a\a'}(\h)$ and $c^{\a\a''}(\h)$ are different frame-angles. Next, for any $i=0,1,2$ we define the subset $\mathscr{I}_{k,k',i}^{\a,\a'}$ of the saddles in $\cP_{is,i}$ that are obtained from $\h\in\cP_{is,i}$ during the sliding of the bar $B^{\a'}(\h)$ around the frame-angle $c^{\a'\a}(\h)$. More precisely,
\be{defpistrenino}
\ba{ll}
\mathscr{I}_{k,k',i}^{\a,\a'}&:=\{\h\in\cP_{is,i}:|r^\a(\h)|=k-1, |r^{\a'}(\h)|=k'-k+1, k'\leq1+||r^{\a}(\h)||,|c^{\a'\a}(\h)|=1,\\
&(r^{\a}(\h)\cup c^{\a'\a}(\h))\cap\h_{cl}=r^{\a,1}_{cl}\dot{\cup}r^{\a,2}_{cl} \hbox{ with } d(r^{\a,1}_{cl},r^{\a,2}_{cl})=2\}, \quad i=0,1,2,
\ea
\ee
\noindent
where $\a,\a'\in\{n,s,w,e\}$ is such that $c^{\a'\a}(\h)$ is a frame-angle, $r_{cl}^{\a,1}$, $r_{cl}^{\a,2}$ are two disjoint connected components in $r^{\a}(\h)\cup c^{\a'\a}(\h)$ and $k'=2,...,l_c$, $k=2,...,k'$ (see for example the configuration in Figure \ref{fig:figiso}(c), which is in $\mathscr{I}_{l_c-2,l_c-1,1}^{e,s}$). The index $i$ in (\ref{defpistrenino}) has to be different from 3 because if the system is in $\mathscr{I}_{3}^{\a}$, then it is not possible to complete the sliding of a bar around the frame-angle. Note that the conditions in (\ref{defpistrenino}) guarantee that these configurations are obtained during a sliding of a bar around a frame-angle, that is identified by the indeces $\a$ and $\a'$. Moreover, the index $k'$ denotes the length of the bar that we are sliding. The index $k$ counts the number of particles that are in $r^{\a}(\h)\cup c^{\a'\a}(\h)$ during the sliding and can be less or equal than $l_c$, but for some values of $k$ the set $\mathscr{I}_{k,k',i}^{\a,\a'}$ can be empty. Our notation does not distinguish whether $\mathscr{I}_{k,k',i}^{\a,\a'}$  is empty or not in order to avoid the presence of an additional index.

Furthermore, for any $i=-1,0,1,2$ we define the subset $\mathscr{I}_i^{\a,\a'}$ of the saddles in $\cP_{is,i}^{fp}$ as the last saddle at the end of a path that describes the sliding of a bar around a frame-angle, i.e., the saddle where the last particle of the bar (protuberance) is detached. (See configuration (12) in the path described in Figure \ref{fig:columntorow}, imposing $U_1=U_2=U$ because we are in the isotropic case). More precisely,
\be{gatetrenino}
\ba{ll}
\mathscr{I}_i^{\a,\a'}:=\{\h\in\cP_{is,i}^{fp}&:\exists \ \h'\in\mathscr{I}_{2,k',i}^{\a,\a'} \text{ such that } \h \text{ is obtained from } \h' \text{ by removing the }\\
&\text{non monotonicity, then removing the protuberance and }\\
&\text{possibly moving the free particle at zero cost}\}, \quad i=-1,0,1,2,
\ea
\ee

\noindent
where $\a,\a'\in\{n,s,w,e\}$. Again, the indices $\a$ and $\a'$ identify the frame-angle with respect to which the sliding of the bar takes place. Note that the configuration in Figure \ref{fig:figiso}(a) is in $\mathscr{I}_0^{s,w}\cup\mathscr{I}_0^{e,n}$ and the configuration in Figure \ref{fig:figiso}(b) is in $\mathscr{I}_{-1}^{w,s}\cup\mathscr{I}_{-1}^{w,n}\cup\mathscr{I}_{-1}^{e,n}\cup\mathscr{I}_{1}^{n,e}$.

In the discussion below \cite[Theorem 1.4.3]{BHN}, the authors state that the full identification of the set $\cG_{is}(\vuoto,\pieno)$ is not known. The following result fills this gap. 

\bt{giso}{\rm{(Union of minimal gates for isotropic interactions).}} 
We obtain the following description for $\cG_{is}(\vuoto,\pieno)$:
\be{}
\cG_{is}(\vuoto,\pieno)=\cigeo\cup\displaystyle\bigcup_{ i=0}^{3}\bigcup_{\a}\mathscr{I}_i^{\a}\cup\bigcup_{i=0}^{2}\bigcup_{\a,\a'}\bigcup_{k,k'}\mathscr{I}_{k,k',i}^{\a,\a'}\cup\bigcup_{i=-1}^{2}\bigcup_{\a,\a'}\mathscr{I}_i^{\a,\a'}
\ee
\et

\noindent
We refer to Section \ref{propiso} for the proof of the main Theorem \ref{giso}.

In \cite[Theorem 1.4.3(i)]{BHN} the authors show that in $\cS_{is}(\vuoto,\pieno)$ there are unessential saddles, also called dead-ends, without fully identifying them, while in Corollary \ref{corstrategy} and Proposition \ref{selle3} we identify three types of unessential saddles. Moreover, in Proposition \ref{c*contenuto} we prove that $\cC_{is}^*$ is contained in $\cG(\vuoto,\pieno)$, which contradicts what is said in the discussion below \cite[Theorem 1.4.3]{BHN}.

\subsection{Gate for weakly anisotropic interactions}
\label{wani}

In this Section we impose $U_1\neq U_2$ and $U_1<2U_2-2\e$ in (\ref{hamilt}) (recall (\ref{defepsilon}) for the definition of $\e$), i.e., we consider weakly anisotropic interactions between nearest neighboring sites. We will consider $0<\e\ll U_2$, where $\ll$ means suffiently smaller; for instance $\e\leq\frac{U_2}{100}$ is enough. Many interesting quantities that follow have lower index $wa$ to remind the reader that they refer to weakly anisotropic interactions. In order to state our main results for the gates in the weakly anisotropic regime we need the following definitions. Let
\be{defbarl}
\bar l:=\Bigg\lceil{U_1-U_2\over U_1+U_2-\D}\Bigg\rceil.
\ee

For any $s> \bar l+2 $, if $s$ has the same parity as
$\bar l$ i.e.,
$[s-\bar l]_2=[0]_2$, then we define the set of {\it $0$-standard
	rectangles} as
$\cR^{0-st}(s):=\cR(\ell_1(s),\ell_2(s))$
with side lengths 
\be{ls0} \ell_1(s):= {s+\bar l\over 2} \qquad \ell_2(s):= {s-\bar
	l\over 2},\qquad \hbox{
	for } [s-\bar l]_2=[0]_2.
\ee
If $s$ has the same parity as $\bar l-1$ i.e.,  $[s-\bar
l]_2=[1]_2$, we define the set of {\it $1$-standard rectangles} to be
$\cR^{1-st}(s):=\cR(\ell_1(s),\ell_2(s))$
with side lengths
\be{ls1}
\ell_1(s):={s+\bar l-1\over 2}, \qquad \ell_2(s):={s-\bar l +1\over 2} \qquad
\hbox{ for } [s-\bar l]_2=[1]_2.
\ee

\noindent
For this value of $s$ we define the set of {\it quasi-standard}
rectangles as $\cR^{q-st}(s):=\cR(\ell_1(s)+1,\ell_2(s)-1)$. Finally, we set
\be{standardgen} \cR^{st}(s):=\left\{\ba{ll}
\cR^{0-st}(s)  &\mbox{if } [s-\bar l]_2= [0]_2\\
\cR^{1-st}(s)  &\mbox{if } [s-\bar l]_2= [1]_2. \ea\right. \ee

\noindent
We define the critical horizontal length $l_1^*$ and the critical vertical length $l_2^*$ as 
\be{critinteri}
l_1^*:=\left\lceil{U_1 \over U_1+U_2-\D}\right\rceil, \quad \quad l_2^*:=\left\lceil {U_2 \over U_1+U_2 -\D}\right\rceil,
\ee

\noindent
where $\lceil \;\rceil$ denotes the integer part plus 1. We set the critical value of $s_{wa}$ as
\be{defweaks^*} 
s_{wa}^*:=l_1^*+l_2^*-1.
\ee

\noindent
We need the following definitions.

\bd{dweak}
\bi
\item[(a)] Define $\bar\cQ_{wa}$ as the set of configurations having one cluster anywhere in $\L_0$ consisting of a $(l_1^*-1)\times l_2^*$ rectangle with a single protuberance attached to one of the shortest sides. Similarly, we define $\widetilde\cQ_{wa}$ as the set of configurations having one cluster anywhere in $\L_0$ consisting of a $(l_1^*-1)\times l_2^*$ rectangle with a single protuberance attached to one of the longest sides. 
\item[(b)] Define
\be{gammaweak}
\gw:=U_1l_2^*+U_2l_1^*+U_1+U_2+\e l_2^*-\e l_1^*l_2^*-2\e.
\ee

\item[(c)] Define the volume of the clusters in $\bar\cQ_{wa}$ as
\be{defnwa}
n_{wa}^c:=l_2^*(l_1^*-1)+1,
\ee

\noindent
and
\be{defDwa}
\ba{lll}
\bar\cD_{wa}:=\{\h'\in\cV_{n_{wa}^c}| \ \exists \ \h\in\cQ_{wa}: H(\h)=H(\h') \hbox{ and } \Phi_{|\cV_{n_{wa}^c}}(\h,\h')\leq H(\h)+U_1\},\\
\widetilde\cD_{wa}:=\{\h'\in\cV_{n_{wa}^c}| \ \exists \ \h\in\widetilde\cQ_{wa}: H(\h)=H(\h') \hbox{ and } \Phi_{|\cV_{n_{wa}^c}}(\h,\h')\leq H(\h)+U_1\}.
\ea
\ee

\noindent
{\rm Note that the last condition in (\ref{defDwa}) is the same as requiring that $\Phi_{|\cV_{n_{wa}^c}}(\h,\h')<\gw+H(\vuoto)=\gw$. We encourage the reader to consult Proposition \ref{cardweak}, where we give the geometric description of the sets $\bar\cD_{wa}$ and $\widetilde\cD_{wa}$. Roughly speaking, one can think of $\bar\cD_{wa}$ and $\widetilde\cD_{wa}$ as the sets of configurations consisting of a rectangular cluster with four bars attached to its four sides, whose lengths satisfy precise conditions.}

\item[(d)] Define
\be{c*wa}
\cwgeo:=\bar\cD_{wa}^{fp}.
\ee

\ei
\ed

\noindent
The reason why only the set $\bar\cD_{wa}$ is relevant for the set $\cwgeo$ will be clarified later (see Lemma \ref{dtilde}). Note that
\be{}
\ba{lll}
H(\cwgeo)&=&H(\bar\cD_{wa}^{fp})=H(\bar\cD_{wa})+\D=H(\bar\cQ_{wa})+\D\\
&=&U_1l_2^*+U_2(l_1^*-1)-\e l_2^*(l_1^*-1)+2\D-U_1\\
&=&U_1l_2^*+U_2l_1^*+U_1+U_2+\e l_2^*-\e l_1^*l_2^*-2\e\\
&=&\gw.
\ea
\ee

\noindent
See Figure \ref{fig:Pweak} for examples of configurations in $\cwgeo$.

\br{enq}
Note that $H(\bar\cQ_{wa})<H(\widetilde\cQ_{wa})$. Indeed,
\be{energiaQwa}
\ba{lll}
H(\bar\cQ_{wa})=\gw-\D, \\
H(\widetilde\cQ_{wa})=\gw-\D+U_1-U_2.
\ea
\ee
\er

\noindent
The first main result of Section \ref{wani} is a refinement of \cite[Theorem 2]{NOS}.

\setlength{\unitlength}{1.1pt}
\begin{figure}[htp]
	\centering
	\begin{picture}(400,90)(0,30)
		\thinlines
		\qbezier[51](20,0)(100,0)(180,0)
		\qbezier[51](20,90)(100,90)(180,90)
		\qbezier[51](20,0)(20,45)(20,90)
		\qbezier[51](180,0)(180,45)(180,90)
		\thinlines
		\put(20,0){\line(1,0){150}}
		\put(20,90){\line(1,0){150}}
		\put(20,0){\line(0,1){90}}
		\put(170,0){\line(0,1){60}}
		\put(170,70){\line(0,1){20}}
		\put(170,60){\line(1,0){10}}
		\put(180,60){\line(0,1){10}}
		\put(170,70){\line(1,0){10}}
		\put(140,100){\line(1,0){10}}
		\put(150,100){\line(0,1){10}}
		\put(140,100){\line(0,1){10}}
		\put(140,110){\line(1,0){10}}
		\thinlines 
		\put(190,45){$l_2^*$}
		\put(100,-15){$l_1^*$}

		\thinlines
		\qbezier[51](220,0)(300,0)(380,0)
		\qbezier[51](220,90)(300,90)(380,90)
		\qbezier[51](220,0)(220,45)(220,90)
		\qbezier[51](380,0)(380,45)(380,90)
		\thinlines
		\put(230,10){\line(1,0){10}}
		\put(240,10){\line(0,-1){10}}
		\put(240,0){\line(1,0){130}}
		\put(370,0){\line(0,1){20}}
		\put(370,20){\line(1,0){10}}
		\put(380,20){\line(0,1){60}}
		\put(370,80){\line(1,0){10}}
		\put(370,80){\line(0,1){10}}
		\put(230,10){\line(0,1){10}}
		\put(230,20){\line(-1,0){10}}
		\put(220,20){\line(0,1){60}}
		\put(220,80){\line(1,0){20}}
		\put(240,80){\line(0,1){10}}
		\put(240,90){\line(1,0){130}}
		\put(300,100){\line(1,0){10}}
		\put(310,100){\line(0,1){10}}
		\put(300,100){\line(0,1){10}}
		\put(300,110){\line(1,0){10}}
		\put(390,45){$l_2^*$}
		\put(300 ,-15){$l_1^*$}
	\end{picture}
	\vskip 1.5 cm
	\caption{Critical configurations in $\cwgeo$ in the weakly anisotropic case. Moreover, if we remove the free particle, then we obtain on the left a configuration in $\bar\cQ_{wa}$ and on the right a configuration in $\bar\cD_{wa}\setminus\bar\cQ_{wa}$.}
	\label{fig:Pweak}
\end{figure}
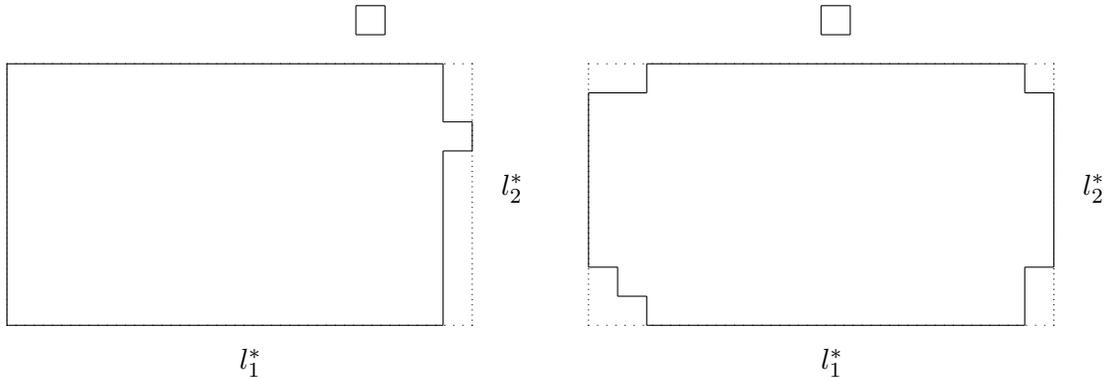

\bt{thgateweak}{\rm{(Gate for weakly anisotropic interactions).}}
The set $\cwgeo$ is a gate for the transition from $\vuoto$ to $\pieno$.
\et

\noindent
We refer to Section \ref{gateweak1} for the proof of Theorem \ref{thgateweak}.

In order to give the result regarding the geometric characterization of $\cG_{wa}(\vuoto,\pieno)$, we need some definitions. For any $i=0,1$ we define $\cP_{wa,i}\subseteq\cS_{wa}(\vuoto,\pieno)$ to consist of configurations with a single cluster and no free particle, and a fixed number of vacancies that is not monotone with circumscribed rectangles obtained from the one of the configurations in $\bar\cD_{wa}$ via increasing by one the horizontal or vertical length. More precisely,
\be{defPwa}
\ba{ll}
\cP_{wa,i}:= &\{\h:\, n(\h)= 0,\, v(\h)=2l_2^*+i(l_1^*-l_2^*)-2,\, g_1'(\h)=i,\,g_2'(\h)=1-i,\,\h_{cl} \hbox{ is}\\
&\hbox{ connected }\hbox{with circumscribed rectangle in }
\cR(l_1^*-i+1,l_2^*+i)\}, \ i=0,1.
\ea
\ee
\noindent 
See Figure \ref{fig:Pweak2} for examples of configurations in $\cP_{wa,0}$ (on the left-hand side) and in $\cP_{wa,1}$ (on the right-hand side). 

\setlength{\unitlength}{1.1pt}
\begin{figure}
	\centering
	\begin{picture}(400,90)(0,30)
		\thinlines
		\qbezier[51](20,0)(100,0)(180,0)
		\qbezier[51](20,90)(100,90)(180,90)
		\qbezier[51](20,0)(20,45)(20,90)
		\qbezier[51](180,0)(180,45)(180,90)
		\thinlines
		\put(30,0){\line(1,0){140}}
		\put(20,90){\line(1,0){30}}
		\put(50,90){\line(0,-1){10}}
		\put(50,80){\line(1,0){10}}
		\put(60,80){\line(0,1){10}}
		\put(60,90){\line(1,0){110}}
		\put(30,0){\line(0,1){70}}
		\put(30,70){\line(-1,0){10}}
		\put(20,70){\line(0,1){20}}
		\put(170,0){\line(0,1){60}}
		\put(170,70){\line(0,1){20}}
		\put(170,60){\line(1,0){10}}
		\put(180,60){\line(0,1){10}}
		\put(170,70){\line(1,0){10}}
		\thinlines 
		\put(190,45){$l_2^*$}
		\put(88,-15){$l_1^*+1$}

		\thinlines
		\qbezier[51](220,0)(295,0)(370,0)
		\qbezier[51](220,100)(295,100)(370,100)
		\qbezier[51](220,0)(220,50)(220,100)
		\qbezier[51](370,0)(370,50)(370,100)
		\thinlines
		\put(230,0){\line(1,0){130}}
		\put(360,0){\line(0,1){20}}
		\put(360,20){\line(1,0){10}}
		\put(370,20){\line(0,1){10}}
		\put(370,30){\line(-1,0){10}}
		\put(360,30){\line(0,1){10}}
		\put(360,40){\line(1,0){10}}
		\put(370,40){\line(0,1){10}}
		\put(370,50){\line(-1,0){10}}
		\put(360,50){\line(0,1){40}}
		\put(360,90){\line(-1,0){100}}
		\put(260,90){\line(0,1){10}}
		\put(260,100){\line(-1,0){10}}
		\put(250,100){\line(0,-1){10}}
		\put(250,90){\line(-1,0){30}}
		\put(220,90){\line(0,-1){80}}
		\put(220,10){\line(1,0){10}}
		\put(230,10){\line(0,-1){10}}
		\put(380,45){$l_2^*+1$}
		\put(300 ,-15){$l_1^*$}
	\end{picture}
	\vskip 1.5 cm
	\caption{Critical configurations in the weakly anisotropic case: on the left hand-side is represented a configuration in $\cP_{wa,0}$ and on the right hand-side a configuration in $\cP_{wa,1}$.}
	\label{fig:Pweak2}
\end{figure}
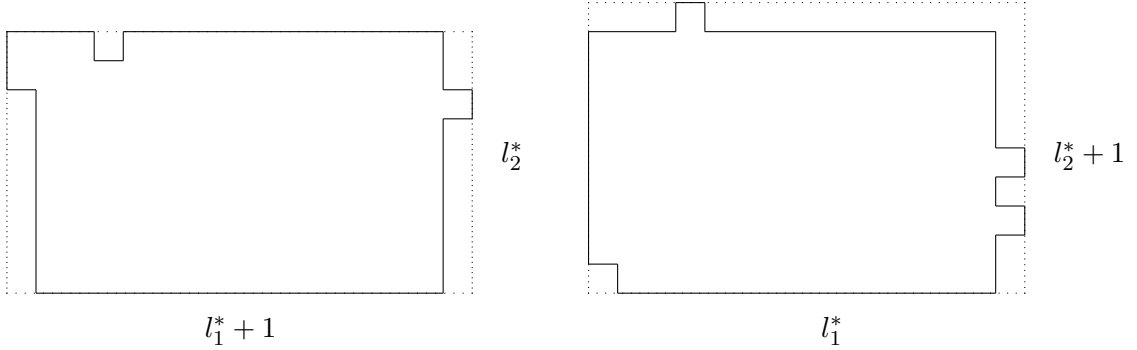

The set $\cG_{wa}(\vuoto,\pieno)$ contains all the configurations that are in the sets defined in (\ref{defPwa}) with the following further conditions. First, we define the subset $\cN_0^{\a'}$ (resp.\ $\cN_1^{\a}$) of the saddles in $\cP_{wa,0}$ (resp.\ $\cP_{wa,1}$) as the set that contains only one occupied unit square in either a vertical (resp.\ horizontal) row or in one of its two adjacent frame-angles. More precisely,

\be{N1}
\cN_0^{\a'}:=\{\h\in\cP_{wa,0}: |r^{\a'}(\h)\cup c^{\a'\bar\a}(\h)\cup c^{\a'\widetilde\a}(\h)|=1\},
\ee

\noindent
for any $\a'\in\{w,e\}$ and $\bar\a,\widetilde\a\in\{n,s\}$ such that $\bar\a\neq\widetilde\a$, and
\be{N2}
\cN_1^{\a}:=\{\h\in\cP_{wa,1}: |r^{\a}(\h)\cup c^{\a\a''}(\h)\cup c^{\a\a'''}(\h)|=1\},
\ee
for any $\a\in\{n,s\}$ and $\a'',\a'''\in\{w,e\}$ such that $\a''\neq\a'''$. Note that in Figure \ref{fig:Pweak2} the configuration on the left-hand side is in $\cN_0^{e}$ and the configuration on the right-hand side is in $\cN_1^{n}$.

Next, we define the subset $\cN_{k,k'}^{\a,\a'}$ of the saddles in $\cP_{wa,0}$ as the set of configurations that are obtained from $\h\in\cP_{wa,0}$ during the sliding of the bar $B^{\a'}(\h)$ around the frame-angle $c^{\a'\a}(\h)$. More precisely,
\be{defpwatrenino}
\ba{ll}
\cN_{k,k'}^{\a,\a'}:=&\{\h\in\cP_{wa,0}:|r^\a(\h)|=k-1, |r^{\a'}(\h)|=k'-k+1, k'\leq1+||r^{\a}(\h)||,\\
&|c^{\a'\a}(\h)|=1,(r^{\a}(\h)\cup c^{\a'\a}(\h))\cap\h_{cl}=r^{\a,1}_{cl}\dot{\cup}r^{\a,2}_{cl} \hbox{ with } d(r^{\a,1}_{cl},r^{\a,2}_{cl})=2\},
\ea
\ee

\noindent
where $\a\in\{n,s\}, \ \a'\in\{w,e\}$, $r_{cl}^{\a,1}$, $r_{cl}^{\a,2}$ are two disjoint connected components in $r^{\a}(\h)\cup c^{\a'\a}(\h)$ and $k'=l_1^*-l_2^*+1,...,l_2^*$, $k=2,...,k'$. Note that the conditions in (\ref{defpwatrenino}) guarantee that these configurations are obtained during a sliding of a bar around a frame-angle, that is identified by the indeces $\a$ and $\a'$. Note that the index $k'$ denotes the length of the bar that we are sliding. The index $k$ counts the number of particles that are in $r^{\a}(\h)\cup c^{\a'\a}(\h)$ during the sliding and can be less or equal than $l_2^*$, but it is possible that for some values of $k$ the set $\cN_{k,k'}^{\a,\a'}$ is empty. Our notation does not distinguish whether $\cN_{k,k'}^{\a,\a'}$ is empty or not in order to avoid the presence of an additional index. 

Now we are able to give the second main result of Section \ref{wani}.

\bt{gweak}{\rm{(Union of minimal gates for weakly anisotropic interactions).}}
We obtain the following description for $\cG_{wa}(\vuoto,\pieno)$:
\be{}
\cG_{wa}(\vuoto,\pieno)=\cwgeo\cup\displaystyle\bigcup_{\a}\bigcup_{\a'}\bigcup_{k'=l_1^*-l_2^*+1}^{l_2^*}\bigcup_{k=2}^{k'}\cN_{k,k'}^{\a,\a'}\cup\displaystyle\bigcup_{\a'}\cN_0^{\a'}\cup\displaystyle\bigcup_{\a}\cN_1^{\a}
\ee
\et

\noindent
We refer to Section \ref{weak1} for the proof of Theorem \ref{gweak}.

\subsection{Main results: sharp asymptotics for weakly anisotropic interactions}
\label{sharpestimates}

If the reader is interested in the sharp asymptotics below can find the proofs and some discussions in Section \ref{sharpasymptotics}.

\bt{sharptimeweak}
There exists a constant $K_{wa}=K_{wa}(\L,l_2^*)$ such that
\be{nuclsharpweak}
\E_{\,\square}(\tau_\blacksquare) = K_{wa}e^{\Gamma_{wa}^*\beta}\,[1+o(1)],
\qquad \b\to\infty.
\ee
Moreover, as $\L\to\Z^2$,
\be{Casymp}
K_{wa}(\L,l_2^*) \ra \frac{1}{4\pi N_{wa}}\,\frac{\log|\L|}{|\L|}
\ee
with
\be{Nid}
N_{wa}=\displaystyle\sum_{k=1}^4\binom{4}{k}\binom{l_2^*+k-2}{2k-1}
\ee
the cardinality of $\bar\cD_{wa}=\bar\cD_{wa}(\L,l_2^*)$ modulo shifts.
\et

\noindent
Theorem \ref{sharptimeweak} investigates the prefactor for the weakly anisotropic case. This analysis for the isotropic case is given in \cite[Theorem 1.4.4]{BHN}, while for the strongly anisotropic case is given in \cite[Theorem 4.12]{BN3}.

\bt{entrunif}
Let $\t_{\cC_{wa}^{*-}}$ be the time just prior $\t_{\cC_{wa}^*}$. Then the entrance distribution of $\cC_{wa}^*$ is uniform, i.e.,
\be{uniforme}
\lim_{\b\ra\infty}\P_{\vuoto}\Big(\h_{\t_{\cC_{wa}^{*-}}}=\h|\t_{\cC_{wa}^*}<\t_{\vuoto}\Big)=\frac{1}{|\bar\cD_{wa}|} \quad \forall \ \h\in\bar\cD_{wa}.
\ee
\et

\br{}
Note that Theorem \ref{entrunif} concerning the uniform entrance distribution in the gate can not be extended to the strongly anisotropic case due to the two possible entrances in the gate (see \cite[Lemma 5.17]{BN3}).
\er

\noindent
We define the {\it mixing time} as
\be{}
t_{mix}(\e):=\min\{n\geq0: \max_{x\in\cX}||P_n(x,\cdot)-\mu(\cdot)||_{TV}\leq\e\},
\ee

\noindent
where $||\nu-\nu'||_{TV}:=\frac{1}{2}\sum_{x\in\cX}|\nu(x)-\nu'(x)|$ for any two probability distributions $\nu,\nu'$ on $\cX$. The {\it spectral gap} of the Markov chain is defined as
\be{}
\r:=\-a^{(2)}
\ee

\noindent
where $1=a^{(1)}>a^{(2)}\geq...\geq a^{|\cX|}\geq-1$ are the eigenvalues of the matrix $(P(x,y))_{x,y\in\cX}$ defined in (\ref{defkaw}).

\bt{gapspettrale}
Let $int\in\{is,wa\}$. For any $\e\in(0,1)$,
\be{}
\displaystyle\lim_{\b\ra\infty}\frac{1}{\b}\log t_{mix}(\e)=\G_{int}^*=\lim_{\b\ra\infty}-\frac{1}{\b}\log \r
\ee

\noindent
Furthermore, there exist two constants $0<c_1\leq c_2<\infty$ independent of $\b$ such that for every $\b>0$,
\be{}
c_1e^{-\b\G_{int}^*}\leq\r\leq c_2e^{-\b\G_{int}^*}
\ee
\et

\noindent
Theorem \ref{gapspettrale} holds also for the strongly anisotropic case (see \cite[Theorem 4.16]{BN3}).

\section{Main results for the simplified model}
\label{simpletheorem}
In this Section we focus on the simplified model described in Section \ref{simplemodel} for $int\in\{is,wa\}$. For the case $int=sa$, i.e., for the strongly anisotropic case, see \cite[Remarks 2.2 and 4.11]{BN3}. The case $int=is$ has been already studied in \cite{HOS}, see \cite[Theorem 1.53]{HOS} for the main result. In this paper we extend that result to $int=wa$ and we derive a result concerning the union of the minimal gates for $int\in\{is,wa\}$. In this Section we set $\t_{A}$ as the first hitting time of the dynamics in $A$ for $A\subseteq\cX^\b$. We refer to Section \ref{simpleproof} for the proof of the two following Theorems.

\bt{simpleiso}
Let $int=is$. Fix $\D\in(\frac{3}{2}U, 2U)$, with $U/(2U-\D)$ not integer. Suppose that $\lim_{\b\ra\infty}\frac{1}{\b}\log|\L_\b|=\infty$. Then the union of minimal gates $\cG_{is}(\vuoto^\b,\pieno^\b)$ for the transition from $\vuoto^\b$ to $\pieno^\b$ is
\be{}
\cG_{is}(\vuoto^\b,\pieno^\b)=\cJ^\b(\cigeo)\cup\displaystyle\bigcup_{ i=0}^{3}\bigcup_{\a}\cJ^\b(\mathscr{I}_i^{\a})\cup\bigcup_{i=0}^{2}\bigcup_{\a,\a'}\bigcup_{k,k'}\cJ^\b(\mathscr{I}_{k,k',i}^{\a,\a'})\cup\bigcup_{i=-1}^{2}\bigcup_{\a,\a'}\cJ^\b(\mathscr{I}_i^{\a,\a'}),
\ee

\noindent
where $\cigeo$, $\mathscr{I}_i^{\a}$, $\mathscr{I}_{k,k',i}^{\a,\a'}$ and $\mathscr{I}_i^{\a,\a'}$ are defined in Definition \ref{critsets}(d), (\ref{I1}), (\ref{defpistrenino}) and (\ref{gatetrenino}) respectively.
\et

\bt{simpleweak}
Let $int=wa$. Fix $\D\in(U_1+\frac{U_2}{2}, U_1+U_2)$, with $U_2/(U_1+U_2-\D)$ not integer, and $U_1<2U_2-2\e$, where $\e$ is defined in (\ref{defepsilon}). Suppose that $\lim_{\b\ra\infty}\frac{1}{\b}\log|\L_\b|=\infty$.
\bi
\item[(a)] Let $\cR^{\leq(l_1,l_2)}$ (resp.\ $\cR^{\geq(l_1,l_2)}$) be the set of configurations whose single contour is a rectangle contained in (resp.\ containing) a rectangle with sides $l_1$ and $l_2$. Then
\be{}
\ba{lll}\quad
\h\in \cJ^\b(\cR^{\le(l_1^*-1,l_2^*-1)})&\Longrightarrow&
\lim\limits_{\b\to\infty}\P_{\nu_\h} (\t_{\vuoto^\b} < \t_{\pieno^\b}) = 1,\\
\quad \h\in \cJ^\b(\cR^{\ge(l_1^*,l_2^*)})&\Longrightarrow&
\lim\limits_{\b\to\infty}\P_{\nu_\h} (\t_{\pieno^\b} < \t_{\vuoto^\b}) = 1.
\ea
\ee
\item[(b)] The set of configurations $\cJ^\b(\cC_{wa}^*)$, with $\cC_{wa}^*$ defined in (\ref{c*wa}), is a gate for the transition from $\vuoto^\b$ to $\pieno^\b$ and there exists $c>0$ such that, for sufficiently large $\b$,
\be{}
\P_{\nu_{\vuoto}}(\t_{\cJ^\b(\cC_{wa}^*)}>\t_{\pieno^\b})\leq e^{-\b c}.
\ee
\item[(c)] For any $\d>0$,
\be{}
\lim_{\b\ra\infty}\P_{\nu_{\vuoto}}(e^{\b(\G_{wa}^*-\d)}<\t_{\pieno^\b}<e^{\b(\G_{wa}^*+\d)})=1.
\ee
\item[(d)] The union of minimal gates $\cG_{wa}(\vuoto^\b,\pieno^\b)$ for the transition from $\vuoto^\b$ to $\pieno^\b$ is
\be{}
\cG_{wa}(\vuoto^\b,\pieno^\b)=\cJ^\b(\cwgeo)\cup\displaystyle\bigcup_{\a}\bigcup_{\a'}\bigcup_{k'=l_1^*-l_2^*+1}^{l_2^*}\bigcup_{k=2}^{k'}\cJ^\b(\cN_{k,k'}^{\a,\a'})\cup\displaystyle\bigcup_{\a'}\cJ^\b(\cN_0^{\a'})\cup\displaystyle\bigcup_{\a}\cJ^\b(\cN_1^{\a}),
\ee
\noindent
where the sets $\cN_0^{\a'}$, $\cN_1^{\a}$ and $\cN_{k,k'}^{\a,\a'}$ are defined in (\ref{N1}), (\ref{N2}) and (\ref{defpwatrenino}) respectively. 
\ei
\et

\br{}
In the simplified model we focus on the local aspects of metastability and nucleation: the removal of the interactions outside $\L_0$ forces the critical droplet to appear inside $\L_0$. In the original model with interaction and exclusion throughout $\L_\b$, if $\liminf_{\b\ra\infty}\frac{1}{\b}\log|\L_\b|$ is large enough, then the decay from the metastable to the stable state is driven by the formation of many droplets far away from the origin, which subsequently grow, coalesce and reach $\L_0$. This is a much harder problem. In \cite{{GHNOS},{GN}} important steps in this direction are achieved. The authors show that, in the limit as the temperature and the particle density tend to zero simultaneously, the gas of Kawasaki particles evolves as a system of “Quasi-Random Walks”, in other words, it is close to an ideal gas where particles have no interaction. In \cite{GHNOS} the authors are able to deal with a large class of initial conditions having no anomalous concentration of particles and with time horizons that are much larger than the typical particle collision time.
\er

\section{Proof of the model-independent propositions}
\label{modproof}
In this Section we give the proof of Propositions \ref{selle1} and \ref{selle2}.

\subsection{Proof of Proposition \ref{selle1}}
\label{proofdep1}

\bpr
We denote by $\s_1,...,\s_j$ the saddles in the statement. We want to prove that these saddles are unessential (see Section \ref{modinddef} point 4 for the definition). Since we can repeat the following argument $j$ times, we may focus on a single configuration $\s_i$. Consider any $\o\in(m\ra\cX^s)_{opt}$ such that $\o\cap\s_i\neq\emptyset$. Since $\cW(m,\cX^s)$ is a gate for the transition from $m$ to $\cX^s$ and $\s_i\in\cS(m,\cX^s)\setminus\cW(m,\cX^s)$ for any $i=1,...,j$, we note that $\{\arg\max_{\o}H\}\setminus\{\s_i\}\neq\emptyset$. Thus our strategy consists in finding $\o'\in(m\ra\cX^s)_{opt}$ such that $\{\arg\max_{\o'}H\}\subseteq\{\arg\max_{\o}H\}\setminus\{\s_i\}$. We analyze separately the two following cases.

\medskip
\noindent
{\bf Case 1.} Suppose that the path $\o$ reaches $\cS(m,\cX^s)$ for the first time in the configuration $\s_i\in\partial\cC_{\cX^s}^{m}(\G)\cap(\cS(m,\cX^s)\setminus(\cW(m,\cX^s)\cup K))$, i.e., there exists the configuration $\s_i$ (as above) such that $\o=(m,...,\s_i,...,\h_1^{(1)},...,\cX^s)$, where $\h_1^{(1)}\in\partial\cC_{\cX^s}^{m}(\G)\cap(\cW(m,\cX^s)\cup K)$. Any such $\o$ can be written in the form 
\be{omega1} \o=(m,\o_1,...,\o_{k_1},\s_i,\o_{k_1+1},...,\o_{k_2},\g_1,\h_1^{(1)},...,\h_{m_1}^{(1)},...,\o_{k_q},...,\o_{k_{q+1}},\g_q,\h_1^{(q)},...,\h_{m_q}^{(q)})\circ\bar\o,
\ee

\noindent
where $\o_1,...,\o_{k_1},\o_{k_1+1},...,\g_1\in\cC_{\cX^s}^{m}(\G)$, $\o_{k_2+1},...,\o_{k_3},\g_2,...,\o_{k_q+1},...,\o_{k_{q+1}},\g_q\in\cX\setminus\cS(m,\cX^s)$ and $\h_i^{(j)}\in\cS(m,\cX^s)$ for all $i=1,...,m$, $j=1,...,q$. At least one among these saddles belongs to $\cW(m,\cX^s)$ and $\bar\o$ is a path that connects $\h_{m_q}^{(q)}$ to $\cX^s$ such that $\max_{\s\in\bar\o}H(\s)<\G+H(m)$. Note that $q$ and $m_1,...,m_q$ could be 1. We want to prove that $\s_i$ is unessential, thus we define a new path 
\be{omega2} 
\o'=(m,\o'_1,...,\o'_h,\g_1,\h_1^{(1)},...,\h_{m_1}^{(1)},...,\o_{k_q},...,\o_{k_{q+1}},\g_q,\h_1^{(q)},...,\h_{m_q}^{(q)})\circ \bar\o,
\ee

\noindent
where $(m,\o'_1,...,\o'_h,\g_1)$ is a path that is contained in $\cC_{\cX^s}^{m}(\G)$ such that its time-reversal exists by \cite[Lemma 2.28]{MNOS} with $\h=\g_1$ and $\cA=m$. We note that the part of $\o'$ after $\g_1$ is the same as in equation (\ref{omega1}), thus $\{\hbox{arg max}_{\o'}H\}=\{\h_1^{(1)},...,\h_{m_1}^{(1)},...,\h_1^{(q)},...,\h_{m_q}^{(q)}\}$ and therefore
\be{formulaselle}
\{\hbox{arg max}_{\o'}H\}\subseteq\{\hbox{arg max}_{\o}H\}\setminus\{\s_i\}, \quad i=1,...,n.
\ee

\noindent
This implies that the saddle $\s_i$ is unessential for any $i=1,...,n$ and thus, using \cite[Theorem 5.1]{MNOS}, $\s_i\in\cS(m,\cX^s)\setminus\cG(m,\cX^s)$.

\medskip
\noindent
{\bf Case 2.} The path $\o$ reaches the set $\cS(m,\cX^s)$ before reaching $\partial\cC_{\cX^s}^{m}(\G)\cap(\cS(m,\cX^s)\setminus(\cW(m,\cX^s)\cup K))$ in $\s_i$. In this case we can bypass the saddle $\s_i$ by arguing in a similar way as in case 1, indeed we can write 
\be{}
\o=(m,\o_1,...,\o_{k_1},\g_1,\h_1^{(1)},...,\h_{m_1}^{(1)},...,\s_i,...,\g_t,\h_1^{(t)},...,\h_{m_t}^{(t)},...,\o_{k_{q+1}},\g_q,\h_1^{(q)},...,\h_{m_q}^{(q)})\circ\bar\o
\ee

\noindent
and define
\be{}
\o'=(m,\o'_1,...,\o'_h,\g_t,\h_1^{(t)},...,\h_{m_t}^{(t)},...,\o_{k_{q+1}},\g_q,\h_1^{(q)},...,\h_{m_q}^{(q)})\circ\bar\o,
\ee

\noindent
where $(m,\o'_1,...,\o'_h,\g_t)$ is is a path that is contained in $\cC_{\cX^s}^{m}(\G)$ such that its time-reversal exists by \cite[Lemma 2.28]{MNOS} with $\h=\g_t$ and $\cA=m$. Thus $\{\hbox{arg max}_{\o'}H\}=\{\h_1^{(t)},...,\h_{m_t}^{(t)},...,\h_1^{(q)},...,\h_{m_q}^{(q)}\}$ and therefore (\ref{formulaselle}) holds.
\epr

\subsection{Proof of Proposition \ref{selle2}}
\label{proofdep2}

\bpr
We denote by $\z_1,...,\z_l$ the saddles in the statement. We want to prove that these saddles are unessential (see Section \ref{modinddef} point 4 for the definition). Since we can repeat the following argument $l$ times, we may focus on a single configuration $\z_i$. Consider any $\o\in(m\ra\cX^s)_{opt}$ such that $\o\cap\z_i\neq\emptyset$. Since $\cW(m,\cX^s)$ is a gate for the transition from $m$ to $\cX^s$ and $\z_i\in\cS(m,\cX^s)\setminus\cW(m,\cX^s)$ for any $i=1,...,j$, we note that $\{\arg\max_{\o}H\}\setminus\{\z_i\}\neq\emptyset$. Thus our strategy consists in finding $\o'\in(m\ra\cX^s)_{opt}$ such that $\{\arg\max_{\o'}H\}\subseteq\{\arg\max_{\o}H\}\setminus\{\z_i\}$. Due to Proposition \ref{selle1}, we can reduce the proof to consider any $\o\in(m\ra\cX^s)_{opt}$ such that the first saddle that is visited is $\h_1^{(1)}\in\partial\cC_{\cX^s}^{m}(\G)\cap(\cW(m,\cX^s)\cup K)$. Note that there exists $\h_{m_q}^{(q)}\in\partial\cC_{m}^{\cX^s}(\G+H(m)-H(\cX^s))\cap(\cW(m,\cX^s)\cup\widetilde{K})$, different from $\z_i$, that can be connected to the set $\cL^G$ via one step of the dynamics. By the model-dependent input (iii) we deduce that $\z_i$ can be reached either after visiting the set $\cL^G$ 
\be{om1} \o=(m,\o_1,...,\o_{k_1},\g_1,\h_1^{(1)},...,\h_{m_1}^{(1)},...,\o_{k_q},...,\o_{k_{q+1}},\g_q,\h_1^{(q)},...,\h_{m_q}^{(q)},\h^G,...,\z_i)\circ\widetilde\o
\ee

\noindent
or directly from $\h_{m_q}^{(q)}$
\be{om2} \o=(m,\o_1,...,\o_{k_1},\g_1,\h_1^{(1)},...,\h_{m_1}^{(1)},...,\o_{k_q},...,\o_{k_{q+1}},\g_q,\h_1^{(q)},...,\h_{m_q}^{(q)},\z_i,...,\bar\h^{G})\circ\bar\o,
\ee

\noindent
where $\o_1,...,\o_{k_1},\g_1\in\cC_{\cX^s}^{m}(\G)$ and $\o_{k_1+1},...,\o_{k_2},\g_2,...,\o_{k_q+1},...,\o_{k_{q+1}}\in\cX\setminus\cS(m,\cX^s)$. Additionally, the configurations $\h_i^{(j)}\in\cS(m,\cX^s)$ for all $i=1,...,m$, $j=1,...,q$ and at least one among these saddles belongs to $\cW(m,\cX^s)$. Moreover, $\h^G$ (resp.\ $\bar\h^G$) is in $\cL^G$ and $\widetilde\o$ (resp.\ $\bar\o$) is a path that connects $\z_i$ (resp.\ $\bar\h^G$) to $\cX^s$. Note that $q$ and $m_1,...,m_q$ could be 1. We want to prove that $\z_i$ is unessential, thus for both $\o$ in (\ref{om1}) and (\ref{om2}) we define a new path
\be{} \o'=(m,\o_1,...,\o_{k_1},\g_1,\h_1^{(1)},...,\h_{m_1}^{(1)},...,\o_{k_q},...,\o_{k_{q+1}},\g_q,\h_1^{(q)},...,\h_{m_q}^{(q)},\h^G)\circ\hat\o, 
\ee

\noindent
where by the model-dependent input (iii)-(a) there exists a path $\hat\o$ that connects $\h^G$ to $\cX^s$ such that
\be{}
\max_{\s\in\hat\o}H(\s)<\G+H(m).
\ee

\noindent
Thus $\{\hbox{arg max}_{\o'}H\}=\{\h_1^{(1)},...,\h_{m_1}^{(1)},...,\h_1^{(q)},...,\h_{m_q}^{(q)}\}$ and therefore
\be{}
\{\hbox{arg max}_{\o'}H\}\subseteq\{\hbox{arg max}_{\o}H\}\setminus\{\z_i\}, \quad i=1,...,n.
\ee

\noindent
This implies that the saddle $\z_i$ is unessential for any $i=1,...,n$ and thus, using \cite[Theorem 5.1]{MNOS}, $\z_i\in\cS(m,\cX^s)\setminus\cG(m,\cX^s)$.
\epr

\section{Useful model-dependent definitions and tools}
\label{dependentdef}
In this Section (and only here) we set 
\be{definizioneqbarsa}
\cQ_{wa}=\bar\cQ_{wa}, \qquad \cD_{wa}=\bar\cD_{wa}
\ee
\noindent
to show the similarites between the results with the isotropic model.

\subsection{Geometric description of the protocritical droplets}
In \cite[Theorem 1.4.1]{BHN} the authors obtain the geometric description of the set $\cD_{is}$ as $\cD_{is}=\bar\cD_{is}\cup\widetilde\cD_{is}$. In this Section we derive the geometric description of the analogous sets for the weakly anistropic models $\bar\cD_{wa}$ and $\widetilde\cD_{wa}$ following the argument proposed in \cite{BHN}. The geometric description of the sets $\bar\cD_{sa}$ and $\widetilde\cD_{sa}$ is given in \cite[Proposition 5.1]{BN3}. Recall definition (\ref{defDwa}).

\bp{cardweak}{\rm{(Geometric description of $\widetilde{\cD}_{wa}$ and $\bar\cD_{wa}$).}}
We obtain the following geometric description of $\bar\cD_{wa}$ and $\widetilde\cD_{wa}$:
\bi
\item[(a)] $\bar\cD_{wa}$ is the set of configurations having one cluster $\h$ anywhere in $\L_0$ consisting of a $(l_1^*-2)\times(l_2^*-2)$ rectangle with four bars $B^\a(\h)$, with $\a\in\{n,w,e,s\}$, attached to its four sides satisfying 
\be{condbardweak2}
1\leq|B^w(\h)|,|B^e(\h)|\leq l_2^*, \qquad l_1^*-l_2^*+1\leq |B^n(\h)|,|B^s(\h)|\leq l_1^*,
\ee

\noindent
and
\be{condbarweak}
\displaystyle\sum_\a|B^\a(\h)|-\displaystyle\sum_{\a\a'\in\{nw,ne,sw,se\}}|c^{\a\a'}(\h)|=2l_1^*+l_2^*-3.
\ee

\item[(b)] $\widetilde\cD_{wa}$ is the set of configurations having one cluster $\h$ anywhere in $\L_0$ consisting of a $(l_1^*-3)\times(l_2^*-1)$ rectangle with four bars $B^\a(\h)$, with $\a\in\{n,w,e,s\}$, attached to its four sides satisfying 
\be{condbardweak1}
1\leq|B^w(\h)|,|B^e(\h)|\leq l_2^*+1, \qquad 1\leq |B^n(\h)|,|B^s(\h)|\leq l_1^*-1,
\ee

\noindent
and
\be{}
\displaystyle\sum_\a|B^\a(\h)|-\displaystyle\sum_{\a\a'\in\{nw,ne,sw,se\}}|c^{\a\a'}(\h)|=l_1^*+2l_2^*-2.
\ee
\ei
\ep

\br{barreweak}
Let $\h\in\bar\cD_{wa}$. 

\bi
\item[(i)] Note that (\ref{condbarweak}) takes into account the number of occupied unit squares in $\partial^-\hbox{CR}(\h)$ due to Remark \ref{sommabarre}. We deduce that at most three frame-angles of CR$(\h)$ can be occupied, otherwise $|\partial^-\hbox{CR}(\h)|=2l_1^*+2l_2^*-4>2l_1^*+l_2^*-3$, which is absurd. 

\item[(ii)] Since $|B^{s}(\h)|+|B^{w}(\h)|\leq l_1^*+l_2^*-4+k-|c^{ne}(\h)|$, we get
\be{} 
|B^{n}(\h)|+|B^{e}(\h)|=2l_1^*+l_2^*-3-(|B^{s}(\h)|+|B^{w}(\h)|)+k\geq l_1^*+1+|c^{ne}(\h)|.
\ee

\noindent
By symmetry, we generalize the inequality above for any $\a\in\{n,s\}$ and $\a'\in\{w,e\}$: we get $|B^{\a}(\h)|+|B^{\a'}(\h)|\geq l_1^*+1+|c^{\a\a'}(\h)|$.
\ei
\er

\begin{proof*}{\bf of Proposition \ref{cardweak}}
	(a) We denote by $\bar\cD_{wa}^{geo}$ the geometric set with the properties specified in point (a) that we introduce to make the argument more clear. The proof will be given in two steps:
	
	\bi
	\item[(i)] $\bar\cD_{wa}^{geo}\subseteq\bar\cD_{wa}$;
	\item [(ii)] $\bar\cD_{wa}^{geo}\supseteq\bar\cD_{wa}$.
	\ei
	
	\noindent
	{\bf Proof of (i).} To prove (i) we must show that for all $\h\in\bar\cD_{wa}^{geo}$,
	
	\bi
	\item[(i1)] $H(\h)=H(\bar\cQ_{wa})$;
	\item[(i2)] $\exists \ \o:\bar\cQ_{wa}\ra\h$, i.e., $\o=(\o_1,...,\o_k=\h)$, such that $\displaystyle\max_i H(\o_i)\leq H(\bar\cQ_{wa})+U_1$, with $|\o_i|=n_{wa}^c$ for all $i=1,...,k$ and $\o_1\in\bar\cQ_{wa}$ (see (\ref{defnwa}) for the definition of $n_{wa}^c$).
	\ei

	{\bf Proof of (i1).} Any $\h\in\bar\cD_{wa}^{geo}$ satisfies $n(\h)=0$, $|C(\h_{cl})|=(l_1^*-2)(l_2^*-2)+2l_1^*+l_2^*-3=n_{wa}^c$, and $g_1(\h)=l_1^*$ and $g_2(\h)=l_2^*$ since the configuration is monotone. Thus by (\ref{Hcont}) we deduce that $H$ is constant on $\bar\cD_{wa}^{geo}$. Since $\bar\cQ_{wa}\subseteq\bar\cD_{wa}^{geo}$, this completes the proof of (i1).
	
	\medskip
	{\bf Proof of (i2).} Consider $\z\in\bar\cQ_{wa}$ and $\h\in\bar\cD_{wa}^{geo}$. Here, without loss of generality, we assume that the protuberance is in $r^w(\z)$. Then we have
	
	\bi
	\item[-] $|B^w(\z)|=1$;
	\item[-] $|B^n(\z)|=|B^s(\z)|=l_1^*-1$;
	\item[-] $|B^e(\z)|=l_2^*$;
	\item[-] $|c^{ne}(\z)|=|c^{se}(\z)|=1$.
	\ei
	
	\noindent
	Using the sliding of a unit square around a frame-angle described in Figure \ref{fig:trenino} (see Definition \ref{movepart}), we move, one by one, $|B^n(\z)|-|B^n(\h)|$ particles from around the frame-angle $c^{nw}(\z)$. After that we move $|B^e(\z)|-|B^e(\h)|+|B^s(\z)|-|B^s(\h)|$ particles around the frame-angle $c^{sw}(\z)$. Finally, we move $|B^e(\z)|-|B^e(\h)|$ particles around the frame-angle $c^{es}(\z)$. The result is the configuration $\h\in\bar\cD_{wa}^{geo}$. This concludes the proof of (i2).
	
	\medskip
	\noindent
	{\bf Proof of (ii).} By (i2), we know that all configurations in $\bar\cD_{wa}^{geo}$ are connected via $U_1$-path to $\bar\cQ_{wa}$. Since $\bar\cQ_{wa}\subseteq\bar\cD_{wa}\cap\bar\cD_{wa}^{geo}$, in order to prove (ii) it suffices to show that following $U_1$-paths it is not possible to exit $\bar\cD_{wa}^{geo}$. We call a path {\it clustering} if all the configurations in the path consist of a single cluster and no free particles. Below we will prove that for any $\h\in\bar\cD_{wa}^{geo}$ and any $\h'$ connected to $\h$ by a clustering $U_1$-path,
	
	\bi
	\item[(A)] $\hbox{CR}(\h')=\hbox{CR}(\h)$;
	\item[(B)] $\h'\supseteq\hbox{CR}^-(\h)$.
	\ei
	
	\medskip
	{\bf Proof of (A).} Starting from any $\h\in\cX$, it is geometrically impossible to modify $\hbox{CR}(\h)$ without detaching a particle, that contradicts the hypotheses of clustering $U_1$-path.
	
	\medskip
	{\bf Proof of (B).} Fix $\h\in\bar\cD_{wa}^{geo}$. The proof is done in two steps.
	
	\medskip
	\noindent
	{\bf 1.} First, we consider clustering $U_1$-paths along which we do not move a particle from $\hbox{CR}^-(\h)$. Along such paths we only encounter configurations in $\bar\cD_{wa}^{geo}$ or those obtained from $\bar\cD_{wa}^{geo}$ by breaking one of the bars in $\partial^-\hbox{CR}(\h)$ into two pieces at cost $U_1$ (resp.\ $U_2$) if the bar is horizontal (resp.\ vertical). This holds because there is no particles outside $\hbox{CR}(\h)$ that can lower the cost. 
	
	If the broken bar is horizontal, then only moves at zero cost are admissible, so any particle can be detached. This implies that the unique way to regain $U_1$ and complete the $U_1$-path is to restore the bar.
	
	If the broken bar is vertical, then the admissible moves are those with cost less or equal than $U_1-U_2$. Again any particle can be detached, indeed its cost is at least $U_1$. The moves at cost $U_2$ are not possible, since $U_1<2U_2-2\e$. Thus the unique way to complete the $U_1$-path is to restore the broken bar. Thus we have proved that $\h'\supseteq\hbox{CR}^-(\h)$.

		\begin{figure}
		\centering
		\begin{tikzpicture}[scale=0.4,transform shape]

			\draw (0,0)--(4,0);
			\draw (0,0)--(0,1);
			\draw (0,1)--(-1,1);
			\draw (-1,1)--(-1,4);
			\draw(-1,4)--(0,4);
			\draw (0,4)--(0,8);
			\draw (0,8)--(1,8);
			\draw(1,8)--(1,9);
			\draw (1,9)--(3,9);
			\draw (3,9)--(3,8);
			\draw (3,8)--(5,8);
			\draw[->] (4,5)--(6.5,5);
			\node at (5,4.4) {{\Huge{$+U_2$}}};

			\draw (8,0)--(12,0);
			\draw (8,0)--(8,1);
			\draw (8,1)--(7,1);
			\draw (7,1)--(7,4);
			\draw(7,4)--(8,4);
			\draw (8,4)--(8,7);
			\draw (8,7)--(9,7);
			\draw (9,7)--(9,8);
			\draw (9,8)--(8,8);
			\draw(8,8)--(8,9);
			\draw (8,9)--(11,9);
			\draw (11,9)--(11,8);
			\draw (11,8)--(13,8);
			\draw[->] (12,5)--(14.5,5);
			\node at (13.3,4.5) {{\Huge{0}}};

			\draw (16,0)--(20,0);
			\draw (16,0)--(16,1);
			\draw (16,1)--(15,1);
			\draw (15,1)--(15,4);
			\draw(15,4)--(16,4);
			\draw (16,4)--(16,6);
			\draw (16,6)--(17,6);
			\draw (17,6)--(17,7);
			\draw (17,7)--(16,7);
			\draw(16,7)--(16,9);
			\draw (16,9)--(19,9);
			\draw (19,9)--(19,8);
			\draw (19,8)--(21,8);
			\draw[->] (20,5)--(22.5,5);
			\node at (21.3,4.5) {{\Huge{0}}};

			\draw (24,0)--(28,0);
			\draw (24,0)--(24,1);
			\draw (24,1)--(23,1);
			\draw (23,1)--(23,4);
			\draw(23,4)--(24,4);
			\draw (24,4)--(24,5);
			\draw (24,5)--(25,5);
			\draw (25,5)--(25,6);
			\draw (25,6)--(24,6);
			\draw(24,6)--(24,9);
			\draw (24,9)--(27,9);
			\draw (27,9)--(27,8);
			\draw (27,8)--(29,8);
			\draw[->] (28,5)--(30.5,5);
			\node at (29.3,4.5) {{\Huge{0}}};

			\draw (32,0)--(36,0);
			\draw (32,0)--(32,1);
			\draw (32,1)--(31,1);
			\draw (31,1)--(31,4);
			\draw(31,4)--(33,4);
			\draw (33,4)--(33,5);
			\draw (33,5)--(32,5);
			\draw (32,5)--(32,9);
			\draw (32,9)--(35,9);
			\draw (35,9)--(35,8);
			\draw (35,8)--(37,8);

		\end{tikzpicture}
		
		\vskip 0 cm
		\caption{Creation and motion of the recess at cost 0.}
		\label{fig:recess}
	\end{figure}
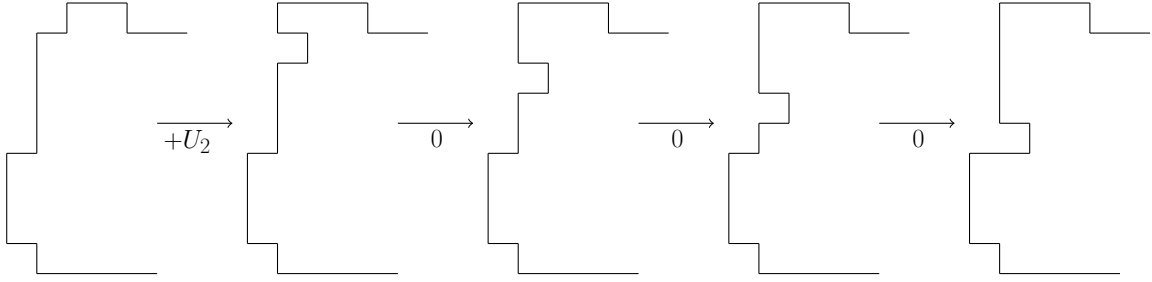
	
	\medskip
	\noindent
	{\bf 2.}  Consider now clustering $U_1$-paths along which we move a particle from a corner of $\hbox{CR}^-(\h)$. It is not allowed to move at cost $U_1+U_2$, because it exceeds $U_1$, thus the overshoot $U_2$ must be regained by letting the particle slide next to a bar that is attached to a side of $\hbox{CR}^-(\h)$ (see Figure \ref{fig:recess}). If the particle moves vertically (resp.\ horizontally), we regain $U_1$ (resp.\ $U_2$). Since there are never two bars attached to the same side, we can at most regain $U_1$, thus it is not possible to move a particle from $\hbox{CR}^-(\h)$ other than from a corner. From now on, since $U_2<2U_1-2\e$, only moves at cost at most zero are admissible. There are no protuberances present anymore, because only the configurations in $\bar\cQ_{wa}$ have a protuberance. Thus no particle outside $\hbox{CR}^-(\h)$ can move, except those just detached from $\hbox{CR}^-(\h)$. These particles can move back, in which case we return to the same cofiguration $\h$ (see Figure \ref{fig:recess}). In fact, all possible moves at zero cost consist in moving the recess just created in $\hbox{CR}^-(\h)$ along the same side of $\hbox{CR}^-(\h)$, until it reaches the top of the bar, after which it cannot advance anymore at zero cost (see Figure \ref{fig:recess}). All these moves do not change the energy, except the last one that returns the particle to its original position and regains $U_1$. This concludes the proof of (B).
	
	From (A), we deduce that $\hbox{CR}(\h')=\cR(l_1^*,l_2^*)$. From (A) and (B), we deduce that the number of particles that are in $\partial^-\hbox{CR}(\h)$ is equal to the number of particles that are in $\partial^-\hbox{CR}(\h')$, thus (\ref{condbarweak}), $1\leq|B^w(\h')|,|B^e(\h')|\leq l_2^*$ and $1\leq |B^n(\h')|,|B^s(\h')|\leq l_1^*$ hold. In order to prove that following clustering $U_1$-paths it is not possible to exit $\bar\cD_{wa}^{geo}$, we have to prove the lower bound in (\ref{condbardweak2}) for the lengths $|B^n(\h')|$ and $|B^s(\h')|$. We set
	\be{} 
	k=\displaystyle\sum_{\a\a'\in\{nw,ne,sw,se\}}|c^{\a\a'}(\h')|. 
	\ee
	
	\noindent
	Since $|B^w(\h')|+|B^e(\h')|\leq2l_2^*-4+k$, by (\ref{condbarweak}) we get
	\be{contobarreweak}
	|B^n(\h')|+|B^s(\h')|=2l_1^*+l_2^*-3-(|B^w(\h')|+|B^e(\h')|)+k\geq2l_1^*-l_2^*+1.
	\ee
	
	\noindent
	Since $|B^s(\h')|\leq l_1^*$, (\ref{contobarreweak}) implies
	\be{}
	|B^n(\h')|\geq2l_1^*-l_2^*+1-|B^s(\h')|\geq l_1^*-l_2^*+1.
	\ee
	
	\noindent
	By symmetry we can similarly argue for the length $|B^s(\h')|$. This implies that following $U_1$-paths it is not possible to exit $\bar\cD_{wa}^{geo}$.
	The argument goes as follows. Detaching a particle costs at least $U_1+U_2$ unless the particle is a protuberance, in which case the cost is $U_1$. The only configurations in $\bar\cD_{wa}^{geo}$ having a protuberance are those in $\bar\cQ_{wa}$. If we detach the protuberance from a configuration in $\bar\cQ_{wa}$, then we obtain a $(l_1^*-1)\times l_2^*$ rectangle with a free particle. Since in the sequel only moves at zero cost are allowed, it is only possible to move the free particle. Since in a $U_1$-path the particle number is conserved, the only way to regain $U_1$ and complete the $U_1$-path is to reattach the free particle to a vertical side of the rectangle, thus return to $\bar\cQ_{wa}$. This implies that for any $\h\in\bar\cD_{wa}^{geo}$ and any $\h'$ connected to $\h$ by a $U_1$-path we must have that $\h'\in\bar\cD_{wa}^{geo}$. This concludes the proof.
	
	(b) The proof is analogue to the one in (a).
\end{proof*}

\subsection{Definitions}
\label{sitigood}

Since we are considering the isotropic and weakly anisotropic models and some properties are in common with the strongly anisotropic model, we choose the lower index $int\in\{is,wa,sa\}$ to make clear in the notation which of the three models we are referring to. We set
\be{L*}
L_{int}^*:=
\begin{cases}
	L-l_c & \hbox{if } int=is, \\
	L-l_2^* & \hbox{if } int=wa.
\end{cases}
\ee

\noindent
Let $int\in\{is,wa\}$. For $\h\in\cC_{int}^*$, we associate $(\hat\h,x)$ with $\hat\h\in\cD_{int}$ protocritical droplet and $x\in\L$ the position of the free particle. We denote by $\cC_{int}^G(\hat\h)$ (resp.\ $\cC_{int}^B(\hat\h)$) the configurations that can be reached from $(\hat\h,x)$ by a path that moves the free particle towards the cluster and attaches the particle in $\partial^-CR(\hat\h)$ (resp.\ $\partial^+CR(\hat\h)$). In Figure \ref{fig:sitiGB} on the left-hand side we depict explicitly the good and bad sites for a specific $\hat\h$. Let 
\be{sitiGB}
\cC_{int}^G=\displaystyle\bigcup_{\hat\h\in\cD_{int}}\cC_{int}^G(\hat\h), \qquad \cC_{int}^B=\displaystyle\bigcup_{\hat\h\in\cD_{int}}\cC_{int}^B(\hat\h).
\ee

\noindent
For $\h\in\cC_{int}^*$, let $\hat\h\in\cD_{int}$ be the configuration obtained from $\h$ by removing the free particle. For $A\subseteq\L$ and $x\in\L$, recall that $d(x,A)$ denotes the lattice distance between $x$ and $A$. As in \cite[Section 3.5]{BHN}, we need the following definitions.

\bd{cornis}
Let $\L_4$ be $\L$ without its four frame-angles. We define, recursively, 
\be{}
B_1(\hat\h):=\{x\in\L_4| \ x\notin\hat\h, \ d(x,\hat\h)=1\}
\ee
and
\be{} 
\ba{lll} 
B_2(\hat\h):=\{x\in\L_4| \ x\notin\hat\h, \ d(x,B_1(\hat\h))=1\},  \\
\bar B_2(\hat\h):=B_2(\hat\h),
\ea
\ee
and
\be{}
\ba{lll}
B_3(\hat\h):=\{x\in\L_4| \ x\notin B_1(\hat\h), \ d(x,B_2(\hat\h))=1\},  \\
\bar B_3(\hat\h):=B_3(\hat\h)\cup \{\bar B_2(\hat\h)\cap \partial^-\L_4\},
\ea
\ee

\setlength{\unitlength}{0.2cm}

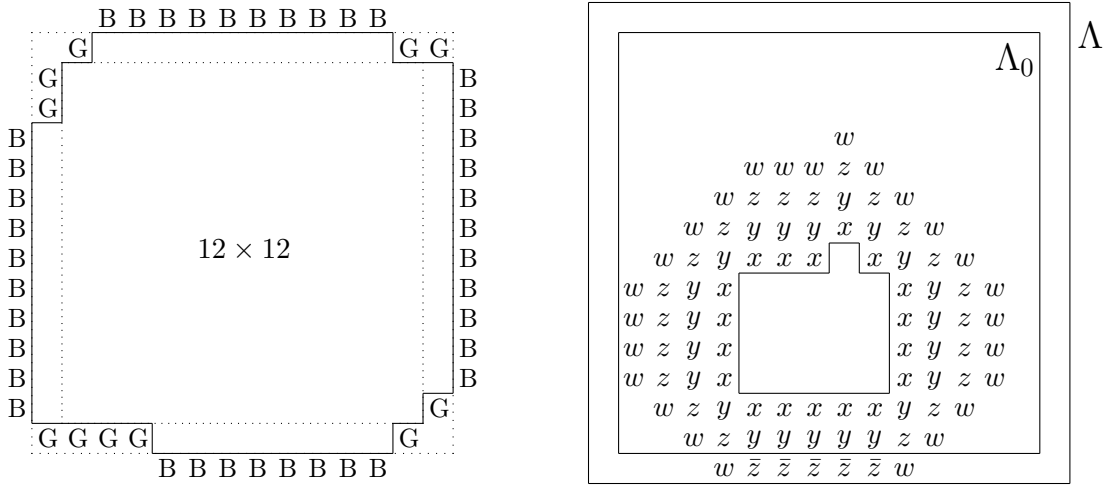
\begin{figure}
	\begin{picture}(-2,30)(-6,0)
		
		
		
		\thinlines
		\qbezier[51](0,0)(14,0)(28,0)
		\qbezier[51](28,0)(28,14)(28,28)
		\qbezier[51](0,28)(14,28)(28,28)
		\qbezier[51](0,0)(0,14)(0,28)
		
		\thinlines
		\qbezier[51](2,2)(13,2)(26,2)
		\qbezier[51](26,2)(26,13)(26,26)
		\qbezier[51](2,26)(13,26)(26,26)
		\qbezier[51](2,2)(2,13)(2,26)
		
		\put(0,2){\line(1,0){8}}
		\put(8,0){\line(1,0){16}}
		\put(24,2){\line(1,0){2}}
		\put(26,4){\line(1,0){2}}
		\put(24,26){\line(1,0){4}}
		\put(4,28){\line(1,0){20}}
		\put(2,26){\line(1,0){2}}
		\put(0,22){\line(1,0){2}}
		
		\put(8,0){\line(0,1){2}}
		\put(24,0){\line(0,1){2}}
		\put(26,2){\line(0,1){2}}
		\put(28,4){\line(0,1){22}}
		\put(24,26){\line(0,1){2}}
		\put(4,26){\line(0,1){2}}
		\put(2,22){\line(0,1){4}}
		\put(0,2){\line(0,1){20}}
		
		\put(0.4,0.4){\small G}
		\put(2.4,0.4){\small G}
		\put(4.4,0.4){\small G}
		\put(6.4,0.4){\small G}
		
		\put(24.4,0.4){\small G}
		\put(26.4,2.4){\small G}
		\put(24.4,26.4){\small G}
		\put(26.4,26.4){\small G}
		\put(0.4,22.4){\small G}
		\put(0.4,24.4){\small G} 
		\put(2.4,26.4){\small G}
		
		\put(8.4,-1.6){\small B}
		\put(10.4,-1.6){\small B}
		\put(12.4,-1.6){\small B}
		\put(14.4,-1.6){\small B}
		\put(16.4,-1.6){\small B}
		\put(18.4,-1.6){\small B}
		\put(20.4,-1.6){\small B}
		\put(22.4,-1.6){\small B}
		
		\put(28.4,4.4){\small B}
		\put(28.4,6.4){\small B}
		\put(28.4,8.4){\small B}
		\put(28.4,10.4){\small B}
		\put(28.4,12.4){\small B}
		\put(28.4,14.4){\small B}
		\put(28.4,16.4){\small B}
		\put(28.4,18.4){\small B}
		\put(28.4,20.4){\small B}
		\put(28.4,22.4){\small B}
		\put(28.4,24.4){\small B}
		\put(4.4,28.4){\small B}
		\put(6.4,28.4){\small B}
		\put(8.4,28.4){\small B}
		\put(10.4,28.4){\small B}
		\put(12.4,28.4){\small B}
		\put(14.4,28.4){\small B}
		\put(16.4,28.4){\small B}
		\put(18.4,28.4){\small B}
		\put(20.4,28.4){\small B}
		\put(22.4,28.4){\small B}
		\put(-1.6,2.4){\small B}
		\put(-1.6,4.4){\small B}
		\put(-1.6,6.4){\small B}
		\put(-1.6,8.4){\small B}
		\put(-1.6,10.4){\small B}
		\put(-1.6,12.4){\small B}
		\put(-1.6,14.4){\small B}
		\put(-1.6,16.4){\small B}
		\put(-1.6,18.4){\small B}
		\put(-1.6,20.4){\small B}
		
		\put(11,13){$12 \times 12$}

		\put(37,-2){\line(1,0){32}}
		\put(37,-2){\line(0,1){32}}
		\put(69,-2){\line(0,1){32}}
		\put(37,30){\line(1,0){32}}
		\put(39,0){\line(1,0){28}}
		\put(39,0){\line(0,1){28}}
		\put(39,28){\line(1,0){28}}
		\put(67,28){\line(0,-1){28}}
		\put(69.5,27){\Large{$\L$}}
		\put(64,25.5){\Large{$\L_0$}}
		\put(47,4){\line(1,0){10}}
		\put(47,4){\line(0,1){8}}
		\put(47,12){\line(1,0){6}}
		\put(53,12){\line(0,1){2}}
		\put(53,14){\line(1,0){2}}
		\put(55,12){\line(0,1){2}}
		\put(55,12){\line(1,0){2}}
		\put(57,12){\line(0,-1){8}}
		\put(57.5,10.5){$x$}
		\put(57.5,8.5){$x$}
		\put(57.5,6.5){$x$}
		\put(57.5,4.5){$x$}
		\put(55.5,2.5){$x$}
		\put(53.5,2.5){$x$}
		\put(51.5,2.5){$x$}
		\put(49.5,2.5){$x$}
		\put(47.5,2.5){$x$}
		\put(45.5,4.5){$x$}
		\put(45.5,6.5){$x$}
		\put(45.5,8.5){$x$}
		\put(45.5,10.5){$x$}
		\put(47.5,12.5){$x$}
		\put(49.5,12.5){$x$}
		\put(51.5,12.5){$x$}
		\put(53.5,14.5){$x$}
		\put(55.5,12.5){$x$}
		\put(57.5,12.6){$y$}
		\put(59.5,10.6){$y$}
		\put(59.5,8.6){$y$}
		\put(59.5,6.6){$y$}
		\put(59.5,4.6){$y$}
		\put(57.5,2.6){$y$}
		\put(55.5,0.7){$y$}
		\put(53.5,0.7){$y$}
		\put(51.5,0.7){$y$}
		\put(49.5,0.7){$y$}
		\put(47.5,0.7){$y$}
		\put(45.5,2.6){$y$}
		\put(43.5,4.6){$y$}
		\put(43.5,6.6){$y$}
		\put(43.5,8.6){$y$}
		\put(43.5,10.6){$y$}
		\put(45.5,12.6){$y$}
		\put(47.5,14.5){$y$}
		\put(49.5,14.5){$y$}
		\put(51.5,14.5){$y$}
		\put(53.5,16.5){$y$}
		\put(55.5,14.5){$y$}
		\put(57.5,14.5){$z$}
		\put(59.5,12.5){$z$}
		\put(61.5,10.5){$z$}
		\put(61.5,8.5){$z$}
		\put(61.5,6.5){$z$}
		\put(61.5,4.5){$z$}
		\put(59.5,2.5){$z$}
		\put(57.5,0.5){$z$}
		\put(55.5,-1.5){$\bar z$}
		\put(53.5,-1.5){$\bar z$}
		\put(51.5,-1.5){$\bar z$}
		\put(49.5,-1.5){$\bar z$}
		\put(47.5,-1.5){$\bar z$}
		\put(45.5,0.5){$z$}
		\put(43.5,2.5){$z$}
		\put(41.5,4.5){$z$}
		\put(41.5,6.5){$z$}
		\put(41.5,8.5){$z$}
		\put(41.5,10.5){$z$}
		\put(43.5,12.5){$z$}
		\put(45.5,14.5){$z$}
		\put(47.5,16.5){$z$}
		\put(49.5,16.5){$z$}
		\put(51.5,16.5){$z$}
		\put(53.5,18.5){$z$}
		\put(55.5,16.5){$z$}
		\put(57.3,16.5){$w$}
		\put(59.3,14.5){$w$}
		\put(61.3,12.5){$w$}
		\put(63.3,10.5){$w$}
		\put(63.3,8.5){$w$}
		\put(63.3,6.5){$w$}
		\put(63.3,4.5){$w$}
		\put(61.3,2.5){$w$}
		\put(59.3,0.5){$w$}
		\put(57.3,-1.5){$w$}
		\put(45.3,-1.5){$w$}
		\put(43.3,0.5){$w$}
		\put(41.3,2.5){$w$}
		\put(39.3,4.5){$w$}
		\put(39.3,6.5){$w$}
		\put(39.3,8.5){$w$}
		\put(39.3,10.5){$w$}
		\put(41.3,12.5){$w$}
		\put(43.3,14.5){$w$}
		\put(45.3,16.5){$w$}
		\put(47.3,18.5){$w$}
		\put(49.3,18.5){$w$}
		\put(51.3,18.5){$w$}
		\put(53.3,20.5){$w$}
		\put(55.3,18.5){$w$}

	\end{picture}
	\vskip 0.5 cm
	\caption{On the left-hand side we represent good sites ($G$) and bad sites ($B$) for $\ell_c=14$ in the isotropic case. On the right-hand side we depict with $x$ the sites in $B_1(\hat\h)$, with $y$ the sites in $\bar{B}_2(\hat\h)$, with $z$ and $\bar z$ the sites in $\bar{B}_3(\hat\h)$ and with $\bar z$ and $w$ the sites in $\bar{B}_4(\hat\h)$.}
	\label{fig:sitiGB} 
\end{figure}

and for $i=4,5,...,L_{int}^*$
\be{}
\ba{lll}
B_i(\hat\h):=\{x\in\L_4| \ x\notin B_{i-2}(\hat\h), \ d(x,B_{i-1}(\hat\h))=1\},  \\
\bar B_i(\hat\h):=B_i(\hat\h)\cup \{\bar B_{i-1}(\hat\h)\cap \partial^-\L_4\}.
\ea
\ee
\ed

In words, $B_1(\bar\h)$ is the ring of sites in $\L_4$ at distance $1$ from $\hat\h$, while $\bar B_i(\hat\h)$ is the ring of sites in $\L_4$ at distance $i$ from $\hat\h$ union all the sites in $\partial^-\L_4$ at distance $1<j<i$ from $\hat\h$ ($i=2,3,...,L_{int}^*$) (see Figure \ref{fig:sitiGB} on the right-hand side). Note that, depending on the location of $\hat\h$ in $\L$, the $\bar B_i(\hat\h)$ coincide for large enough $i$. The maximal number of rings is $L_{int}^*$.  

Now we need to introduce specific sets that will be crucial later on.

\bd{c*i}
Let $int\in\{is,wa\}$. We define
\be{}
\cC_{int}^*(i):=\{(\hat\h,x): \ \hat\h\in\cD_{int}, \ x\in\bar B_i(\hat\h)\}, \quad i=2,3,...,L_{int}^*.
\ee
\ed

\noindent
First, note that the sets $\cC_{int}^*(i)$ are not disjoint. 

\br{}
From the definitions of the sets $\cigeo$ and $\cwgeo$, for any $int\in\{is,wa\}$ we deduce that
\be{}
\cC_{int}^*=\displaystyle\bigcup_{i=2}^{L_{int}^*}\cC_{int}^*(i).
\ee
\er

\noindent
For this discussion in the case $int=sa$ we refer to \cite[Section 5.2]{BN3}.

\subsection{Useful lemmas for the gates}
\label{lemmi3modelli}
In this Section we give some useful lemmas that help us to characterize the gates.

\subsubsection{Lemmas valid for the two models}

Here we state Lemma \ref{trenini} for the case $int\in\{is,wa\}$, but it holds also for the case $int=sa$ (see \cite[Lemma 5.6]{BN3}).

\bl{trenini}
Let $int\in\{is,wa\}$. Starting from $\cC_{int}^*\setminus\cQ_{int}^{fp}$, if the free particle is attached to a bad site obtaining $\h^B\in\cC^B_{int}$, the only transitions that does not exceed the energy $\G_{int}^*$ are either detaching the protuberance, or a sequence of $1$-translations of a bar or slidings of a bar around a frame-angle. Moreover, we get:

\bi
\item[(i)] if it is possible to slide a bar around a frame-angle, then the saddles that are crossed are essential; 
\item[(ii)] if it is not possible to slide a bar around a frame-angle, then the path must come back to the starting configuration and the saddles that are crossed are unessential.
\ei
\el

\bpr
Let $int\in\{is,wa\}$. Let $\h^B\in\cC^B_{int}$ the configuration obtained by attaching the free particle as a protuberance to a bar, thus $H(\h^B)=\gi-U$ if $int=is$ and either $H(\h^B)=\G_{int}^*-U_1$ or $H(\h^B)=\G_{int}^*-U_2$ if $int=wa$. Note that it is impossible to move particles in $\partial^-\hbox{CR}(\h^B)$ before further lowering the energy, since this move costs at least $2U$ if $int=is$ and $U_1+U_2$ if $int=wa$. Moreover, it is impossible to create a new particle before further lowering the energy, since this move costs $\D$. On the other hand there are no moves available to lower the energy. If the protuberance is detached, then the energy reaches the value $\G_{int}^*$. Analyzing motions of particles along the border of the droplet (both sequence of $1$-translations of a bar and sliding around a frame-angle), if $int=is$ the energy raises by $U$ at the first step, it is constant in the following steps but the last, when it decreases by $U$. If $int=wa$, the energy raises by either $U_1$ or $U_2$ at the first step, it is constant in the following steps but the last, when it decreases by $U_1$ or $U_2$ respectively. Thus these are admissible moves.

First, we prove (i). Let $\x_1^{(e)},...,\x_{m}^{(e)}\notin\cC_{int}^*$ the saddles visited during the sliding of a bar around a frame-angle. We want to prove that these saddles are essential (see Section \ref{modinddef} point 4 for the definition). Since we can repeat the following argument $m$ times, we may focus on a single configuration $\x_i^{(e)}$. Since $\cC_{int}^*$ is a gate for the transition and $\x_i^{(e)}\in\cS_{int}(\vuoto,\pieno)\setminus\cC_{int}^*$ for any $i=1,...,m$, we note that a path $\o\in(\vuoto\ra\pieno)_{opt}$ such that $\{\arg\max_{\o}H\}=\{\x_i^{(e)}\}$ does not exist. Thus our strategy consists in finding a path $\o\in(\vuoto\ra\pieno)_{opt}$ such that for any $\o'\in(\vuoto\ra\pieno)_{opt}$ 
\be{condess}
\o\cap\x_i^{(e)}\neq\emptyset \hbox{ and } \{\hbox{arg max}_{\o'}H\}\nsubseteq\{\hbox{arg max}_{\o}H\}\setminus\{\x_i^{(e)}\}, \quad i=1,...,m.
\ee

\begin{figure}
	\centering
	\begin{tikzpicture}[scale=0.5,transform shape]
		
		\draw (0,0) rectangle (16,16);
		\draw(-1,-1) rectangle (17,17);
		\node at (17.5,16) {\Huge{$\L$}};
		\node at (15.4,15) {\Huge{$\L_0$}};
		\draw (5,1)--(11,1);
		\draw (11,1)--(11,2);
		\draw(11,2)--(12,2);
		\draw(12,2)--(12,6);
		\draw(12,6)--(11,6);
		\draw(11,6)--(11,8);
		\draw(11,8)--(10,8);
		\draw(10,8)--(10,9);
		\draw(10,9)--(3,9);
		\draw(3,9)--(3,4);
		\draw(3,4)--(4,4);
		\draw(4,4)--(4,2);
		\draw (4,2)--(5,2);
		\draw (5,2)--(5,1);
		\draw [grigio,fill=grigio] (3,8)--(10,8)--(10,9)--(3,9)--(3,8);
		\draw [grigio,fill=grigio] (4,2) rectangle (11,8);
		\draw [grigio,fill=grigio] (3,8) rectangle (10,9);
		\draw [grigio,fill=grigio] (3,4) rectangle (4,9);
		\draw [grigio,fill=grigio] (11,2) rectangle (12,6);
		\draw [grigio,fill=grigio] (5,1) rectangle (11,2);
		\draw[dashed] (4,2) rectangle (11,8);
		\draw (10,8)--(10,9);
		\draw (3,9)--(10,9);
		\draw (3,9)--(3,4);
		\draw (3,4)--(4,4);
		\draw (4,4)--(4,2);
		\draw (4,2)--(5,2);
		\draw (5,2)--(5,1);
		\draw (5,1)--(11,1);
		\draw (11,1)--(11,2);
		\draw (11,2)--(12,2);
		\draw (12,2)--(12,6);
		\draw (12,6)--(11,6);
		\draw (11,6)--(11,8);
		\draw (11,8)--(10,8);
		\draw[dashed] (5,9) rectangle (6,10);
		\node at (5.5,9.5) {$\bullet$};
		\draw[fill=grigio] (5,15) rectangle (6,16);
		\draw[->](8.5,9.5)--(5.7,9.5);
		\node at (9.5,9.5) {\Huge{$(i,j)$}};
		\node at (5.5,15.5) {$\bullet$};
		\draw[->](8.5,15.5)--(5.7,15.5);
		\node at (10.5,15.5) {\Huge{$(i,L-1)$}};
		\draw[dashed] (5,14) rectangle (6,15);
		\node at (5.5,14.5) {$\bullet$};
		\draw[dashed] (5,13) rectangle (6,14);
		\node at (5.5,13.5) {$\bullet$};
		\draw[dashed] (5,12) rectangle (6,13);
		\node at (5.5,12.5) {$\bullet$};
		\draw[dashed] (5,11) rectangle (6,12);
		\node at (5.5,11.5) {$\bullet$};
		\draw[dashed] (5,10) rectangle (6,11);
		\node at (5.5,10.5) {$\bullet$};
		\draw[->](8.5,10.5)--(5.7,10.5);
		\node at (10.3,10.5) {\Huge{$(i,j+1)$}};
		\draw[->](8.5,14.5)--(5.7,14.5);
		\node at (10.5,14.5) {\Huge{$(i,L-2)$}};
		
		\path
		(4.7,15.5) edge[bend right=70,->] node [left] {} (4.7,14.5)
		(4.7,14.4) edge[bend right=70,->] node [left] {} (4.7,13.5)
		(4.7,13.4) edge[bend right=70,->] node [left] {} (4.7,12.5)
		(4.7,12.4) edge[bend right=70,->] node [left] {} (4.7,11.5)
		(4.7,11.4) edge[bend right=70,->] node [left] {} (4.7,10.5)
		(4.7,10.4) edge[bend right=70,->] node [left] {} (4.7,9.5);

	\end{tikzpicture}
	
	\vskip 0 cm
	\caption{Here we depict the configuration $\h_1$ that consists of a cluster $\h\in\cD_{int}$ union a free particle, in grey, that is in position $(i,L-1)$. The dotted unit squares represent the following positions of the free particle that moves as represented by the arrows on the left, until the particle is attached to the cluster in position $(i,j)$. The latter is the configuration $\h^B$.}
	\label{fig:camminoij}
\end{figure}
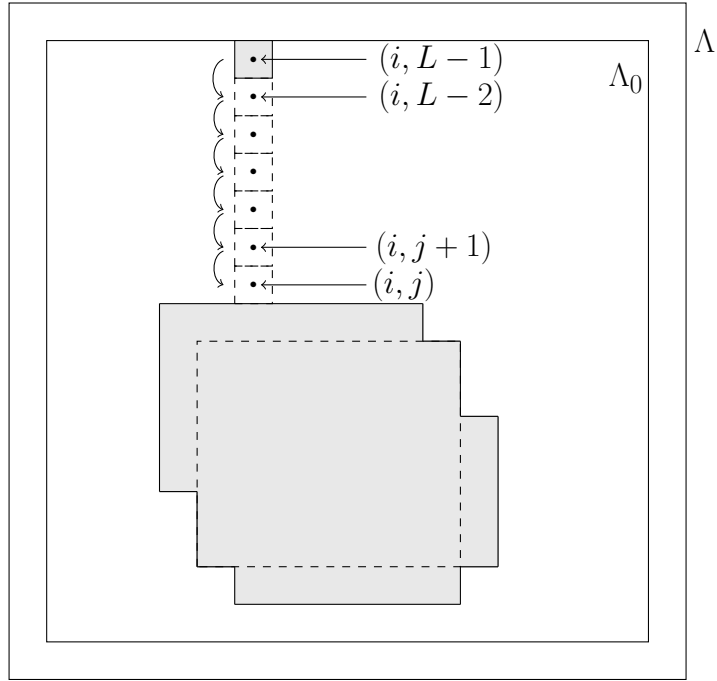

\noindent
Let $\h^B$ be the union of a cluster $\h\in\cD_{int}$ and a protuberance attached to one of its bars in a site with coordinates $(i,j)$. Without loss of generality assume that the bar is $B^n(\h)$ and that $|c^{wn}(\h)|=1$, otherwise a sequence of $1$-translations of the bars $B^n(\h)$ and $B^w(\h)$ can take place before creating the free particle in order to obtain $|c^{wn}(\h)|=1$. Note that during these translations the path does not cross any saddle. We define the specific path $\o$ of the strategy above as
\be{} 
\o=(\vuoto,\o_1,...,\o_k,\h,\h_1,...,\h_{L-j-2},\h^B,\x_1^{(e)},...,\x_m^{(e)})\circ\bar\o,
\ee

\noindent
where $\o_1,...,\o_k\in\cC_{\pieno}^{\vuoto}(\G_{int}^*)$, $\h\in\cD_{int}$, $\h_1\in\cC_{int}^*(L-j-1),...,\h_{L-j-2}\in\cC_{int}^*(2)$, $\h^B\in\cC_{int}^B(\hat\h_{L-j-2})$ (see Figure \ref{fig:camminoij} for a picture of this situation) and $\bar\o$ is a path that connects $\x_m^{(e)}$ to $\pieno$ such that $\max_{\s\in\bar\o}H(\s)\leq\G_{int}^*$. Now we show that for any $\o'$ the condition (\ref{condess}) is satisfied. If $\o'$ passes through the configuration $\x_i^{(e)}$, 
$ \{\hbox{arg max}_{\o'}H\}\supseteq \{\x_i^{(e)}\}$, thus  (\ref{condess}) is satisfied. Therefore we can assume that $\o'\cap\x_i^{(e)}=\emptyset$. If $\o'$ crosses the set $\cS_{int}(\vuoto,\pieno)$ through a configuration $\widetilde\h$ such that $\widetilde\h\cap\o=\emptyset$, then the condition (\ref{condess}) holds. In the sequel $\o'$ visits the configurations $\h_1,...,\h_{L-j-2}\in\cC_{int}^*$. Starting from $\h_{L-j-2}$, there are four allowed directions for moving the free particle. If we move it in the direction of the cluster (south in Figure \ref{fig:camminoij}), we deduce that the path $\o'$ visits the configuration $\h^B$. For the other three choices, the free particle still remains free after the move, indeed by construction of the path $\o$, starting from $\h_{L-j-2}$ it is not possible to reach the set $\cC_{int}^G$ via one step of the dynamics. Thus the path $\o'$ can visit either a saddle not already visited by $\o$ (west or east in Figure \ref{fig:camminoij}) or a saddle that has been already visited by $\o$ (north in Figure \ref{fig:camminoij}). In the first case, we obtain that (\ref{condess}) is satisfied. In the latter case, we can iterate this argument and, since $\o'$ goes from $\vuoto$ to $\pieno$, we can assume that the path $\o'$ visits the configuration $\h^B\in\cC_{int}^B(\hat\h_{L-j-2})$. From now on, starting from $\h^B$, there are two possible scenarios:
\bi
\item[I.] $\o'$ activates the same sliding of a bar around a frame-angle as $\o$;
\item[II.] $\o'$ activates a sliding of a bar around a frame-angle different from $\o$.
\ei

In case I, since $\o'\cap\x_i^{(e)}=\emptyset$, the sliding of a bar around a frame-angle has been stopped before hitting $\x_i^{(e)}$. Thus we can assume that $\o'$ comes back to $\h^B$, otherwise the energy exceeds $\G_{int}^*$. Since the path $\o'$ must reach $\pieno$, starting from $\h^B$ the protuberance is detached and in the sequel is attached in another site. Thus $\o'$ reaches a saddle that is not visited by $\o$. This implies that (\ref{condess}) is satisfied. 

In case II, when the path $\o'$ initiates the sliding of a bar around a frame-angle, it reaches at the first step a saddle that it is not visited by $\o$, thus the condition (\ref{condess}) is satisfied. Therefore the unique possibility is not to start this sliding and thus the path $\o'$ must come back to $\h^B$, since it has to reach $\pieno$. From now on, as before, the path $\o'$ has to detach the protuberance that in the sequel is attached in another site, thus $\o'$ visits a saddle that is not visited by $\o$. This implies that (\ref{condess}) is satisfied. Thus we have proved that the saddle $\x_i^{(e)}$ is essential for any $i=1,...,m$.

Finally, we prove (ii). By assumptions we know that it is not possible to complete a sliding of a bar around a frame-angle and thus this sliding must stop. Let $\x_1^{(ne)},...,\x_n^{(ne)}$ the saddles that are visited during this motion, we want to prove that these saddles are unessential (see Section \ref{modinddef} point 4 for the definition). Since we can repeat the following argument $n$ times, we may focus on a single configuration $\x_i^{(ne)}$. Consider any $\o\in(\vuoto\ra\pieno)_{opt}$ such that $\o\cap\x_i^{(ne)}\neq\emptyset$. Since $\cC_{int}^*$ is a gate for the transition from $\vuoto$ to $\pieno$ and $\x_i^{(ne)}\in\cS_{int}(\vuoto,\pieno)\setminus\cC_{int}^*$ for any $i=1,...,n$, we note that $\{\arg\max_{\o}H\}\setminus\{\x_i^{(ne)}\}\neq\emptyset$. Thus our strategy consists in finding $\o'\in(\vuoto\ra\pieno)_{opt}$ such that $\{\arg\max_{\o'}H\}\subseteq\{\arg\max_{\o}H\}\setminus\{\x_i^{(ne)}\}$. Starting from $\x_i^{(ne)}$, the unique admissible moves in order to not exceed $\G_{int}^*$ are the time-reversal of the previous moves. This implies that the path must come back to the starting configuration $\h^B$. Thus we can write $\o=(\vuoto,\o_1,...,\o_k,\g_1,...,\g_l,\h,\h^B,\x_1^{(ne)},...,\x_i^{(ne)},..,\x_1^{(ne)},\h^B)\circ\bar\o$, where $\o_1,...,\o_k\in\cC_{\pieno}^{\vuoto}(\G_{int}^*)$, $\g_1,...,\g_l\in\cX\setminus\cC_{\vuoto}^{\pieno}(\G_{int}^*-H(\cX^s))$ such that $H(\g_i)\leq\G_{int}^*$ for any $i=1,...,l$ and $\h\in\cC_{int}^*$, $\h^B\in\cC_{int}^B(\hat\h)$ and $\bar\o$ is a path that connects $\h^B$ to $\pieno$ such that $\max_{\s\in\bar\o}H(\s)\leq\G_{int}^*$. For this path we define a new path $\o'=(\vuoto,\o_1,...,\o_k,\g_1,...,\g_l,\h,\h^B)\circ\bar\o$. Thus we deduce that $\{\hbox{arg max}_{\o'}H\}\subseteq\{\hbox{arg max}_{\o}H\}\setminus\{\x_i^{(ne)}\}$, which implies that the saddle $\x_i^{(ne)}$ is unessential for any $i=1,..,n$.
\epr

\subsubsection{Lemmas valid for the weakly anisotropic model}
\label{lemmiani}

In this Section we give some lemmas for the weakly anisotropic model, i.e., $int=wa$, that will be useful later on. We postone the proof to Appendix \ref{app*}. Lemmas \ref{coltorow} and \ref{dtilde} are valid also in the case $int=sa$ (see \cite[Lemma 5.8]{BN3} and \cite[Lemma 5.9]{BN3} respectively), while Lemma \ref{trasless} has a corresponding version for the case $int=sa$ (see \cite[Lemma 5.7]{BN3}).

\bl{trasless}
Starting from $\h^B\in\cC_{wa}^B$, the saddles obtained by a $1$-translation of a bar are essential and in $\cN_{0}^{\a'}\cup\cN_1^{\a}$. Moreover, all the saddles in $\cN_{0}^{\a'}\cup\cN_1^{\a}$ can be obtained from this $\h^B$ by a $1$-translation of a bar. 
\el

\bl{coltorow}
Starting from a configuration $\h\in\cC_{wa}^*$, it is not possible to slide a vertical bar around a frame-angle without exceeding the energy $\G_{wa}^*$.
\el

\noindent
With the following Lemma we can justify the definition of $\cwgeo$ given in (\ref{c*wa}).

\bl{dtilde}
Starting from $\widetilde\cD_{wa}$, the dynamics either passes through $\bar\cD_{wa}$ or it is not possible that a free particle is created without exceeding the energy level $\G_{wa}^*$.
\el

\subsection{Model-dependent strategy}
\label{S6.4}
Our goal is to characterize the union of all the minimal gates for isotropic and weakly anisotropic interactions. To this end, due to \cite[Theorem 5.1]{MNOS}, we will characterize all the essential saddles for the transition from the metastable to the stable state. In this Section we apply the model-independent strategy explained in Section \ref{genstrategy} in order to identify some unessential saddles. Let $int\in\{is,wa\}$. We apply (\ref{ciclo}) both for $\s=\vuoto$, $\cA=\{\pieno\}$ and $\G=\G_{int}^*$ defining $\cC_{\pieno}^{\vuoto}(\G_{int}^*)$, and for $\s=\pieno$, $\cA=\{\vuoto\}$ and $\G=\G_{int}^*-H(\pieno)$ defining $\cC_{\vuoto}^{\pieno}(\G_{int}^*-H(\pieno))$. We chose this notation in order to emphasize the dependence on $\G_{int}^*$. First, we prove the required model-dependent inputs (iii)-(a) and (iii)-(b) in Section \ref{genstrategy} (see Proposition \ref{pathstrong}(i) and Proposition \ref{pathstrong}(ii)). Second, by \cite[Theorem 1.3.3(iii)]{BHN} for $int=is$ and Theorem \ref{thgateweak} for $int=wa$, we know that $\cC_{int}^*$ is a gate for the transition from $\vuoto$ to $\pieno$ for $int\in\{is,wa\}$. Thus we apply the model-independent strategy explained in Section \ref{genstrategy} to Kawasaki dynamics by taking $m=\vuoto$, $\cX^s=\{\pieno\}$, $\cW(m,\cX^s)=\cC_{int}^*$, $\cL^B=\cC_{int}^B$ and $\cL^G=\cC_{int}^G$. In Proposition \ref{c*contenuto} we prove that $\cC_{int}^*\subseteq\cG_{int}(\vuoto,\pieno)$, that allows us to study the essentiality only of the saddles that are not in $\cC_{int}^*$.

In order to apply Propositions \ref{selle1} and \ref{selle2}, we need to characterize the sets $K_{int}$ and $\widetilde{K}_{int}$ (see (\ref{defK}) and (\ref{defKtilde}) respectively for the definitions) for our models. This is done in Proposition \ref{Kvuoto}. Due to this result, our strategy consists in partitioning the saddles that are not in $\cC_{int}^*$ in three types: the saddles that are in the boundary of $\cC_{\pieno}^{\vuoto}(\G_{int}^*)$, i.e., $\s\in\partial\cC_{\pieno}^{\vuoto}(\G_{int}^*)\cap(\cS_{int}(\vuoto,\pieno)\setminus\cC_{int}^*)$, the saddles that are in the boundary of $\cC_{\vuoto}^{\pieno}(\G_{int}^*-H(\pieno))$ and not in $\widetilde{K}_{int}$, i.e., $\z\in\partial\cC_{\vuoto}^{\pieno}(\G_{int}^*-H(\pieno))\cap(\cS_{int}(\vuoto,\pieno)\setminus(\cC_{int}^*\cup\widetilde{K}_{int}))$, and the remaining saddles $\x\in
\cS_{int}(\vuoto,\pieno)\setminus(\partial\cC_{\pieno}^{\vuoto}(\G_{int}^*)\cup(\partial\cC_{\vuoto}^{\pieno}(\G_{int}^*-H(\pieno))\setminus\widetilde{K}_{int})\cup\cC_{int}^*)$. By Propositions \ref{selle1} and \ref{selle2}, we obtain Corollary \ref{corstrategy} that states that the saddles of the first and second types are respectively unessential. In Proposition \ref{selle3} for $int\in\{is,wa\}$ we highlight some of the saddles of type three that are unessential.

We need to distinguish the analysis for $int\in\{is,wa\}$ and $int=sa$ due to the different entrance in $\cC_{int}^*$ for $int\in\{is,wa\}$ and $int=sa$ (see Lemma \ref{entrataweak} and \cite[Lemma 5.17]{BN3} respectively). For the case $int=sa$ this strategy is presented in \cite[Section 5.4]{BN3}.

Finally, we identify the essential saddles of the third type in Proposition \ref{isoselless} for the isotropic interactions and in Proposition \ref{weakselless} for the weakly anisotropic interactions.

\subsubsection{Main Propositions}
In this Subsection we give the main results for our model-dependent strategy. We refer to Subsection \ref{proofproposizionidep} for the proof of these propositions.

The next proposition shows that when the dynamics reaches $\cC_{int}^G$ it has gone “over the hill”, while when it reaches $\cC_{int}^B$ the energy has to increase again to the level $\G_{int}^*$ to visit $\vuoto$ or $\pieno$. An analogue version for $int=is$ is proven in \cite[Proposition 2.3.9]{BHN}, while here we extend that result to $int=wa$ following a similar argument. Note that this result holds also in the case $int=sa$, see \cite[Proposition 5.10]{BN3}.

\bp{pathstrong}
Let $int\in\{is,wa\}$.
\item[(i)]  If $\h\in\cC_{int}^G$, then there exists a path $\o:\h\ra\pieno$ such that $\max_{\z\in\o}H(\z)<\G_{int}^*$. 
\item[(ii)] If $\h\in\cC_{int}^B$, then there are no $\o:\h\ra\vuoto$ or $\o:\h\ra\pieno$ such that $\max_{\z\in\o}H(\z)<\G_{int}^*$.
\ep

\bp{c*contenuto}
Let $int\in\{is,wa\}$, then $\cC_{int}^*\subseteq\cG_{int}(\vuoto,\pieno)$.
\ep

\noindent
Proposition \ref{c*contenuto} holds also in the case $int=sa$ (see \cite[Proposition 5.11]{BN3}).

\bp{Kvuoto}
Let $int\in\{is,wa\}$, then 
\bi
\item[(i)] $K_{int}=\emptyset$; 
\item[(ii)] For the set $\widetilde{K}_{int}$ we obtain the following characterization:
\bi
\item[(a)]$\widetilde{K}_{is}\cap\partial\cC_{\vuoto}^{\pieno}(\G_{is}^*-H(\pieno))=\{\h\in\bigcup_{i}\bigcup_{\a,\a'}\mathscr{I}_{i}^{\a,\a'}\setminus\cigeo: \text{it is possible to attach the free} \\ \text{\qquad \qquad \qquad \qquad  \qquad \quad particle in }  \partial^-\text{CR}(\h) \text{ via one step of the dynamics}\}=:\bar{\mathscr{I}}_{is}$;
\item[(b)]$\widetilde{K}_{wa}\cap\partial\cC_{\vuoto}^{\pieno}(\G_{wa}^*-H(\pieno))=\bigcup_{\a}\bigcup_{\a'}\bigcup_{k'}\cN_{2,k'}^{\a,\a'}$.
\ei
\ei
\ep

\noindent
For the corresponding result of Proposition \ref{Kvuoto} for $int=sa$, see \cite[Proposition 5.12]{BN3}.

\bc{corstrategy}
Let $int\in\{is,wa\}$.
\bi
\item[(i)]The saddles of the first type $\s\in\partial\cC_{\pieno}^{\vuoto}(\G_{int}^*)\cap(\cS_{int}(\vuoto,\pieno)\setminus\cC_{int}^*)$ are unessential;
\item[(ii)] The saddles of the second type
\bi
\item[(a)] $\z\in\partial\cC_{\vuoto}^{\pieno}(\G_{is}^*-H(\pieno))\cap(\cS_{is}(\vuoto,\pieno)\setminus(\cC_{is}^*\cup\bar{\mathscr{I}}_{is}))$ are unessential;
\item[(b)] $\z\in\partial\cC_{\vuoto}^{\pieno}(\G_{wa}^*-H(\pieno))\cap(\cS_{wa}(\vuoto,\pieno)\setminus(\cC_{wa}^*\cup\bigcup_{\a}\bigcup_{\a'}\bigcup_{k'}\cN_{2,k'}^{\a,\a'}))$ are unessential.
\ei
\ei
\ec

\bpr
Combining Propositions \ref{selle1}, \ref{selle2} and \ref{Kvuoto} we get the claim.
\epr

\bp{selle3}
Let $int\in\{is,wa\}$. Any saddle $\x$ that is neither in $\cC_{int}^*$, nor in the boundary of the cycle $\cC_{\pieno}^{\vuoto}(\G_{int}^*)$, nor in $\partial\cC_{\vuoto}^{\pieno}(\G_{int}^*-H(\pieno))\setminus\widetilde{K}_{int}$, i.e.,  $\x\in\cS_{int}(\vuoto,\pieno)\setminus(\partial\cC_{\pieno}^{\vuoto}(\G_{int}^*)\cup(\partial\cC_{\vuoto}^{\pieno}(\G_{int}^*-H(\pieno))\setminus\widetilde{K}_{int})\cup\cC_{int}^*)$, such that $\t_{\x}<\t_{\cC^B_{int}}$ is unessential. Therefore it is not in $\cG_{int}(\vuoto,\pieno)$.
\ep

\noindent
For the corresponding result of Proposition \ref{selle3} for $int=sa$, we refer to \cite[Proposition 5.14]{BN3}.

\subsubsection{Useful Lemmas for the model-dependent strategy}
\label{lemminuovi}
In this Subsection we give some useful lemmas about the entrance in the gate and the minimality of the sets $\cC_{int}^*(i)$ with $i=3,...,L_{int}^*$ for $int\in\{is,wa\}$. We stress that the behavior for $int\in\{is,wa\}$ is very different from that observed for $int=sa$, indeed we note that the weakly anisotropic model has some characteristics similar to the isotropic and some similar to the strongly anisotropic model. For the corresponding results obtained in the case $int=sa$ we refer to \cite[Subsection 5.4.2]{BN3}. Recall Definition \ref{critsets} for the definitions of $\cQ_{is}$, $\cD_{is}$ and $\cigeo$, and (\ref{definizioneqbarsa}) and (\ref{c*wa}) for the corrisponding definition for $int=wa$. The next lemma generalizes \cite[Proposition 2.3.8]{BHN}, proved for $int=is$, to the case $int=wa$ following similar arguments. In the case $int=sa$, this result is given in \cite[Lemma 5.15]{BN3}.

\bl{motion}
Let $int\in\{is,wa\}$.

\bi
\item[(i)] Starting from $\cC_{int}^*\setminus\cQ_{int}^{fp}$, the only transitions that do not raise the energy are motions of the free particle in the region where the free particle is at lattice distance $\geq3$ from the protocritical droplet.
\item[(ii)] Starting from $\cQ_{int}^{fp}$, the only transitions that do not raise the energy are motions of the free particle in the region where the free particle is at lattice distance $\geq3$ from the protocritical droplet and motions of the protuberance along the side of the rectangle where it is attached. When the lattice distance is 2, either the free particle can be attached to the protocritical droplet or the protuberance can be detached from the protocritical droplet and attached to the free particle, to form a rectangle plus a dimer. From the latter configuration the only transition that does not raise the energy is the reverse move.
\item[(iii)] Starting from $\cC_{int}^*$, the only configurations that can be reached by a path that lowers the energy and does not decrease the particle number, are those where the free particle is attached to the protocritical droplet.
\ei

\el

\bpr
The proof is analogue to the one reported in \cite{BHN} for $int=is$: the path we consider is the same as in the isotropic case, but the energy of the moves is different. Indeed the energy decreases by $U$ if $int=is$ and by either  $U_1$ or $U_2$ depending whether it is attached to the vertical or horizontal side respectively. 
\epr

The next lemma investigates how the entrance in $\cC_{int}^*$ occurs when $int\in\{is,wa\}$. This result for $int=is$ is proven in \cite[Proposition 2.3.7]{BHN}, while here we extend that result for $int=wa$. We encourage the reader to inspect the difference between lemma \ref{entrataweak} and \cite[Lemma 5.17]{BN3}, where the peculiar entrance in the gate for the strongly anisotropic case is analyzed.

\bl{entrataweak}
Let $int\in\{is,wa\}$. Any $\o\in(\vuoto\ra\pieno)_{opt}$ passes first through $\cQ_{int}$, then possibly through $\cD_{int}\setminus\cQ_{int}$, and finally through $\cC_{int}^*$.
\el

\noindent
We postpone the proof of Lemma \ref{entrataweak} in Appendix \ref{app*}. We refer to Subsection \ref{prooflemmiapp} for the proof of the remaining lemmas. The next Lemma  is proven for $int=is$ in \cite[Section 3.5]{BHN}, while here we extend that result for $int=wa$. Recall the definition of minimal gate given in Section \ref{modinddef} point 4. In Lemma \ref{mingate} we extend to $int=wa$ the statement in \cite[eq.\ (3.5.5)]{BHN} for $\cC_{is}^*(i)$ with $i=3,...,L_{is}^*$. Concerning $\cC_{is}^*(2)$, in Lemma \ref{C*2ess}, we correct the statement in \cite[eq.\ (3.5.5)]{BHN} by replacing the minimality of the gate $\cC_{is}^*(2)$ with the sentence ``$\cC_{is}^*(2)$ is composed by essential saddles". We stress that this correction does not effect the results where the statement was used in \cite{BHN}. Moreover, in Lemma \ref{C*2ess}, we prove that the saddles in $\cC_{int}^*(2)$ are essential also for $int=wa$. The result for $int=sa$ is given in \cite[Lemma 5.16]{BN3}.

\bl{mingate}
Let $int\in\{is,wa\}$, then $\cC_{int}^*(i)$ is a minimal gate for any $i=3,...,L_{int}^*$.
\el

\br{}
In the case $int=sa$, the statement of Lemma \ref{mingate} does not hold. A different result is derived in \cite[Lemma 5.18]{BN3}.
\er

\bl{C*2ess}
Let $int\in\{is,wa\}$. The saddles in $\cC_{int}^*(2)$ are essential.
\el

\subsubsection{Proof of Propositions}
\label{proofproposizionidep}

\begin{proof*}{\bf of Proposition \ref{pathstrong}}
	The case $int=is$ is proven in \cite[Proposition 2.3.9]{BHN}, thus we consider $int=wa$.
	
	(i) If $\h\in\cC_{wa}^G$, then its energy is either $\G_{wa}^*-U_1-U_2$ or $\G_{wa}^*-U_1$ (resp.\ $\G_{wa}^*-U_2$), depending on whether the attached particle is in a corner or is a protuberance on a vertical (resp.\ horizontal) side. In the latter case we can move the particle at no cost and gain an extra $-U_2$ (resp.\ $-U_1$) when it has become a corner. After that it is possible to create a new particle and attach it, which leads to energy $\G_{wa}^*-U_1-U_2-(U_1+U_2-\D)<\G_{wa}^*$. We can continue in this way, filling up all the sites in $\partial^-\hbox{CR}(\h)$. Now we can proceed along the reference path for the nucleation constructed in \cite[Section 3.2]{NOS} until the path reaches $\pieno$. We have exhibited a path $\o$ such that $\max_{\s\in\o}H(\s)<\G_{wa}^*$.
	
	(ii) If $\h\in\cC_{wa}^B$, then $H(\h)=\G_{wa}^*-U_1$ (resp.\ $H(\h)=\Gamma_{wa}^*-U_2$) if the protuberance has been attached to a vertical (resp.\ horizontal) side. As long as the energy does not exceed $\G_{wa}^*$, it is impossible to create a new particle before further lowering the energy. But there are no moves available to lower the energy. As a consequence the unique admissible moves are those where the last particle that was attached is moving along the side at zero cost or detaching again, or start a sliding of a bar around a frame-angle (see the explanation in the third case). In the first case we obtain a configuration that is analogue to $\h\in\cC_{wa}^B$ and therefore we can iterate the argument by taking this configuration as $\h^B$ for a finite number of steps, since the path has to reach $\pieno$. In the second case, we obtain a configuration that is in $\cC_{wa}^*$, thus the path has to reach again the energy value $\G_{wa}^*$. In the third case, we justify separately when the sliding is at cost $U_1$ or at cost $U_2$. If $H(\h)=\G_{wa}^*-U_1$, the only admissible motions along the border of the droplet that do not exceed $\G_{wa}^*$ are those at cost $U_1$, since the unique possibility is to move the particle in a frame-angle in such a way that it connects to the protuberance, otherwise all the other moves have cost at least $U_1+U_2$. Similarly, by symmetry we deduce that if $H(\h)=\G_{wa}^*-U_2$, then the only admissible move is starting a sliding of a bar around a frame-angle at cost $U_2$. In both cases the energy returns to $\G_{wa}^*$, which concludes the proof.
\end{proof*}

\begin{proof*}{\bf of Proposition \ref{c*contenuto}}
	Let $int\in\{is,wa\}$. By Lemma \ref{C*2ess} we know that the saddles in $\cC_{int}^*(2)$ are essential and thus are in $\cG_{int}(\vuoto,\pieno)$ due to \cite[Theorem 5.1]{MNOS}. Moreover, by Lemma \ref{mingate} we know that the set $\cC_{int}^*(i)$ is a minimal gate for any $i=3,...,L_{int}^*$, thus
	\be{} 
	\cC_{int}^*=\cC_{int}^*(2)\cup\bigcup_{i=3}^{L_{int}^*}\cC_{int}^*(i)\subseteq\cG_{int}(\vuoto,\pieno).
	\ee
\end{proof*}

\begin{proof*}{\bf of Proposition \ref{Kvuoto}}
	Let $int\in\{is,wa\}$. 
	
	\noindent
	(i) To prove that $K_{int}=\emptyset$ we argue by contradiction. Let $\bar\h\in K_{int}$, thus there exist $\h\in\cC_{int}^*$ and $\o=\o_1\circ\o_2$ from $\h$ to $\vuoto$ with the properties described in (\ref{defK}), where $\circ$ denotes the composition of two paths. We know that $\h$ is composed by the union of a protocritical droplet in $\cD_{int}$ and a free particle. Since $\o_1\cap\cC_{int}^*=\{\h\}$, we note that the free particle must be in $\L^-$, otherwise the free particle has to cross at least $\L^-$ and $\partial\L^-$, the latter in the configuration $\h'\in\cC_{int}^*$, with $\h'\neq\h$, which contradicts the conditions in (\ref{defK}). Therefore, starting from $\h$, by the optimality of the path we deduce that the unique admissible move is to remove the free particle. The configuration that is obtained in this way is in $\cD_{int}$, that belongs to $\cC_{\pieno}^{\vuoto}(\G_{int}^*)$, which is absurd since (\ref{defK}) requires that $\o_1\cap\cC_{\pieno}^{\vuoto}(\G_{int}^*)=\emptyset$. Thus it is not possible to find $\o_1$ and $\o_2$, therefore $K_{int}=\emptyset$. 
	
	\medskip
	\noindent
	(ii) Let $\bar\h\in\widetilde{K}_{int}\cap\partial\cC_{\vuoto}^{\pieno}(\G_{int}^*-H(\pieno))$ for any $int\in\{is,wa\}$. By the definition of the set $\widetilde{K}_{int}$ we know that there exist $\h\in\cC_{int}^*$ and $\o=\o_1\circ\o_2$ from $\h$ to $\pieno$ with the properties described in (\ref{defKtilde}). We know that $\h$ is composed by the union of a protocritical droplet $\hat\h\in\cD_{int}$ and a free particle. Since $\o_1\cap\cC_{int}^*=\{\h\}$, we note that $\h\in\cC_{int}^*(2)$, otherwise the free particle has to cross at least $\bar{B}_2(\hat\h)$ and $\bar{B}_3(\hat\h)$, the latter in the configuration $\h'\in\cC_{int}^*$, with $\h'\neq\h$, which contradicts the conditions in (\ref{defKtilde}). Therefore, starting from $\h$, by the optimality of the path we deduce that the unique admissible move is to attach the free particle to the cluster. If $\bar\h$ is obtained from $\h$ by attaching the free particle in a good site giving rise to a configuration in $\cC_{int}^G(\hat\h)$, by Proposition \ref{pathstrong}(i) we know that $\o_1\cap\cC_{\vuoto}^{\pieno}(\G_{int}^*-H(\pieno))\neq\emptyset$, that contradicts (\ref{defKtilde}), thus t is not possible to find $\o_1$ and $\o_2$, therefore $\bar\h\notin\widetilde{K}_{int}$, which is in contradiction with the assumption. 
	
	Assume now that $\bar\h$ is obtained from $\h$ by attaching the free particle in a bad site giving rise to a configuration in $\cC_{int}^B(\hat\h)$. If $\h\in\cQ_{int}^{fp}$, then by Lemma \ref{motion}(ii) the unique admissible move is the reverse one, thus we may assume that $\h\in\cC_{int}^*\setminus\cQ_{int}^{fp}$ and that the path does not go back to $\h$, otherwise we can iterate this argument for a finite number of steps since the path has to reach $\pieno$. Starting from $\h$, by Lemma \ref{trenini} we know that $\bar\h$ is obtained either via a sequence of $1$-translations of a bar or via a sliding of a bar around a frame-angle. If a sequence of $1$-translations takes place, by the optimality of the path we deduce that the unique possibility is either detaching the protuberance or sliding a bar around a frame-angle. In the first case the configuration that is obtained is in $\cC_{int}^*$ and thus $\bar\h\notin\widetilde{K}_{int}$, which contradicts the assumption. Now we analyze separately the case in which a sliding of a bar around a frame-angle takes place for $int=is$ and $int=wa$.

	If $int=is$, the configurations visited by the path $\o$ during this sliding are $\bar\h_1,...,\bar\h_m\in\mathscr{I}_{k,k',i}^{\a,\a'}$ for some $\a,\a'\in\{n,s,w,e\}$, $k'=2,...,l_c$ and $k=2,...,k'$, while the last configuration is a saddle $\widetilde\h\in\mathscr{I}_{i}^{\a,\a'}$ when the last particle of the bar is detached. Thus $\{\bar\h_1,...,\bar\h_m\}\cap\partial\cC_{\vuoto}^{\pieno}(\G_{int}^*-H(\pieno))=\emptyset$. Starting from $\widetilde\h$, the free particle can move and be attached in $\partial^+\hbox{CR}(\widetilde\h)$ and another sliding of a bar around a frame-angle can take place. If this is the case, as proved above the saddles visited during this motion are not in $\partial\cC_{\vuoto}^{\pieno}(\G_{is}^*-H(\pieno))$ except the last configuration $\widetilde{\widetilde\h}$ visited during the sliding, that is a saddle, if $\P(\widetilde{\widetilde\h},\cC^G_{is})>0$. Thus the unique possibility to have $\bar\h\in\widetilde{K}_{is}\cap\partial\cC_{\vuoto}^{\pieno}(\G_{is}^*-H(\pieno))$ is that $\bar\h\in\mathscr{I}_{i}^{\a,\a'}\setminus\cigeo$ and it is possible to attach the free particle in $\partial^-\hbox{CR}(\bar\h)$ via one step of the dynamics. Taking the union over all $i\in\{-1,0,1,2\}$ and $\a,\a'\in\{n,s,w,e\}$ we get the claim.

	If $int=wa$, by (\ref{condtrenino}), Proposition \ref{cardweak}(a) and Lemma \ref{coltorow} we deduce that the unique possibility to slide a bar around a frame-angle is that the bar is horizontal and it has length between $l_1^*-l_2^*+1$ and $l_2^*-1$. Thus the configurations visited by the path $\o$ during this sliding are $\bar\h_1,...,\bar\h_m\in\cN_{k,k'}^{\a,\a'}$ for some $\a\in\{n,s\}$, $\a'\in\{w,e\}$, $k'=l_1^*-l_2^*+1,...,l_2^*-1$ and $k=2,...,k'$, while the last configuration $\widetilde\h$ obtained when the last particle of the bar is detached is not a saddle. Note that $\widetilde\h$ is not in the set $\cB$ defined in \cite[eq.\ (3.64)]{NOS}, since $s(\widetilde\h)=s^*_{wa}+1$ and $v(\widetilde\h)=2l_2^*-l_1^*-2<p_{min}(\widetilde\h)-1=l_2^*-1$. Thus by \cite[Proposition 11]{NOS} we know that $\widetilde\h\in\cC_{\vuoto}^{\pieno}(\G_{wa}^*-H(\pieno))$. This implies that only the configuration $\bar\h_m$, that belongs to $\cN_{2,k'}^{\a,\a'}$, is in $\widetilde{K}_{wa}\cap\partial\cC_{\vuoto}^{\pieno}(\G_{wa}^*-H(\pieno))$. Taking the union of $\cN_{2,k'}^{\a,\a'}$ over all $\a\in\{n,s\}$ and $\a'\in\{w,e\}$, we get the claim.
\end{proof*}

\begin{proof*}{\bf of Proposition \ref{selle3}}
	Let $int\in\{is,wa\}$. We denote by $\x_1,...,\x_n$ the saddles in the statement. We want to prove that these saddles are unessential (see Section \ref{modinddef} point 4 for the definition). Since we can repeat the following argument $n$ times, we may focus on a single configuration $\x_i$. Consider any $\o\in(\vuoto\ra\pieno)_{opt}$ such that $\o\cap\x_i\neq\emptyset$. By hypotheses, we have to analyze only the case in which the path $\o$ reaches the saddle $\x_i$ before reaching $\cC_{int}^B$. Since $\cC_{int}^*$ is a gate for the transition and $\x_i\in\cS_{int}(\vuoto,\pieno)\setminus\cC_{int}^*$, we note that $\{\arg\max_{\o}H\}\setminus\{\x_i\}\neq\emptyset$. Thus our strategy consists in finding $\o'\in(\vuoto\ra\pieno)_{opt}$ such that $\{\arg\max_{\o'}H\}\subseteq\{\arg\max_{\o}H\}\setminus\{\x_i\}$.
	
	First, assume that $\o$ reaches the saddle $\x_i$ before reaching $\cC_{int}^G$ and thus $\x_i$ must be obtained by a configuration $\h\in\cC_{int}^*$ without attaching the free particle. In particular, Lemmas \ref{motion}(ii) and \ref{entrataweak} imply that the only possibility is that $\h$ is composed by the union of a cluster $\hat\h\in\cQ_{int}$ and a free particle at distance 2 from the cluster. Moreover, $\x_i$ is either the union of a quasi-square $(l_c-1)\times l_c$ with a dimer if $int=is$, or the union of a rectangle $(l_1^*-1)\times l_2^*$ with an horizontal dimer if $int=wa$. Thus starting from $\x_i$, by Lemma \ref{motion}(ii) we know that the only transition that does not raise the energy is the reverse move giving rise to the configuration $\h$. Thus by Lemma \ref{entrataweak} we can write 
	\be{} 
	\o=(\vuoto,\o_1,...,\o_{k_1},\g_1,\h_1^{(1)},...,\h_{m_1}^{(1)},...,\o_{k_q},\g_q,\h_1^{(q)},...,\h_{m_q}^{(q)},\h,\x_i,\h)\circ\bar\o ,
	\ee
	
	\noindent
	where $\o_1,...,\o_{k_1},...,\o_{k_q}\in\cC_{\pieno}^{\vuoto}(\G_{int}^*)$, $\g_1,...,\g_q\in\cD_{int}$, $\h_1^{(1)},...,\h_{m_1}^{(1)},...,\h_1^{(q)},...,\h_{m_q}^{(q)}\in\cC_{int}^*$ and $\bar\o$ is a path that connects $\h$ to $\pieno$ such that $\max_{\s\in\bar\o}H(\s)\leq\G_{int}^*$. We define a new path 
	\be{} 
	\o'=(\vuoto,\o_1,...,\o_{k_1},\g_1,\h_1^{(1)},...,\h_{m_1}^{(1)},...,\o_{k_q},...,\o_{k_{q+1}},\g_q,\h_1^{(q)},...,\h_{m_q}^{(q)},\h)\circ\bar\o.
	\ee
	
	\noindent
	Thus $\{\arg\max_{\o'}H\}=\{\h_1^{(1)},...,\h_{m_1}^{(1)},...,\h_1^{(q)},...,\h_{m_q}^{(q)},\h\}\cup\{\arg\max_{\bar\o}H\}$ and therefore 
	\be{}
	\{\hbox{arg max}_{\o'}H\}\subseteq\{\hbox{arg max}_\o H\}\setminus\{\x_i\}, \quad i=1,...,n.
	\ee
	
	\noindent
	This implies that the saddle $\x_i$ is unessential for any $i=1,...,n$ and thus, using \cite[Theorem 5.1]{MNOS}, $\x_i\in\cS_{int}(\vuoto,\pieno)\setminus\cG_{int}(\vuoto,\pieno)$.
	
	Finally, if the path $\o$ reaches the saddle $\x_i$ after reaching $\cC_{int}^G$ in the configuration $\h^G$, we can write
	\be{}
	\o=(\vuoto,\o_1,...,\o_{k_1},\g_1,\h_1^{(1)},...,\h_{m_1}^{(1)},...,\h_1^{(q)},...,\h_{m_q}^{(q)},\h^G,...,\x_i,...,\pieno)
	\ee
	
	\noindent
	and define
	\be{}
	\o'=(\vuoto,\o_1,...,\o_{k_1},\g_1,\h_1^{(1)},...,\h_{m_1}^{(1)},...,\h_1^{(q)},...,\h_{m_q}^{(q)},\h^G)\circ\widetilde\o,
	\ee
	\noindent
	where $\widetilde\o$ is a path such that $\max_{\s\in\widetilde\o}H(\s)<\G_{int}^*$. This path exists by Proposition \ref{pathstrong}(i). It easy to check that the saddle $\x_i$ is unessential for any $i=1,...,n$ and thus, using \cite[Theorem 5.1]{MNOS}, $\x_i\in\cS_{int}(\vuoto,\pieno)\setminus\cG_{int}(\vuoto,\pieno)$.
\end{proof*}

\subsubsection{Proof of Lemmas}
\label{prooflemmiapp}

\begin{proof*}{\bf of Lemma \ref{mingate}}
	The case $int=is$ is proven in \cite[Section 3.5]{BHN}, thus we consider $int=wa$. First, we prove that $\cC_{wa}^*(i)$ is a gate. By Lemma \ref{entrataweak} we know that any $\o\in(\vuoto\ra\pieno)_{opt}$ enters $\cC_{wa}^*$ through a configuration of the form $(\hat\h,z)$, with $\hat\h\in\bar\cD_{wa}$ a protocritical droplet (by (eq.)) and $z$ the site occupied by the free particle. Note that either $z\in B_i(\hat\h)$ if $d(\partial^-\L_4,\hat\h)>i$ or $z\in\bar B_i(\hat\h)$ if $d(\partial^-\L_4,\hat\h)\leq i$, thus $\cC_{wa}^*(i)$ is a gate. 
	
	Now we prove that $\cC_{wa}^*(i)$ is a minimal gate. For any $\h\in\cC_{wa}^*(i)$, our strategy consists in proving that $\cC_{wa}^*(i)\setminus\{\h\}$ is not a gate by defining a path $\o\in(\vuoto\ra\pieno)_{opt}$ such that $\o\cap(\cC_{wa}^*(i)\setminus\{\h\})=\emptyset$. For the following the reader can visualize the path described using Figure \ref{fig:sitiGB}. We take an arbitrary path starting from $\vuoto$ and that enters $\cC_{wa}^*(i)$ in $\h=(\hat\h,z)$. Then the path proceeds by moving the free particle from $z$ to $\hat\h$ such that, the distance between the free particle and $\hat\h$ at the first step is strictly decreasing, and at the later steps is not increasing. Finally the free particle is attached in a site $x\in\partial^-\hbox{CR}(\hat\h)$ giving rise to a configuration in $\cC_{wa}^G(\hat\h)$. From this configuration, the path proceeds towards $\pieno$ as the one in Proposition \ref{pathstrong}(i). Since the constructed $\o\in(\vuoto\ra\pieno)_{opt}$ and $\o\cap\cC_{wa}^*(i)=\{\h\}$, the proof is concluded.
\end{proof*}

\begin{proof*}{\bf Lemma \ref{C*2ess}}
	Let $\x_1,...,\x_n$ the saddles in $\cC_{int}^*(2)$, we want to prove that these saddles are essential (see Section \ref{modinddef} point 4 for the definition). Since we can repeat the following argument $n$ times, we may focus on a single configuration $\x_i$. We note that a path $\o\in(\vuoto\ra\pieno)_{opt}$ such that $\{\arg\max_{\o}H\}=\{\x_i\}$ does not exist, thus our strategy consists in finding a path $\o\in(\vuoto\ra\pieno)_{opt}$ such that for any $\o'\in(\vuoto\ra\pieno)_{opt}$ 
	\be{eqC*2ess}
	\o\cap\x_i\neq\emptyset \hbox{ and } \{\hbox{arg max}_{\o'}H\}\nsubseteq\{\hbox{arg max}_{\o}H\}\setminus\{\x_i\}, \quad i=1,...,n.
	\ee
	\noindent
	Let $\x_i$ be the union of a cluster $\h\in\cD_{int}$ and a free particle in a site with coordinates $(i,j)$ at lattice distance 2 from the cluster. We define the specific path $\o$ of the strategy above as
	\be{}
	\o=(\vuoto,\o_1,...,\o_k,\h,\h_1,...,\h_{L-j-3},\x_i,\h^{(1)},...,\h^{(k)},\h^G)\circ\bar\o,
	\ee
	\noindent
	where $\o_1,...,\o_k\in\cC_{\pieno}^{\vuoto}(\G_{int}^*)$, $\h\in\cD_{int}$, $\h_1\in\cC_{int}^*(L-j-1),...,\h_{L-j-3}\in\cC_{int}^*(3)$, $\h^{(1)},...,\h^{(k)}\in\cC_{int}^*$, $\h^G\in\cC_{int}^G(\hat\h^{(k)})$ and $\bar\o$ a path that connects $\h^G$ to $\pieno$ such that $\max_{\s\in\bar\o}H(\s)<\G_{int}^*$. Note that the part of the path $\o$ from $\x_i$ to $\h^{(k)}$ is constructed by moving the free particle at zero cost from $(i,j)$ to a good site depicted on the left-hand side of Figure \ref{fig:sitiGB}, so that we obtain a configuration $\h^G$. Moreover, the path $\bar\o$ exists by Lemma \ref{pathstrong}(i). Now we show that for any $\o'$ the condition (\ref{eqC*2ess}) is satisfied. If $\o'$ passes through the configuration $\x_i$, then $\{\arg \max_{\o'}H\}\supseteq\{\x_i\}$, thus (\ref{eqC*2ess}) is satisfied. Therefore we can assume that $\o'\cap\x_i=\emptyset$. If $\o'$ crosses the set $\cS_{int}(\vuoto,\pieno)$ through a configuration $\widetilde\h$ such that $\o\cap\widetilde\h=\emptyset$, then the condition (\ref{eqC*2ess}) holds. Thus we can reduce our analysis to $\o'$ that visits all the configurations $\h_1,...,\h_{L-j-3}\in\cC_{int}^*$. Starting from $\h_{L-j-3}\in\cC_{int}^*(3)$, there are four allowed directions for moving the free particle. The move can not be in the direction of the cluster, indeed in that case the path $\o'$ visits $\x_i\in\cC_{int}^*(2)$. Concerning the other three choices, we have two cases. In the first case, the path $\o'$ visits a saddle not already present in $\o$, thus (\ref{eqC*2ess}) is satisfied. In the second case, the path $\o'$ visits a saddle that has been already visited by $\o$, thus we can iterate this argument for a finite number of steps, since the path $\o'$ has to reach $\pieno$. Thus we have proved that the saddle $\x_i$ is essential for any $i=1,...,n$. 
\end{proof*}

\section{Proof of the main Theorem \ref{giso}: isotropic case}
\label{proofiso}
In this Section we give the proof of the main Theorem \ref{giso} by analyzing the geometry of the set $\cG_{is}(\vuoto,\pieno)$.

\subsection{Main Proposition}
\label{propiso}

In this Section we give the proof of the main Theorem \ref{giso}, emphasizing the saddles for the transition from $\vuoto$ to $\pieno$ that are essential and the ones that are not. We want to investigate in more detail the saddles $\x\in\cS_{is}(\vuoto,\pieno)\setminus(\partial\cC_{\pieno}^{\vuoto}(\gi)\cup\partial\cC_{\vuoto}^{\pieno}(\gi-H(\pieno))\cup\cigeo)$ visited after crossing the set $\cC^B_{is}$.

\bp{isoselless}
Any saddle $\x$ that is neither in $\cC_{is}^*$, nor in the boundary of the cycle $\cC_{\pieno}^{\vuoto}(\G_{is}^*)$, nor in $\partial\cC_{\vuoto}^{\pieno}(\G_{is}^*-H(\pieno))\setminus\widetilde{K}_{is}$, such that $\t_{\x}\geq\t_{\cC_{is}^B}$ can be essential or not. For those essential, we obtain the following characterization:
\be{}
\ba{ll}
\cG_{is}(\vuoto,\pieno)\cap\cS_{is}(\vuoto,\pieno)\setminus(\partial\cC_{\pieno}^{\vuoto}(\G_{is}^*)\cup(\partial\cC_{\vuoto}^{\pieno}(\G_{is}^*-H(\pieno))\setminus\widetilde{K}_{is})\cup\cC_{is}^*)=\qquad\qquad\qquad\\
\quad\qquad\qquad\qquad\qquad\qquad\qquad\qquad\qquad=\displaystyle\bigcup_{ i=0}^{3}\bigcup_{\a}\mathscr{I}_i^{\a}\cup\bigcup_{i=0}^{2}\bigcup_{\a,\a'}\bigcup_{k,k'}\mathscr{I}_{k,k',i}^{\a,\a'}\cup\displaystyle\bigcup_{i=-1}^{2}\bigcup_{\a,\a'}\mathscr{I}_i^{\a,\a'}
\ea
\ee
\ep

\noindent
We refer to Section \ref{dim4} for the proof of the Proposition \ref{isoselless}.

\medskip
\begin{proof*}{\bf of the main Theorem \ref{giso}} 
	By Corollary \ref{corstrategy} we know that the saddles of the first and second type, defined in Definitions \ref{sellesigma} and \ref{sellezeta} respectively, are unessential. By Propositions \ref{selle3} and \ref{isoselless} we have the characterization of the essential saddles of the third type in Section \ref{S6.4}. Use Proposition \ref{c*contenuto} to get the claim.
\end{proof*}

\subsection{Proof of Proposition \ref{isoselless}}
\label{dim4}
We recall Definitions \ref{translation} and \ref{trenino} for the definitions of the $1$-translation of a bar and for the sliding of a bar around a frame-angle respectively and that $d(\cdot,\cdot)$ denotes the lattice distance. In order to prove Proposition \ref{isoselless} we need the following lemma.

\bl{isotrasless}
\bi
\item[(i)] Starting from $\h\in\mathscr{I}_0^{\a,\a'}$ (resp.\ $\h\in\mathscr{I}_{-1}^{\a,\a'}$), if the free particle is attached in $\partial^+\hbox{CR}(\h)$ obtaining the configuration $\h'$, then the following saddles obtained via a $1$-translation of any bar are essential and in $\mathscr{I}_0^{\a}\cup\mathscr{I}_1^{\a}$ (resp.\ in $\mathscr{I}_0^{\a}$). Moreover, all the saddles in $\mathscr{I}_0^{\a}\cup\mathscr{I}_1^{\a}$ can be obtained from a $\h\in\mathscr{I}_{-1}^{\a,\a'}\cup\mathscr{I}_0^{\a,\a'}$ via a $1$-translation of a bar. In particular, starting from $\h\in\cigeo$, if the free particle is attached in a bad site obtaining $\h^B\in\cC_{is}^B$, then the following saddles obtained via a $1$-translation of any bar are essential. These saddles are in $\mathscr{I}_0^{\a}$ if $\hat\h\in\bar\cD_{is}$ and in $\mathscr{I}_0^{\a}\cup\mathscr{I}_1^{\a}$ if $\hat\h\in\widetilde\cD_{is}$. 
\item[(ii)] Starting from $\h\in\mathscr{I}_{1}^{\a,\a'}$, if the free particle is attached in $\partial^+\hbox{CR}(\h)$ obtaining the configuration $\h'$, then the following saddles obtained via a $1$-translation of any bar are essential and in $\mathscr{I}_1^{\a}\cup\mathscr{I}_2^{\a}$. Moreover, all the saddles in $\mathscr{I}_1^{\a}\cup\mathscr{I}_2^{\a}$ can be obtained from a $\h\in\mathscr{I}_{1}^{\a,\a'}$ via a $1$-translation of a bar.
\item[(iii)] Starting from $\h\in\mathscr{I}_{2}^{\a,\a'}$, if the free particle is attached in $\partial^+\hbox{CR}(\h)$ obtaining the configuration $\h'$, then the following saddles obtained via a $1$-translation of any bar are essential and in $\mathscr{I}_2^{\a}\cup\mathscr{I}_3^{\a}$. Moreover, all the saddles in $\mathscr{I}_2^{\a}\cup\mathscr{I}_3^{\a}$ can be obtained from a $\h\in\mathscr{I}_{2}^{\a,\a'}$ via a $1$-translation of a bar.
\ei
\el

\noindent
The proof of the lemma is postponed to Section \ref{dimlemmiso}.

\bigskip
\begin{proof*}{\bf of Proposition \ref{isoselless}}
	Consider a configuration $\h\in\cigeo(2)$ such that $\h=(\hat\h,x)$, with $\hat\h\in\cD_{is}$ and $x$ the site of the free particle such that $d(\hat\h,x)=2$. By hypotheses we have that the free particle is attached in a bad site obtaining a configuration $\h'\in\cC_{is}^B$ (see Figure \ref{fig:fig1}(a)-(b) for a possible pair of configuration $(\h,\h')$). Due to \cite[Theorem 5.1]{MNOS}, our strategy consists in characterizing the essential saddles that could be visited after attaching the free particle in a bad site. We consider the following cases:
	
	\begin{description} 
		\item[Case 1.] $\hat\h\in\bar\cD_{is}$; 
		\item[Case 2.] $\hat\h\in\widetilde\cD_{is}$.
	\end{description}

	Note that from case 1 one can go to the other cases and viceversa, but since the path has to reach $\pieno$ this back and forth must end in a finite number of steps.

	\medskip
	\noindent
	{\bf Case 1.} Let $\hat\h\in\bar\cD_{is}$, thus by \cite[Theorem 1.4.1]{BHN} we know that $\hat\h$ consists in an $(l_c-2)\times(l_c-2)$ square with four bars $B^{\alpha}(\h)$, with $\alpha\in\{n,e,w,s\}$, attached to its four sides satisfying 
	\be{}
	1\leq|B^{\a}(\h)|\leq l_c, \qquad \quad \displaystyle\sum_\a |B^\a(\h)|-\displaystyle\sum_{\a\a'\in\{nw,ne,sw,se\}}|c^{\a\a'}(\h)|=3l_c-3. 
	\ee
	First, note that at most three frame-angles in $\partial^-\hbox{CR}(\hat\h)$ can be occupied, otherwise $|\partial^-\hbox{CR}(\hat\h)|=4l_c-4>3l_c-3$, which is absurd. Thus we consider separately the following cases:
	
	\begin{itemize}
		\item[A.] three frame-angles in $\partial^-\hbox{CR}(\hat\h)$ are occupied;
		\item[B.] two frame-angles in $\partial^-\hbox{CR}(\hat\h)$ are occupied;
		\item[C.] one frame-angle in $\partial^-\hbox{CR}(\hat\h)$ is occupied;
		\item[D.] no frame-angle in $\partial^-\hbox{CR}(\hat\h)$ is occupied.
	\end{itemize}

	\setlength{\unitlength}{1.1pt}
	\begin{figure}
		\begin{picture}(380,30)(-30,30)
			\thinlines 
			\qbezier[20](30,0)(60,0)(90,0)
			\qbezier[20](30,0)(30,30)(30,60) 
			\qbezier[20](30,60)(60,60)(90,60)
			\qbezier[20](90,0)(90,30)(90,60)
			\put(90,35){\line(1,0){5}}
			\put(95,35){\line(0,1){30}}
			\put(95,65){\line(-1,0){70}}
			\put(25,65){\line(0,-1){70}}
			\put(25,-5){\line(1,0){40}}
			\put(65,-5){\line(0,1){5}}
			\put(65,0){\line(1,0){25}}
			\put(90,0){\line(0,1){35}}
			\put(85,70){\line(1,0){5}}
			\put(85,75){\line(1,0){5}}
			\put(85,70){\line(0,1){5}}
			\put(90,70){\line(0,1){5}}
			\put(50,-27){(a)}
			\put(33,30){\tiny{$(l_c-2)\times(l_c-2)$}}
			
			\thinlines 
			\qbezier[20](160,0)(190,0)(220,0)
			\qbezier[20](160,0)(160,30)(160,60) 
			\qbezier[20](160,60)(190,60)(220,60)
			\qbezier[20](220,0)(220,30)(220,60)
			\put(220,35){\line(1,0){5}}
			\put(225,35){\line(0,1){30}}
			\put(225,65){\line(-1,0){5}}
			\put(220,65){\line(0,1){5}}
			\put(220,70){\line(-1,0){5}}
			\put(215,70){\line(0,-1){5}}
			\put(215,65){\line(-1,0){60}}
			\put(155,65){\line(0,-1){70}}
			\put(155,-5){\line(1,0){40}}
			\put(195,-5){\line(0,1){5}}
			\put(195,0){\line(1,0){25}}
			\put(220,0){\line(0,1){35}}
			\put(180,-27){(b)}
			\put(163,30){\tiny{$(l_c-2)\times(l_c-2)$}}
			
			\thinlines
			\qbezier[20](290,0)(320,0)(345,0)
			\qbezier[20](290,0)(290,30)(290,65) 
			\qbezier[20](290,65)(320,65)(345,65)
			\qbezier[20](345,0)(345,30)(345,65)
			\put(345,0){\line(1,0){5}}
			\put(350,0){\line(0,1){70}}
			\put(350,70){\line(-1,0){30}}
			\put(320,70){\line(0,-1){5}}
			\put(320,65){\line(-1,0){35}}
			\put(285,65){\line(0,-1){70}}
			\put(285,-5){\line(1,0){40}}
			\put(325,-5){\line(0,1){5}}
			\put(325,0){\line(1,0){20}}
			\put(355,65){\line(1,0){5}}
			\put(360,65){\line(0,1){5}}
			\put(360,70){\line(-1,0){5}}
			\put(355,70){\line(0,-1){5}}
			\put(310,-27){(c)}
			\put(291,30){\tiny{$(l_c-3)\times(l_c-1)$}}
		\end{picture}
		\vskip 2.cm
		\caption{Here we depict a possible configuration $\h$ in (a), its corresponding $\h'$ in (b) and the configuration obtained from $\h'$ after the sliding of the bar $B^e(\h)$ around the frame-angle $c^{en}(\h')$ in (c) for the case 1A.}
		\label{fig:fig1}
	\end{figure}
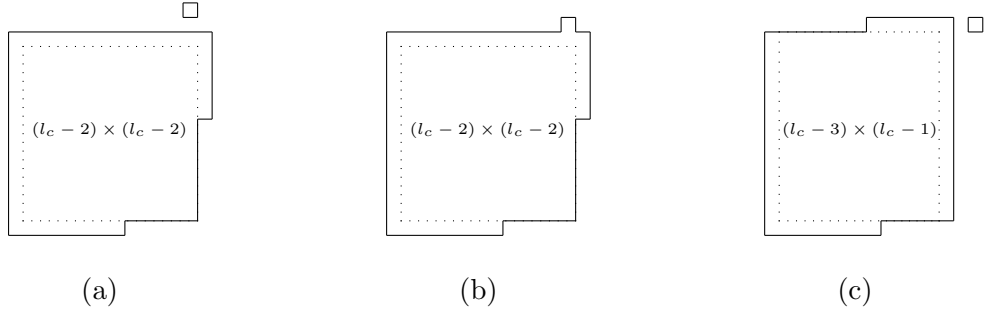

	\medskip
	{\bf Case 1A.} Without loss of generality we consider $\h$ as in Figure \ref{fig:fig1}(a). If we are considering the case in which a sequence of $1$-translations of a bar is possible and takes place, then by Lemma \ref{isotrasless}(i) the saddles that are crossed are essential and in $\mathscr{I}_{0}^{\a}$. If a sequence of $1$-translations of a bar takes place in such a way that the last configuration has at most two occupied frame-angles, then we can reduce our proof to the cases B, C and D below. Thus we are left to analyze the case in which there is the activation of a sliding of a bar around a frame-angle. Consider again, for example, $\h$ as in Figure \ref{fig:fig1}(a). If the free particle is attached to the bar $B^e(\h)$, then it is not possible to slide the bar $B^n(\h)$ around the frame-angle $c^{ne}(\h')$, since the condition (\ref{condtrenino}) is not satisfied. Thus by Lemma \ref{trenini}(ii) we know that the saddles that could be crossed are unessential. If the free particle is attached to the bar $B^s(\h)$, we conclude as before. If the free particle  is attached to the bar $B^n(\h)$ (see Figure \ref{fig:fig1}(b)), then it is not possible to slide the bar $B^w(\h)$ around the frame-angle $c^{wn}(\h')$ because the condition (\ref{condtrenino}) is not satisfied and thus we can conclude as before. The unique possibility is to slide the bar $B^e(\h)$ around the frame-angle $c^{en}(\h')$ if $|B^e(\h)|<|B^n(\h)|$, otherwise (\ref{condtrenino}) is not satisfied. The saddles that are possibly visited (except the last one) are in $\mathscr{I}_{k,k',0}^{\a,\a'}$ and by Lemma \ref{trenini}(i) they are essential. The last configuration visited during this sliding of a bar is depicted in Figure \ref{fig:fig1}(c). It belongs to $\cigeo$, indeed the cluster is in $\widetilde\cD_{is}$ and therefore the saddles that could be crossed starting from it will be investigated in case 2. If the free particle is attached to the bar $B^w(\h)$, we conclude in a similar way as before. This concludes case 1A.
	
	\medskip
	{\bf Case 1B.} We consider separately the following subcases:
	
	\bi
	\item[(i)] the two occupied frame-angles are $c^{\a\a'}(\h)$ and $c^{\a''\a'''}(\h)$, with all the indeces $\a,\a',\a'',\a'''$ different between each other (see Figure \ref{fig:fig2}(a));
	\item[(ii)] the two occupied frame-angles are $c^{\a\a'}(\h)$ and $c^{\a'\a''}(\h)$, with $\a\neq\a''$ (see Figure \ref{fig:fig3}(a)).
	\ei

	\setlength{\unitlength}{1.1pt}
	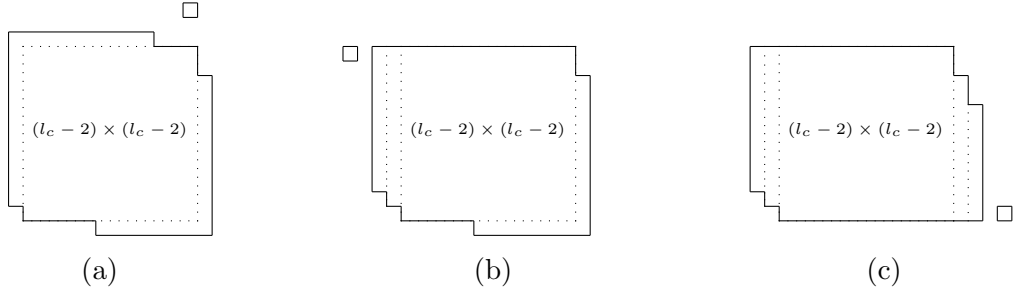
\begin{figure}
		\begin{picture}(380,30)(-30,30)
			\thinlines 
			\qbezier[20](30,0)(60,0)(90,0)
			\qbezier[20](30,0)(30,30)(30,60) 
			\qbezier[20](30,60)(60,60)(90,60)
			\qbezier[20](90,0)(90,30)(90,60)
			\put(90,60){\line(-1,0){15}}
			\put(75,60){\line(0,1){5}}
			\put(75,65){\line(-1,0){50}}
			\put(25,65){\line(0,-1){60}}
			\put(25,5){\line(1,0){5}}
			\put(30,5){\line(0,-1){5}}
			\put(30,0){\line(1,0){25}}
			\put(55,-5){\line(0,1){5}}
			\put(55,-5){\line(1,0){40}}
			\put(95,-5){\line(0,1){55}}
			\put(95,50){\line(-1,0){5}}
			\put(90,50){\line(0,1){10}}
			\put(85,70){\line(1,0){5}}
			\put(85,75){\line(1,0){5}}
			\put(85,70){\line(0,1){5}}
			\put(90,70){\line(0,1){5}}
			\put(50,-20){(a)}
			\put(33,30){\tiny{$(l_c-2)\times(l_c-2)$}}
			
			\thinlines 
			\qbezier[20](160,0)(190,0)(220,0)
			\qbezier[20](160,0)(160,30)(160,60) 
			\qbezier[20](160,60)(190,60)(220,60)
			\qbezier[20](220,0)(220,30)(220,60)
			\qbezier[20](155,60)(155,27.5)(155,5)
			\put(220,60){\line(-1,0){70}}
			\put(150,60){\line(0,-1){50}}
			\put(150,10){\line(1,0){5}}	
			\put(155,10){\line(0,-1){5}}
			\put(155,5){\line(1,0){5}}
			\put(160,5){\line(0,-1){5}}
			\put(160,0){\line(1,0){25}}
			\put(185,-5){\line(0,1){5}}
			\put(185,-5){\line(1,0){40}}
			\put(225,-5){\line(0,1){55}}
			\put(225,50){\line(-1,0){5}}
			\put(220,50){\line(0,1){10}}
			\put(140,55){\line(1,0){5}}
			\put(140,60){\line(1,0){5}}
			\put(145,55){\line(0,1){5}}
			\put(140,55){\line(0,1){5}}
			\put(185,-20){(b)}
			\put(163,30){\tiny{$(l_c-2)\times(l_c-2)$}}
			
			\thinlines
			\qbezier[20](290,0)(320,0)(350,0)
			\qbezier[20](290,0)(290,30)(290,60) 
			\qbezier[20](290,60)(320,60)(350,60)
			\qbezier[20](350,0)(350,30)(350,60)
			\qbezier[20](285,60)(285,27.5)(285,5)
			\qbezier[14](355,0)(355,20)(355,40)
			\put(350,60){\line(-1,0){70}}
			\put(280,60){\line(0,-1){50}}
			\put(280,10){\line(1,0){5}}	
			\put(285,10){\line(0,-1){5}}
			\put(285,5){\line(1,0){5}}
			\put(290,5){\line(0,-1){5}}
			\put(290,0){\line(1,0){35}}
			\put(325,0){\line(1,0){35}}
			\put(360,0){\line(0,1){40}}
			\put(360,40){\line(-1,0){5}}
			\put(355,40){\line(0,1){10}}
			\put(355,50){\line(-1,0){5}}
			\put(350,50){\line(0,1){10}}
			\put(365,0){\line(1,0){5}}
			\put(365,5){\line(1,0){5}}
			\put(370,0){\line(0,1){5}}
			\put(365,0){\line(0,1){5}}
			\put(320,-20){(c)}
			\put(293,30){\tiny{$(l_c-2)\times(l_c-2)$}}
		\end{picture}
		\vskip 2. cm
		\caption{Case 1B(i): we depict a possible starting configuration $\h\in\cigeo$ in (a), the configuration $\widetilde\h$ obtained from $\h$ after the sliding of the bar $B^n(\h)$ around the frame-angle $c^{nw}(\h')$ in (b) and the configuration $\bar\h$ obtained from $\widetilde\h$ after the sliding of the bar $B^s(\widetilde\h)$ around the frame-angle $c^{se}(\h'')$ in (c).}
		\label{fig:fig2}
	\end{figure}

	\medskip
	{\bf Case 1B(i).} Without loss of generality we consider $\h$ as in Figure \ref{fig:fig2}(a). If we are considering the case in which a sequence of $1$-translations of a bar is possible and takes place, then by Lemma \ref{isotrasless}(i) the saddles that are crossed are essential and they are in $\mathscr{I}_{0}^{\a}$. If at least one bar is full, it is possible to activate a sequence of $1$-translations of a bar in order to obtain either two occupied frame-angles with a bar in common or three occupied frame-angles. For example, in Figure \ref{fig:fig2}(a), if the bar $B^s(\h)$ is full, one could attach the free particle to $B^e(\h)$ and translate the bar $B^w(\h)$ in order to have the frame-angle $c^{sw}(\h)$ occupied. In both situations the saddles visited up to this point are essential by Lemma \ref{isotrasless}(i), while the saddles that follow are analyzed in case 1B(ii) and 1A respectively. Thus we can reduce our proof to the case in which there is no translation of a bar and therefore we need to consider only the sliding of a bar around a frame-angle. We may assume without loss of generality that $|B^n(\h)|<|B^w(\h)|$ and $|B^s(\h)|<|B^e(\h)|$, indeed the other cases can be treated with the same argument. By Lemma \ref{trenini} to obtain essential saddles there are one of the following possibilities: attach the free particle to the bar $B^w(\h)$ (resp.\ $B^e(\h)$) and then slide the bar $B^n(\h)$ (resp.\ $B^s(\h)$) around the frame-angle $c^{nw}(\h')$ (resp.\ $c^{se}(\h')$). Assume first that the free particle is attached to $B^w(\h)$. By Lemma \ref{trenini}(i) the saddles that are possibly visited are essential and, except the last one, they are in $\mathscr{I}_{k,k',0}^{\a,\a'}$. The last configuration visited during this sliding of a bar is $\widetilde\h\in\mathscr{I}_{0}^{\a,\a'}$ and it is depicted in Figure \ref{fig:fig2}(b). Starting from $\widetilde\h$, 
	the unique possibility to visit essential saddles is to attach a free particle in $\partial^+\hbox{CR}(\widetilde\h)$ and then either activate a sequence of $1$-translations of bars or slide a bar around a frame-angle. In the first case, by Lemma \ref{isotrasless}(i) the saddles that are possibly visited are essential and in $\mathscr{I}_{0}^{\a}\cup\mathscr{I}_1^{\a}$. In the latter case, the unique possibility is to attach the free particle to the bar $B^e(\widetilde\h)$ obtaining a configuration $\h''$, and then slide $B^s(\widetilde\h)$ around the frame-angle $c^{se}(\h'')$. By Lemma \ref{trenini}(i) the saddles that are possibly visited are essential and, except the last one, they are in $\mathscr{I}_{k,k',1}^{\a,\a'}$. The last configuration visited during this sliding of a bar is in $\mathscr{I}_{1}^{\a,\a'}$ and it is depicted in Figure \ref{fig:fig2}(c). Starting from this configuration it is impossible to slide a bar around any frame-angle, thus by Lemma \ref{trenini}(ii) the saddles that possibly will be crossed if the sliding of a bar is initiated are unessential. If a sequence of $1$-translations of bars takes place, by Lemma \ref{isotrasless}(ii) the saddles that could be crossed are essential and in $\mathscr{I}_1^{\a}\cup\mathscr{I}_2^{\a}$.
	
	Note that if $|B^w(\h)|<|B^n(\h)|$ and/or $|B^e(\h)|<|B^s(\h)|$ a similar argument can be used. This concludes case 1B(i).

	\setlength{\unitlength}{1.1pt}
	\begin{figure}
		\begin{picture}(380,40)(-30,40)
			\thinlines 
			\qbezier[20](30,0)(60,0)(90,0)
			\qbezier[20](30,0)(30,30)(30,60) 
			\qbezier[20](30,60)(60,60)(90,60)
			\qbezier[20](90,0)(90,30)(90,60)
			\put(95,65){\line(-1,0){55}}
			\put(40,65){\line(-1,0){15}}
			\put(25,65){\line(0,-1){35}}
			\put(25,30){\line(1,0){5}}
			\put(30,30){\line(0,-1){30}}
			\put(30,0){\line(1,0){5}}
			\put(35,0){\line(0,-1){5}}
			\put(35,-5){\line(1,0){55}}
			\put(90,-5){\line(0,1){5}}
			\put(90,0){\line(0,1){20}}
			\put(95,20){\line(-1,0){5}}
			\put(95,20){\line(0,1){45}}
			\put(85,70){\line(1,0){5}}
			\put(85,75){\line(1,0){5}}
			\put(85,70){\line(0,1){5}}
			\put(90,70){\line(0,1){5}}
			\put(55,-18){(a)}
			\put(33,30){\tiny{$(l_c-2)\times(l_c-2)$}}
			
			\thinlines 
			\qbezier[20](160,0)(190,0)(215,0)
			\qbezier[20](160,0)(160,30)(160,65) 
			\qbezier[20](160,65)(190,65)(215,65)
			\qbezier[20](215,0)(215,30)(215,65)
			\put(220,65){\line(-1,0){30}}
			\put(190,65){\line(0,1){5}}
			\put(190,70){\line(-1,0){35}}
			\put(155,70){\line(0,-1){70}}
			\put(155,0){\line(1,0){5}} 
			\put(160,0){\line(0,-1){5}}
			\put(160,-5){\line(1,0){55}}
			\put(215,-5){\line(0,1){5}}
			\put(220,65){\line(0,-1){45}}
			\put(215,20){\line(1,0){5}}
			\put(215,20){\line(0,-1){20}}
			\put(145,65){\line(1,0){5}}
			\put(150,65){\line(0,1){5}}
			\put(150,70){\line(-1,0){5}}
			\put(145,70){\line(0,-1){5}}
			\put(183,-18){(b)}
			\put(161,30){\tiny{$(l_c-3)\times(l_c-1)$}}
			
			\thinlines
			\qbezier[20](290,0)(320,0)(345,0)
			\qbezier[20](290,0)(290,30)(290,65) 
			\qbezier[20](290,65)(320,65)(345,65)
			\qbezier[20](345,0)(345,30)(345,65)
			\put(350,70){\line(-1,0){45}}
			\put(305,70){\line(0,-1){5}}
			\put(305,65){\line(-1,0){20}}
			\put(285,65){\line(0,-1){35}}
			\put(285,30){\line(1,0){5}}
			\put(290,30){\line(0,-1){30}}
			\put(290,0){\line(1,0){5}} 
			\put(295,0){\line(0,-1){5}}
			\put(295,-5){\line(1,0){55}}
			\put(350,-5){\line(0,1){5}}
			\put(350,70){\line(0,-1){70}}
			\put(355,65){\line(1,0){5}}
			\put(360,65){\line(0,1){5}}
			\put(360,70){\line(-1,0){5}}
			\put(355,70){\line(0,-1){5}}
			\put(313,-18){(c)}
			\put(291,30){\tiny{$(l_c-3)\times(l_c-1)$}}
		\end{picture}
		\vskip 2. cm
		\caption{Case 1B(ii): we depict a possible starting configuration $\h\in\cigeo$ in (a), the configuration $\widetilde\h$ obtained from $\h$ after the sliding of the bar $B^w(\h)$ around the frame-angle $c^{wn}(\h')$ in (b) and the configuration $\bar\h$ obtained from $\h$ after the sliding of the bar $B^e(\h)$ around the frame-angle $c^{en}(\h'')$ in (c).}
		\label{fig:fig3}
	\end{figure}
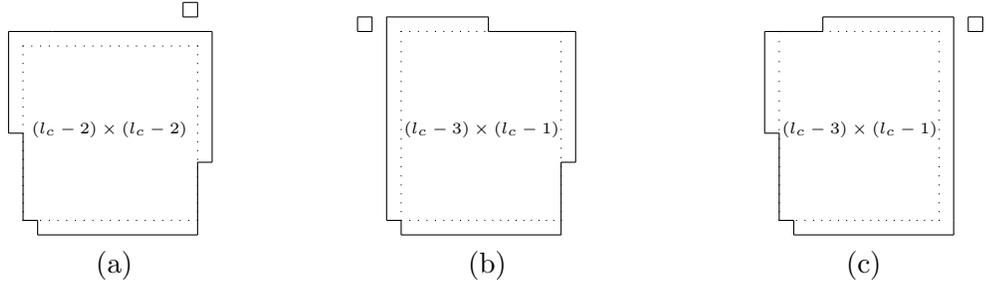

	\medskip
	{\bf Case 1B(ii).} Without loss of generality we consider $\h$ as in Figure \ref{fig:fig3}(a). If we are considering the case in which a sequence of $1$-translations of a bar is possible and takes place, then by Lemma \ref{isotrasless}(i) the saddles that are crossed are essential and they are in $\mathscr{I}_{0}^{\a}$. If one bar among $B^w(\h)$ and $B^e(\h)$ is full, it is possible to activate a sequence of  $1$-translations of the bar $B^s(\h)$ in order to have three occupied frame-angles. This situation has already been analyzed in case 1A. Thus we can reduce our proof to the case in which there is no translation of a bar and therefore we need to consider only the sliding of a bar around a frame-angle. If the free particle is attached to the bar $B^s(\h)$, since it is not possible to slide a bar around any frame-angle, by Lemma \ref{trenini}(ii) we know that the saddles that could be crossed are unessential. If the free particle is attached to one bar among $B^w(\h)$ and $B^e(\h)$, then it is not possible to complete the sliding of the bar $B^n(\h)$ around the frame-angle $c^{nw}(\h')$ or $c^{ne}(\h')$. Thus by Lemma \ref{trenini}(ii) the saddles that could be crossed are unessential. If the free particle is attached to the bar $B^n(\h)$, then it is possible to slide the bar $B^w(\h)$ or $B^e(\h)$ around the frame-angle $c^{wn}(\h')$ or $c^{en}(\h')$ respectively. Thus by Lemma \ref{trenini}(i) we know that the saddles that could be crossed are essential and they are in $\mathscr{I}_{k,k',0}^{\a,\a'}$, except the last one that is in $\cigeo$ (see Figure \ref{fig:fig3}(b)-(c)). Hence the saddles that could be crossed starting from such configuration will be analyzed in case 2. This concludes case 1B(ii).

	\setlength{\unitlength}{1.1pt}
	\begin{figure}
		\begin{picture}(380,40)(-30,40)
			\thinlines 
			\qbezier[20](30,0)(60,0)(90,0)
			\qbezier[20](30,0)(30,30)(30,60) 
			\qbezier[20](30,60)(60,60)(90,60)
			\qbezier[20](90,0)(90,30)(90,60)
			\put(75,65){\line(-1,0){50}}
			\put(75,65){\line(0,-1){5}}
			\put(75,60){\line(1,0){15}}
			\put(25,65){\line(0,-1){30}}
			\put(25,35){\line(1,0){5}}
			\put(30,35){\line(0,-1){40}}
			\put(30,-5){\line(1,0){60}}
			\put(90,-5){\line(0,1){5}}
			\put(90,0){\line(1,0){5}}
			\put(95,0){\line(0,1){60}}
			\put(95,60){\line(-1,0){5}}
			\put(85,70){\line(1,0){5}}
			\put(85,75){\line(1,0){5}}
			\put(85,70){\line(0,1){5}}
			\put(90,70){\line(0,1){5}}
			\put(55,-18){(a)}
			\put(33,30){\tiny{$(l_c-2)\times(l_c-2)$}}
			
			\thinlines 
			\qbezier[20](160,0)(190,0)(220,0)
			\qbezier[20](160,0)(160,30)(160,60) 
			\qbezier[20](160,60)(190,60)(220,60)
			\qbezier[20](220,0)(220,30)(220,60)
			\qbezier[9](160,65)(175,65)(190,65)
			\put(205,65){\line(0,-1){5}}
			\put(205,65){\line(-1,0){15}}
			\put(190,65){\line(0,1){5}}
			\put(190,70){\line(-1,0){30}}
			\put(160,70){\line(0,-1){25}}
			\put(205,60){\line(1,0){15}}
			\put(160,45){\line(0,-1){50}}
			\put(160,-5){\line(1,0){60}}
			\put(220,-5){\line(0,1){5}}
			\put(220,0){\line(1,0){5}}
			\put(225,0){\line(0,1){60}}
			\put(225,60){\line(-1,0){5}}
			\put(150,65){\line(1,0){5}}
			\put(150,70){\line(1,0){5}}
			\put(150,65){\line(0,1){5}}
			\put(155,65){\line(0,1){5}}
			\put(185,-18){(b)}
			\put(163,30){\tiny{$(l_c-2)\times(l_c-2)$}}
			
			\thinlines
			\qbezier[20](290,0)(320,0)(350,0)
			\qbezier[20](290,0)(290,30)(290,60) 
			\qbezier[20](290,60)(320,60)(350,60)
			\qbezier[20](350,0)(350,30)(350,60)
			\put(340,65){\line(-1,0){50}}
			\put(340,65){\line(0,-1){5}}
			\put(340,60){\line(1,0){10}}
			\put(290,65){\line(0,-1){10}}
			\put(285,55){\line(1,0){5}}
			\put(285,55){\line(0,-1){50}}
			\put(285,5){\line(1,0){5}}
			\put(290,5){\line(0,-1){5}}
			\put(290,0){\line(1,0){10}}
			\put(300,-5){\line(0,1){5}}
			\put(300,-5){\line(1,0){50}}
			\put(350,-5){\line(0,1){15}}
			\put(350,10){\line(1,0){5}}
			\put(355,10){\line(0,1){45}}
			\put(355,55){\line(-1,0){5}}
			\put(355,55){\line(-1,0){5}}
			\put(350,55){\line(0,1){5}}
			\put(345,70){\line(1,0){5}}
			\put(345,75){\line(1,0){5}}
			\put(345,70){\line(0,1){5}}
			\put(350,70){\line(0,1){5}}
			\put(315,-18){(c)}
			\put(293,30){\tiny{$(l_c-2)\times(l_c-2)$}}

		\end{picture}
		\vskip 2. cm
		\caption{Case 1C: in (a) we depict a possible starting configuration $\h\in\cigeo$ and in (b) the configuration $\widetilde\h$ obtained from $\h$ after the sliding of the bar $B^w(\h)$ around the frame-angle $c^{wn}(\h')$. Case 1D: in (c) we depict a possible starting configuration in $\cigeo$.}
		\label{fig:fig4}
	\end{figure}

	\medskip
	{\bf Case 1C.} Without loss of generality we consider $\h$ as in Figure \ref{fig:fig4}(a). If we are a considering the case in which a sequence of $1$-translations of a bar is possible and takes place, then by Lemma \ref{isotrasless}(i) the saddles that are crossed are in $\mathscr{I}_{0}^{\a}$. Starting from this configuration it is possible to obtain two occupied frame-angles: this situation has been already analyzed in case 1B. Thus we can reduce our proof to the case in which there is no $1$-translation of a bar and therefore there is the activation of a sliding of a bar around a frame-angle. If the free particle is attached to the bar $B^e(\h)$ or $B^s(\h)$, since it is not possible to complete any sliding of a bar around a frame-angle at cost $U$, by Lemma \ref{trenini}(ii) we know that the saddles that could be crossed are unessential. If $|B^w(\h)|<|B^n(\h)|$ and the free particle is attached to the bar $B^n(\h)$, then it is possible to slide the bar $B^w(\h)$ around the frame-angle $c^{wn}(\h')$. Thus by Lemma \ref{trenini}(i) the saddles that could be crossed are essential and, except the last one, they are in $\mathscr{I}_{k,k',0}^{\a,\a'}$. Note that the last configuration is in $\mathscr{I}_{0}^{\a,\a'}$ (see Figure \ref{fig:fig4}(b)). Starting from such a configuration, since  it is not possible to complete any sliding of a bar around a frame-angle, by Lemma \ref{trenini}(ii) we know that the saddles that could be visited are unessential unless a sequence of $1$-translations of bars takes place. In this case, by Lemma \ref{isotrasless}(ii) the saddles that could be visited are essential and in $\mathscr{I}_0^{\a}\cup\mathscr{I}_1^{\a}$. The case $|B^n(\h)|<|B^w(\h)|$ in which the free particle is attached to the bar $B^w(\h)$ is analogue. This concludes case 1C.

	\medskip
	{\bf Case 1D.} Without loss of generality we consider $\h$ as in Figure \ref{fig:fig4}(c). If we are considering the case in which a sequence of $1$-translations of a bar is possible and takes place, then by Lemma \ref{isotrasless}(i) the saddles that are crossed are in $\mathscr{I}_{0}^{\a}$. Starting from this configuration, it is possible to obtain one or two occupied frame-angles: these situations have been already analyzed in cases 1C and 1B respectively. Thus we can reduce our proof to the case in which there is no $1$-translation of a bar and therefore there is the activation of a sliding of a bar around a frame-angle. If the free particle is attached to one of the bars, since it is not possible to complete any sliding of bar around a frame-angle, by Lemma \ref{trenini}(ii) we know that the saddles that could be crossed are unessential. This concludes case 1D.
	
	\bigskip
	\noindent
	{\bf Case 2.} Let $\hat\h\in\widetilde\cD_{is}$, thus by \cite[Theorem 1.4.1]{BHN} we know that $\hat\h$ consists of an $(l_c-3)\times(l_c-1)$ quasi-square with four bars $B^{\alpha}(\h)$, with $\alpha\in\{n,w,e,s\}$, attached to its four sides satisfying 
	\be{}
	1\leq|B^{\a}(\h)|,|B^{\a'}(\h)|\leq l_c+1, \qquad 1\leq|B^{\a''}(\h)|,|B^{\a'''}(\h)|\leq l_c-1, 
	\ee
	\noindent
	where either $\a,\a'\in\{n,s\}$ and $\a'',\a'''\in\{w,e\}$, or $\a,\a'\in\{w,e\}$ and $\a'',\a'''\in\{n,s\}$, and
	\be{}
	\displaystyle\sum_\a |B^\a(\h)|-\displaystyle\sum_{\a\a'\in\{nw,ne,sw,se\}}|c^{\a\a'}(\h)|=3l_c-2.
	\ee
	
	\noindent
	First, note that at most three frame-angles in $\partial^-\hbox{CR}(\hat\h)$ can be occupied, otherwise $|\partial^-\hbox{CR}(\hat\h)|=4l_c-4>3l_c-2$, which is absurd. By hypotheses we have that the free particle is attached in a bad site obtaining a configuration $\h'\in\cC_{is}^B$. We consider separately the following cases:
	
	\begin{itemize}
		\item[A.] three frame-angles in $\partial^-\hbox{CR}(\h)$ are occupied;
		\item[B.] two frame-angles in $\partial^-\hbox{CR}(\h)$ are occupied;
		\item[C.] one frame-angle in $\partial^-\hbox{CR}(\h)$ is occupied;
		\item[D.] no frame-angle in $\partial^-\hbox{CR}(\h)$ is occupied.
	\end{itemize}
	
	The argument used in case 2 is analogue to that used above in case 1 and is discussed in details in Appendix \ref{app}.
\end{proof*}

\subsection{Proof of Lemma \ref{isotrasless}}
\label{dimlemmiso}

\bpr
(i) Note that $H(\h')=\gi-U$, thus it is possible to translate bars at cost $U$ with cost $\leq\gi$. These saddles are in $\mathscr{I}_0^{\a}$ if $\h\in\mathscr{I}_{-1}^{\a,\a'}$ and in $\mathscr{I}_0^{\a}\cup\mathscr{I}_1^{\a}$ if $\h\in\mathscr{I}_0^{\a,\a'}$. To conclude, all the configurations in $\mathscr{I}_0^{\a}\cup\mathscr{I}_1^{\a}$ can be obtained from a configuration $\h\in\mathscr{I}_{-1}^{\a,\a'}\cup\mathscr{I}_0^{\a,\a'}$ via a $1$-translation of a bar. Since $\bar\cD_{is}^{fp}\subseteq\mathscr{I}_{-1}^{\a,\a'}$ and $\widetilde\cD_{is}^{fp}\subseteq\mathscr{I}_0^{\a,\a'}$, we get the particular case in which $\h'=\h^B\in\cC_{is}^B$ as claimed.

It remains to prove that the saddles in $\mathscr{I}_0^{\a}\cup\mathscr{I}_1^{\a}$ are essential. Let $\x_1,...,\x_m\in\mathscr{I}_0^{\a}\cup\mathscr{I}_1^{\a}$ the saddles visited during a $1$-translation of a bar. We want to prove that these saddles are essential (see Section \ref{modinddef} point 4 for the definition). Since we can repeat the following argument $m$ times, we may focus on a single configuration $\x_i$. Since $\cC_{is}^*$ is a gate for the transition and $\x_i\in\cS_{is}(\vuoto,\pieno)\setminus\cC_{is}^*$, we note that a path $\o\in(\vuoto\ra\pieno)_{opt}$ such that $\{\arg\max_{\o}H\}=\{\x_i\}$ does not exist. Thus our strategy consists in finding a path $\o\in(\vuoto\ra\pieno)_{opt}$ such that for any $\o'\in(\vuoto\ra\pieno)_{opt}$
\be{isocondess}
\o\cap\x_i\neq\emptyset, \ \{\hbox{arg max}_{\o'}H\}\nsubseteq\{\hbox{arg max}_{\o}H\}\setminus\{\x_i\}, \qquad i=1,...,m.
\ee

\begin{figure}
	\centering
	\begin{tikzpicture}[scale=0.45,transform shape]
		\draw (0,0) rectangle (16,16);
		\draw(-1,-1) rectangle (17,17);
		\node at (17.5,16) {\Huge{$\L$}};
		\node at (15.4,15) {\Huge{$\L_0$}};
		\draw (5,1)--(11,1);
		\draw (11,1)--(11,2);
		\draw(11,2)--(12,2);
		\draw(12,2)--(12,6);
		\draw(12,6)--(11,6);
		\draw(11,6)--(11,8);
		\draw(11,8)--(10,8);
		\draw(10,8)--(10,9);
		\draw(4,4)--(4,2);
		\draw (4,2)--(5,2);
		\draw (5,2)--(5,1);
		\draw [grigio,fill=grigio] (4,2) rectangle (11,8);
		\draw [grigio,fill=grigio] (4,8) rectangle (10,9);
		\draw [grigio,fill=grigio] (11,2) rectangle (12,6);
		\draw [grigio,fill=grigio] (5,1) rectangle (11,2);
		\draw[dashed] (4,2) rectangle (11,8);
		\draw (10,8)--(10,9);
		\draw [grigio,fill=grigio] (4,9) rectangle (9,10);
		\draw[dashed](10,9)--(4,9);
		\draw (4,9)--(4,10);
		\draw (4,10)--(9,10);
		\draw (9,10)--(9,9);
		\draw(9,9)--(10,9);
		\draw (4,9)--(4,4);
		\draw (4,4)--(4,2);
		\draw (4,2)--(5,2);
		\draw (5,2)--(5,1);
		\draw (5,1)--(11,1);
		\draw (11,1)--(11,2);
		\draw (11,2)--(12,2);
		\draw (12,2)--(12,6);
		\draw (12,6)--(11,6);
		\draw (11,6)--(11,8);
		\draw (11,8)--(10,8);
		\draw[dashed, fill=grigio] (2,9) rectangle (3,10);
		\node at (2.5,9.5) {$\bullet$};
		\draw[dashed] (2,10) rectangle (3,11);
		\node at (2.5,10.5) {$\bullet$};
		\draw[dashed] (2,11) rectangle (3,12);
		\node at (2.5,11.5) {$\bullet$};
		\draw[dashed] (3,11) rectangle (4,12);
		\node at (3.5,11.5) {$\bullet$};
		\draw[dashed] (4,11) rectangle (5,12);
		\node at (4.5,11.5) {$\bullet$};
		\draw[dashed] (5,11) rectangle (6,12);
		\node at (5.5,11.5) {$\bullet$};
		\draw[dashed] (6,11) rectangle (7,12);
		\node at (6.5,11.5) {$\bullet$};
		\draw[dashed] (7,11) rectangle (8,12);
		\node at (7.5,11.5) {$\bullet$};
		\draw[dashed] (8,11) rectangle (9,12);
		\node at (8.5,11.5) {$\bullet$};
		\draw[dashed] (9,11) rectangle (10,12);
		\node at (9.5,11.5) {$\bullet$};
		\draw[dashed] (10,11) rectangle (11,12);
		\node at (10.5,11.5) {$\bullet$};
		\draw[dashed] (11,11) rectangle (12,12);
		\node at (11.5,11.5) {$\bullet$};
		\draw[dashed] (12,11) rectangle (13,12);
		\node at (12.5,11.5) {$\bullet$};
		\draw[dashed] (13,11) rectangle (14,12);
		\node at (13.5,11.5) {$\bullet$};
		\draw[dashed] (13,10) rectangle (14,11);
		\node at (13.5,10.5) {$\bullet$};
		\draw[dashed] (13,9) rectangle (14,10);
		\node at (13.5,9.5) {$\bullet$};
		\draw[dashed] (13,8) rectangle (14,9);
		\node at (13.5,8.5) {$\bullet$};
		\draw[dashed] (13,7) rectangle (14,8);
		\node at (13.5,7.5) {$\bullet$};
		\draw[dashed] (13,6) rectangle (14,7);
		\node at (13.5,6.5) {$\bullet$};
		\draw[dashed] (13,5) rectangle (14,6);
		\node at (13.5,5.5) {$\bullet$};
		\draw[dashed] (13,4) rectangle (14,5);
		\node at (13.5,4.5) {$\bullet$};
		\draw[dashed] (12,4) rectangle (13,5);
		\node at (12.5,4.5) {$\bullet$};
		
		\draw[->](12,7)--(12.5,4.7);
		\node at (12,7.5) {\Huge{$(i,j)$}};
		\draw[->](18,2.5)--(13.7,4.5);
		\node at (19.7,2.5) {\Huge{$(i+1,j)$}};

		\path
		(1.7,9.5) edge[bend left=70,->] node [left] {} (1.7,10.5)
		(1.7,10.7) edge[bend left=70,->] node [left] {} (1.7,11.5)
		(2.5,12.2) edge[bend left=70,->] node [left] {} (3.5,12.2)
		(3.7,12.2) edge[bend left=70,->] node [left] {} (4.5,12.2)
		(4.7,12.2) edge[bend left=70,->] node [left] {} (5.5,12.2)
		(5.7,12.2) edge[bend left=70,->] node [left] {} (6.5,12.2)
		(6.7,12.2) edge[bend left=70,->] node [left] {} (7.5,12.2)
		(7.7,12.2) edge[bend left=70,->] node [left] {} (8.5,12.2)
		(8.7,12.2) edge[bend left=70,->] node [left] {} (9.5,12.2)
		(9.7,12.2) edge[bend left=70,->] node [left] {} (10.5,12.2)
		(10.7,12.2) edge[bend left=70,->] node [left] {} (11.5,12.2)
		(11.7,12.2) edge[bend left=70,->] node [left] {} (12.5,12.2)
		(12.7,12.2) edge[bend left=70,->] node [left] {} (13.5,12.2)
		(14.2,11.5) edge[bend left=70,->] node [left] {} (14.2,10.5)
		(14.2,10.3) edge[bend left=70,->] node [left] {} (14.2,9.5)
		(14.2,9.3) edge[bend left=70,->] node [left] {} (14.2,8.5)
		(14.2,8.3) edge[bend left=70,->] node [left] {} (14.2,7.5)
		(14.2,7.3) edge[bend left=70,->] node [left] {} (14.2,6.5)
		(14.2,6.3) edge[bend left=70,->] node [left] {} (14.2,5.5)
		(14.2,5.3) edge[bend left=70,->] node [left] {} (14.2,4.5)
		(13.5,3.8) edge[bend left=70,->] node [left] {} (12.5,3.8);
	\end{tikzpicture}
	\vskip 0 cm
	\caption{Here we depict in grey the configuration $\h_1$ that consists of the union of a cluster $C(\h_1)$ and a free particle. The dotted unit squares represent the following positions of the free particle that moves as represented by the arrows on the left, north, east and south, until the particle is attached to the cluster in position $(i,j)$. The latter is the configuration $\h'$. To simplify the exposition we chose to use the path in the external frame of $\partial^+\hbox{CR}(C(\h_1))$.}
	\label{fig:camminobarranord}
\end{figure}
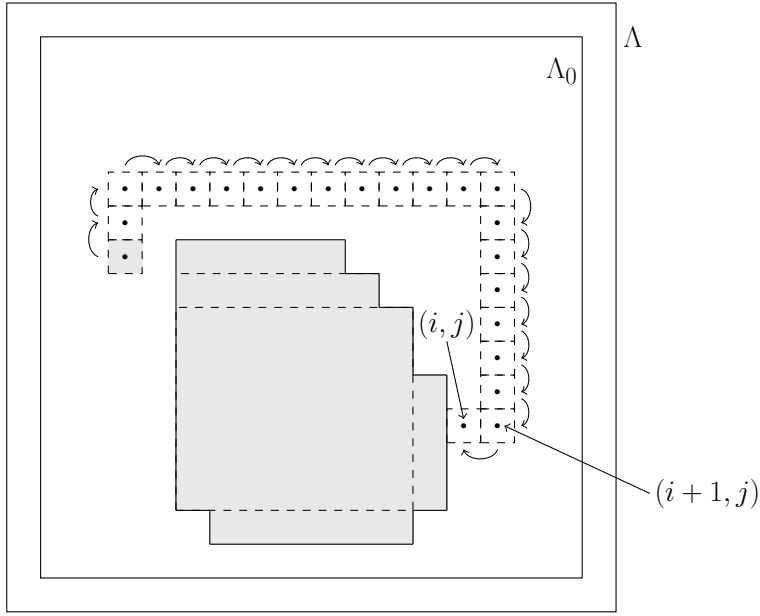

\noindent
Let $\h'$ be the configuration with a single cluster obtained as union of a cluster $C(\h_1)$ and a protuberance attached to one of its bars in a site with coordinates $(i,j)$ such that $H(\h_1)=\gi$ and $\h_1$ is the configuration with cluster $C(\h_1)$ and one free particle. More precisely, by (\ref{gatetrenino}) the configuration $\h_1$ is obtained by the configuration $\h_0$ via a sliding of the bar $B^{\a}(\h_0)$ around the frame-angle $c^{\a\a'}(\h_0)$. Without loss of generality we assume that that the protuberance of $\h'$ is attached to the bar $B^{e}(C(\h_1))$. We define the specific path $\o$ of the strategy above as
\be{} 
\o=(\vuoto,\o_1,...,\o_k,\g_1,...,\h_0,...\g_s,\h_1,...,\h_k,\h',\x_1,...,\x_m)\circ\bar\o,
\ee

\noindent
where $\o_1,...,\o_k\in\cC_{\pieno}^{\vuoto}(\gi)$, $\g_1,...,\g_s\in\cX$ such that $H(\g_i)\leq\gi$ for any $i=1,...,s$ and $\bar\o$ is a path that connects $\x_m$ to $\pieno$ such that $\max_{\s\in\bar\o}H(\s)\leq\G_{int}^*$. Moreover, the saddles $\h_1,...,\h_k$ are composed by the union of the cluster $C(\h_1)$ and a free particle such that the free particle is at distance $\geq3$ from a site in $\partial^-\hbox{CR}(\h_1)$ (see Figure \ref{fig:camminobarranord} for a picture of this situation in the case $\a=w$ and $\a'=n$). In particular, consider that $\h_1$ has the free particle in $\bar{B}_{2}(\h_1)$ and $\h_k$ in the site with coordinates $(i+1,j)\in\bar{B}_{2}(\h_1)$ using the assumption that the protuberance of $\h'$ is attached to $B^{e}(\h_1)$.

Now we show that for any $\o'$ the condition (\ref{isocondess}) is satisfied. If $\o'$ passes through the configuration $\x_i$, $\{\hbox{arg max}_{\o'}H\}\supseteq\{\x_i\}$, thus (\ref{isocondess}) is satisfied. Therefore we can assume that $\o'\cap\x_i=\emptyset$.  If $\o'$ crosses the set $\cS_{is}(\vuoto,\pieno)$ through a configuration $\h''$ such that $\h''\cap\o=\emptyset$, then the condition (\ref{isocondess}) holds. Thus we can reduce our analysis to $\o'$ that visits all the saddles $\h_1,...,\h_k$. Starting from $\h_k$, there are four allowed directions for moving the free particle. If we move it in the direction of the cluster (west in Figure \ref{fig:camminobarranord}), we deduce that the path $\o'$ visits the configuration $\h'$. For the other three choices, the free particle still remains free after the move, indeed by construction of the path $\o$, starting from $\h_k$ it is not possible to attach the free particle in $\partial^-\hbox{CR}(\h_1)$ via one step of the dynamics. Thus the path $\o'$ can visit either a saddle not already visited by $\o$ (south or east in Figure \ref{fig:camminobarranord}) or a saddle that has been already visited by $\o$ (north in Figure \ref{fig:camminobarranord}). In the first case, we obtain that (\ref{isocondess}) is satisfied. In the latter case, we can iterate this argument and, since $\o'$ goes from $\vuoto$ to $\pieno$, we can assume that the path $\o'$ visits the configuration $\h'$. From now on, starting from $\h'$, there are two possible scenarios:
\bi
\item[I.] $\o'$ activates the same $1$-translation of a bar as $\o$;
\item[II.] $\o'$ activates a different $1$-translation of a bar from $\o$.
\ei

In case I, since $\o'\cap\x_i=\emptyset$, we deduce that the $1$-translation of a bar has been stopped before hitting $\x_i$. Thus we can assume that $\o'$ comes back to $\h'$, otherwise the energy exceeds $\gi$. Since the path $\o'$ has to reach $\pieno$, starting from $\h'$ the protuberance is detached and in the sequel is attached in another site. Thus $\o'$ reaches a saddle that is not visited by $\o$. This implies that (\ref{isocondess}) is satisfied. 

In case II, when the path $\o'$ initiates a $1$-translation of a bar different from $\o$, it reaches a saddle that it is not visited by $\o$, thus the condition (\ref{isocondess}) is satisfied. Therefore we deduce that the path does not start this $1$-translation and thus the path $\o'$ must go back to $\h'$ (since it has to reach $\pieno$). From now on, as before, the path $\o'$ has to detach the protuberance and in the sequel it is attached in another site, thus $\o'$ reaches a saddle that is not visited by $\o$. This implies that (\ref{isocondess}) is satisfied. Thus we have proved that the saddle $\x_i$ is essential for any $i=1,...,m$.

\medskip
\noindent
The proof of (ii) and (iii) is analogue to the one done in (i) by modifying the length of the horizontal and vertical sides of $\hbox{CR}(\h)$.
\epr

\section{Proof of the main results: weakly anisotropic case}
\label{proofweak}

\subsection{Proof of the main Theorem \ref{thgateweak}}
\label{gateweak1}

In this Section we give the proof of the main Theorem \ref{thgateweak}. Recall the definitions of standard rectangles given in (\ref{standardgen}). Now we recall the definition of the set $\cP$ given in \cite{NOS} as
\be{}
\ba{ll}
\cP:= \{\h&:\, n(\h)= 1,\, v(\h)=\ell_2(s_{wa}^*)-1,\,
\h_{cl} \hbox{ is connected},\, \hbox{monotone},\\
&\hbox{ with circumscribed rectangle in }
\cR(\ell_1(s_{wa}^*)+1,\ell_2(s_{wa}^*))\}.
\ea
\ee

\noindent
In particular, in order to state that the set $\cwgeo$ is a gate for the transition from $\vuoto$ to $\pieno$, we need the following 

\bl{gateweak}
If $\o\in(\vuoto\ra\pieno)_{opt}$ passes through the set $\cP$, then $\o\cap\cwgeo\neq\emptyset$.
\el

\noindent
Since the proof of Lemma \ref{gateweak} follows a similar argument used in the proof of \cite[Lemma 6.1]{BN3}, we postpone it to Appendix \ref{app2}.

\medskip
\begin{proof*}{\bf of the main Theorem \ref{thgateweak}}
	By \cite[Theorem 2]{NOS}, we know that the set $\cP$ is a gate for the transition from $\vuoto$ to $\pieno$. By Lemma \ref{gateweak} we know that every path $\o\in(\vuoto\ra\pieno)_{opt}$ that crosses $\cP$ has to intersect also $\cwgeo$. This implies that every optimal path $\o$ from $\vuoto$ to $\pieno$ is such that $\o\cap\cwgeo\neq\emptyset$ and thus $\cwgeo\cap\cP\equiv\cwgeo$ is a gate.
\end{proof*}

\subsection{Proof of the main Theorem \ref{gweak}}
\label{weak1}

In this Section we analyze the geometry of the set $\cG_{wa}(\vuoto,\pieno)$ (recall (\ref{defg})). In particular, we give the proof of the main Theorem \ref{gweak} by giving in Proposition \ref{weakselless} the geometric characterization of the essential saddles of the third type that are not in $\cwgeo$ and that are visited after crossing the set $\cC_{wa}^B$.

\bp{weakselless}
Any saddle $\x$ that is neither in $\cC_{wa}^*$, nor in the boundary of the cycle $\cC_{\pieno}^{\vuoto}(\G_{wa}^*)$, nor in $\partial\cC_{\vuoto}^{\pieno}(\G_{wa}^*-H(\pieno))\setminus\widetilde{K}_{wa}$, such that $\t_{\x}\geq\t_{\cC_{wa}^B}$ can be essential or not. For those essential, we obtain the following characterization:
\be{}
\ba{ll}
\cG_{wa}(\vuoto,\pieno)\cap\cS_{wa}(\vuoto,\pieno)\setminus(\partial\cC_{\pieno}^{\vuoto}(\G_{wa}^*)\cup(\partial\cC_{\vuoto}^{\pieno}(\G_{wa}^*-H(\pieno))\setminus\widetilde{K}_{wa})\cup\cC_{wa}^*)=\qquad\qquad\\ \qquad\qquad\qquad\qquad\qquad\qquad\qquad\qquad\qquad\qquad=\displaystyle\bigcup_{\a}\bigcup_{\a'}\bigcup_{k'}\bigcup_{k=2}^{k'}\cN_{k,k'}^{\a,\a'}\cup\displaystyle\bigcup_{\a'}\cN_0^{\a'}\cup\displaystyle\bigcup_{\a}\cN_1^{\a}
\ea
\ee
\ep

\noindent
Since the proof of Proposition \ref{weakselless} is similar to the one of \cite[Proposition 6.3]{BN3}, we postpone it to Appendix \ref{app2}.

\medskip
\begin{proof*}{\bf of Theorem \ref{gweak}}
	By Corollary \ref{corstrategy} we know that the saddles of the first and second type, defined in Definitions \ref{sellesigma} and \ref{sellezeta} respectively, are unessential. By Propositions \ref{selle3} and \ref{weakselless} we have the characterization of the essential saddles of the third type in Section \ref{S6.4}. We use Proposition \ref{c*contenuto} to get the claim.
\end{proof*}

\section{Proof of the sharp asymptotics}
\label{sharpasymptotics}

In this Section, following the approach initiated in \cite{BHN} for Kawasaki dynamics and developed in \cite[Chapters 16 and 18]{BH}, we provide a comparison with some model-independent results given in \cite[Chapter 16]{BH}.

\subsection{Model-independent results for the prefactor}
\label{modindipdisc}

	 In order to give the results, we need some definitions. In \cite{BH} the authors let the protocritical and critical sets as $\mathscr{P}^*(m,s)$ and $\mathscr{C}^*(m,s)$ respectively. Since they differ from our notation, we refer to them as $\mathscr{P}_{PTA}^*(m,s)$ and $\mathscr{C}_{PTA}^*(m,s)$. Given $\x,\x'\in\cX$, we set $\x\sim\x'$ if the two configurations can be obtained from each other via an allowed move.

\bd{setlibro}{\cite[Definition 16.3]{BH}}
Let 
\be{}
\G^*=\Phi(m,s)-H(m).
\ee

\noindent
Then $(\mathscr{P}_{PTA}^*(m,s),\mathscr{C}_{PTA}^*(m,s))$ is the maximal subset of $\cX\times\cX$ such that:
\bi
\item[(1)] $\forall \ \x\in\mathscr{P}_{PTA}^*(m,s)\ \exists \ \x'\in\mathscr{C}_{PTA}^*(m,s): \ \x\sim\x'$ and $\forall \ \x'\in\mathscr{C}_{PTA}^*(m,s)\ \exists \ \x\in\mathscr{P}_{PTA}^*(m,s): \ \x'\sim\x$;
\item[(2)] $\forall \ \x\in\mathscr{P}_{PTA}^*(m,s), \ \Phi(\x,m)<\Phi(\x,s)$;
\item[(3)] $\forall \ \x'\in\mathscr{C}_{PTA}^*(m,s)\ \exists \ \g:\x'\ra s$ such that $\displaystyle\max_{\z\in\g}H(\z)-H(m)\leq\G^*,  \g\cap\{\z\in\cX:  \Phi(\z,m)<\Phi(\z,s)\}=\emptyset$.
\ei
\ed

\noindent
Now we abbreviate $\mathscr{P}_{PTA}^*=\mathscr{P}_{PTA}^*(m,s)$ and  $\mathscr{C}_{PTA}^*=\mathscr{C}_{PTA}^*(m,s)$. In \cite{BH} the following results (Theorems \ref{th16.4'} and \ref{th16.5'}) are proved subject to the two hypotheses
\bi
\item[(H1)] $\cX^{m}=\{m\}$ and $\cX^{s}=\{s\}$;
\item[(H2)] $\x'\ra|\{\x\in\mathscr{P}_{PTA}^*: \ \x\sim\x'\}|$ is constant on $\mathscr{C}_{PTA}^*$.
\ei

\bt{th16.4'}{\cite[Theorem 16.4]{BH}}
\bi
\item[(a)] $\displaystyle\lim_{\b\ra\infty}\P_m(\t_{\mathscr{C}_{PTA}^*}<\t_s|\t_s<\t_m)=1$;
\item[(b)] $\displaystyle\lim_{\b\ra\infty}\P_m(\x_{\t_{\mathscr{C}_{PTA}^*}}=\chi)=\frac{1}{|\mathscr{C}_{PTA}^*|}$ for all $\chi\in\mathscr{C}_{PTA}^*$.
\ei
\et

\bt{th16.5'}{\cite[Theorem 16.5]{BH}}
There exists a constant $K\in(0,\infty)$ such that
\be{}
\displaystyle\lim_{\b\ra\infty}e^{-\b\G^*}\E_{m}(\t_s)=K.
\ee
\et

\noindent
Concerning Theorem \ref{th16.5'}, we add some comments below about the general strategy developed in \cite[Subsection 16.3.2]{BH}. A key role is played by the \emph{Dirichlet form}
\be{Diri}
\cE_\b(h) = \frac12 \sum_{\h,\h'\in\cX} \mu_\b(\h)c_\b(\h,\h')[h(\h)-h(\h')]^2,
\qquad h\colon\,\cX \to [0,1],
\ee
where $\mu_\b$ is the Gibbs measure defined in \cite[eq.\ (16.1.1)]{BH} and $c_\b$ is the kernel of transition rates defined in \cite[eq.\ (16.1.2)]{BH}. Given two non-empty disjoint sets $\cA,\cB\subseteq\cX$, the \emph{capacity} of the pair $\cA,\cB$ is defined by
\be{Dirimin}
\CAPA_\b(\cA,\cB) = \min_{{h\colon\,\cX \to [0,1]} \atop 
	{h|_\cA\equiv 1,\,h|_\cB\equiv 0}} \cE_\b(h),
\ee
where $h|_\cA\equiv 1$ means that $h(\h)=1$ for all $\h\in \cA$ and $h|_\cB\equiv 0$
means that $h(\h)=0$ for all $\h\in \cB$. The right-hand side of (\ref{Dirimin})
has a unique minimizer $h_{\cA,\cB}^*$, called the \emph{equilibrium potential}
of the pair $\cA,\cB$, given by
\be{mini}
h_{\cA,\cB}^*(\h) = \P_\h(\t_\cA<\t_\cB), \qquad \h\in\cX\setminus (\cA\cup\cB)
\ee
This is the solution of the equation
\be{soleq}
\begin{aligned}
	(c_\b h)(\eta) &= 0, \qquad  \eta\in \cX\setminus (\cA\cup\cB),\\ 
	h(\eta)&=1, \qquad  \eta\in \cA,\\
	h(\eta)&=0,\qquad  \eta\in \cB.\\ 
\end{aligned}
\ee
Moreover,
\be{caprep}
\CAPA_\b(\cA,\cB) = \sum_{\h\in\cA} \mu_\b(\h)\,c_\b(\h,\cX\setminus\h)
\,\P_\h(\t_\cB<\t_\cA)
\ee
with $c_\b(\h,\cX\setminus\h)=\sum_{\h'\in\cX\setminus\h} c_\b(\h,\h')$ the rate of 
moving out of $\h$. This rate enters because $\t_\cA$ is the first hitting time of 
$\cA$ after the initial configuration is left. Note that from 
(\ref{Diri}--\ref{Dirimin}) it follows that
\be{capsym}
\CAPA_\b(\cA,\cB)=\CAPA_\b(\cB,\cA).
\ee

\noindent
First, we introduce a graph representation of the configuration space. View $\cX$ as a graph whose vertices are configurations and whose edges connect communicating configurations. Let
\begin{itemize}
	\item[-] 
	$\cX^*$ be the subgraph of $\cX$ obtained by removing all vertices $\eta$ with $H(\eta)>\G^*+H(m)$ and all edges incident to these vertices;
	\item[-] 
	$\cX^{**}$ be the subgraph of $\cX^*$ obtained by removing all vertices $\eta$ with $H(\eta)=\G^*+H(m)$ and all edges incident to these vertices; 
	\item[-]
	$\cX^{meta}$ and $\cX^{stab}$ be the connected components of $\cX^{**}$ containing $m$ and $s$ respectively.
\end{itemize}

\noindent
Moreover, we consider the set
\be{X**part}
\cX^{**}\setminus(\cX^{meta}\cup\cX^{stab})=\bigcup_{i=1}^I \cX(i),
\ee

\noindent
where each $\cX(i)$ is a \emph{well} in $\cS(m,s)$, i.e., a set of communicating configurations with energy $<\G^*+H(m)$ but with communication height $\G^*+H(m)$ towards both $m$ and $s$. Among all the wells $\cX(i)$, we can highlight the wells $\cZ_j^{m}$ (resp.\ $\cZ_j^{s}$) of the unessential saddles of the first (resp. second) type $\s_j$ (resp.\ $\z_j$) (see Definitions \ref{sellesigma} and \ref{sellezeta}, and Propositions \ref{selle1} and \ref{selle2}). In particular, these wells can be defined as follows. 

\bd{vallenoness}
We define
\bi
\item[1)] $\cZ_j^m\subset\cX^{**}$, $j=1,...,J_m$, is a connected set such that, for all $\h\in\cZ_j^m$, $\Phi(m,\h)=\Phi(s,\h)$ and any path $\o:\h\ra s$ must be such that $\o\cap\cX^{meta}\neq\emptyset$;
\item[2)] $\cZ_j^s\subset\cX^{**}$, $j=1,...,J_s$, is a connected set such that, for all $\h\in\cZ_j^m$, $\Phi(s,\h)=\Phi(m,\h)$ and any path $\o:\h\ra m$ must be such that $\o\cap\cX^{stab}\neq\emptyset$.
\ei
\ed

\bp{vallesigma}
If $\cZ_j^m\neq\emptyset$, $\cX(i)\equiv\cZ_j^m$ if and only if $\cZ_j^m$ is a connected component in $\cX^{**}\setminus(\cX^{meta}\cup\cX^{stab})$ such that there exists a saddle of the first type $\s_j$ that communicates via one step with a configuration in $\cZ_j^m$. 
\ep

\bpr
First, assume that $\cX(i)\equiv\cZ_j^m$ for some $i\in\{1,...,I\}$ and $j\in\{1,...,J_m\}$. It is clear that $\cZ_j^m$ is a connected component in $\cX^{**}\setminus(\cX^{meta}\cup\cX^{stab})$, indeed any configuration $\h\in\cX^{meta}\cup\cX^{stab}$ has $\Phi(\h,m)\neq\Phi(\h,s)$. Furthermore, by definition of $\cZ_j^m$, we note that $\Phi(s,\h)=\G^*+H(m)$ and any path $\o:\h\ra s$ must be such that $\o\cap\cX^{meta}\neq\emptyset$. This implies that there exists $\s_j\in\o\cap\cS(s,m)$ that is an unessential saddle of the first type (see Proposition \ref{selle1}) such that it communicates via one step with a configuration in $\cZ_j^m$.

Conversely, assume that any fixed $\cZ_j^m$ is a connected component in $\cX^{**}\setminus(\cX^{meta}\cup\cX^{stab})$ such that there exists a saddle of the first type $\s_j$ that communicates via one step with a configuration in $\cZ_j^m$. Thus we deduce that $\cX(i)\equiv\cZ_j^m$ for some $i\in\{1,...,I\}$, indeed for all $\h\in\cZ_j^m$ it hold $\Phi(\h,m)=\Phi(\h,s)=\G^*+H(m)$ and $H(\h)<\G^*+H(m)$ since $\cZ_j^m\subset\cX^{**}$.
\epr

\bp{vallezeta}
If $\cZ_j^s\neq\emptyset$, $\cX(i)\equiv\cZ_j^s$ if and only if $\cZ_j^s$ is a connected component in $\cX^{**}\setminus(\cX^{meta}\cup\cX^{stab})$ such that there exists a saddle of the second type $\z_j$ that communicates via one step with a configuration in $\cZ_j^s$.
\ep

\bpr
The proof is analogue to the proof of Proposition \ref{vallesigma} by replacing $\cZ_j^m$ with $\cZ_j^s$ and ``$\o:\h\ra s$ such that $\o\cap\cX^{meta}\neq\emptyset$" with ``$\o:\h\ra m$ such that $\o\cap\cX^{stab}\neq\emptyset$".
\epr

 First, note that the unessential saddles of the first and second type are not in the set $\mathscr{C}_{PTA}^*$. Indeed, on the one hand, the saddles $\{\s_j\}_{j=1}^{J_m}$ do not verify the condition (3) in Definition \ref{setlibro}, because every optimal path that connects any fixed $\s_j$ to $s$ passes through $\cX^{meta}$. On the other hand the saddles $\{\z_j\}_{j=1}^{J_s}$ do not verify conditions (1) and (2) in Definition \ref{setlibro}, because they communicate only with configurations that are not in $\mathscr{P}_{PTA}^*$. By \cite[eq.\ (16.3.4)]{BH} and \cite[Lemma 16.16]{BH}, we know that $h$ is constant on each wells, but for the wells $\cZ_j^{m}$ and $\cZ_j^{s}$ we compute this constant in Lemma \ref{hselleiness}, indeed \cite[Lemma 16.15]{BH} can be extended also for these sets together with the unessential saddles of the first and second type.

\bl{hselleiness}
Recall Definitions \ref{sellesigma} and \ref{sellezeta} for the definition of the saddles $\s_j$ of the first type and $\z_j$ second type respectively and Definition \ref{vallenoness} for the definition of the wells $\cZ_j^m$ and $\cZ_j^s$. As $\b\ra\infty$,
\be{}
Z_\b\CAPA_\b(m,s)=[1+o(1)]\Theta e^{-\b\G^*},
\ee

\noindent
with
\be{varform}
\Theta=\displaystyle\min_{c_1,...,c_{\bar{I}}}\min_{h:\cX^*\ra[0,1] \atop h_{|\cX^{meta}_I}=1,\ h_{|\cX^{stab}_{II}}=0, \ h_{|\cX(i)}=c_i,i=1,...,\bar{I}}\frac{1}{2}\sum_{\x,\x'\in\cX^*}\mathbbm{1}_{\{\x\sim\x'\}} [h(\x)-h(\x')]^2,
\ee
\noindent
where
\be{}
\cX^{meta}_I:=\displaystyle\cX^{meta}\cup\bigcup_{j=1}^{J_m}(\{\s_j\}\cup\cZ_j^{m}), \qquad \cX^{stab}_{II}:=\displaystyle\cX^{stab}\cup\bigcup_{j=1}^{J_s}(\{\z_j\}\cup\cZ_j^{s}).
\ee
\noindent
and
$\cX(i)$, $i=1,...,\bar{I}$, are all the wells of the transition except $\bigcup_{j=1}^{J_m}\cZ_j^{m}$ and $\bigcup_{j=1}^{J_s}\cZ_j^{s}$. 
\el

\bpr
The statement is similar to the one of \cite[Lemma 16.17]{BH}, but the difference is in the variational formula for $\Theta$. More precisely, comparing (\ref{varform}) with \cite[eq.\ (16.3.11)]{BH}, the proof is analogue to the one done for \cite[Lemma 16.17]{BH}, but we have to prove that the function $h$ is constant equal to 1 (resp.\ 0) on $\bigcup_{j=1}^{J_m}(\{\s_j\}\cup\cZ_j^{m})$ (resp.\ $\bigcup_{j=1}^{J_s}(\{\z_j\}\cup\cZ_j^{s})$).

Fix any saddle of the first type $\s_j$. 
By \cite[Lemma 16.16]{BH} we set $h(\h)=c_j$ for any $\h\in\cZ_j^m$, $h(\h)=c_k$ for any $\h\in\cZ_k^m$ with $k\neq j$, and $h(\s_j)=\bar{c}_j$. By definition of saddles of the first type, note that $\s_j$ communicates only with configurations either in $\cX^{meta}$, or in $\cZ_j^m$ or in $\cZ_k^m$ with $k\neq j$. Thus the contribution to (\ref{varform}) of the saddle of the first type $\s_j$ is
\be{contononess}
\displaystyle\sum_{\x\in\cX^{meta}}|\s_j\sim\x|(1-\bar{c}_j)^2+\sum_{\x\in\cZ_j^m}|\s_j\sim\x|(\bar{c}_j-c_j)^2+\sum_{\x\in\cZ_k^m \atop k\neq j}|\s_j\sim\x|(c_k-\bar{c}_j)^2
\ee

\noindent
The cases $\cZ_j^m=\emptyset$ (resp.\ $\cZ_k^m=\emptyset$) or there is no $\x\in\cZ_j^m$ (resp.\ $\x\in\cZ_k^m$) such that $\s_j\sim\x$ correspond to the situation in which either there is not the well $\cZ_j^m$ (resp.\ $\cZ_k^m$ for any $k\neq j$) or the dynamics does not allow the communication via one step from $\s_j$ and $\cZ_j^m$ (resp.\ $\cZ_k^m$). In the first situation we get $\bar{c}_j=1$ and in the second one we get $\bar{c}_j=c_j=1$. Otherwise, since the quantity in (\ref{contononess}) is greater or equal than zero, the minimum w.r.t.\ $c_j$, $c_k$ and $\bar{c}_j$ of (\ref{contononess}) is obtained for $c_j=c_k=\bar{c}_j=1$, thus $h$ is constant equal to 1 on $\bigcup_{j=1}^{J_m}(\{\s_j\}\cup\cZ_j^{m})$. If we consider all the possible transitions of this type and use (\ref{contononess}), we obtain a contribution to (\ref{varform}) equal to
\be{noness}
\displaystyle\sum_{j=1}^{J_m}\Bigg(\sum_{\x\in\cX^{meta}}|\s_j\sim\x|(1-\bar{c}_j)^2+\sum_{\x\in\cZ_j^m}|\s_j\sim\x|(\bar{c}_j-c_j)^2+\sum_{\x\in\cZ_k^m \atop k\neq j}|\s_j\sim\x|(c_k-\bar{c}_j)^2\Bigg)
\ee

\noindent
Again, the minimum of (\ref{noness}) with respect to $c_1,...,c_{J_m}$ and $\bar{c}_1,...,\bar{c}_{J_m}$ is obtained for $c_i=\bar{c}_i=1$ for any $i=1,...,J_m$.

Similarly, we deduce that $h$ is constant equal to 0 on $\bigcup_{j=1}^{J_s}(\{\z_j\}\cup\cZ_j^{s})$. Indeed, if we fix any saddle of the second type $\z_j$, by \cite[Lemma 16.16]{BH} we set $h(\h)=c_j$ for any $\h\in\cZ_j^s$, $h(\h)=c_k$ for any $\h\in\cZ_k^s$ with $k\neq j$, and $h(\z_j)=\bar{c}_j$. By definition of saddles of the second type, note that $\z_j$ communicates only with configurations either in $\cX^{stab}$, or in $\cZ_j^s$ or in $\cZ_k^s$ with $k\neq j$. Thus the contribution to (\ref{varform}) of the saddle of the second type $\z_j$ is
\be{contononess2}
\displaystyle\sum_{\x\in\cX^{stab}}|\z_j\sim\x|\bar{c}_j^2+\sum_{\x\in\cZ_j^s}|\z_j\sim\x|(\bar{c}_j-c_j)^2+\sum_{\x\in\cZ_k^s \atop k\neq j}|\z_j\sim\x|(c_k-\bar{c}_j)^2.
\ee

\noindent
The cases $\cZ_j^s=\emptyset$ (resp.\ $\cZ_k^s=\emptyset$) or there is no $\x\in\cZ_j^s$ (resp.\ $\x\in\cZ_k^s$) such that $\z_j\sim\x$ correspond to the situation in which either there is not the well $\cZ_j^s$ (resp.\ $\cZ_k^s$ for any $k\neq j$) or the dynamics does not allow the communication via one step from $\z_j$ and $\cZ_j^s$ (resp.\ $\cZ_k^s$). In the first situation we get $\bar{c}_j=0$ and in the second one we get $\bar{c}_j=c_j=0$. Otherwise, since the quantity in (\ref{contononess2}) is greater or equal than zero, the minimum w.r.t.\ $c_j$, $c_k$ and $\bar{c}_j$ of (\ref{contononess2}) is obtained for $c_j=c_k=\bar{c}_j=0$, thus $h$ is constant equal to 0 on $\bigcup_{j=1}^{J_s}(\{\z_j\}\cup\cZ_j^{s})$. If we consider all the possible transitions of this type, we obtain a contribute to (\ref{varform}) equal to
\be{noness2}
\displaystyle\sum_{j=1}^{J_s}\Bigg(\sum_{\x\in\cX^{stab}}|\z_j\sim\x|\bar{c}_j^2+\sum_{\x\in\cZ_j^m}|\z_j\sim\x|(\bar{c}_j-c_j)^2+\sum_{\x\in\cZ_k^m \atop k\neq j}|\z_j\sim\x|(c_k-\bar{c}_j)^2\Bigg)
\ee

\noindent
Again, the minimum of (\ref{noness2}) with respect to $c_1,...,c_{J_s}$ and $\bar{c}_1,...,\bar{c}_{J_s}$ is obtained for $c_i=\bar{c}_i=0$ for any $i=1,...,J_s$. Therefore formula (16.3.15) in the proof of \cite[Lemma 16.17]{BH} should be modified as
\be{defh}
h=
\begin{cases}
	1 &\hbox{on }  \displaystyle\cX^{meta}\cup\bigcup_{j=1}^{J_m}(\{\s_j\}\cup\cZ_j^{m}), \\
	0 &\hbox{on } \displaystyle\cX^{stab}\cup\bigcup_{j=1}^{J_s}(\{\z_j\}\cup\cZ_j^{s}),\\
	c_i &\hbox{on } \cX(i), i=1,...,\bar{I},
\end{cases}
\ee

\noindent
where $\cX(i)$, $i=1,...,\bar{I}$, are all the wells of the transition except $\bigcup_{j=1}^{J_m}\cZ_j^{m}$ and $\bigcup_{j=1}^{J_s}\cZ_j^{s}$. We get the claim.
\epr

\br{}
Lemma \ref{hselleiness} implies that also the unessential saddles $\s_j$ and $\z_j$ have to be considered in the estimate of the prefactor. However, since $h(\s_j)=1$ and $h(\z_j)=0$ for any $j$, the transitions that involve these unessential saddles do not contribute numerically to the computation of $K$. The variational formula for $\Theta$ in (\ref{varform}) is non-trivial because it depends on the geometry of all the wells $\cX(i)$, $i=1,...,I$, and on the form of the function $h$ on the configurations in $\cX^{*}\setminus\cX^{**}$, namely the saddle configurations. These two steps are the model-dependent keys to compute the prefactor $K=1/\Theta$.
\er

\br{confronto}
For the hexagonal Ising model that evolves under Glauber dynamics, the estimate of the upper bound of the prefactor, given in \cite[Section 6.1]{AJNT}, is done setting the equilibrium potential $h$, for that specific model, according to our discussion (see (\ref{defh}) and \cite[eq.\ (6.19)]{AJNT})). The authors analyzed the unessential saddles of first and second type that they used for the definition of $h$ together with their valleys. Indeed they define the sets $\bigcup_{j=1}^{J_A}\cN_j^A$ and $\bigcup_{j=1}^{J_B}\cN_j^B$ in \cite[Section 6.1]{AJNT} that coincide, in their model, with $\bigcup_{j=1}^{J_m}\cZ_j^{m}$ (resp.\ $\bigcup_{j=1}^{J_s}\cZ_j^{s}$) by Proposition \ref{vallesigma} (resp.\ Proposition \ref{vallezeta}). In \cite[Section 6.1]{AJNT} an explicit example of saddle $\s_i$ is given.
\er

\br{}
In \cite[Section 7]{BGNneg2021} the authors study the estimate of the prefactor for the $q$-state Potts model that evolves under Glauber dynamics using the above discussion and Lemma \ref{hselleiness}. In \cite[Lemma 7.4(b)]{BGNneg2021} the authors identify geometrically unessential saddles of the first type (see \cite[Figure 18(b)]{BGNneg2021}) and describe their wells in the proof. In \cite[Lemma 7.4(c)]{BGNneg2021} the authors identify geometrically one unessential saddle of the second type (see \cite[Figure 19]{BGNneg2021}). Choosing $q=2$, this Lemma gives the same results for the standard Ising model. We refer also to \cite[Chapter 17]{BH}, where the authors compute the prefactor in \cite[Theorem 17.4]{BH} using \cite[Lemma 16.17]{BH} and some model-dependent properties without identifying the unessential saddles.
\er

\br{}
For the Ising model with strongly anisotropic interactions that evolves under Kawasaki dynamics, the estimate of the prefactor is given in \cite[Theorem 4.12]{BN3} according to the above discussion and Lemma \ref{hselleiness}. For the strongly anisotropic case the authors are able to obtain a sharp estimate for $K_{sa}$ in \cite[eq.\ (4.27)]{BN3}. Nevertheless, the asymptotic behavior of the prefactor $K_{sa}$ as $\L\ra\Z^2$ is the same as $K_{wa}$ (see Theorem \ref{sharptimeweak}) and $K_{is}$ (see \cite[Theorem 1.4.5]{BHN}). Moreover, in \cite[Figure 14]{BN3} an example of unessential saddle of the second type is given. 
\er

\subsection{Isotropic and weakly anisotropic cases}
\label{kawisoandweak}

\noindent
In this Subsection we consider $int\in\{is,wa\}$. For our model $\cX^{m}_{int}=\{\vuoto\}$ and $\cX^{s}_{int}=\{\pieno\}$, thus (H1) holds and $\G^*=\Phi(\vuoto,\pieno)-H(\vuoto)=\G^*_{int}$. Moreover, $\mathscr{P}_{PTA}^*(\vuoto,\pieno)=\cD_{int}$ and $\mathscr{C}_{PTA}^*(\vuoto,\pieno)$ is the union of all the configurations that are composed by a cluster in $\cD_{int}$ and a free particle in $\partial^-\Lambda$. Therefore it is clear that $\cC_{int}^*\neq\mathscr{C}_{PTA}^*(\vuoto,\pieno)$. Note that (H2) follows from Lemma \ref{entrataweak}. Here we abbreviate $\mathscr{P}_{PTA}^*=\mathscr{P}_{PTA}^*(\vuoto,\pieno)$ and  $\mathscr{C}_{PTA}^*=\mathscr{C}_{PTA}^*(\vuoto,\pieno)$. Note that $\cX^{meta}=\cC_{\pieno}^{\vuoto}(\G_{int}^*)$ and $\cX^{stab}=\cC^{\pieno}_{\vuoto}(\G_{int}^*-H(\pieno))$. Recall Definition \ref{vallenoness} for the definition of the wells $\cZ_{int,j}^{\vuoto}$ and $\cZ_{int,j}^{\pieno}$. Concerning Theorem \ref{th16.5'} for our cases, by \cite[Lemma 3.3.2]{BHN} we know that $h$ is constant on each wells for the isotropic case and this holds also for the weakly anisotropic case. Thanks to the model-independent discussion given in Section \ref{modindipdisc} and Lemma \ref{hselleiness}, for $int\in\{is,wa\}$ formula (\ref{defh}) becomes
\be{defhkaw}
h=
\begin{cases}
	1 &\hbox{on }  \displaystyle\cC_{\pieno}^{\vuoto}(\G^*_{int})\cup\bigcup_{j=1}^{J_{\vuoto}}(\{\s_{int,j}\}\cup\cZ_{int,j}^{\vuoto}), \\
	0 &\hbox{on } \displaystyle\cC^{\pieno}_{\vuoto}(\G^*_{int}-H(\pieno))\cup\bigcup_{j=1}^{J_{\pieno}}(\{\z_{int,j}\}\cup\cZ_{int,j}^{\pieno}),\\
	c_i &\hbox{on } \cX_{int}(i), i=1,...,\bar{I},
\end{cases}
\ee

\noindent
where $\cX_{int}(i)$, $i=1,...,\bar{I}$, are all the wells of the transition except $\bigcup_{j=1}^{J_{\vuoto}}\cZ_{int,j}^{\vuoto}$ and $\bigcup_{j=1}^{J_{\pieno}}\cZ_{int,j}^{\pieno}$. Note that, for $int=is$, formula (3.3.20) in the proof of \cite[Proposition 3.3.3]{BHN} should be modified accordingly to (\ref{defhkaw}). 

\subsection{Proof of Theorem \ref{sharptimeweak}}

In this Section we give the proof of Theorem \ref{sharptimeweak}, that represents sharp asymptotics for the weakly anisotropic case. To this end, we need a theorem that summarizes what we have shown so far. It represents the analogue version of \cite[Theorem 2.3.10]{BHN} valid for the isotropic case.

\bt{graph}
Let $int\in\{wa,sa\}$. The following conditions hold:

\bi
\item[(i)] $\cC_{\pieno}^{\vuoto}(\G_{int}^*)\neq\cC_{\vuoto}^\pieno(\G_{int}^*-H(\pieno))$, so $\cC_{\pieno}^\vuoto(\G_{int}^*)$ and $\cC_{\vuoto}^\pieno(\G_{int}^*-H(\pieno))$ are disconnected in $\cX^{**}_{int}$;
\item[(ii)] $\bar\cD_{int}\subseteq\cC_{\pieno}^\vuoto(\G_{int}^*)$, $\cC_{int}^G\subseteq\cC_{\vuoto}^\pieno(\G_{int}^*-H(\pieno))$ and $\cC_{int}^B\subseteq\cX^{**}_{int}\setminus(\cC_{\pieno}^\vuoto(\G_{int}^*)\cup\cC_{\vuoto}^\pieno(\G_{int}^*-H(\pieno)))$.
\ei
\et

\noindent
In \cite{BEGK} metastability is defined in terms of properties of capacities, namely:

\bd{metastableset}\cite[Definition 3.2.1.]{BHN},\cite[Definition 8.2]{BH}
Consider a family of Markov chains, indexed by $\b$, on a finite state space 
$\cX$. A set $\cM \subseteq \cX$ is called PTA-metastable if 
\be{metaset}
\lim_{\b\to\infty}\,\, \frac
{\max_{\h \notin \cM} \mu_\b(\h)[\CAPA_\b(\h,\cM)]^{-1}}
{\min_{\h \in \cM} \mu_\b(\h)[\CAPA_\b(\h,\cM\setminus \h)]^{-1}}
= 0.
\ee
\ed

\noindent 
By \cite[Lemma 3.2.2]{BHN} we know that the set $\{\vuoto,\pieno\}$ is PTA-metastable in the sense of Definition \ref{metastableset}. To obtain our sharp estimate of $\E_{\vuoto}(\pieno)$, we need the following

\bp{EDir}\cite[Proposition 3.2.3]{BHN}
\begin{equation*}
	\E_{\,\square}(\t_\blacksquare) = \frac{1}{Z_\b\,\CAPA_\b(\square,\blacksquare)}
	\,[1+o(1)]
\end{equation*}

\noindent
as $\b\to\infty$.
\ep

\noindent
We follow the general strategy outlined in \cite{{BM},{B04},{BH},{BHN}}:

\begin{itemize}
	\item[--] 
	Note that all terms in the Dirichlet form in (\ref{Diri}) involving configurations 
	$\eta$ with $H(\eta)>\Gamma^*_{wa}$, i.e., $\eta\in\cX_{wa}\setminus\cX^*_{wa}$, contribute 
	at most $Ce^{-(\Gamma^*_{int}+\d)\b}$ for some $\d>0$ and can be neglected. 
	Thus, effectively we can replace $\cX_{wa}$ by $\cX^*_{wa}$.
	\item[--] 
	Show that $h^*_{\vuoto,\pieno}=O(e^{-\d\b})$ on $\cC_{\vuoto}^\pieno(\G_{wa}^*-H(\pieno))$ and 
	$h^*_{\vuoto,\pieno}=1-O(e^{-\d\b})$ on $\cC_{\pieno}^\vuoto(\G_{wa}^*)$ for some $\d>0$.
	\item[--] 
	Prove sharp upper and lower bounds for $h^*_{\square,\blacksquare}$ 
	on $\cX^*_{wa}\setminus(\cC_{\pieno}^\square(\G_{wa}^*)\cup\cC_{\vuoto}^\blacksquare(\G_{wa}^*-H(\pieno)))$ in terms of a 
	variational problem involving only the vertices and the bonds on and 
	incident to $\cX^*_{wa}\setminus(\cC_{\pieno}^\square(\G_{wa}^*)\cup\cC_{\vuoto}^\blacksquare(\G_{wa}^*-H(\pieno)))$.
\end{itemize}

\noindent
Note that
\be{Adefs}
\begin{aligned}
	\cC_{\pieno}^\square(\G_{wa}^*)&= \{\h\in\cX^*_{wa}\colon\,\Phi(\h,\square)<\Phi(\h,\blacksquare)\},\\
	\cC_{\vuoto}^\blacksquare(\G_{wa}^*-H(\pieno))&= \{\h\in\cX^*_{wa}\colon\,\Phi(\h,\blacksquare)<\Phi(\h,\square)\}.
\end{aligned}
\ee
The guiding idea behind the sharp estimate of $Z_\b\,\CAPA_\b(\square,\blacksquare)$ 
is that $h_{\square,\blacksquare}^*$ is exponentially close to 1 on $\cX_{\pieno}^\square(\G_{wa}^*)$ 
and exponentially close to 0 on $\cX_{\vuoto}^\blacksquare(\G_{wa}^*-H(\pieno))$. By \cite[Lemma 3.3.1]{BHN} we know that $h^*_{\square,\blacksquare}$ is trivial on $\cC_{\pieno}^\square(\G_{wa}^*)\cup\cC_{\vuoto}^{\blacksquare}(\G_{wa}^*-H(\pieno))$, thus it remains to understand what $h^*_{\square,\blacksquare}$
looks like on the set
\be{A=def}
\cX^*_{wa}\setminus(\cC_{\pieno}^\square(\G_{wa}^*)\cup\cC_{\vuoto}^\blacksquare(\G_{wa}^*-H(\pieno)))= \{\h\in\cX^*_{wa}\colon\,\Phi(\h,\square)=\Phi(\h,\blacksquare)\},
\ee
which separates $\cC_{\pieno}^\square(\G_{wa}^*)$ and $\cC_{\vuoto}^\blacksquare(\G_{wa}^*-H(\pieno))$ and contains $\cS_{wa}(\square,\blacksquare)$. Before doing so, we first show that $h^*_{\square,\blacksquare}$ is also trivial on $\cX^{**}_{wa}\setminus(\cC_{\pieno}^\square(\G_{wa}^*)\cup\cC_{\vuoto}^\blacksquare(\G_{wa}^*-H(\pieno)))$. This set
can be partitioned into maximally connected components,
\be{X**part}
\cX^{**}_{wa}\setminus(\cC_{\pieno}^\square(\G_{wa}^*)\cup\cC_{\vuoto}^\blacksquare(\G_{wa}^*-H(\pieno))) = \bigcup_{i=1}^I \cX_{wa}(i),
\ee
where each $\cX_{wa}(i)$ is a \emph{well} in $\cS_{wa}(\square,\blacksquare)$, i.e., a set of communicating configurations with energy $<\G^*_{wa}$ but with 
communication height $\G^*_{wa}$ towards both $\square$ and $\blacksquare$. By \cite[Lemma 3.3.2]{BHN} we know that $h^*_{\square,\blacksquare}$ is close 
to a constant on each of these wells. By Proposition \ref{pathstrong}(ii) we know that for each $\hat\h\in\bar\cD_{wa}$ the four bars of bad sites in $\partial^+\hbox{CR}(\hat\h)$ form a well. These are not the only wells, but \cite[Lemma 3.3.2]{BHN} shows that we not need care too much about wells anyway. For the weakly anisotropic case, \emph{only the transitions in and out of the wells contribute to the Dirichlet form} at the order we are after, not those inside the wells.

Using Lemma \ref{hselleiness}, \cite[Proposition 3.3.4]{BHN} and \cite[Lemma 3.4.1]{BHN}, in order to prove Theorem \ref{sharptimeweak} it remains to count the cardinality of $\bar\cD_{wa}$ modulo shifts. We denote by $N_{wa}$ that quantity.

\bp{Nweak}
\begin{equation*}
	N_{wa}=\displaystyle\sum_{k=1}^4\binom{4}{k}\binom{l_2^*+k-2}{2k-1}.
\end{equation*}
\ep

\bpr
We have to count the number of different shapes of the clusters in $\bar\cD_{wa}$. We do this by counting in how many ways $l_2^*-1$ particles can be removed from the four bars of a $l_1^*\times l_2^*$ rectangle starting from the corners. We split the counting according to the number $k=1,2,3,4$ of corners from which particles are removed. The number of ways in which we can choose $k$ corners is $\binom{4}{k}$. After we have removed the particles at these corners, we need to remove $l_2^*-1-k$ more particles frome either side of each corner. The number of ways in which this can be done is
\be{}	
\ba{lll}
|\{(m_1,...,m_{2k})\in\N_0^{2k}:\ m_1+...+m_{2k}=l_2^*-1-k\}| \\
\qquad=|\{(m_1,...,m_{2k})\in\N^{2k}:\ m_1+...+m_{2k}=l_2^*-1+k\}| \\
\qquad=\displaystyle\binom{l_2^*+k-2}{2k-1}.
\ea
\ee

\noindent
Thus we get the claim.
\epr

\subsection{Proof of Theorem \ref{entrunif}}

In this Section we give the proof of Theorem \ref{entrunif}, that represents the uniform entrance distribution for $int=wa$. Let $\partial^-\cC_{wa}^*$ be those configurations in $\cC_{wa}^*$ where the free particle is in $\partial^-\L$. Following the same argument used in \cite{BHN} for the isotropic regime, since $\bar\cD_{wa}\subseteq\cC_{\pieno}^{\vuoto}(\gw)$ by Theorem \ref{graph}(ii), it follows from \cite[Lemma 3.3.1]{BHN} and $\cC_{wa}^*\subseteq\cS_{wa}(\vuoto,\pieno)$ that
\be{bdvar1}
\min_{\h'\in\bar\cD_{wa}} h^*_{\vuoto,\partial^-\cC_{wa}^*}(\eta')
\geq 1-Ce^{-\d\b},
\ee

\noindent
where 
\be{h*rels}
h^*_{\vuoto,\partial^-\cC_{wa}^*}(\h')
= \left\{\begin{array}{ll}
	0 &\mbox{if } \h'\in\partial^-\cC_{wa}^*,\\
	\P_{\h'}(\t_\square<\t_{\partial^-\cC_{wa}^*}) &\mbox{otherwise}.
\end{array}
\right.
\ee 

\noindent
Moreover, letting $\partial^{--}\cC_{wa}^*$ be the set of configurations obtained from $\partial^-\cC_{wa}^*$ by moving the free particle from $\partial^-\L$ to $\partial^{--}\L=\partial^-(\L^-)$, we deduce that
\be{bdvar2}
\max_{\h'\in\partial^{--}\cC_{wa}^*} h^*_{\vuoto,\partial^-\cC_{wa}^*}(\eta')\leq Ce^{-\d\b}.
\ee

\noindent
From now on, following the argument proposed in \cite{BHN} we are able to prove the assertion in (\ref{uniforme}).

\subsection{Proof of Theorem \ref{gapspettrale}}

Thanks to \cite[Lemma 3.6]{BNZ}, we deduce that for our model the quantity $\widetilde\G(B)$, with $B\subsetneq\cX$, defined in \cite[eq.\ (21)]{BNZ} is such that $\widetilde\G(\cX\setminus\{\pieno\})=\G_{int}^*$, with $int\in\{is,wa\}$. Thus Theorem \ref{gapspettrale} follows by the following proposition.

\bp{}{\cite[Proposition 3.24]{BNZ}}
For any $\e\in(0,1)$ and any $s\in\cX^s$
\be{}
\displaystyle\lim_{\b\ra\infty}\frac{1}{\b}\log t_{\b}^{mix}(\e)=\widetilde\G(\cX\setminus\{s\})=\lim_{\b\ra\infty}-\frac{1}{\b}\log \r_{\b}
\ee

\noindent
Furthermore, there exist two constants $0<c_1\leq c_2<\infty$ independent of $\b$ such that for every $\b>0$
\be{}
c_1e^{-\b\widetilde\G(\cX\setminus\{s\})}\leq\r_{\b}\leq c_2e^{-\b\widetilde\G(\cX\setminus\{s\})}
\ee
\ep

\section{Proof of the main results for the simplified model}
\label{simpleproof}

In this Section we provide the proof of results in Section \ref{simpletheorem} that are extensions to the simplified model (described in Section \ref{simplemodel}) of the local model. In the case $int=is$, Theorem \ref{simpleiso} follows from \cite{HOS} and Theorem \ref{giso}. Now we consider $int\in\{wa,sa\}$ and we follow the strategy proposed in \cite{HOS} for $int=is$. In \cite[Section 2]{HOS} the authors give several large deviation estimates concerning exponential clocks, that hold also for the anisotropic cases. In \cite[Section 3]{HOS} the authors give several large deviation estimates concerning random walks. All these results are valid for the cases $int\in\{wa,sa\}$ without changes except for \cite[Proposition 3.13]{HOS}, in which we have to replace $U$ with $U_1$. The recurrence property for the anisotropic simplified model is obtained with similar arguments carried out in \cite[Section 6]{HOS}. To this end, we modify the definition of the set $\bar\cX_2$ given in \cite[eq.\ (5.8)]{HOS} by replacing $U$ with $U_1$. Therefore also the definition of the set $\cX_2$ given in \cite[eq.\ (6.1)]{HOS} should be modified accordingly. Thus, if we define for the anisotropic model $T_1=e^{0\b}$, $T_2=e^{U_1\b}$ and $T_3=e^{\D\b}$, \cite[Proposition 6.2]{HOS} holds also for the anisotropic cases. Concerning the reduction, we follow the strategy proposed in \cite[Section 7]{HOS}. In particular, we have to study the behavior of the gas and its interaction with the dynamics in the box $\L$. There are two classes of gas particles with different behavior: particles that have been in $\L^\b\setminus\L$ for a long time (say of order $T_3$), which we call green particles; and particles that exit from $\L$ and afterwards return to $\L$ in a short time (say of order 1), which we call red particles. The effect of green (resp.\ red) particles is studied in \cite[Section 7.6]{HOS} (resp.\ \cite[Section 7.7]{HOS}) and can be extended to the cases $int\in\{wa,sa\}$ by modifying the times $T_1$, $T_2$ and $T_3$, and the sets $\cX_2$ and $\bar\cX_2$ as above. In the case $int=wa$, from this discussion and \cite[Theorem 3]{NOS} (resp.\ Theorem \ref{thgateweak}) for the local model, Theorem \ref{simpleweak}(a) (resp.\ Theorem \ref{simpleweak}(b)) follows. Analogously, using \cite[Theorem 1]{NOS} (resp.\ Theorem \ref{gweak}) for the local model, Theorem \ref{simpleweak}(c) (resp.\ Theorem \ref{simpleweak}(d)) holds.

\appendix
\section{Appendix}
\label{appendice}

\subsection{Additional material for Section \ref{dependentdef}}
\label{app*}

\begin{proof*}{\bf of Lemma \ref{trasless}}
	Let $int=wa$. Note that $H(\h^B)=\G_{wa}^*-U_2$ (resp.\ $H(\h^B)=\G_{wa}^*-U_1$) if the free particle has been attached to an horizontal (resp.\ vertical) bar. In the first case, in order to avoid exceeding the energy value $\G_{wa}^*$ it is possible to translate only the vertical bars. These saddles are in $\cN_{1}^{\a}$. In the latter case, it is possible to translate both vertical and horizontal bars. If the translated bar is horizontal, the saddles that are crossed are in $\cN_{0}^{\a'}$. If the translated bar is vertical, the configurations obtained do not reach the level $\G_{wa}^*$, thus they are not saddles. To conclude, all the configurations in $\cN_{0}^{\a'}\cup\cN_{1}^{\a}$ can be obtained from a configuration $\h^B$ via a $1$-translation of a bar.
	
	It remains to prove that the saddles in $\cN_{0}^{\a'}\cup\cN_{1}^{\a}$ are essential. This part of the proof is analogue to the corresponding one done for $int=is$ in Lemma \ref{isotrasless}.
\end{proof*}

\begin{proof*}{\bf of Lemma \ref{coltorow}}
	Let $int=wa$. Since $H(\h)=\G_{wa}^*$, it is possible to activate a sliding of a vertical bar around a frame-angle only after lowering the energy. Thus the only admissible move is to attach the free particle to an horizontal bar, since we want to slide a vertical bar. When this happens, the energy reaches the value $\G_{wa}^*-U_2$. Now the only possibility to slide the bar is to activate a $1$-translation of the vertical bar at cost $U_2$, thus the subsequent moves must be at zero cost, until the last one that costs $-U_2$. This implies that the last configuration has energy $\G_{wa}^*-U_2$. If a $1$-translation of a horizontal bar is activated, the energy increases by $U_1$ and thus it reaches the value $\G_{wa}^*-U_2+U_1>\G_{wa}^*$, which is in contradiction with the optimality of the path.
\end{proof*}

\begin{proof*}{\bf of Lemma \ref{dtilde}}
	By Proposition \ref{cardweak}(a) we know the geometric description of $\widetilde\cD_{wa}$. Starting from a configuration $\h\in\widetilde\cD_{wa}$, since $H(\h)=\G_{wa}^*-\D+U_1-U_2$, by optimality of the path, it is possible to create a free particle only after lowering the energy. This is possible only if $\h\in\widetilde\cQ_{wa}$, where it is possible to detach the protuberance and reattach it to a vertical side, thus we obtain a configuration in $\bar\cQ_{wa}$. This concludes the proof.
\end{proof*}

\noindent
We give the proof of Lemma \ref{entrataweak} here in Appendix \ref{app*}, since it is similar to the proof of \cite[Lemma 5.17]{BN3} that concerns the strongly anisotropic interactions. 

\begin{proof*}{\bf of Lemma \ref{entrataweak}}
	The proof is analogue to the one of \cite[Lemma 5.17]{BN3} done for the strongly anisotropic interactions. The difference is only in case (iii), indeed, considering the time-reversal of the path $\o$ from $\h\in\cwgeo$ to $\vuoto$, if a sliding of a bar around a frame-angle takes place at cost $U_1$, the configuration $\o_{\bar{k}}$ that we obtain does not belong to the set $\cB$ defined in \cite[eq.\ (3.64)]{NOS}, because $s(\o_{\bar{k}})=s^*_{wa}+1$ and $v(\o_{\bar{k}})=2l_2^*-l_1^*-2<p_{min}(\o_{\bar{k}})-1=l_2^*-1$. Thus by \cite[Proposition 11]{NOS} we know that the time-reversal of the path $\o$ visits a configuration $\bar\s\in\cC_{wa}^*$. Thus we can iterate the argument by taking this configuration as $\h$ and the iteration involves a finite number of steps since $\o$ has to reach $\vuoto$. This concludes the proof.
\end{proof*}

\subsection{Additional material for Section \ref{proofiso}}
\label{app}
We give explicit argument to complete the proof of Proposition \ref{isoselless}, considering the cases that were left in Section \ref{dim4}, since the proofs are analogue to the ones discussed in that Section. Due to \cite[Theorem 5.1]{MNOS}, our strategy consists in characterizing the essential saddles that could be visited after attaching the free particle in a bad site.

\medskip
{\bf Case 2A.} Without loss of generality we consider $\h$ as in Figure \ref{fig:figA1}(a). If we are considering the case in which a sequence of $1$-translations of a bar is possible and takes place, then by Lemma \ref{isotrasless}(i) the saddles that are crossed are essential and in $\mathscr{I}_{0}^{\a}\cup\mathscr{I}_{1}^{\a}$. If a sequence of $1$-translations of a bar takes place so that the last configuration has at most two occupied frame-angles and it belongs to the cases 2B, 2C and 2D treated below. Thus we are left to analyze the case in which there is the activation of a sliding of a bar around a frame-angle. 

\setlength{\unitlength}{1.1pt}
\begin{figure}
	\begin{picture}(380,30)(-30,30)
		\thinlines 
		\qbezier[20](30,0)(57.5,0)(85,0)
		\qbezier[20](30,0)(30,32.5)(30,65) 
		\qbezier[20](30,65)(57.5,65)(85,65)
		\qbezier[20](85,0)(85,32.5)(85,65)
		\put(85,35){\line(1,0){5}}
		\put(90,35){\line(0,1){35}}
		\put(90,70){\line(-1,0){65}}
		\put(25,70){\line(0,-1){75}}
		\put(25,-5){\line(1,0){40}}
		\put(65,-5){\line(0,1){5}}
		\put(65,0){\line(1,0){20}}
		\put(85,0){\line(0,1){35}}
		\put(85,75){\line(1,0){5}}
		\put(85,80){\line(1,0){5}}
		\put(85,75){\line(0,1){5}}
		\put(90,75){\line(0,1){5}}
		\put(50,-20){(a)}
		\put(31,30){\tiny{$(l_c-3)\times(l_c-1)$}}
		
		\thinlines 
		\qbezier[20](160,0)(185,0)(210,0)
		\qbezier[20](160,0)(160,35)(160,70) 
		\qbezier[20](160,70)(185,70)(210,70)
		\qbezier[20](210,0)(210,35)(210,70)
		\put(215,0){\line(0,1){75}}
		\put(215,75){\line(-1,0){35}}
		\put(180,75){\line(0,-1){5}}
		\put(180,70){\line(-1,0){25}}
		\put(155,70){\line(0,-1){75}}
		\put(155,-5){\line(1,0){40}}
		\put(195,-5){\line(0,1){5}}
		\put(195,0){\line(1,0){20}}
		\put(220,65){\line(1,0){5}}
		\put(225,65){\line(0,1){5}}
		\put(225,70){\line(-1,0){5}}
		\put(220,70){\line(0,-1){5}}
		\put(180,-20){(b)}
		\put(166,30){\tiny{$(l_c-4)\times l_c$}}
		
		\thinlines
		\qbezier[20](290,0)(320,0)(350,0)
		\qbezier[20](290,0)(290,30)(290,60) 
		\qbezier[20](290,60)(320,60)(350,60)
		\qbezier[20](350,0)(350,30)(350,60)
		\put(350,-5){\line(0,1){35}}
		\put(350,30){\line(1,0){5}}
		\put(355,30){\line(0,1){35}}
		\put(355,65){\line(-1,0){65}}
		\put(285,-5){\line(1,0){65}}
		\put(285,-5){\line(0,1){40}}
		\put(285,35){\line(1,0){5}}
		\put(290,35){\line(0,1){30}}
		\put(310,-20){(c)}
		\put(293,30){\tiny{$(l_c-2)\times(l_c-2)$}}
	\end{picture}
	\vskip 2. cm
	\caption{Case 2A: in (a) we depict a possible starting configuration $\h\in\cigeo$, in (b) the configuration  $\widetilde\h$ obtained from $\h$ after the sliding of the bar $B^e(\h)$ around the frame-angle $c^{en}(\h')$ and in (c) the configuration $\widetilde\h$ obtained from $\h$ after the sliding of the bar $B^s(\h)$ around the frame-angle $c^{sw}(\h')$.}
	\label{fig:figA1}
\end{figure}
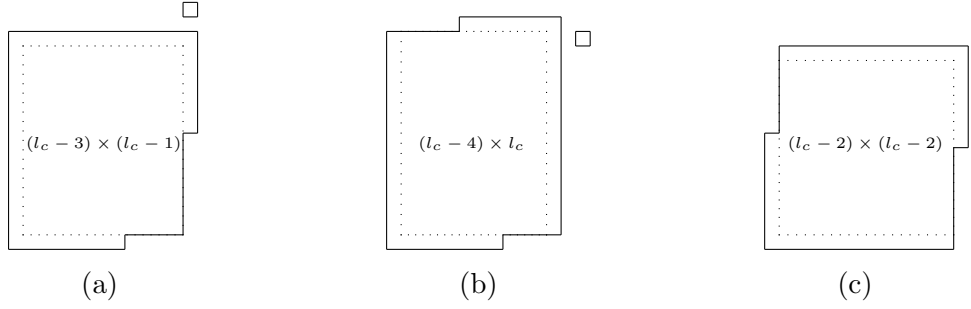

If the free particle is attached to the bar $B^e(\h)$, then it is not possible to complete the sliding of the bar $B^n(\h)$ around the frame-angle $c^{ne}(\h')$, since the condition (\ref{condtrenino}) is not satisfied. Thus by Lemma \ref{trenini}(ii) we know that the saddles that could be crossed in this attempt are unessential. If the free particle is attached to the bar $B^s(\h)$, we conclude as before.

If the free particle is attached to the bar $B^n(\h)$, then it is not possible to complete the sliding of the bar $B^w(\h)$ around the frame-angle $c^{wn}(\h')$ because the condition (\ref{condtrenino}) is not satisfied and thus we can conclude as before. If $|B^e(\h)|<|B^n(\h)|$ the condition (\ref{condtrenino}) is satisfied, thus it is possible to slide the bar $B^e(\h)$ around the frame-angle $c^{en}(\h')$. By Lemma \ref{trenini}(i) the saddles that are possibly visited are essential and, except the last one, they are in $\mathscr{I}_{k,k',1}^{\a,\a'}$. The last configuration visited during this sliding of a bar is $\widetilde\h\in\mathscr{I}_{1}^{\a,\a'}$, has a free particle and it is depicted in Figure \ref{fig:figA1} (b). Starting from such a configuration, if a sequence of $1$-translations of a bar is possible and takes place, then by Lemma \ref{isotrasless}(ii) the saddles that could be crossed are essential and in $\mathscr{I}_{1}^{\a}\cup\mathscr{I}_{2}^{\a}$. If this free particle is now attached to $\partial^-\hbox{CR}(\widetilde\h)$, then the path reaches $\cC_{\pieno}^{\vuoto}(\gi)$ and the saddles are those obtained up to this point. Otherwise, if the free particle is now attached to $B^w(\widetilde\h)$ obtaining the configuration $\h''$, and if a sliding of the bar $B^s(\widetilde\h)$ around the frame-angle $c^{sw}(\h'')$ takes place, by Lemma \ref{trenini}(i) the saddles that could be crossed are essential and, except the last one, are in $\mathscr{I}_{k,k',1}^{\a,\a'}$. The last configuration is in $\cC_{is}^*$, because the cluster is in $\widetilde\cD_{is}$. If this sliding of a bar does not take place, the path $\o$ has to go back to a configuration in $\mathscr{I}_{1}^{\a,\a'}$ and the saddles that could be crossed are already considered. From now on, we can iterarate this argument for a finite number of steps since the path has to reach $\pieno$.

If the free particle is first attached to the bar $B^w(\h)$, we argue in a similar way as before. Indeed, if $|B^s(\h)|<|B^w(\h)|$, it is possible to slide the bar $B^s(\h)$ around the frame-angle $c^{sw}(\h')$. By Lemma \ref{trenini}(i) the saddles that are possibly visited are essential and, except the last one, are in $\mathscr{I}_{k,k',0}^{\a,\a'}$. The last configuration visited during this sliding of a bar is in $\cigeo$ (see Figure \ref{fig:figA1}(c)) and it belongs to case 1A. This concludes case 2A.

\setlength{\unitlength}{1.6pt}
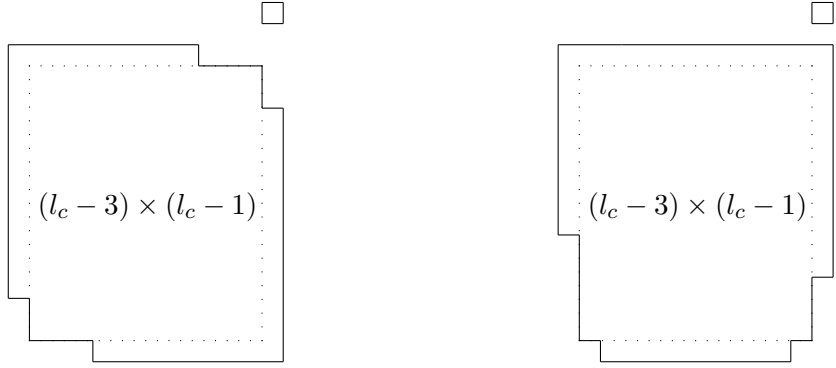
\begin{figure}
	\begin{picture}(380,60)(-30,40)
		\thinlines 
		\qbezier[20](30,20)(57.5,20)(85,20)
		\qbezier[20](30,20)(30,52.5)(30,85) 
		\qbezier[20](30,85)(57.5,85)(85,85)
		\qbezier[20](85,20)(85,52.5)(85,85)
		\put(85,85){\line(-1,0){15}}
		\put(70,85){\line(0,1){5}}
		\put(70,90){\line(-1,0){45}}
		\put(25,90){\line(0,-1){60}}
		\put(25,30){\line(1,0){5}}
		\put(30,30){\line(0,-1){10}}
		\put(30,20){\line(1,0){15}}
		\put(45,15){\line(0,1){5}}
		\put(45,15){\line(1,0){45}}
		\put(90,15){\line(0,1){60}}
		\put(90,75){\line(-1,0){5}}
		\put(85,75){\line(0,1){10}}
		\put(85,95){\line(1,0){5}}
		\put(85,100){\line(1,0){5}}
		\put(85,95){\line(0,1){5}}
		\put(90,95){\line(0,1){5}}
		\put(32,50){$(l_c-3)\times(l_c-1)$}

		\thinlines
		\qbezier[20](160,20)(187.5,20)(215,20)
		\qbezier[20](160,20)(160,52.5)(160,85) 
		\qbezier[20](160,85)(187.5,85)(215,85)
		\qbezier[20](215,20)(215,52.5)(215,85)
		\put(220,90){\line(-1,0){50}}
		\put(170,90){\line(-1,0){15}}
		\put(155,90){\line(0,-1){45}}
		\put(155,45){\line(1,0){5}}
		\put(160,45){\line(0,-1){25}}
		\put(160,20){\line(1,0){5}}
		\put(165,15){\line(0,1){5}}
		\put(165,15){\line(1,0){45}}
		\put(210,15){\line(0,1){5}}
		\put(210,20){\line(1,0){5}}
		\put(215,20){\line(0,1){15}}
		\put(220,35){\line(-1,0){5}}
		\put(220,35){\line(0,1){55}}
		\put(215,95){\line(1,0){5}}
		\put(215,100){\line(1,0){5}}
		\put(215,95){\line(0,1){5}}
		\put(220,95){\line(0,1){5}}
		\put(162,50){$(l_c-3)\times(l_c-1)$}
	\end{picture}
	\vskip 1.4 cm
	\caption{On the left-hand side we depict a possible starting configuration $\h\in\cigeo$ for the case 2B(i) and on the right-hand side a possible starting configuration $\h\in\cigeo$ for the case 2B(ii).}
	\label{fig:figA2}
\end{figure}

\medskip
{\bf Case 2B.} We consider separately the following subcases:

\bi
\item[(i)] the two occupied frame-angles are $c^{\a\a'}(\h)$ and $c^{\a''\a'''}(\h)$, with all the indices $\a,\a',\a'',\a'''$ different between each other (see Figure \ref{fig:figA2} on the left-hand side);
\item[(ii)] the two occupied frame-angles are $c^{\a\a'}(\h)$ and $c^{\a'\a''}(\h)$, with $\a\neq\a''$ (see Figure \ref{fig:figA2} on the right-hand side).
\ei

\medskip
{\bf Case 2B(i).} Without loss of generality we consider $\h$ as in Figure \ref{fig:figA2} on the left-hand side. If we are considering the case in which a sequence of $1$-translations of a bar is possible and takes place, then by Lemma \ref{isotrasless}(i) the saddles that are crossed are essential and they are in $\mathscr{I}_{0}^{\a}\cup\mathscr{I}_{1}^{\a}$. If at least one bar is full and a sequence of $1$-translations of a bar takes place, it is possible to obtain a configuration either with two occupied frame-angles with a bar in common or with three occupied frame-angles. The first case will been analyzed in case 2B(ii) and the latter one has been already analyzed in case 2A. Thus we can reduce our proof to the case in which there is no translation of bars and therefore there is the activation of a sliding of a bar around a frame-angle. First, assume that $|B^n(\h)|<|B^w(\h)|$ and $|B^s(\h)|<|B^e(\h)|$. By Lemma \ref{trenini} the only two possibilities to obtain essential saddles is to attach the free particle to the bar $B^w(\h)$ or $B^e(\h)$ and then slide the bar $B^n(\h)$ around the frame-angle $c^{nw}(\h')$ or $B^s(\h)$ around $c^{se}(\h')$ respectively. If the free particle is first attached to $B^w(\h)$, by Lemma \ref{isotrasless}(i) the saddles that are possibly visited are essential and, except the last one, are in $\mathscr{I}_{k,k',0}^{\a,\a'}$. The last configuration is $\widetilde\h\in\mathscr{I}_{-1}^{\a,\a'}$. Starting from $\widetilde\h$, it is possible to attach the free particle to the bar $B^e(\widetilde\h)$ obtaining a configuration $\h''$, and then slide the bar $B^s(\widetilde\h)$ around the frame-angle $c^{se}(\h'')$. Thus by Lemma \ref{trenini}(i) we know that these saddles are essential and, except the last one, are in $\mathscr{I}_{k,k',0}^{\a,\a'}$. The last configuration is $\bar\h\in\mathscr{I}_{0}^{\a,\a'}$. Starting from $\bar\h$, it is not possible to complete any sliding of a bar around a frame-angle and thus by Lemma \ref{trenini}(ii) we conclude that the saddles that will be possibly crossed are unessential unless a sequence of $1$-translations of bars takes place. In this case, by Lemma \ref{isotrasless}(i) the saddles that could be crossed are essential and in $\mathscr{I}_0^{\a}\cup\mathscr{I}_1^{\a}$. If the free particle is first attached to $B^e(\h)$, we conclude similarly.

Assume now that $|B^w(\h)|<|B^n(\h)|$ and $|B^e(\h)|<|B^s(\h)|$: we argue in the same way as before. If the free particle is first attached to the bar $B^n(\h)$, the essential saddles that could be crossed are in $\mathscr{I}_{k,k',1}^{\a,\a'}$ and the last one is $\widetilde\h\in\mathscr{I}_{1}^{\a,\a'}$ and has a free particle. Again, starting from $\widetilde\h$, if the free particle is attached to the bar $B^s(\widetilde\h)$ obtaining the configuration $\h''$, by Lemma \ref{trenini}(i) we know that the saddles that will be possibly crossed are essential and, except the last one, are in $\mathscr{I}_{k,k',2}^{\a,\a'}$. The last configuration is $\bar\h\in\mathscr{I}_{2}^{\a,\a'}$. If the free particle is first attached to the bar $B^s(\h)$, we conclude as above. Starting from $\bar\h$, it is not possible to complete any sliding of a bar around a frame-angle and thus by Lemma \ref{trenini}(ii) the saddles that could be crossed are unessential unless a sequence of $1$-translations of bars takes place. In the latter case, by Lemma \ref{isotrasless}(iii) the saddles that could be crossed are essential and in $\mathscr{I}_2^{\a}\cup\mathscr{I}_3^{\a}$.

The cases $|B^w(\h)|<|B^n(\h)|$, $|B^s(\h)|<|B^e(\h)|$ and $|B^n(\h)|<|B^w(\h)|$, $|B^e(\h)|<|B^s(\h)|$ can be treated with the same argument as the previous ones and the essential saddles encountered have been already considered. This concludes case 2B(i).

\medskip
{\bf Case 2B(ii).} Without loss of generality we consider $\h$ as on the right-hand side in Figure \ref{fig:figA2}. If we are considering the case in which a sequence of $1$-translations of a bar is possible and takes place, then by Lemma \ref{isotrasless}(i) the saddles that are crossed are essential and they are in $\mathscr{I}_{0}^{\a}\cup\mathscr{I}_{1}^{\a}$. If one bar among $B^w(\h)$ and $B^e(\h)$ is full, then it is possible that a sequence of $1$-translations of the bar $B^s(\h)$ takes place in order to have three occupied frame-angles. This situation has already been analyzed in case 2A. Thus we can reduce our proof to the case in which there is no translation of bars and therefore there is the activation of a sliding of a bar around a frame-angle. If the free particle is attached to the bar $B^s(\h)$, it is not possible to complete any sliding of a bar around a frame-angle at cost $U$ and by Lemma \ref{trenini}(ii) we know that the saddles that could be crossed are unessential. If the free particle is attached to one bar among $B^w(\h)$ and $B^e(\h)$, since (\ref{condtrenino}) is not satisfied it is not possible to complete the sliding of the bar $B^n(\h)$ around the frame-angle $c^{nw}(\h')$ and $c^{ne}(\h')$ respectively. Thus by Lemma \ref{trenini}(ii) the saddles that will be possibly crossed are unessential. If the free particle is attached to the bar $B^n(\h)$, then it is possible to complete the sliding of the bar $B^w(\h)$ or $B^e(\h)$ around the frame-angle $c^{wn}(\h')$ or $c^{en}(\h')$ respectively. Thus by Lemma \ref{trenini}(i) we know that the saddles that will be crossed are essential and, except the last one, are in $\mathscr{I}_{k,k',1}^{\a,\a'}$. The last configuration is $\widetilde\h\in\mathscr{I}_{1}^{\a,\a'}$ that has a free particle. If this free particle is now attached to $\partial^-\hbox{CR}(\widetilde\h)$, then the path reaches $\cC_{\pieno}^{\vuoto}(\gi)$ and the saddles are those obtained up to this point. Otherwise, since condition (\ref{condtrenino}) is not satisfied, it is not possible to complete any sliding of a bar around a frame-angle, thus by Lemma \ref{trenini}(ii) the saddles that will be possibly crossed during the sliding of a bar are unessential. On the other hand, if a sequence of $1$-translations of a bar is possible and takes place, then, starting from $\widetilde\h$, by Lemma \ref{isotrasless}(ii) the saddles that could be crossed are essential and in $\mathscr{I}_{1}^{\a}\cup\mathscr{I}_{2}^{\a}$. This concludes case 2B(ii).

\setlength{\unitlength}{1.6pt}
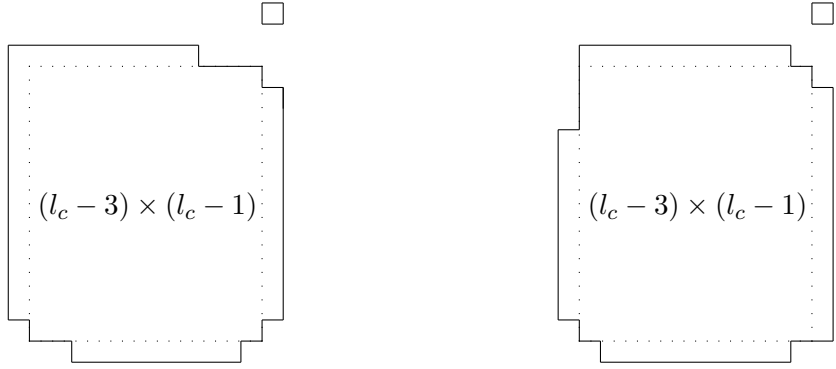
\begin{figure}
	\begin{picture}(380,40)(-30,40)
		\thinlines 
		\qbezier[20](30,20)(57.5,20)(85,20)
		\qbezier[20](30,20)(30,52.5)(30,85) 
		\qbezier[20](30,85)(57.5,85)(85,85)
		\qbezier[20](85,20)(85,52.5)(85,85)
		\put(70,90){\line(-1,0){45}}
		\put(70,90){\line(0,-1){5}}
		\put(70,85){\line(1,0){15}}
		\put(25,90){\line(0,-1){65}}
		\put(25,25){\line(1,0){5}}
		\put(30,25){\line(0,-1){5}}
		\put(30,20){\line(1,0){10}}
		\put(40,15){\line(0,1){5}}
		\put(40,15){\line(1,0){40}}
		\put(80,15){\line(0,1){5}}
		\put(80,20){\line(1,0){5}}
		\put(85,20){\line(0,1){5}}
		\put(90,25){\line(-1,0){5}}
		\put(90,25){\line(0,1){55}}
		\put(90,80){\line(-1,0){5}}
		\put(85,80){\line(0,1){5}}
		\put(90,75){\line(0,1){5}}
		\put(85,95){\line(1,0){5}}
		\put(85,100){\line(1,0){5}}
		\put(85,95){\line(0,1){5}}
		\put(90,95){\line(0,1){5}}
		\put(32,50){$(l_c-3)\times(l_c-1)$}
		
		\thinlines
		\qbezier[20](160,20)(187.5,20)(215,20)
		\qbezier[20](160,20)(160,52.5)(160,85) 
		\qbezier[20](160,85)(187.5,85)(215,85)
		\qbezier[20](215,20)(215,52.5)(215,85)
		\put(210,90){\line(-1,0){50}}
		\put(210,90){\line(0,-1){5}}
		\put(210,85){\line(1,0){5}}
		\put(160,90){\line(0,-1){20}}
		\put(155,70){\line(1,0){5}}
		\put(155,70){\line(0,-1){45}}
		\put(155,25){\line(1,0){5}}
		\put(160,25){\line(0,-1){5}}
		\put(160,20){\line(1,0){5}}
		\put(165,15){\line(0,1){5}}
		\put(165,15){\line(1,0){45}}
		\put(210,15){\line(0,1){5}}
		\put(210,20){\line(1,0){10}}
		\put(220,20){\line(0,1){60}}
		\put(220,80){\line(-1,0){5}}
		\put(215,80){\line(0,1){5}}
		\put(215,95){\line(1,0){5}}
		\put(215,100){\line(1,0){5}}
		\put(215,95){\line(0,1){5}}
		\put(220,95){\line(0,1){5}}
		\put(162,50){$(l_c-3)\times(l_c-1)$}

	\end{picture}
	\vskip 1.4 cm
	\caption{On the left-hand side we depict a possible starting configuration $\h\in\cigeo$ for the case 2C and on the right-hand side a possible starting configuration $\h\in\cigeo$ for the case 2D.}
	\label{fig:figA3}
\end{figure}

\medskip
{\bf Case 2C.} Without loss of generality we consider $\h$ as in Figure \ref{fig:figA3} on the left-hand side. If we are considering the case in which a sequence of $1$-translations of a bar is possible and takes place, then by Lemma \ref{isotrasless}(i) the saddles that are crossed are in $\mathscr{I}_{0}^{\a}\cup\mathscr{I}_{1}^{\a}$. Thus we can reduce our proof to the case in which there is no $1$-translation of a bar and therefore there is only the activation of a sliding of a bar around a frame-angle. Starting from this configuration it is possible to obtain two occupied frame-angles via a sequence of $1$-translations of a bar: this situation has been already analyzed in case 2B. If the free particle is attached to the bar $B^e(\h)$ or $B^s(\h)$, since it is not possible to complete any sliding of bar around a frame-angle at cost $U$, by Lemma \ref{trenini}(ii) we know that the saddles that could be crossed are unessential. If $|B^w(\h)|<|B^n(\h)|$ (resp.\ $|B^n(\h)|<|B^w(\h)|$) and the free particle is attached to the bar $B^n(\h)$ (resp.\ $B^w(\h)$), then it is possible to complete a sliding of the bar $B^w(\h)$ (resp.\ $B^n(\h)$) around the frame-angle $c^{wn}(\h')$ (resp.\ $c^{nw}(\h')$). Thus by Lemma \ref{trenini}(i) the saddles that could be crossed are essential and in $\mathscr{I}_{k,k',1}^{\a,\a'}$ (resp.\ $\mathscr{I}_{k,k',0}^{\a,\a'}$), except the last one that is $\widetilde\h\in\mathscr{I}_{1}^{\a,\a'}$ (resp.\ $\widetilde\h\in\mathscr{I}_{-1}^{\a,\a'}$) with a free particle. If this free particle is now attached to $\partial^-\hbox{CR}(\widetilde\h)$, then the path reaches $\cC_{\pieno}^{\vuoto}(\gi)$ and the saddles are those obtained up to this point. Otherwise, condition (\ref{condtrenino}) is not satisfied, thus it is not possible to complete any sliding of a bar around a frame-angle. Therefore by Lemma \ref{trenini}(ii) we know that the saddles that could be visited are unessential unless a sequence of $1$-translations of a bar is possible and takes place. In this case, by Lemma \ref{isotrasless}(ii) (resp.\ Lemma \ref{isotrasless}(i)) the saddles that could be crossed are essential and in $\mathscr{I}_{1}^{\a}\cup\mathscr{I}_{2}^{\a}$ (resp.\ $\mathscr{I}_{0}^{\a}$). This concludes case 2C.

\medskip
{\bf Case 2D.} Without loss of generality we consider $\h$ as in Figure \ref{fig:figA3} on the right-hand side. If we are considering the case in which a sequence of $1$-translations of a bar is possible and takes place, then by Lemma \ref{isotrasless}(i) the saddles that are crossed are essential and in $\mathscr{I}_{0}^{\a}\cup\mathscr{I}_{1}^{\a}$. Thus we can reduce our proof to the case in which there is no $1$-translation of a bar and therefore there is only the activation of a sliding of a bar around a frame-angle. Starting from this configuration, it is possible to obtain one or two occupied frame-angles via a sequence of $1$-translations: these situations have been already analyzed in cases 2C and 2B respectively. If the free particle is attached to one of the bars, since it is not possible to complete any sliding of bar around a frame-angle at cost $U$, by Lemma \ref{trenini}(ii) we know that the saddles that could be crossed are unessential. This concludes case 2D.

\subsection{Additional material for Section \ref{proofweak}}
\label{app2}

\medskip
\begin{proof*}{\bf of Lemma \ref{gateweak}}
Consider $\o\in(\vuoto\ra\pieno)_{opt}$. If $\o\cap\cwgeo\neq\emptyset$, we get the claim. Thus we can reduce our analysis to the case in which the path $\o$ reaches the set $\cP$ in a configuration $\h\in\cP\setminus\cwgeo$. We set $\o=(\vuoto,\o_1,...,\o_k,\h)\circ\bar\o$, where $\bar\o$ is a path that connects $\h$ to $\pieno$ such that $\max_{\s\in\o}H(\s)\leq\gw$. We are interested in the time-reversal of the path $\o$. Since $\h\in\cP\setminus\cwgeo$, we know that it is composed by the union of a cluster $\hbox{CR}^-(\h)=\cR(l_1^*-2,l_2^*-2)$, such that at least one frame-angle of $\hbox{CR}^-(\h)$ is empty, a free particle and four bars attached to the four sides of $\hbox{CR}^-(\h)$ in such a way that $\h$ contains $n^c_{wa}+1$ particles (see (\ref{defnwa}) for the definition of $n_{wa}^c$). For the entire proof we refer to Figure \ref{fig:figesempiostrong}, where for $int=wa$ the horizontal and vertical lengths have to be changed to $l_1^*$ and $l_2^*$ respectively. Suppose that $\hbox{CR}^-(\h)$ contains $x$ empty frame-angles, with $1\leq x\leq4$, (see Figure \ref{fig:figesempiostrong}(a) to visualize the configuration $\h$ in the case $x=1$). Since $H(\h)=\gw$, the move from $\h$ to $\o_k$ must have a non-positive cost and thus the unique admissible moves are:

\bi
\item[(i)] either moving the free particle at zero cost;
\item[(ii)] or removing the free particle;
\item[(iii)] or attaching the free particle at cost $-U_1$ (see Figure \ref{fig:figesempiostrong}(b)) or $-U_2$, or $-U_1-U_2$.
\ei

\noindent
{\bf Case (i).} In this case the configuration $\o_{k}$ is analogue to $\h$ and therefore we can iterate this argument by taking this configuration as $\h$.

\medskip
\noindent
{\bf Case (ii).} In this case $H(\o_k)=\gw-\D$. We may assume that the configuration $\o_{k-1}$ is not obtained by $\o_{k}$ via adding a free particle, otherwise $\o_{k-1}$ is analogue to $\h$ and thus we can iterate the argument by taking this configuration as $\h$. By the optimality of the path, again considering the time-reversal, we deduce that the unique admissible move to obtain $\o_{k-1}$ from $\o_{k}$ is breaking a horizontal (resp.\ vertical) bar at cost $U_1$ (resp.\ $U_2$). Thus it is possible that either a sequence of $1$-translations of a bar or a sliding of a bar around a frame-angle takes place. In the first case, we obtain a configuration that is analogue to $\o_{k-1}$ and thus we can iterate the argument for a finite number of steps, since the path has to reach $\vuoto$. In the latter case, by Remark \ref{barreweak}(ii) we deduce that the condition (\ref{condtrenino}) is not satisfied and therefore it is not possible to complete any sliding of a bar around a frame-angle. This implies that the unique admissible moves are the reverse ones, thus we obtain a configuration that is analogue to $\o_{k-1}$ and therefore we can iterate the argument for a finite number of steps, since the path has to reach $\vuoto$. In this way we can reduce ourselves to consider the case (iii).

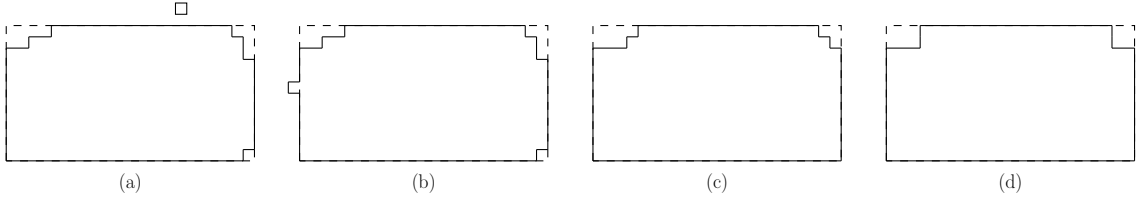
\begin{figure}
	\centering
	\begin{tikzpicture}[scale=0.3,transform shape]

		\draw (16,2.5)--(16,3);
		\draw (16,3) -- (16.5,3);
		\draw (16.5,3) -- (16.5,7);
		\draw (16.5,7) -- (16,7);
		\draw (16,7) -- (16,8);
		\draw (16,8)--(15.5,8);
		\draw(15.5,8)--(15.5,8.5);
		\draw (15.5,8.5)--(7.5,8.5);
		\draw (7.5,8.5)--(7.5,8);
		\draw [dashed] (5.5,2.5) rectangle (16.5,8.5);
		\draw (5.5,2.5) -- (16,2.5);
		\draw (5.5,2.5)--(5.5,7.5);
		\draw (5.5,7.5)--(6.5,7.5);
		\draw (7.5,8)--(6.5,8);
		\draw (6.5,8)--(6.5,7.5);
		\draw (13,9) rectangle (13.5,9.5);
		\node at (11,1.5){\Huge{(a)}};

		\draw (29,2.5)--(29,3);
		\draw (29,3) -- (29.5,3);
		\draw (29.5,3) -- (29.5,7);
		\draw (29.5,7) -- (29,7);
		\draw (29,7) -- (29,8);
		\draw (29,8)--(28.5,8);
		\draw(28.5,8)--(28.5,8.5);
		\draw (28.5,8.5)--(20.5,8.5);
		\draw (20.5,8.5)--(20.5,8);
		\draw [dashed] (18.5,2.5) rectangle (29.5,8.5);
		\draw (18.5,2.5) -- (29,2.5);
		\draw (18.5,2.5)--(18.5,5.5);
		\draw (18.5,5.5)--(18,5.5);
		\draw (18,5.5)--(18,6);
		\draw (18,6)--(18.5,6);
		\draw (18.5,6)--(18.5,7.5);
		\draw (18.5,7.5)--(19.5,7.5);
		\draw (20.5,8)--(19.5,8);
		\draw (19.5,8)--(19.5,7.5);
		\node at (24,1.5){\Huge{(b)}};

		\draw (42.5,2.5)--(42.5,3);
		\draw (42,2.5) -- (42.5,2.5);
		\draw (42.5,3) -- (42.5,7.5);
		\draw (42.5,7.5) -- (42,7.5);
		\draw (42,7.5) -- (42,8);
		\draw (42,8)--(41.5,8);
		\draw(41.5,8)--(41.5,8.5);
		\draw (41.5,8.5)--(33.5,8.5);
		\draw (33.5,8.5)--(33.5,8);
		\draw [dashed] (31.5,2.5) rectangle (42.5,8.5);
		\draw (31.5,2.5) -- (42,2.5);
		\draw (31.5,2.5)--(31.5,7.5);
		\draw (31.5,7.5)--(33,7.5);
		\draw (33.5,8)--(33,8);
		\draw (33,8)--(33,7.5);
		\node at (37,1.5){\Huge{(c)}};

		\draw (55.5,2.5)--(55.5,3);
		\draw (55,2.5) -- (55.5,2.5);
		\draw (55.5,3) -- (55.5,7.5);
		\draw (55.5,7.5) -- (54.5,7.5);
		\draw (54.5,7.5) -- (54.5,8);
		\draw(54.5,8)--(54.5,8.5);
		\draw (54.5,8.5)--(46,8.5);
		\draw (46,8.5)--(46,8);
		\draw [dashed] (44.5,2.5) rectangle (55.5,8.5);
		\draw (44.5,2.5) -- (55,2.5);
		\draw (44.5,2.5)--(44.5,7.5);
		\draw (44.5,7.5)--(46,7.5);
		\draw (46,8)--(46,7.5);
		\node at (50,1.5){\Huge{(d)}};

	\end{tikzpicture}
	
	\vskip 0 cm
	\caption{Here we depict in (a) the configuration $\h$; in (b) the configuration obtained by $\h$ by attaching the free particle at cost $-U_1$ to $B^w(\h)$; in (c) the configuration $\h'$ obtained from $\h$ by attaching the free particle to $c^{se}(\h)$ and then detaching the particle in $c^{nw}(\hbox{CR}^-(\h))$ and attach it to $B^e(\h)$, and in (d) the configuration $\h''$ obtained from $\h'$ by detaching the particle in $c^{ne}(\hbox{CR}^-(\h'))$ attach it to $B^n(\h')$.}
	\label{fig:figesempiostrong}
\end{figure}

\medskip
\noindent
{\bf Case (iii).} (a) We consider the case where from $\h$, again considering the time-reversal, we attach a particle at cost $-U_1$ in $\partial^+\hbox{CR}(\h)$ giving rise to the configuration $\o_k$, i.e., $H(\o_k)=\gw-U_1$ (see Figure \ref{fig:figesempiostrong}(b)). Thus it is possible that either a sequence of $1$-translations of a bar or a sliding of a bar around a frame-angle takes place. In the first case, we obtain a configuration that is analogue to $\o_{k}$ and thus we can iterate the argument for a finite number of steps, since the path has to reach $\vuoto$. In the latter case, by Remark \ref{barreweak}(ii) we deduce that the condition (\ref{condtrenino}) is not satisfied and therefore it is not possible to complete any sliding of a bar around a frame-angle. This implies that the unique admissible moves are the reverse ones, thus we obtain a configuration that is analogue to $\o_{k}$ and therefore we can iterate the argument for a finite number of steps, since the path has to reach $\vuoto$.

(b) We consider the case where from $\h$, again considering the time-reversal, we attach a particle at cost $-U_2$ in $\partial^+\hbox{CR}(\h)$ giving rise to the configuration $\o_k$, i.e., $H(\o_k)=\gw-U_2$. We argue in a similar way as above.

(c) We consider the case where from $\h$, again considering the time-reversal, we attach a particle at cost $-U_1-U_2$ in $\partial^-\hbox{CR}(\h)$ giving rise to the configuration $\o_k$, i.e., $H(\o_k)=\gw-U_1-U_2$. Thus it is possible either to have a sequence of $1$-translations of a bar, or to have a sliding of a bar around a frame-angle, or to detach a particle at cost $U_1+U_2$. In the first two possibilities, analogously to what has been discussed previously in (a) and (b), the unique admissible moves are the reverse ones and therefore we conclude as above. In the latter possibility, we have that either $\o_{k-1}$ is obtained from $\o_k$ by detaching a particle from a bar at cost $U_1+U_2$ or from a corner of $\h$ that is in $\hbox{CR}^-(\h)$. In the first case, the particle can be attached to an empty frame-angle of $\hbox{CR}^-(\h)$ and we can repeat these steps at most $x-1$ times (if $x\geq2$), that implies that there exists $\bar{k}<k-1$ such that $\o_{\bar{k}}$ is composed by the union of a free particle and a rectangle $\cR(l_1^*-2,l_2^*-2)$ with four bars attached to its four sides in such a way $\o_{\bar{k}}$ contains $n_{wa}^c+1$ particles, namely $\o_{\bar{k}}\in\cwgeo$. In the second case, we may assume that the detached particle is attached to a bar in $\partial^-\hbox{CR}(\h)$ giving rise to a configuration $\h'$ (see Figure \ref{fig:figesempiostrong}(c)), otherwise we obtain a configuration that is analogue to $\h$. Starting from $\h'$, similarly we obtain $\h''$ (see Figure \ref{fig:figesempiostrong}(d)) if $\h'$ has a corner in $\hbox{CR}^-(\h')$. If this is the case, we can proceed in a similar way until we obtain a configuration $\h'''$ that has no corner in $\hbox{CR}^-(\h''')$. Starting from $\h'''$, by the optimality of the path we deduce that the unique admissible moves are the reverse ones and therefore the path goes back to $\h$. This concludes the proof.
\end{proof*}

\medskip
\begin{proof*}{\bf of Proposition \ref{weakselless}}
Consider a configuration $\h\in\cwgeo(2)$ such that $\h=(\hat\h,x)$, with $\hat\h\in\bar\cD_{wa}$ and $d(\hat\h,x)=2$. By Proposition \ref{cardweak}(a) we deduce that that $\hat\h$ consists of an $(l_1^*-2)\times(l_2^*-2)$ rectangle with four bars $B^\a$, with $\a\in\{n,s,w,e\}$, attached to its four sides satisfying
\be{}
1\leq|B^{w}(\h)|,|B^e(\h)|\leq l_2^*, \quad l_1^*-l_2^*+1\leq|B^{n}(\h)|,|B^s(\h)|\leq l_1^*, 
\ee

\noindent
and
\be{} 
\displaystyle\sum_\a |B^\a(\h)|-\displaystyle\sum_{\a\a'\in\{nw,ne,sw,se\}}|c^{\a\a'}(\h)|=2l_1^*+l_2^*-3.
\ee

\noindent
Assume that the free particle is attached in a bad site obtaining a configuration $\h'\in\cC_{wa}^B$. Due to \cite[Theorem 5.1]{MNOS} and Proposition \ref{selle3}, our strategy consists in characterizing the essential saddles that could be visited after attaching the free particle in a bad site. By Remark \ref{barreweak}(i) we consider separately the following cases:

\begin{itemize}
	\item[A.] three frame-angles of CR$(\hat\h)$ are occupied;
	\item[B.] two frame-angles of CR$(\hat\h)$ are occupied;
	\item[C.] one frame-angle of CR$(\hat\h)$ is occupied;
	\item[D.] no frame-angle of CR$(\hat\h)$ is occupied.
\end{itemize}

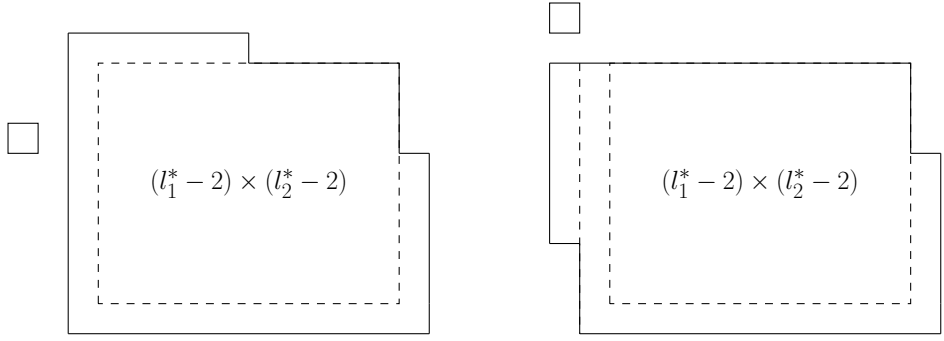
\begin{figure}
	\centering
	\begin{tikzpicture}[scale=0.4,transform shape]

		\draw [dashed] (6,3) rectangle (16,11);
		\node at (11,7) {\Huge{$(l_1^*-2)\times(l_2^*-2)$}};
		\draw(11,11)--(11,12);
		\draw(11,12)--(5,12);
		\draw(5,12)--(5,2);
		\draw(5,2)--(17,2);
		\draw(17,2)--(17,3);
		\draw(17,3)--(17,8);
		\draw(17,8)--(16,8);
		\draw(16,8)--(16,11);
		\draw(16,11)--(11,11);
		\draw(3,8) rectangle (4,9);

		\draw [dashed] (23,3) rectangle (33,11);
		\node at (28,7) {\Huge{$(l_1^*-2)\times(l_2^*-2)$}};
		\draw(28,11)--(21,11);
		\draw[dashed] (22,11)--(22,2);
		\draw(21,11)--(21,5);
		\draw(21,5)--(22,5);
		\draw(22,5)--(22,2);
		\draw(22,2)--(34,2);
		\draw(34,2)--(34,3);
		\draw(34,3)--(34,8);
		\draw(34,8)--(33,8);
		\draw(33,8)--(33,11);
		\draw(33,11)--(28,11);
		\draw(21,12) rectangle (22,13);

	\end{tikzpicture}
	
	\vskip 0 cm
	\caption{Case A: on the left-hand side we represent a possible starting configuration $\h\in\cwgeo$ and on the right-hand side the configuration $\widetilde\h$ obtained from $\h$ after the sliding of the bar $B^n(\h)$ around the frame-angle $c^{nw}(\h')$.}
	\label{fig:figA14}
\end{figure}

Note that from case A one can go to the other cases and viceversa, but since the path has to reach $\pieno$ this back and forth must end in a finite number of steps.

\noindent
{\bf Case A.} Without loss of generality we consider $\h$ as in Figure \ref{fig:figA14} on the left-hand side. If we are considering the case in which a $1$-translation of a bar is possible and takes place, then by Lemma \ref{trasless} the saddles that are crossed are essential and in $\cN_0^{\a'}\cup\cN_1^{\a}$. If a sequence of $1$-translations of a bar takes place in such a way that the last configuration has at most two occupied frame-angles, then the saddles that could be visited starting from such a configuration will be analyzed in cases B, C and D. Thus we are left to analyze the case in which there is the activation of a sliding of a bar around a frame-angle. In the following we quickly exclude the cases in which the particle is attached to $B^n(\h)$, $B^s(\h)$ or $B^e(\h)$ and then explain the more interesting case in which it is attached to $B^w(\h)$ giving rise to Figure \ref{fig:figA14} on the right-hand side. If the free particle is attached to the bar $B^n(\h)$ (resp.\ $B^s(\h)$), by Lemma \ref{coltorow} we know that it is not possible to complete the sliding of the bar $B^w(\h)$ (resp.\ $B^e(\h)$) around the frame-angle $c^{wn}(\h')$ (resp.\ $c^{es}(\h')$). If the free particle is attached to the bar $B^e(\h)$ or $B^w(\h)$, then it is not possible to slide the bar $B^s(\h)$ around the frame-angle $c^{se}(\h')$ or $c^{sw}(\h')$ respectively, since (\ref{condtrenino}) is not satisfied. In the last two cases by Lemma \ref{trenini}(ii), we know that the saddles that are visited are unessential. This implies that the unique possibility to activate and complete a sliding of a bar around a frame-angle is attaching the free particle to the bar $B^w(\h)$, then sliding the bar $B^n(\h)$ around the frame-angle $c^{nw}(\h')$ when $|B^n(\h)|<|B^w(\h)|$, otherwise (\ref{condtrenino}) is not satisifed. The saddles that are possibly visited by the sliding path are in $\cN_{k,k'}^{\a,\a'}$ (see definition (\ref{defpwatrenino})) except the last one, thus by Lemma \ref{trenini}(i) they are essential. The last configuration visited during this sliding of a bar is depicted in Figure \ref{fig:figA14} on the right-hand side. This configuration has energy $\gw-U_1+U_2$ and therefore it is not a saddle and is in $\cC_{\vuoto}^{\pieno}(\gw-H(\pieno))$. By Propositions \ref{selle2} and \ref{Kvuoto}(ii)-(b), the latter implies that the saddles that could be visited are either unessential or in $\cN_{2,k'}^{\a,\a'}$ and therefore the case A is concluded.

\begin{figure}
	\centering
	\begin{tikzpicture}[scale=0.32,transform shape]

		\draw [dashed] (6,3) rectangle (16,11);
		\node at (11,7) {\Huge{$(l_1^*-2)\times(l_2^*-2)$}};
		\draw(11,11)--(11,12);
		\draw(11,12)--(5,12);
		\draw(5,12)--(5,4);
		\draw(5,4)--(6,4);
		\draw(6,4)--(6,3);
		\draw(6,3)--(7,3);
		\draw(7,3)--(7,2);
		\draw(7,2)--(17,2);
		\draw(17,2)--(17,11);
		\draw(17,11)--(16,11);
		\draw(16,11)--(11,11);
		\draw(3,8) rectangle (4,9);
		\node at (12,0.5) {\Huge{(a)}};

		\draw [dashed] (23,3) rectangle (33,11);
		\node at (28,7) {\Huge{$(l_1^*-2)\times(l_2^*-2)$}};
		\draw[dashed](22,11)--(22,4);
		\draw(22,4)--(23,4);
		\draw(23,4)--(23,3);
		\draw(23,3)--(24,3);
		\draw(24,3)--(24,2);
		\draw(24,2)--(34,2);
		\draw(34,2)--(34,11);
		\draw(34,11)--(33,11);
		\draw(33,11)--(21,11);
		\draw(21,11)--(21,5);
		\draw(21,5)--(22,5);
		\draw(22,5)--(22,4);
		\draw(21,12) rectangle (22,13);
		\node at (29,0.5) {\Huge{(b)}};
		
		\draw [dashed] (40,3) rectangle (50,11);
		\node at (45,7) {\Huge{$(l_1^*-2)\times(l_2^*-2)$}};
		\draw(51,11)--(51,12);
		\draw(51,11)--(51,7);
		\draw(51,7)--(50,7);
		\draw(51,12)--(39,12);
		\draw(39,12)--(39,4);
		\draw(39,4)--(40,4);
		\draw(40,4)--(40,3);
		\draw(40,3)--(41,3);
		\draw(41,3)--(41,2);
		\draw(41,2)--(49,2);
		\draw(49,2)--(49,3);
		\draw(49,3)--(50,3);
		\draw(50,3)--(50,7);
		\draw(37,8) rectangle (38,9);
		\node at (46,0.5) {\Huge{(c)}};
		
	\end{tikzpicture}
	
	\vskip -0.2 cm
	\caption{Case B(i): in (a) we depict a possible starting configuration $\h\in\cwgeo$ and in (b) the configuration $\widetilde\h$ obtained from $\h$ after the sliding of the bar $B^n(\h)$ around the frame-angle $c^{nw}(\h')$. Case B(ii): in (c) we depict a possible starting configuration $\h\in\cwgeo$.}
	\label{fig:figA5}
\end{figure}

\medskip
\noindent 
{\bf Case B.} If we are considering the case in which a $1$-translation of a bar is possible and takes place, then by Lemma \ref{trasless} the saddles that are crossed are essential and in $\cN_0^{\a'}\cup\cN_1^{\a}$. We consider separately the following subcases:

\bi
\item[(i)] The two occupied frame-angles are $c^{\a\a'}(\h)$ and $c^{\a''\a'''}(\h)$, with all the indices $\a,\a',\a'',\a'''$ different between each other (see Figure \ref{fig:figA5}(a));
\item[(ii)] The two occupied frame-angles are $c^{\a\a'}(\h)$ and $c^{\a'\a''}(\h)$, with $\a'\in\{n,s\}$ and $\a\neq\a''$ (see Figure \ref{fig:figA5}(c));
\item[(iii)] The two occupied frame-angles are $c^{\a\a'}(\h)$ and $c^{\a'\a''}(\h)$, with $\a'\in\{e,w\}$ and $\a\neq\a''$ (see Figure \ref{fig:figA6} on the left-hand side).
\ei

{\bf Case B(i).} Without loss of generality we consider $\h$ as in Figure \ref{fig:figA5}(a). We can reduce our proof to the case in which there is no translation of a bar and therefore there is the activation of a sliding of a bar around a frame-angle. If the free particle is attached to the bar $B^n(\h)$ (resp.\ $B^s(\h)$), by Lemma \ref{coltorow} we know that it is not possible to complete the sliding of the bar $B^w(\h)$ (resp.\ $B^e(\h)$) around the frame-angle $c^{wn}(\h')$ (resp.\ $c^{es}(\h')$). By Lemma \ref{trenini}(ii), this implies that the saddles that could be crossed are unessential. Note that if the free particle is attached to the bar $B^n(\h)$ (resp.\ $B^s(\h)$), it is not possible to slide the bar $B^e(\h)$ (resp.\ $B^w(\h)$) by definition. If the free particle is attached to the bar $B^w(\h)$ (resp.\ $B^e(\h)$) it is possible to slide the bar $B^n(\h)$ (resp.\ $B^s(\h)$) around the frame-angle $c^{nw}(\h')$ (resp.\ $c^{se}(\h')$) when $|B^n(\h)|<|B^w(\h)|$ (resp. $|B^s(\h)|<|B^e(\h)|$), otherwise (\ref{condtrenino}) is not satisifed. The saddles that are possibly visited by the sliding path are in $\cN_{k,k'}^{\a,\a'}$ except the last one, thus by Lemma \ref{trenini}(i) they are essential. The last configuration visited during this sliding of a bar is depicted in Figure \ref{fig:figA5}(b). This configuration has energy $\gw-U_1+U_2$ and therefore it is not a saddle and is in $\cC_{\vuoto}^{\pieno}(\gw-H(\pieno))$. By Propositions \ref{selle2} and \ref{Kvuoto}(ii)-(b), the latter implies that the saddles that could be visited are either unessential or in $\cN_{2,k'}^{\a,\a'}$ and therefore the case B(i) is concluded.

{\bf Case B(ii).} Without loss of generality we consider $\h$ as in Figure \ref{fig:figA5}(c). If one bar among $B^w(\h)$ and $B^e(\h)$ is full, it is possible to translate $B^s(\h)$ in order to have three occupied frame-angles. This situation has already been analyzed in case A. Thus we can reduce our proof to the case in which there is no translation of a bar and therefore there is the activation of a sliding of a bar around a frame-angle. If the free particle is attached to the bar $B^n(\h)$ or $B^s(\h)$, by Lemma \ref{coltorow} we know that it is not possible to complete the sliding of a vertical bar around any frame-angle. If the free particle is attached to the bar $B^w(\h)$ or $B^e(\h)$, since the bar $B^n(\h)$ is full, we deduce that (\ref{condtrenino}) is not satisfied. This implies that it is not possible to slide the bar $B^n(\h)$ around the frame-angle $c^{nw}(\h')$ and $c^{ne}(\h')$. In the last two cases by Lemma \ref{trenini}(ii) we know that the saddles that are visited are unessential. This concludes case B(ii).

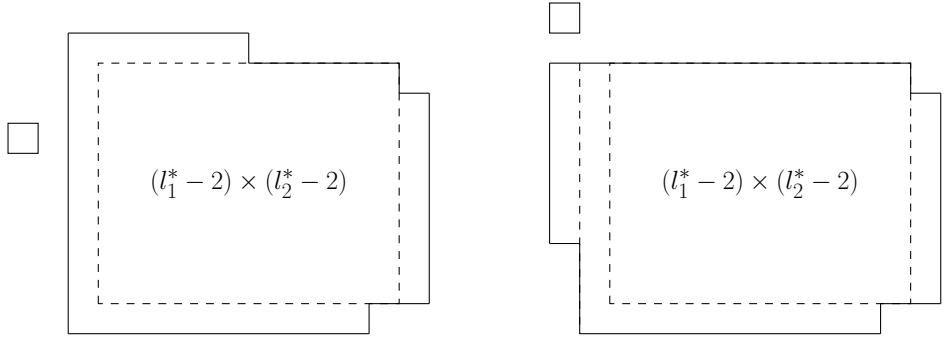
\begin{figure}
	\centering
	\begin{tikzpicture}[scale=0.4,transform shape]

		\draw [dashed] (6,3) rectangle (16,11);
		\node at (11,7) {\Huge{$(l_1^*-2)\times(l_2^*-2)$}};
		\draw(11,11)--(11,12);
		\draw(11,12)--(5,12);
		\draw(5,12)--(5,2);
		\draw(5,2)--(15,2);
		\draw(15,2)--(15,3);
		\draw(15,3)--(17,3);
		\draw(17,3)--(17,10);
		\draw(17,10)--(16,10);
		\draw(16,10)--(16,11);
		\draw(16,11)--(11,11);
		\draw(3,8) rectangle (4,9);

		\draw [dashed] (23,3) rectangle (33,11);
		\node at (28,7) {\Huge{$(l_1^*-2)\times(l_2^*-2)$}};
		\draw(28,11)--(21,11);
		\draw[dashed] (22,11)--(22,2);
		\draw(21,11)--(21,5);
		\draw(21,5)--(22,5);
		\draw(22,5)--(22,2);
		\draw(22,2)--(32,2);
		\draw(32,2)--(32,3);
		\draw(32,3)--(34,3);
		\draw(34,3)--(34,10);
		\draw(34,10)--(33,10);
		\draw(33,10)--(33,11);
		\draw(33,11)--(28,11);
		\draw(21,12) rectangle (22,13);
		
	\end{tikzpicture}
	
	\vskip -0.2 cm
	\caption{Case B(iii): on the left-hand side we depict a possible starting configuration $\h\in\cwgeo$ and on the right-hand side the configuration $\widetilde\h$ obtained from $\h$ after the sliding of the bar $B^n(\h)$ around the frame-angle $c^{nw}(\h')$.}
	\label{fig:figA6}
\end{figure}

{\bf Case B(iii).} Without loss of generality we consider $\h$ as in Figure \ref{fig:figA6} on the left-hand side. If the bar $B^n(\h)$ (resp.\ $B^s(\h)$) is full, it is possible to translate $B^e(\h)$ to occupy the frame-angle $c^{ne}(\h')$ (resp.\ $c^{se}(\h')$). This situation has already been analyzed in case A. Otherwise, it is possible to translate a bar with one occupied frame-angle in order to have two occupied frame-angles in such a way that they have no bar in common. This situation has already been analyzed in case B(i). Thus we can reduce our proof to the case in which there is no translation of a bar and therefore there is the activation of a sliding of a bar around a frame-angle. If the free particle is attached to the bar $B^n(\h)$ (resp.\ $B^s(\h)$), by Lemma \ref{coltorow} we know that it is not possible to complete the sliding of the bar $B^w(\h)$ around the frame-angle $c^{wn}(\h')$ (resp.\ $c^{ws}(\h')$). If the free particle is attached to the bar $B^e(\h)$, we deduce that (\ref{condcorner}) is not satisfied. In the last two cases by Lemma \ref{trenini}(ii) we know that the saddles that are visited are unessential. If the free particle is attached to the bar $B^w(\h)$, it is possible to slide the bar $B^n(\h)$ (resp.\ $B^s(\h)$) around the frame-angle $c^{nw}(\h')$ (resp.\ $c^{sw}(\h')$) when $|B^n(\h)|<|B^w(\h)|$ (resp.\  $|B^s(\h)|<|B^w(\h)|$), otherwise (\ref{condtrenino}) is not satisifed. The saddles that are possibly visited by the sliding path are in $\cN_{k,k'}^{\a,\a'}$ except the last one, thus by Lemma \ref{trenini}(i) they are essential. The last configuration visited during the sliding of the bar $B^n(\h)$ around the frame-angle $c^{nw}(\h')$ is depicted in Figure \ref{fig:figA6} on the right-hand side. This configuration has energy $\gw-U_1+U_2$ and therefore it is not a saddle and is in $\cC_{\vuoto}^{\pieno}(\gw-H(\pieno))$. By Propositions \ref{selle2} and and \ref{Kvuoto}(ii)-(b), the latter implies that the saddles that could be visited are either unessential or in $\cN_{2,k'}^{\a,\a'}$ and therefore the case B(iii) is concluded.

\medskip
\noindent
{\bf Case C.} Without loss of generality we consider $\h$ as in Figure \ref{fig:figA7}(a). If we are considering the case in which a $1$-translation of a bar is possible and takes place, then by Lemma \ref{trasless} the saddles that are crossed are essential and in $\cN_0^{\a'}\cup\cN_1^{\a}$. Starting from this configuration it is possible to obtain two occupied frame-angles: this situation has been already analyzed in Case B. Thus we can reduce our proof to the case in which there is no translation of a bar and therefore there is the activation of a sliding of a bar around a frame-angle. If the free particle is attached to the bar $B^n(\h)$ (resp.\ $B^s(\h)$), by Lemma \ref{coltorow} we know that it is not possible to complete the sliding of the bar $B^w(\h)$ around the frame-angle $c^{wn}(\h')$ (resp.\ $c^{ws}(\h')$). If the free particle is attached to the bar $B^e(\h)$, we deduce that (\ref{condcorner}) is not satisfied. In the last two cases by Lemma \ref{trenini}(ii) we know that the saddles that are visited are unessential. If the free particle is attached to the bar $B^w(\h)$, it is possible to slide the bar $B^n(\h)$ around the frame-angle $c^{nw}(\h')$ when $|B^n(\h)|<|B^w(\h)|$, otherwise (\ref{condtrenino}) is not satisifed. The saddles that are possibly visited by the sliding path are in $\cN_{k,k'}^{\a,\a'}$  except the last one, thus by Lemma \ref{trenini}(i) they are essential. The last configuration visited during this sliding of a bar is depicted in Figure \ref{fig:figA7}(b). This configuration has energy $\gw-U_1+U_2$ and therefore it is not a saddle and is in $\cC_{\vuoto}^{\pieno}(\gw-H(\pieno))$. By Propositions \ref{selle2} and \ref{Kvuoto}(ii)-(b), the latter implies that the saddles that could be visited are either unessential or in $\cN_{2,k'}^{\a,\a'}$ and therefore the case C is concluded.

\begin{figure}
	\centering
	\begin{tikzpicture}[scale=0.32,transform shape]

		\draw [dashed] (6,3) rectangle (16,11);
		\node at (11,7) {\Huge{$(l_1^*-2)\times(l_2^*-2)$}};
		\draw(12,11)--(12,12);
		\draw(12,12)--(5,12);
		\draw(5,12)--(5,3);
		\draw(5,3)--(6,3);
		\draw(6,3)--(7,3);
		\draw(7,3)--(7,2);
		\draw(7,2)--(16,2);
		\draw(16,2)--(16,3);
		\draw(16,3)--(17,3);
		\draw(17,3)--(17,10);
		\draw(17,10)--(16,10);
		\draw(16,10)--(16,11);
		\draw(16,11)--(12,11);
		\draw(3,8) rectangle (4,9);
		\node at (12,0.5) {\Huge{(a)}};

		\draw [dashed] (23,3) rectangle (33,11);
		\node at (28,7) {\Huge{$(l_1^*-2)\times(l_2^*-2)$}};
		\draw[dashed](22,11)--(22,3);
		\draw(22,3)--(23,3);
		\draw(23,3)--(24,3);
		\draw(24,3)--(24,2);
		\draw(24,2)--(33,2);
		\draw(33,2)--(33,3);
		\draw(33,3)--(34,3);
		\draw(34,3)--(34,10);
		\draw(34,10)--(33,10);
		\draw(33,10)--(33,11);
		\draw(33,11)--(21,11);
		\draw(21,11)--(21,4);
		\draw(21,4)--(22,4);
		\draw(22,4)--(22,3);
		\draw(21,12) rectangle (22,13);
		\node at (29,0.5) {\Huge{(b)}};
		
		\draw [dashed] (40,3) rectangle (50,11);
		\node at (45,7) {\Huge{$(l_1^*-2)\times(l_2^*-2)$}};
		\draw(51,11)--(50,11);
		\draw(50,11)--(50,12);
		\draw(51,11)--(51,4);
		\draw(51,4)--(50,4);
		\draw(50,12)--(41,12);
		\draw(41,12)--(41,11);
		\draw(41,11)--(40,11);
		\draw(40,11)--(40,10);
		\draw(40,10)--(39,10);
		\draw(39,10)--(39,4);
		\draw(39,4)--(40,4);
		\draw(40,4)--(40,3);
		\draw(40,3)--(41,3);
		\draw(41,3)--(41,2);
		\draw(41,2)--(50,2);
		\draw(50,2)--(50,3);
		\draw(50,3)--(50,4);
		\draw(37,8) rectangle (38,9);
		
	\end{tikzpicture}
	
	\vskip -0.2 cm
	\caption{Case C: in (a) we depict a possible starting configuration $\h\in\cwgeo$ and in (b) the configuration $\widetilde\h$ obtained from $\h$ after the sliding of the bar $B^n(\h)$ around the frame-angle $c^{nw}(\h')$. Case D: in (c) we depict a possible starting configuration $\h\in\cwgeo$.}
	\label{fig:figA7}
\end{figure}

\medskip
\noindent
{\bf Case D.} Without loss of generality we consider $\h$ as in Figure \ref{fig:figA7}(c). If we are considering the case in which a $1$-translation of a bar is possible and takes place, then by Lemma \ref{trasless} the saddles that are crossed are essential and in $\cN_0^{\a'}\cup\cN_1^{\a}$. Starting from this configuration it is possible to obtain one or two occupied frame-angles: these situations have been already analyzed in Cases C and B respectively. Thus we can reduce our proof to the case in which there is no translation of a bar and therefore there is the activation of a sliding of a bar around a frame-angle. If the free particle is attached to the bar $B^n(\h)$ (resp.\ $B^s(\h)$), by Lemma \ref{coltorow} we know that it is not possible to complete the sliding of the bar $B^w(\h)$ around the frame-angle $c^{wn}(\h')$ (resp.\ $c^{ws}(\h')$). If the free particle is attached to the bar $B^w(\h)$ or $B^e(\h)$, we deduce that (\ref{condcorner}) is not satisfied. In the last two cases by Lemma \ref{trenini}(ii) we know that the saddles that are visited are unessential. This concludes case D.
\end{proof*}

\end{document}